\def\Santtu{Santtu Söderholm}

\def\Maryam{Maryam Samavaki}

\def\Sampsa{Sampsa Pursiainen}

\def\myauthors{\Santtu,\Maryam,\Sampsa}

\def\mytitle{A Complete-Electrode-Model-Based Forward Approach for Transcranial Temporal Interference Stimulation with Linearization: A Numerical Simulation Study}

\def\mysubject{A scientific paper on simulating tTIS via a finite element method.}

\def\mykeywords{%
Transcranial Temporal Interference Stimulation (tTIS);%
Non-Invasive Brain Stimulation;%
Deep Brain Stimulation;%
Complete Electrode ModeL (CEM);%
Linear Programming%
}

\DocumentMetadata{
    pdfversion=1.7,
    pdfstandard=A-2b,
    lang=en,
}

\documentclass[5p,times,nopreprintline]{elsarticle}

\usepackage{subcaption}
\usepackage{mathtools} 
\usepackage{hyperref}
\usepackage{xcolor}
\usepackage[format=plain, labelfont={bf}]{caption}
\usepackage{bm}
\usepackage{lineno}
\usepackage{siunitx}
\usepackage{multirow}
\usepackage{booktabs}

\usepackage{tikz}
\usepackage[outline]{contour} 
\usetikzlibrary{angles,quotes} 
\contourlength{1.2pt}

\tikzset{>=latex} 

\colorlet{veccol}{green!70!black}
\colorlet{vcol}{green!70!black}
\colorlet{xcol}{blue!85!black}
\colorlet{projcol}{red!80}
\colorlet{initcol}{black!60}
\colorlet{diffcol}{black!100}
\colorlet{unitcol}{xcol!60!black!85}
\colorlet{myblue}{blue!70!black}
\colorlet{myred}{red!90!black}
\colorlet{mypurple}{blue!50!red!80!black!80}

\tikzstyle{vector}=[->,very thick,xcol]

\numberwithin{figure}{section}

\numberwithin{table}{section}

\numberwithin{equation}{section}

\NewDocumentCommand\numSetCmd{m}{\mathbb{#1}}

\NewDocumentCommand\Rset{}{\numSetCmd{R}}
\NewDocumentCommand\Cset{}{\numSetCmd{C}}

\NewDocumentCommand\vectorcmd{m}{\mathbf{#1}}

\NewDocumentCommand\matrixcmd{m}{\mathbf{#1}}

\NewDocumentCommand\transferMat{}{\matrixcmd{T}}

\NewDocumentCommand\schurMat{}{\matrixcmd{S}}

\NewDocumentCommand\stiffMat{}{\matrixcmd{A}}

\NewDocumentCommand\electrodeImpedanceMat{}{\matrixcmd{C}}

\NewDocumentCommand\resistanceMat{}{\matrixcmd{R}}

\NewDocumentCommand\massMat{}{\matrixcmd{M}}

\NewDocumentCommand\idMat{}{\matrixcmd{I}}

\NewDocumentCommand\inverse{m}{#1^{-1}}

\NewDocumentCommand\leadFieldMat{}{\matrixcmd{L}}

\NewDocumentCommand\deltaLLinRef{}{\Delta\leadFieldMat_{\mathrm{ref}}^{\mathrm{lin}}}

\NewDocumentCommand\biCGStab{}{\text{biCGStab}}

\NewDocumentCommand\PCG{}{\text{PCG}}

\NewDocumentCommand\realPart{}{\mathfrak{R}}

\NewDocumentCommand\imagPart{m}{\mathfrak{I}#1}

\NewDocumentCommand\Zell{o}{
    \IfNoValueTF{#1}{
        {Z_\ell}
    }{
        {Z_{#1}}
    }
}

\NewDocumentCommand\Zvec{}{\vectorcmd{Z}}

\NewDocumentCommand\jacobianOf{smm}{%
    \IfBooleanTF{#1}{%
        \paren*{\pder{#2}{#3}}%
    }{%
        \pder{#2}{#3}%
    }%
}

\NewDocumentCommand\transpose{m}{#1^{\mathsf{T}}}

\NewDocumentCommand\ctranspose{m}{#1^{\mathsf{H}}}

\NewDocumentCommand\beatF{}{\Delta f}

\DeclarePairedDelimiter\areaDelim{\lvert}{\rvert}

\NewDocumentCommand\areaOf{m}{\areaDelim*{#1}}

\NewDocumentCommand\electrode{O{\ell}}{e_{#1}}

\NewDocumentCommand\frequency{}{f}

\NewDocumentCommand\angFreq{}{\omega}

\NewDocumentCommand\electricField{}{\vectorcmd{E}}

\NewDocumentCommand\Nof{m}{N_{\mathrm{#1}}}

\NewDocumentCommand\azimuth{}{\alpha}

\NewDocumentCommand\elevation{}{\beta}

\NewDocumentCommand\bavec{m}{\vec{\mathbf{#1}}}

\NewDocumentCommand\admittivity{}{{\bm\gamma}}

\NewDocumentCommand\conductivity{}{{\bm\sigma}}

\NewDocumentCommand\permittivity{O{}}{{\bm\varepsilon_{\mathrm{#1}}}}

\NewDocumentCommand\relativepermittivity{}{{\bm\varepsilon_r}}

\NewDocumentCommand\basisFn{}{\psi}
\NewDocumentCommand\testFn{}{v}

\NewDocumentCommand\conj{m}{\overline{#1}}

\NewDocumentCommand\domain{}{\mathbf\Omega}
\NewDocumentCommand\domainboundary{}{\partial\domain}

\NewDocumentCommand\surface{}{S}

\NewDocumentCommand\normalVec{}{\vec{\vectorcmd{n}}}

\NewDocumentCommand\current{}{I}

\NewDocumentCommand\potential{}{u}
\NewDocumentCommand\potentialVec{}{\vectorcmd{u}}
\NewDocumentCommand\voltage{}{{U}}

\NewDocumentCommand\electrodeCurrentMat{}{\matrixcmd{B}}
\NewDocumentCommand\potentialLoss{}{U}
\NewDocumentCommand\potentialLossVec{}{\mathbf{\mathcal{U}}}

\NewDocumentCommand\sobolevSpace{}{\mathrm{H}}

\NewDocumentCommand\diff{}{\mathrm{d}}

\NewDocumentCommand\der{mm}{\frac{\diff{#1}}{\diff{#2}}}
\NewDocumentCommand\pder{mm}{\frac{\partial{#1}}{\partial{#2}}}

\NewDocumentCommand\refVal{m}{#1_{\mathrm{ref}}}
\NewDocumentCommand\iniVal{m}{#1_{\mathrm{ini}}}

\NewDocumentCommand\linVal{m}{#1_{\mathrm{lin}}}

\NewDocumentCommand\iniR{}{\iniVal\resistanceMat}
\NewDocumentCommand\refR{}{\refVal\resistanceMat}
\NewDocumentCommand\linR{}{\linVal\resistanceMat}

\NewDocumentCommand\iniL{}{\iniVal\leadFieldMat}
\NewDocumentCommand\refL{}{\refVal\leadFieldMat}
\NewDocumentCommand\linL{}{\linVal\leadFieldMat}

\NewDocumentCommand\refDiff{m}{\Delta#1_{\mathrm{ref}}}

\NewDocumentCommand\linDiff{m}{\Delta#1_{\mathrm{lin}}}

\NewDocumentCommand\capacitance{}{C}

\NewDocumentCommand\geometricConstant{}{G}

\NewDocumentCommand\unitOf{m}{\brackets*{#1}}

\NewDocumentCommand\contactResistance{}{R_{\mathrm{c}}}
\NewDocumentCommand\doubleLayerResistance{}{R_{\mathrm{d}}}
\NewDocumentCommand\doubleLayerCapacitance{}{C_{\mathrm{d}}}

\NewDocumentCommand\currentPattern{}{\vectorcmd{i}}

\NewDocumentCommand\position{}{\vectorcmd{x}}

\DeclarePairedDelimiter\abs{\lvert}{\rvert}
\DeclarePairedDelimiter\paren{(}{)}
\DeclarePairedDelimiter\args{(}{)}
\DeclarePairedDelimiter\set{\{}{\}}
\DeclarePairedDelimiter\norm{\lVert}{\rVert}

\DeclarePairedDelimiter{\brackets}{[}{]}

\NewDocumentCommand\eell{}{{e_\ell}}

\NewDocumentCommand\restrictionOf{mm}{\paren*{#1\;\middle|\,#2}}

\DeclareMathOperator\vecspan{span}

\DeclareMathOperator\iu{i}

\DeclareMathOperator\diag{diag}

\NewDocumentCommand\labeledQuantity{mO{}}{{#1}_{\mathrm{#2}}}

\NewDocumentCommand\quantile{O{}}{{q_{\mathrm{#1}}}}

\NewDocumentCommand\volumeCurrentDensity{}{{\matrixcmd{J}_{\mathrm{v}}}}

\DeclareMathOperator{\interferenceField}{IF}

\NewDocumentCommand\errorDiscrepancy{}{\delta}

\NewDocumentCommand\volumeCurrent{O{}}{{\matrixcmd{I}_{\mathrm{v}}}}

\NewDocumentCommand\iniJv{}{\volumeCurrentDensity_{\mathrm{ini}}}

\NewDocumentCommand\refJv{}{\volumeCurrentDensity_{\mathrm{ref}}}

\NewDocumentCommand\linJv{}{\volumeCurrentDensity_{\mathrm{lin}}}

\NewDocumentCommand\tetrahedron{}{\vectorcmd{t}}

\NewDocumentCommand\relDelta{}{{\Delta_{\mathrm{rel}}}}

\NewDocumentCommand\divergence{}{\nabla\cdot}

\DeclareMathOperator{\support}{supp}

\NewDocumentCommand\conductanceDensity{}{\matrixcmd{\Sigma}}

\NewDocumentCommand\gradient{}{\nabla}

\DeclareMathOperator\dBT{dBT}

\NewDocumentCommand\logten{}{\log_{10}}

\NewDocumentCommand\someField{}{g}

\usepackage{url}

\hypersetup{
    colorlinks=true,
    linkcolor=.,
    anchorcolor=.,
    citecolor=.,
    filecolor=.,
    menucolor=.,
    runcolor=.,
    urlcolor=blue,
    pdftitle={\mytitle},
    pdfauthor={\myauthors},
    pdfsubject={\mysubject},
    pdfkeywords={\mykeywords},
}

\NewDocumentCommand\CtoR{}{$\Cset$-to-$\Rset$}

\begin{document}

\begin{frontmatter}

\title{\mytitle}

\author[tau]{\Santtu\corref{cor1}}

\address[tau]{Computing Sciences Unit, Faculty of Information Technology and Communication Sciences, P.O. Box 1001, FI-33014 Tampere University, Tampere, Finland}

\address[uef]{Department of Physics and Mathematics, Faculty of Science, Forestry and Technology, University of Eastern Finland, P.O. Box 111, FI-80101 Joensuu, Finland}

\cortext[cor1]{Corresponding author at: Tietotalo building, Korkeakoulunkatu 1, Tampere, 33720, FI}

\ead{santtu.soderholm@tuni.fi}

\author[uef]{\Maryam}

\author[tau]{\Sampsa}

\begin{abstract}

\paragraph{Background and Objective} Transcranial temporal
interference stimulation~(tTIS) is a promising non-invasive
brain stimulation technique in which interference between
electrical current fields extends the possibilities of
electrical brain stimulation.  This study aims to establish
and evaluate the complete electrode model (CEM), a set of
boundary conditions incorporating electrode impedance and
contact patch, as a forward simulation technique for tTIS
and investigate linearized CEM as a surrogate.

\paragraph{Methods} The electric potential distribution
in the head was simulated using the finite element method
(FEM) combined with the CEM. A frequency-dependent lead
field formulation and its linearization with respect
to electrode impedance were implemented to reflect
the capacitive properties of electrode interfaces. The
simulations were conducted in a realistic multi-compartment
head model, where the accuracy of the surrogate was
evaluated against the full nonlinear model across
two different changes in contact impedance from the
initial $\contactResistance = \qty{270}{\ohm}$ to
$\contactResistance\in\set{\num{1270},\num{5270}}\,\si\ohm$.

\paragraph{Results} The CEM-based forward simulation
successfully reproduced the volumetric stimulating
fields induced by tTIS. The linearized CEM model closely
matched the full nonlinear model within a predefined peak
signal-to-noise ratio (PSNR) threshold for relative error.
Both models exhibited the highest sensitivity near the
focal region of tTIS.

\paragraph{Conclusions} The results demonstrate that
linearization with respect to frequency can yield a
computationally efficient approximation of the lead
field, retaining sufficient accuracy within a practical
PSNR threshold. The surrogate captures key sensitivity
characteristics of the full model, particularly near the
focal region, highlighting its potential for tasks where
frequent updates of stimulation parameters are required.

\end{abstract}

\begin{keyword}
Transcranial Temporal Interference Stimulation (tTIS);
Non-Invasive Brain Stimulation;
Deep Brain Stimulation;
Complete Electrode ModeL (CEM);
Linear Programming

\end{keyword}

\end{frontmatter}

\section{Introduction}%
\label{sec:intro}

The ability to precisely stimulate neurons at depth
within the brain using electrical methods offers new
possibilities for treating brain disorders. This simulation
study concerns the mathematical and computational forward
modelling aspects of non-invasive transcranial temporal
interference stimulation (tTIS) \cite{GROSSMAN-2017,
RAMPERSAD_computational_2019, Wessel_2021, Herrmann_2021,
tTIS_2022}, which through the interference of the
stimulating fields can extend the possibilities of
classical transcranial electrical stimulation (tES)
\cite{herrmann2013transcranial,Nitsche_tDCS_2008}
and  deep brain stimulation (DBS)
\cite{lozano2019deep,BUTSON-2005,BUTSON-2006,anderson2018optimized}.
In tTIS, two or more pure sine wave currents
with high frequencies $f_1$ and $f_2$ are delivered through
contact electrodes placed on the skin to the subject's
head, where they generate an interference pattern within
the brain tissue. When these high-frequency signals
intersect within the brain, they interfere with each other,
creating a beat frequency \(\beatF\), which is represented
by the difference between the original frequencies
($\Delta\frequency=\abs{\frequency_1-\frequency_2}$)
\cite{RAMPERSAD_computational_2019}. This beat frequency
can match specific brain wave frequencies and create
a low-frequency envelope by modulating the amplitude
at the target stimulation frequency of a particular
brain function, such as electroencephalogram (EEG)
rhythms. An interference envelope can reach deep-layer
neurons in the brain, providing a non-invasive
way to influence deep-brain activity without
affecting overlying structures~\cite{GROSSMAN-2017,
acerbo2022focal,lee2016investigational}. The effectiveness
of electrical brain stimulation is significantly based on
the mapping of the lead field $\leadFieldMat: \mathbb{R}^L
\times \mathbb{C}^L \to \mathbb{C}^N$, a forward model
that describes how the volumetric current field is
distributed in response to the applied stimulation. Due
to the alternating nature of the stimulation currents,
\(\leadFieldMat\) is frequency dependent, adding an
additional layer of complexity to the optimization problem
of selecting the best frequencies and amplitudes for the
stimulation currents~\cite{Hutcheon_2000, Grossman_2018}.

In earlier studies~\cite{Prieto_optimal_tES_2024,
GALAZPRIETO-2022, Agsten-2018}, we introduced a
lead field mapping for tES based on the finite element
method (FEM), incorporating the complete electrode model
(CEM) boundary conditions~\cite{Somersalo1992,Agsten-2018},
together with a numerical implementation written for the
open Matlab toolbox Zeffiro interface~\cite{he2020zeffiro},
and  explored advanced optimization techniques (L1L1) for
multi-channel tES~\cite{GALAZPRIETO-2022}, demonstrating
its potential to enhance the precision and effectiveness
of brain stimulation methods. This study introduces a
mathematical CEM-based lead field model for interference
stimulation.

We also introduce and evaluate a surrogate
lead field model to enable, e.g., coupling the
forward model with  optimization techniques
\cite{fernandez2020unification,GALAZPRIETO-2022}, which
might require gradient-based or repetitive evaluation
of the lead field mapping. With this surrogate, we
aim to establish a clear and predictable relationship
between changes in contact impedances. This simplified
model approximates outcomes with reduced computational
effort, in situations, in which a rough correction of the
field dynamics can be considered sufficient, e.g., due
to impedance changes observed during an experiment. We
consider linearization as a potentially effective surrogate
strategy to exhibit the non-linearity of the capacitive
effects in biological tissues and electrode interfaces:  it
provides an  alternative compared to repetitive exact lead
field evaluation and has proven out to be  tolerant towards
latent non-linearity effects induced by the variable
frequency in electrical impedance tomography (EIT) studies
\cite{Boverman_linearization_2008}.

An important aspect of our work involves evaluating the
precision of the surrogate in capturing the true behavior
of the system. Towards this end, we perturbed the contact
impedance $\contactResistance$ of a single electrode by
\qty{1000}{\ohm} and \qty{5000}{\ohm} with a base stimulation
frequency $\frequency=\qty{1000}{\hertz}$ and a beat
frequency of $\beatF=\qty{10}{\hertz}$~\cite{GROSSMAN-2017}
between the $2$ stimulating electrode pairs, and observed
how the reference and linearized mappings deviated from
each other, where the reference field represents the true
behavior of the system.

\section{Methods}%
\label{sec:methods}

\subsection{Conductivity and permittivity distribution}

We base our forward model on the assumption that
the input currents are alternating, and hence the
electric field and its parameters are considered to
be complex-valued. For a given conductive object, the
admittivity  is a frequency-dependent  complex scalar
valued function  \begin{equation} \admittivity(\position,
\angFreq)=\admittivity_{ij}(\position, \angFreq)
\end{equation} for all $i, j=1, \ldots, N$, a second-rank
tensor, in which the positive real part represents the
electrical conductivity $\conductivity$ at a frequency
$\angFreq$ of the applied current, while the imaginary part
is produced by multiplying the permittivity of the tissue
with the frequency of the applied alternating current.
Its unit is $\unitOf\admittivity=\si{\siemens\per\meter}$.
Mathematically, this can be  expressed
as~\cite{Bisegna_2008, waveform_2022,wang-etal-2024}
\begin{equation}\label{eq:admittivity}
    \admittivity(\position, \angFreq)
    =
    \conductivity(\position, \angFreq)
    +
    \iu
    \angFreq
    \,
    \permittivity[0]
    \,
    \permittivity[r]
    (\position, \angFreq)
    \,,
\end{equation}
for all $\position\in\domain\subset\Rset^3$.
In the second term, \(\permittivity[0] \,
\) (\qty{8.854e-12}{\farad\per\meter}) and
$\relativepermittivity$  denote vacuum and relative
permittivity, respectively. The second term is a
function of the angular frequency $\angFreq$ with $\iu$
being the imaginary unit, that satisfies $\iu^2=-1$.
The relation between the capacitance $\capacitance$ of
an object and its $\permittivity$ is usually given by
$\capacitance = \geometricConstant\permittivity$, where
$\geometricConstant$ is a geometric scaling factor with
$\unitOf\geometricConstant = \si\meter$.

This change in current flow results in variations
in the voltage measurements made at the surface
electrodes. EIT employs this recorded boundary
voltage values to reconstruct images of the complex
conductivity distribution within the imaging object
\cite{Metherall_1996, Holder_2005}. These boundary
voltages are computed to numerically simulate
measured data (Appendix \ref{app:BConditions}). The
second-rank tensor $\conductivity$ ($\mathrm{Sm}^{-1}$)
is symmetric, positive-definite, and represents the
electrical conductivity of the tissue. The imaginary
term $\angFreq\,\permittivity(\position, \angFreq)$ is
conneced to the capacitance, as described above. Here,
$\angFreq\in(\angFreq_1,\cdots, \angFreq_L)$, where
$L$ is the number of electrodes and $\angFreq_\ell$
denotes the angular frequency of the injected current
through electrode $\electrode[\ell]$, with a typical
\(\angFreq_\ell\approx\qty{1}{\kilo\hertz}\) for
tTIS~\cite{GROSSMAN-2017}.

\subsection{Complete Electrode Model}%
\label{app:BConditions}

The primary partial differential equation (PDE)
that describes the electric potential within the
three-dimensional (3D) smooth bounded head model
$\domain\subset\mathbb{R}^3$, where $\domain$
is a simply connected open region with a smooth
boundary $\partial\domain$. The model incorporates
a set of electrodes $\electrode[\ell]$ attached to
$\partial\domain$, each with a surface contact area
$\areaOf{\electrode}$ and potential $\potential_{\ell}$.
The current flowing through the $\ell$-th electrode
is indicated by $\current_\ell$. A current $\current$
is injected at a low frequency with a chosen
pair of coupled electrodes, $\electrode[a]$ and
$\electrode[b]$, to generate potential over the
domain $\domain$. Then the resulting complex potential
distribution $u=(u_1, \ldots, u_N)$ satisfies the governing
Poisson-type elliptic partial differential equation (PDE)
given by~\cite{GALAZPRIETO-2022, dougherty-etal-2014}:
\begin{equation}
    \divergence
    (
        \admittivity(\position, \angFreq)
        \nabla
        \potential(\position, \angFreq)
    )
    =
    0
    \quad\quad\quad\hbox{with}\quad\position\in\domain\,,
    \label{poisson_com_2}
\end{equation}
which follows from the quasi-static formulation of
Maxwell's equations. Then the resulting scalar potential
distribution $\potential$ on each electrode satisfies~\cite{dougherty-etal-2014}
\begin{equation}
    (
        \admittivity(\position, \angFreq_\ell)
        \nabla
        \potential
        (\position, \angFreq_\ell)
    )
    \cdot
    \normalVec
    =
    \current_\ell
    \quad\quad\quad \hbox{with}\quad\position\in\electrode[\ell]\,,
    \label{poisson_com_c}
\end{equation}
for all $\ell=1, \cdots, L$ and positions $\position$.
The frequency-dependent admittivity $\admittivity$
($\mathrm{Sm}^{-1}$) is described by a complex tensor,
while each electrode potential ${ \potential_\ell}$ with
$\unitOf\potential = \si\volt$ for all $\ell=1, \ldots, L$
is a complex scalar field, and $\angFreq_\ell$ denotes the
angular frequency of the applied current.

If the capacitive effects are considered to be negligible,
i.e., $\permittivity(\position, \angFreq) \approx {\bf
0}$, the admittivity $\admittivity$ can be approximated
by the real valued conductivity $\conductivity$. In CEM,
the complex boundary conditions for the complete electrode
model are presented as follows~\cite{Agsten-2018}:
\begin{align}
\label{CEM1}
    0
    &=
     \admittivity\,\pder{ u}{\bavec n}(\position)
    &
    \text{for } \position \in\domainboundary\setminus\cup_{\ell=1}^L \electrode[\ell]\,,
    \\
    \label{CEM2}
   \current_\ell
    &=
    \int_{\electrode[\ell]}{ \admittivity\,\pder{ u}{\bavec n}(\position) \,\diff\surface}
    &
    \text{for } \ell \!= \!1, \ldots,  L\,,
    \\
    \label{CEM3}
    {U_\ell}
    &=
    \potential(\position)
    +
    \tilde\Zell
    \admittivity
    \pder{ u}{\bavec n}(\position)
    &
    \text{ for }  \position \in   \electrode[\ell], \, \ell \!=\! 1,  \ldots,  L
    \,,
\end{align}
where the effective contact impedance and potential mean
and voltage loss across an electrode contact surface are
given, respectively, by
\begin{align}
\label{CEM4}
   \tilde\Zell
    &=
    \Zell
    \areaOf{ \electrode[\ell]}
    \quad\text{and}
    \\
    \label{CEM5}
    {\voltage_\ell}
    &=
    \frac 1 { \areaOf { \electrode[\ell] }  }
    \int_{\electrode[\ell]}{u \, \diff\surface}
    +
    \Zell\current_\ell
    \,,
\end{align}
where \(\areaOf\electrode =
\int_{\electrode}\diff\surface\) is the contact area of
an electrode. Furthermore, to guarantee the existence and
uniqueness of the solution, the following two constraints
on the injected currents and the measured voltages are
required due to the conservation of charge and the choice
of ground~\cite{Somersalo1992}:
\begin{align*}
    \sum_{\ell=1}^L \current_\ell=0\,,
\quad\quad
    \sum_{\ell=1}^L {\voltage_\ell}=0\,.
\end{align*}
The process of determining the potential distribution
$\potentialVec$ within the domain $\domain$ and the
voltages $\voltage_\ell$ on the electrodes, based on the
given admittivity $\admittivity$ and boundary conditions,
is known as the forward problem.

\subsection{Lead field matrix and potential field}%
\label{app:Lead_Field}

The resistance matrix $\resistanceMat$
~\cite{GALAZPRIETO-2022} with
$\unitOf\resistanceMat=\si\ohm$, described in
\ref{app:res_matrix}, provides a mapping  $\potentialVec
= \resistanceMat\currentPattern$ from the input
currents $\currentPattern$ to the discretized potential
field $\potentialVec$ inside the head $\domain$.
Further, we define the mapping from the stimulation
current pattern $\currentPattern$ to the volume
current densities $\volumeCurrentDensity$ in the
form of a lead field matrix $\leadFieldMat$, with
$\unitOf\leadFieldMat=\si{\per\meter\squared}$, defined as
\begin{equation}
\label{eq:L_Matrix}
    {\leadFieldMat}
    =
    \conductanceDensity
    \resistanceMat
    \,.
\end{equation}
Here the \emph{conductance
density} $\conductanceDensity$ with
$\unitOf\conductanceDensity=\si{\siemens\per\meter\squared}
$ is defined in the domain $\domain=\cup_{i}
\tetrahedron_i$ constructed from tetrahedra
$\tetrahedron_i$ via
\begin{equation}
     \conductanceDensity_{i,h}
     =
     \begin{cases}
        -
        \admittivity_i
        \gradient
        \basisFn_h
        ,
        & \hbox{if}
        \quad
        \support\{\basisFn_h\} \cap\tetrahedron_i\neq\emptyset \\
        0 ,
        & \text{otherwise} \,.
    \end{cases}
\end{equation}
Of note is that in general, the lead field
$\leadFieldMat$ contains both real and imaginary
components $\realPart\leadFieldMat^{(k)}$ and
$\imagPart\leadFieldMat^{(k)}$. This is due to the
complex nature of the admittivity $\admittivity$ and the
electrode impedances $\Zell$. Formula (\ref{eq:L_Matrix})
can be considered as the \textit{forward mapping},
producing the volume current density distribution
$\volumeCurrentDensity = \leadFieldMat\currentPattern$
in the domain, in the process of optimizing the current
pattern $\currentPattern$.

\subsection{Linearization of \(\resistanceMat\) as a function of impedance}%
\label{sec:linearization}

 Following from the complete electrode model, the efficacy
and safety of electrical stimulation is influenced by
tissue conductivity and electrode impedance, omitting
tissues' capacitive effects, which we assume to be minor
at kilohertz frequency range. Modifying the impedance of
each electrode, often reached by changing the frequency
of alternating voltage used for stimulation, impacts
the gradient of   electrode potential. By understanding
how impedance $\Zell$ affects the resistance matrix
\(\resistanceMat\), we can optimize stimulation parameters
to achieve desired outcomes while minimizing adverse
effects. Using linearization, the resistance matrix
\(\resistanceMat=\resistanceMat\paren\Zvec\) is updated
based on changes in a subset of the impedances $\Zvec$
as follows:
\begin{equation}\label{eq:linR}
    \resistanceMat
    \paren*{\Zvec+\diff\Zvec}
    =
    \resistanceMat\paren\Zvec
    +
    \sum_\ell
    \pder{\resistanceMat\paren\Zell}\Zell
    \diff\Zell\,,
\end{equation}
for selected \(\ell\). To compute \(\diff\Zell\) for each
chosen \(\ell\), we first update the contact impedance
$\contactResistance$ of an electrode by applying a
perturbation \(\diff\contactResistance\) via
\begin{equation}
    \contactResistance{}_{,\mathrm{new}}
    =
    \contactResistance{}_{,\mathrm{old}}
    +
    \diff\contactResistance
    \,.
\end{equation}
A new impedance at the updated frequency $f_{\mathrm{new}}$
is computed using the typical impedance relation for
a double-layer electrode, with a parallel double layer
capacitor with resistance \(\doubleLayerResistance\)
and capacitance \(\doubleLayerCapacitance\) in parallel,
and a contact resistance \(\contactResistance\) in
series~\cite{BUTSON-2005, electrode-modelling-2017}:
\begin{equation}
\label{eq:impedance-of-f}
    \Zell_{,\mathrm{new}}
    =
    \frac
        1
        {
            \frac 1 \doubleLayerResistance
            -
            \frac 1 {
                \frac
                {1\iu}
                {2\pi f\doubleLayerCapacitance}
            }
        }
    +
    \contactResistance{}_{,\mathrm{new}}\,.
\end{equation}
Once the new impedances have been computed, changes or
perturbations in impedance follow the equation
\begin{equation}\label{eq:impedance-perturbation}
\diff\Zell=\Zell_{,\mathrm{new}} - \Zell_{,\mathrm{old}}\,.
\end{equation}

\begin{figure}[!h]
    \centering
    \includegraphics[width=0.95\linewidth]{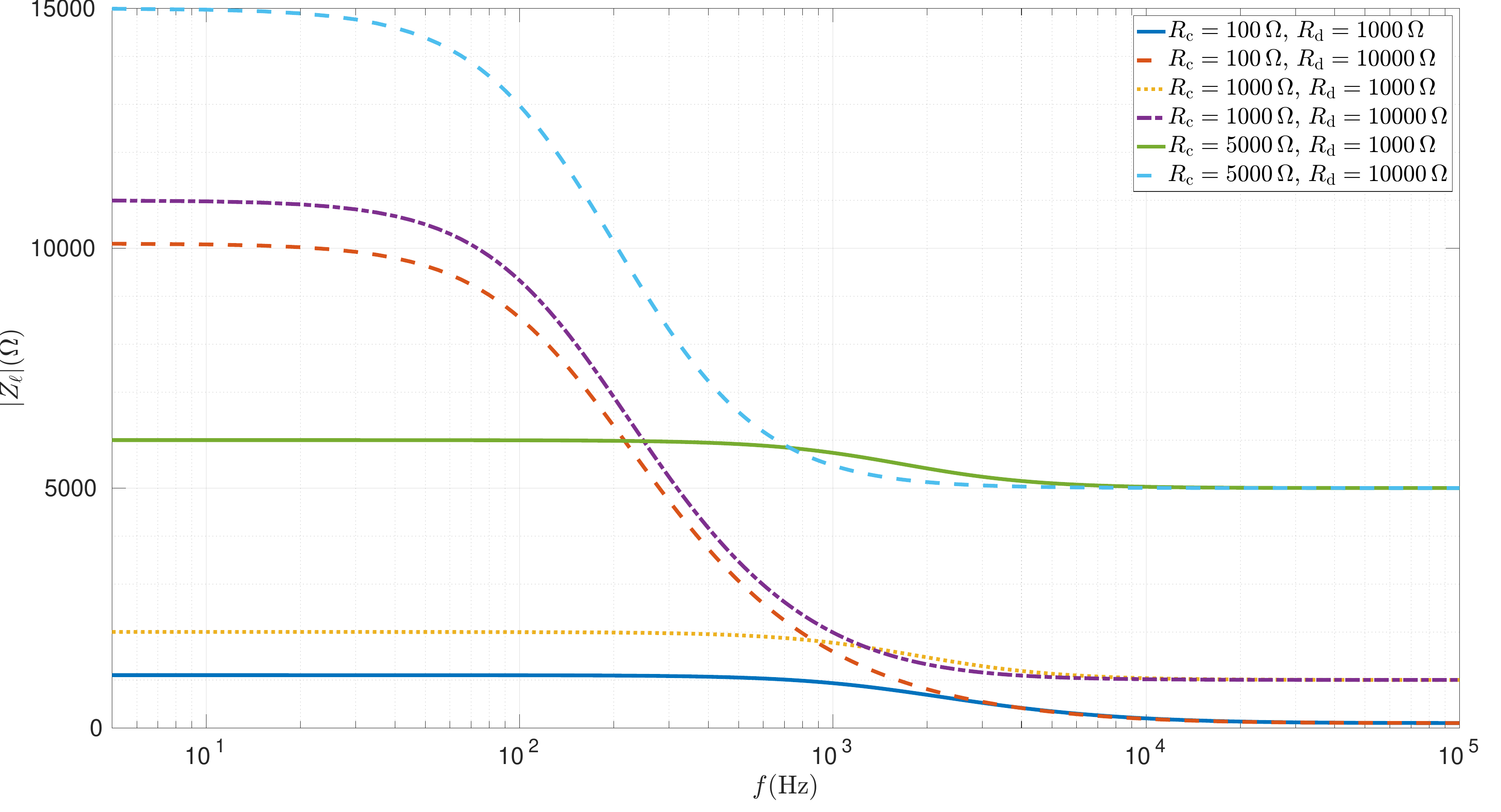}
    \caption{Expected behaviour of \(\abs{\Zell}\) as function of stimulation frequency~\cite{electrode-modelling-2017}, computed via \eqref{eq:impedance-of-f}, with \mbox{\(\contactResistance\in\set{\num{100},\num{1000},\num{10000}}\,\si\ohm\)}, \mbox{\(\doubleLayerResistance\in\set{\num{1000},\num{10000}}\,\si\ohm\)} and \(\doubleLayerCapacitance=\qty{0.1}{\micro\farad}\).}
    \label{fig:impedance-of-f}
\end{figure}

The behaviour of the impedance of such double layer
electrodes for selected \(\contactResistance\)
and \(\doubleLayerResistance\) can be seen in
Figure~\ref{fig:impedance-of-f}. At a large frequency
limit, the values of \(\abs{\Zell}\) approach
\(\contactResistance\), whereas the low-frequency limit is
\(\contactResistance + \doubleLayerResistance\), as can be
seen from \eqref{eq:impedance-of-f}.

The utilized value for the double layer resistance
was \(\doubleLayerResistance=\qty{10}{\kilo\ohm}\),
whereas capacitance was chosen to be
\(\doubleLayerCapacitance=\qty{0.1}{\micro\farad}\)
~\cite{BUTSON-2005}. Sensible initial values
of \(\contactResistance\) can be
expected to float in the neighborhood of
\num{100}--\qty{460}{\ohm}~\cite{Journee-2004}, and in our
case it was chosen as \qty{270}{\ohm}.

The total derivative of the resistance matrix
\(\resistanceMat\) above corresponds to a Jacobian
matrix $\iniVal {\jacobianOf*\resistanceMat\Zvec}$,
which is employed to linearize the lead field matrix
${\leadFieldMat} : \mathbb{R}^ N \times \mathbb{C}^L \to
\mathbb{C}^N$. It provides a linear approximation of how
small changes in the impedance vector $\Zvec$ affects
the resistance matrix $\resistanceMat\paren\Zvec$ around
the initial point $\iniVal\Zvec$. The lead field matrix
\(\leadFieldMat=\leadFieldMat\paren\resistanceMat\)
is a function that depends on the impedance vector
$\Zvec\in\Rset^L$, where each component $\Zell$ corresponds
to an impedance of a stimulating electrode. The Jacobian
matrix $\iniVal{\jacobianOf*\resistanceMat\Zvec}$
represents the first-order partial derivatives of the
resistance matrix $\resistanceMat$ with respect to the
components of the impedance vector ${\Zvec}$. It can be
expressed in matrix form as:
\begin{align}
    \iniVal{\jacobianOf*\resistanceMat\Zvec}
    \paren\Zvec
    &=
    \restrictionOf
        {\jacobianOf\resistanceMat\Zvec}
        {\iniVal\Zvec}
    =
    \begin{bmatrix}
        \begin{bmatrix}
        \pder{}{\Zell[1]}
        \resistanceMat
        \end{bmatrix}
        \\
        \vdots
        \\
        \begin{bmatrix}
        \pder{}{\Zell[N]}
        \resistanceMat
        \end{bmatrix}
    \end{bmatrix}
    _{\Zvec=\iniVal\Zvec}
\end{align}
where ${\iniVal\Zvec}=\paren{\Zell[1]_{,\mathsf{ini}},
\ldots,\Zell[L]_{,\mathsf{ini}}}$ and
$\restrictionOf{f}{x}$ denotes the restriction of a
function $f$ to $x$.

\subsection{Interference stimulation}%
\label{sec:interference}

\begin{figure}[!h]
    \centering
    \begin{subfigure}[b]{0.95\linewidth}
        \centering
        \includegraphics[width=\linewidth]{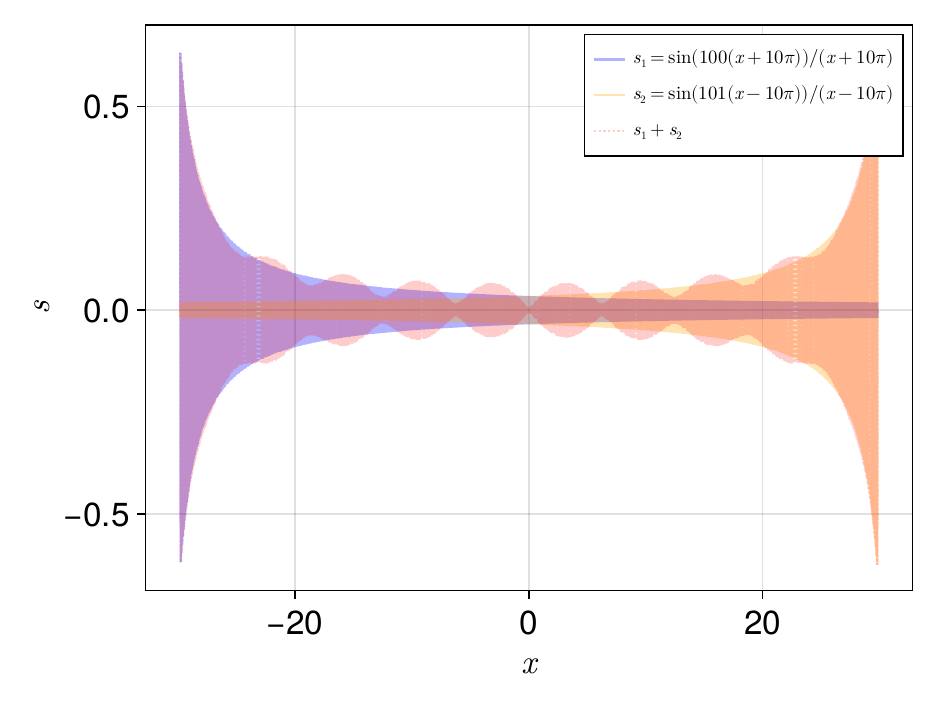}
        \caption{$s_1 + s_2$}
        \label{fig:sum-of-signals}
    \end{subfigure}

    \begin{subfigure}[b]{0.95\linewidth}
        \centering
        \includegraphics[width=\linewidth]{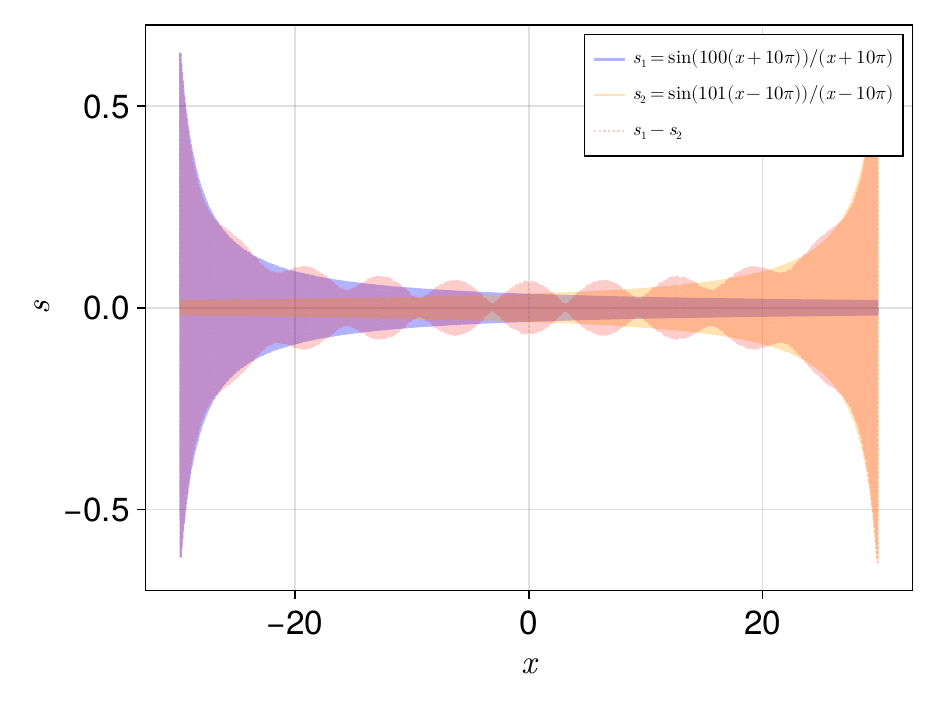}
        \caption{$s_1 - s_2$}
        \label{fig:abs-diff-of-signals}
    \end{subfigure}
    \caption{
        Interference patterns of $2$ sinusoidal
        signals $s_1$ and $s_2$ of roughly equal
        source strength, in a 1-dimensional domain.
        Interference mainly occurs where the envelope
        of each signal is of roughly equal magnitude.
        The difference between summation and subtraction
        mainly introduces a phase shift in the envelope
        of the modulating signal, as can be seen
        from Subfigures \ref{fig:sum-of-signals} and
        \ref{fig:abs-diff-of-signals}. The frequency of the
        stimulating interference signal is considered as
        the envelope of the summed signals $s_1$ and $s_2$
        , seen in red.
    }
    \label{fig:sinusoids}
\end{figure}

While the lead field model of Section~\ref{sec:methods},
equipped with complete electrode boundary
conditions~\cite{Agsten-2018} and extended to complex
electrode impedances $\Zell$ and tissue admittivity
$\admittivity$, seems to work as intended when modelling
the spread of quasi-static electromagnetic radiation from
stimulating electrodes into the brain tissue, additional
considerations need to be taken into account, when
applying the model to tTIS. In order for interference-based
stimulation to work, the stimulating currents need to
be chosen such, that the magnitude of the resulting
electric fields $\electricField_i$ emanating from a set of
electrode pairs have amplitudes $\norm{\electricField_i}$,
that are roughly equal. This is demonstrated in
Figure~\ref{fig:sinusoids}.

If the \(2\) electric fields $\electricField_1 =
\electricField_1(\position)$ and $\electricField_2 =
\electricField_2(\position)$ are unequal in magnitude
at position $\position$, the stronger of the 2 ends up
dominating the other, and there is no clear lower-frequency
envelope formed in the signal $\electricField =
\sum_i\electricField_i$ at $\position$. On the other hand,
if $\norm{\electricField_1}\approx\norm{\electricField_2}$,
we see a clear envelope with a beat frequency $\beatF =
\abs{\frequency_2 - \frequency_1}$~\cite{GROSSMAN-2017}.
This has implications regarding how stimulating electrodes
are to be positioned on the scalp of a patient, and how
much current is to be pushed through them to achieve
a similar electric field intensity at a given location
to achieve optimal stimulation at $\position$. However,
the beat frequency $\beatF$ should be chosen such that
it is equal to the frequency of the type of activity
that one hopes to induce. For $\alpha$ activity, this
is roughly \num{8}--\qty{12}\hertz~\cite{Samaha-2024}.
For $\beta$ activity, $\beatF$ should be roughly
\num{13}--\qty{30}{\hertz}~\cite{Noguchi-2022}.

\subsection{Numerical Simulations}%
\label{sec:Numerical_Simulations}

\begin{figure}[!h]
    \def\subFigWidth{0.45\linewidth}
    \begin{subfigure}[t]{0.45\linewidth}
        \includegraphics[width=\linewidth]{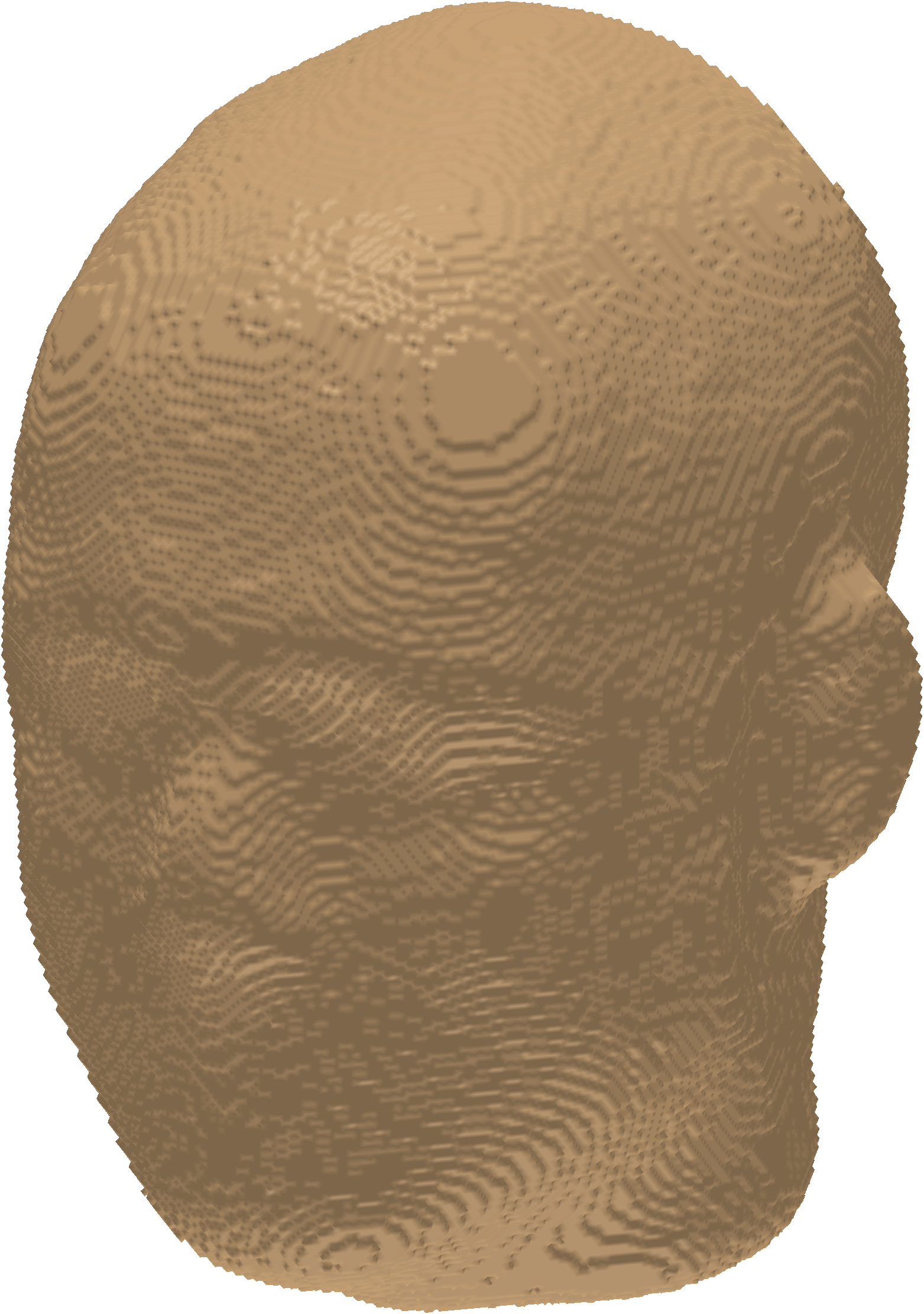}
    \end{subfigure}
    \hfill
    \begin{subfigure}[t]{0.45\linewidth}
        \includegraphics[width=\linewidth]{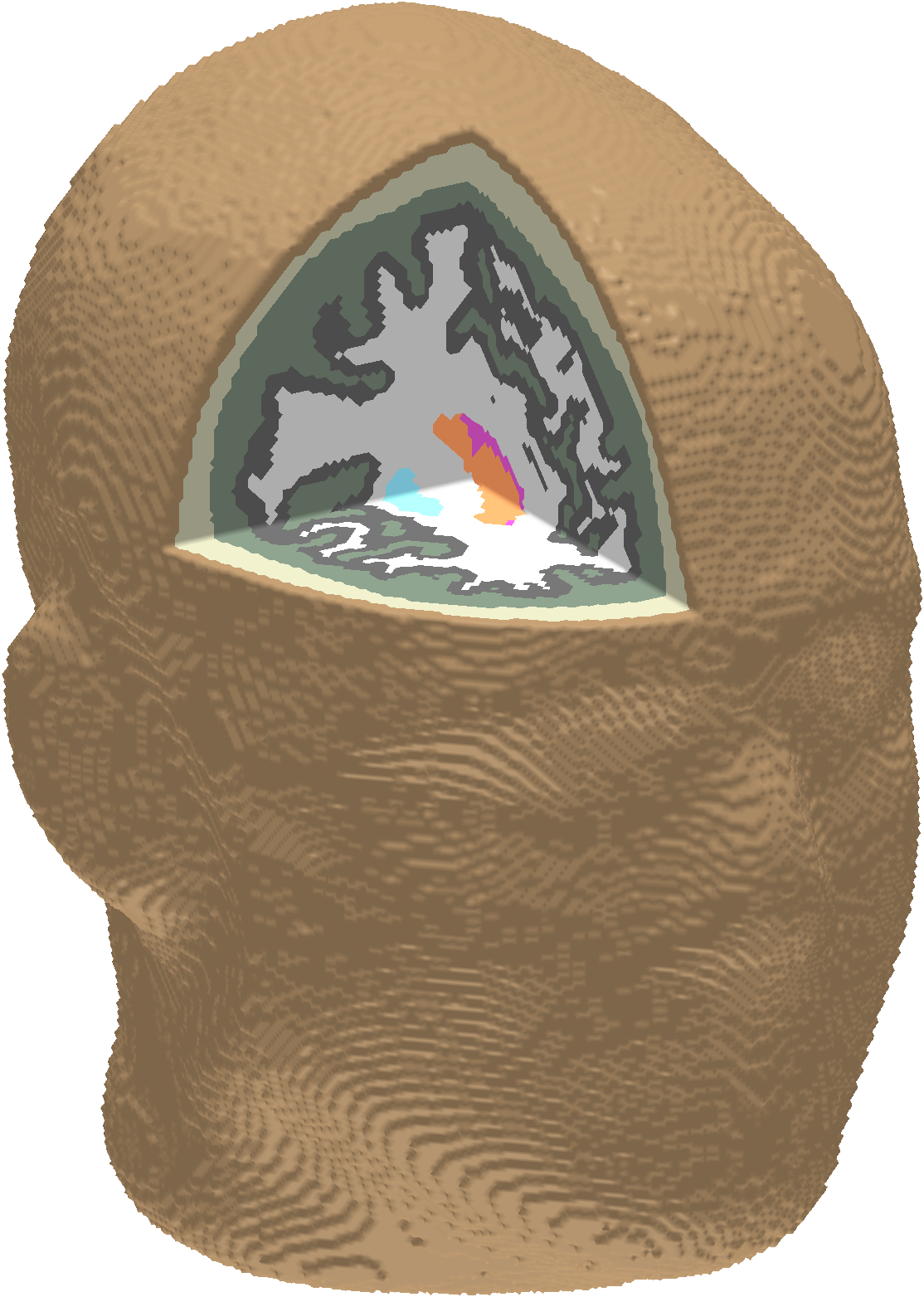}
    \end{subfigure}
    \caption{The volumetric head model used in our experiments.}
    \label{fig:head-model}
\end{figure}

To simulate how tTIS stimulation reaches deeper
regions of the brain, we first generated a volumetric
head model with a \qty{1}{\milli\meter} resolution
seen in Figure~\ref{fig:head-model}, based on the
surface segmentation generated with the FreeSurfer
software suite~\cite{FreeSurferWebSite, Desikan2006},
which itself was based on a publicly available MRI
dataset~\cite{Piastra_2020}.

\begin{figure}
    \centering
    \hfill
    \includegraphics[width=0.4\linewidth]{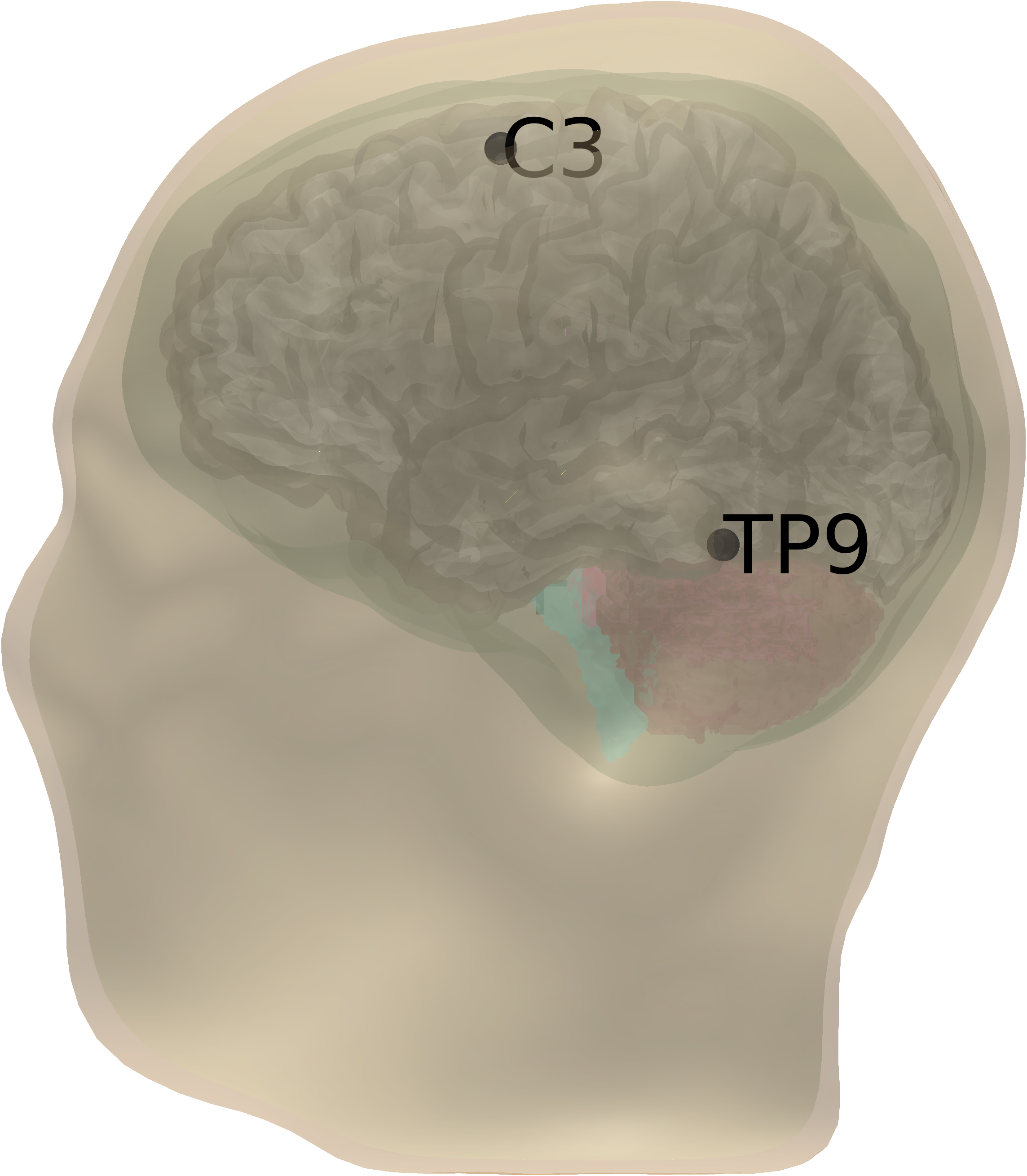}
    \hfill
    \includegraphics[width=0.4\linewidth]{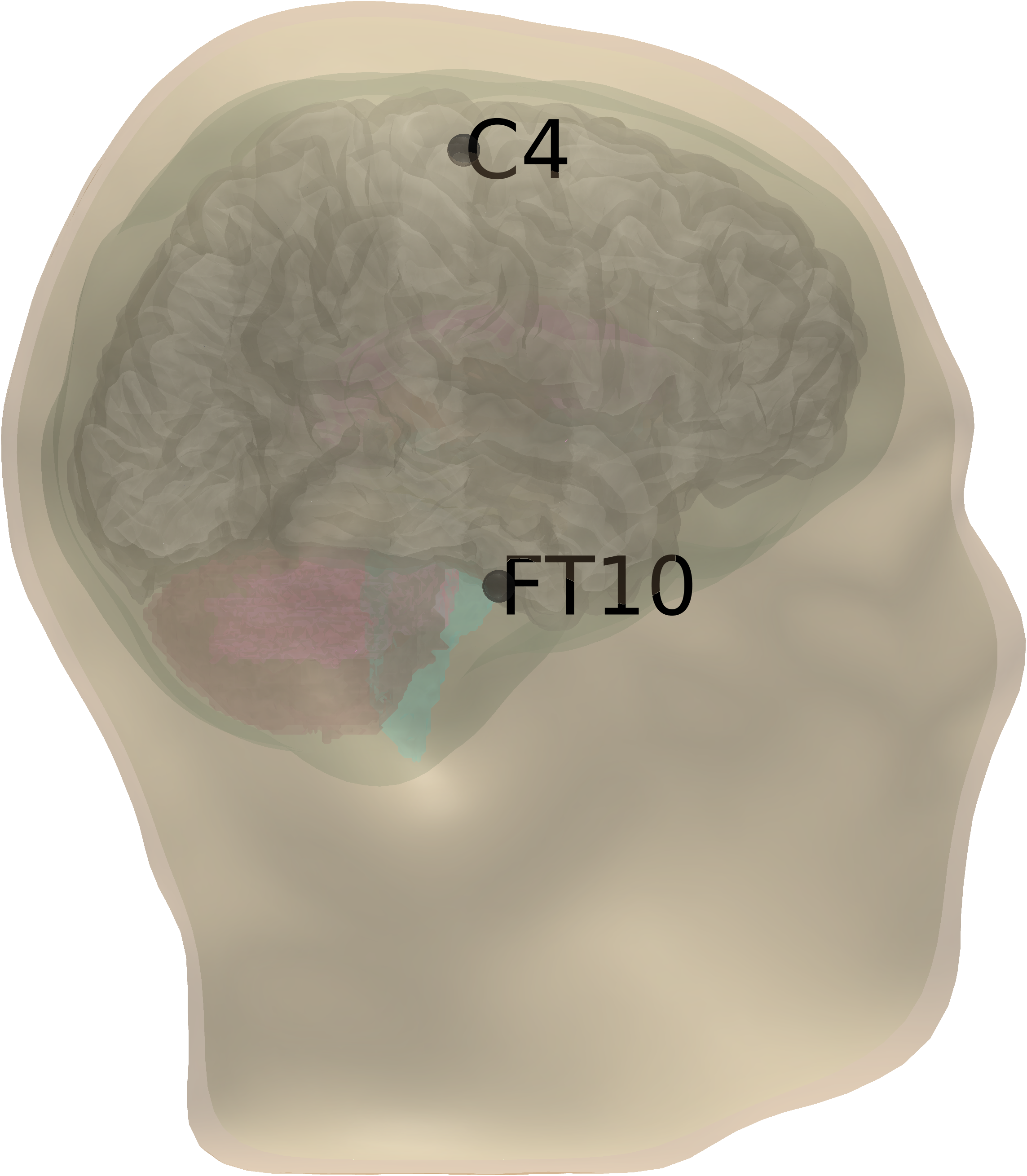}
    \hfill
    \caption{The 4-electrode or 2-voltage-source
    tTIS configuration used in this study. The electrodes
    are labeled according to the standard 10-by-10
    system~\cite{Nuwer-2018}, with electrode pairs
    TP9 + C3, and C4 + FT10 sharing voltage sources.
    The electrode contact surfaces had a diameter of
    $\sim\qty{1.0}{\centi\meter}$.}
    \label{fig:electrode-config}
\end{figure}

With the head model in place, we set up a simulated
experiment, where $2$ tTIS electrode pairs seen in
Figure~\ref{fig:electrode-config} were attached to
the surface of the skin compartment of said head
model. Then a resistance matrix $\iniR$ was computed
according to Section~\ref{app:res_matrix}, followed
by an initial lead field $\iniL$ computation according
to Section~\ref{app:Lead_Field}. The subindex here
refers to an initial matrix, before any contact
impedances $\contactResistance$ had been perturbed.
For a full set of electrode and tissue parameters,
see Tables~\ref{tab:initial-electrode-parameters}
and \ref{tab:conductivity_table}, respectively.
As a computation platform, we used Zeffiro
Interface~\cite{he2020zeffiro,pursiainen_zeffiro_december_2023},
a MATLAB-based open toolbox for forward and
inverse modelling targeting the brain. A Lenovo P620
workstation featuring an AMD Ryzen ThreadRipper PRO 5945WX
processor at \qty{4.1}{\giga\hertz}, 64 GB of DDR4 ECC RAM, and NVIDIA
RTX A5000 GPU with 24 GB of GDDR6, which was applied to
accelerate the Galerkin-FEM-based lead field generation
process~\cite{he2020zeffiro}

\begin{table}[!h]
    \centering
    \caption{Electrode parameter values of
    \eqref{eq:impedance-of-f} for the tTIS experiment. The
    $\contactResistance$ value varied for electrode TP9, and
    was kept at a constant $\qty{270}{\ohm}$ for the other
    electrodes.}
    \begin{tabular}{c|c|c}
        \toprule
        Electrode parameter & Values & Unit\\
        \midrule
        \(\contactResistance\) & \num{270}, \num{1270}, \num{5270} & \si\ohm \\
        \(\doubleLayerResistance\) & \num{10000} & {\si\ohm} \\
        \(\doubleLayerCapacitance\) & \num{1.0e-7} & {\si\farad} \\
        \(\frequency\) & \num{1000} & \si\hertz \\
        \(\beatF\) & \num{10} & \si\hertz \\
        \bottomrule
    \end{tabular}
    \label{tab:initial-electrode-parameters}
\end{table}

\begin{table}[h!]
    \caption{Values of electrical tissue conductivity \(\conductivity\) and relative permittivity \(\relativepermittivity\) at utilized stimulation frequencies \(\frequency\), applied for the skin, skull, cerebrospinal fluid (CSF), grey matter (GM), and white matter (WM), constituting the brain model of this study ~\cite{Gabriel_1996, tissue-property-database}. \label{tab:conductivity_table}}
    \centering
    \label{tab:segmentation}
    \resizebox{\linewidth}{!}{
        \begin{tabular}{clcc}
            \toprule
            \(\frequency\,(\si\hertz)\) & Compartment & \(\conductivity\)~(\si{\siemens\per\meter})~\cite{tissue-property-database} & \(\relativepermittivity\)~\cite{tissue-property-database} \\
            \midrule
            \multirow{5}{*}{\num{100}}
            & Skin & \num{2.0000e-4} &  \num{1135.9} \\
            & Skull & \num{2.0059e-2} & \num{5852.8} \\
            & CSF & \num{2.0} & \num{102} \\
            & Cerebral, cerebellar \& subcortical\textsuperscript{$\ast$} GM & \num{8.9018e-2} & \num{3906100}   \\
            & Cerebral, cerebellar \& subcortical\textsuperscript{$\dag$} WM &  \num{5.8093e-2} & \num{1667700} \\
            \midrule
            \multirow{5}{*}{\num{1000}}
            & Skin & \num{2.0006e-4} &  \num{1135.6} \\
            & Skull & \num{2.0157e-2} & \num{2702.4} \\
            & CSF & \num{2.0} & \num{102} \\
            & Cerebral, cerebellar \& subcortical\textsuperscript{$\ast$} GM & \num{9.8805e-2} & \num{164060}   \\
            & Cerebral, cerebellar \& subcortical\textsuperscript{$\dag$} WM &  \num{6.2574e-2} & \num{69811} \\
            \midrule
            \multirow{5}{*}{\num{10000}}
            & Skin & \num{2.0408e-4} &  \num{1133.6} \\
            & Skull & \num{2.0430e-2} & \num{521.64} \\
            & CSF & \num{2.0} & \num{102} \\
            & Cerebral, cerebellar \& subcortical\textsuperscript{$\ast$} GM & \num{1.1487e-1} & \num{22241}   \\
            & Cerebral, cerebellar \& subcortical\textsuperscript{$\dag$} WM &  \num{6.9481e-2} & \num{12468} \\
            \midrule
            \multirow{5}{*}{\num{100000}}
            & Skin & \num{4.5128e-4} &  \num{1119.2} \\
            & Skull & \num{2.0791e-2} & \num{227.64} \\
            & CSF & \num{2.0} & \num{102} \\
            & Cerebral, cerebellar \& subcortical\textsuperscript{$\ast$} GM & \num{1.3366e-1} & \num{3221.8}   \\
            & Cerebral, cerebellar \& subcortical\textsuperscript{$\dag$} WM &  \num{8.1845e-2} & \num{2107.6} \\
            \bottomrule
        \end{tabular}
    }
    \begin{footnotesize}
    \begin{itemize}
        \item [\(\ast\)] Brainstem, Ventral diencephalon, Amygdala, Thalamus, Caudate, Accumbens, Putamen, Hippocampus, Pallidum.
        \item [\(\dagger\)] Cingulate cortex.
    \end{itemize}
    \end{footnotesize}
\end{table}

With $\iniR$ in place, the contact impedance
$\contactResistance$ of electrode TP9 was then increased
to \qty{1270}{\ohm} and \qty{5270}{\ohm}, and the reference
and linearized resistance matrices $\refR$ and $\linR$
were computed with these new electrode impedances:
the first in accordance with \ref{app:res_matrix} and
the latter according to the linearization presented in
Section~\ref{sec:linearization}. Also, the lead fields
$\refL$ and $\linL$ were produced from these resistance
matrices via multiplication with the conductance
density $\conductanceDensity$, just as with $\iniL$.
With the lead field $\iniL$ in place, a unit current
pattern $\currentPattern=\transpose{(-1,1,-1,1)}$ that
fulfills Kirchhoffs first law was then applied through
the electrodes and mapped to a volume current density
distribution $\iniJv = \iniL\currentPattern$. The
strength of this discretized volume current density field
corresponds to the stimulated activity and was expected
to be the greatest near the electrodes themselves and
in the vicinity of the thalamus, where the distance from
the two sets of electrodes was roughly equal, due to the
reasons seen in Figure~\ref{fig:sinusoids}. Motivated
by Figure~\ref{fig:sinusoids}, for visualizing the
interference effect, we computed the magnitude of the
envelope of two interfering volume current distributions
$\volumeCurrentDensity_1$ and $\volumeCurrentDensity_2$,
emanating from two different electrode
pairs~\cite{GROSSMAN-2017, RAMPERSAD_computational_2019,
MIRZAKHALILI-2020}:
\begin{equation}
    \label{eq:interference-field}
    \interferenceField
    \paren{\volumeCurrentDensity_1,\volumeCurrentDensity_2}
    =
    \abs [\big] {
        \abs*{
            \volumeCurrentDensity_1
            +
            \volumeCurrentDensity_2
        }
        -
        \abs*{
            \volumeCurrentDensity_1
            -
            \volumeCurrentDensity_2
        }
    }
    \,.
\end{equation}
This field should be the strongest, where the
volume current fields $\volumeCurrentDensity_1$ and
$\volumeCurrentDensity_2$ are roughly equal in magnitude
but diminish near the electrodes. In addition to this, the
stabilized relative and absolute differences of the form
\begin{equation}\label{eq:rel-diff}
    \relDelta g
    =
    \frac{
        \abs{
            g_1
            -
            g_2
        }
    }
    {
        \max\paren{
            \abs{g_2},
            \abs{
                \errorDiscrepancy
                \max\abs{g_{2}}
            }
        }
    }
\end{equation}
\begin{equation}\label{eq:abs-diff}
    \Delta g
    \paren{
        g_1,
        g_2
    }
    =
    \abs [\big] {
        \norm{g_1}
        -
        \norm{g_2}
    }
    \,,
\end{equation}
between two volume current fields  ($g =
\volumeCurrentDensity$) and interference current
fields ($g = \interferenceField$) were computed to
see where the two fields differ the most. The inverse
of the stabilization parameter $\errorDiscrepancy$
in \eqref{eq:rel-diff} represents a dynamic range
(\unit\decibel) that determines the maximum contrast
between details in the examination.

The dynamic range $\errorDiscrepancy$ can be associated
with a peak signal-to-noise ratio (PSNR) following
from the model-related uncertainty factors. In the
investigation of the volume and interference current
fields, we choose $\delta = \qty{-35}{\decibel}$
and $\delta = \qty{-18}{\decibel}$, respectively.
The motivation for this choice is explained in the
following section. Thus, this investigation can be
considered valid, where the PSNR of the volume and
interference current is $\geq\qty{35}{\decibel}$ and
$\geq\qty{18}{\decibel}$ (1.78 and 12.59 \% maximum
relative error w.r.t.\ maximum lead field amplitude),
respectively. The measure \eqref{eq:rel-diff} was also
applied to evaluate the difference between solutions
obtained with  reference and linearized lead fields
(Figure~\ref{fig:L-column-vector-plot}) to get a sense
of where the potential distribution changes the most when
the contact resistance $\contactResistance$ is modified at
the electrode TP9. The scale in the related visualizations
is logarithmic, where the visualized field $\someField$ has been
transformed via the equation
\begin{equation}\label{eq:visual-log-transform}
    \begin{aligned}
        \dBT\args{g,\errorDiscrepancy}
        \coloneqq
        &-
        20
        \logten
        \args{
            \max
            \args{
                \min\someField,
                \max\abs\someField\div\errorDiscrepancy
            }
        }
        \\
        &+
        20
        \logten
        \args{
            \max
            \args{
                \someField,
                \max\abs \someField\div\errorDiscrepancy
            }
        }
        \,.
    \end{aligned}
\end{equation}

\begin{figure}[h!]
\centering
\begin{tikzpicture}
  \def\ul{0.52}
  \def\R{3.1}
  \def\ang{26}
  \coordinate (O) at (0,0);
  \coordinate (R) at (1.73,2.57);
  \coordinate (X) at (2.68,1.55);
  \coordinate (Y) at (-5,1.55);
  \draw[vector,mypurple,dashed] (R) -- (X) node[midway, above  right] {$ {\deltaLLinRef}$};
  \draw[vector,dashed] (O) -- (R) node[midway,  left=2] {$\refDiff {\leadFieldMat}$};
  \draw[vector,projcol,dashed]
    (O) -- (X) node[midway, below right] {$\linDiff {\leadFieldMat}$} ;
     \draw[vector]
    (Y) -- (R) node[midway, above left] {$\refVal {\leadFieldMat}$} ;
      \draw[vector,projcol]
    (Y) -- (X) node[midway, above right] {$\linVal {\leadFieldMat}$} ;
  \draw[vector,initcol]
    (Y) -- (O) node[midway, above  right] {$ \iniVal {\leadFieldMat}$} ;
  \draw pic[->,thick,"$\theta$",draw=black,angle radius=26,angle eccentricity=1.3]
    {angle = X--O--R};
\end{tikzpicture}
\caption{Schematic of the column-wise lead field differences analyzed in this study.}
\label{fig:L-column-vector-plot}
\end{figure}

\section{Results}%
\label{sec:results}

\begin{figure}[!h]
    \centering
    \begin{minipage}{0.85\linewidth}
        \def\subFigWidth{0.31\linewidth}
        \begin{subfigure}[b]{\subFigWidth}
            \includegraphics[width=\linewidth]{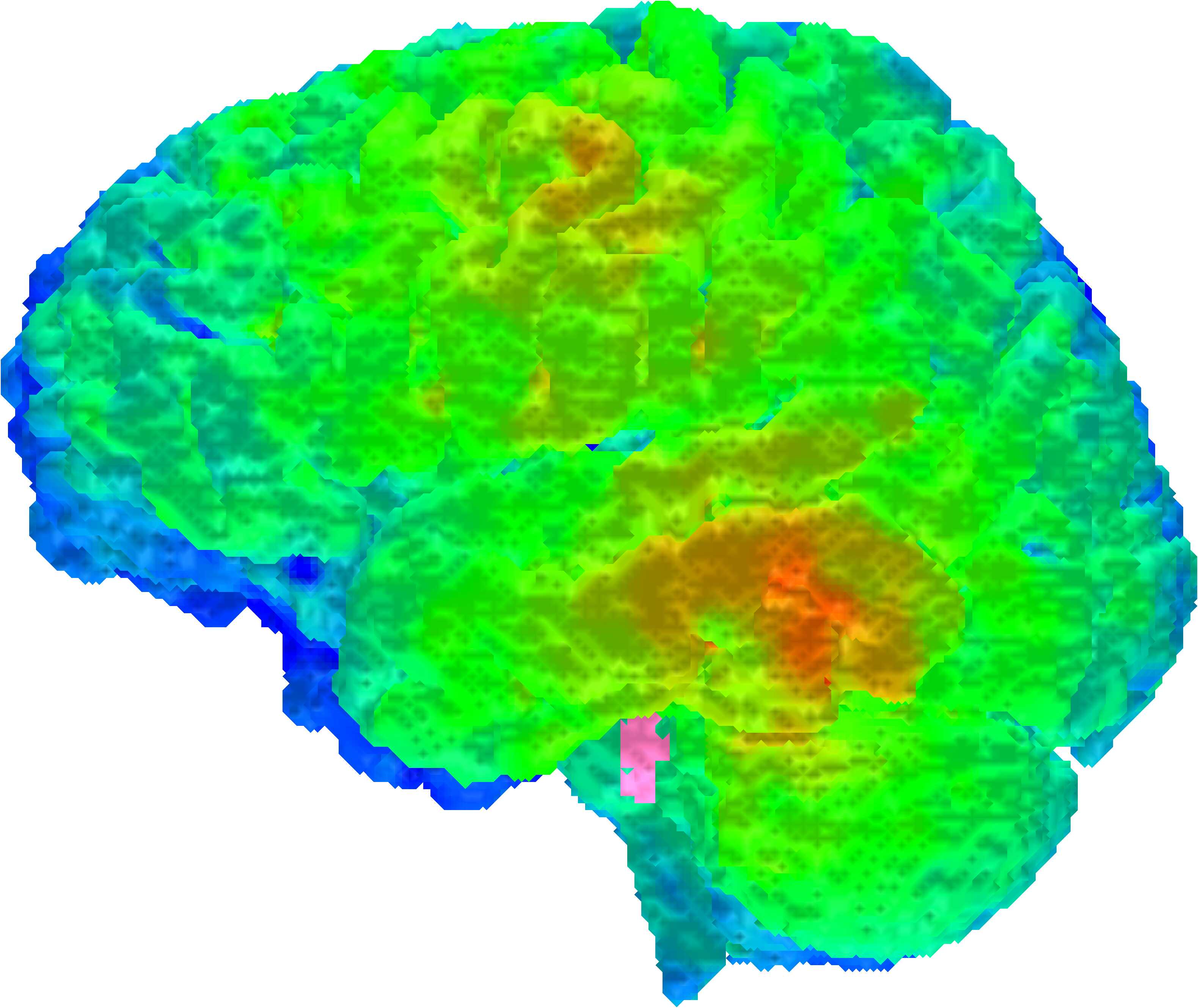}
            \caption{$\azimuth=\qty{-90}{\degree}$, $\elevation=\qty{0}{\degree}$.}
        \end{subfigure}
        \hfill
        \begin{subfigure}[b]{\subFigWidth}
        \includegraphics[width=\linewidth]{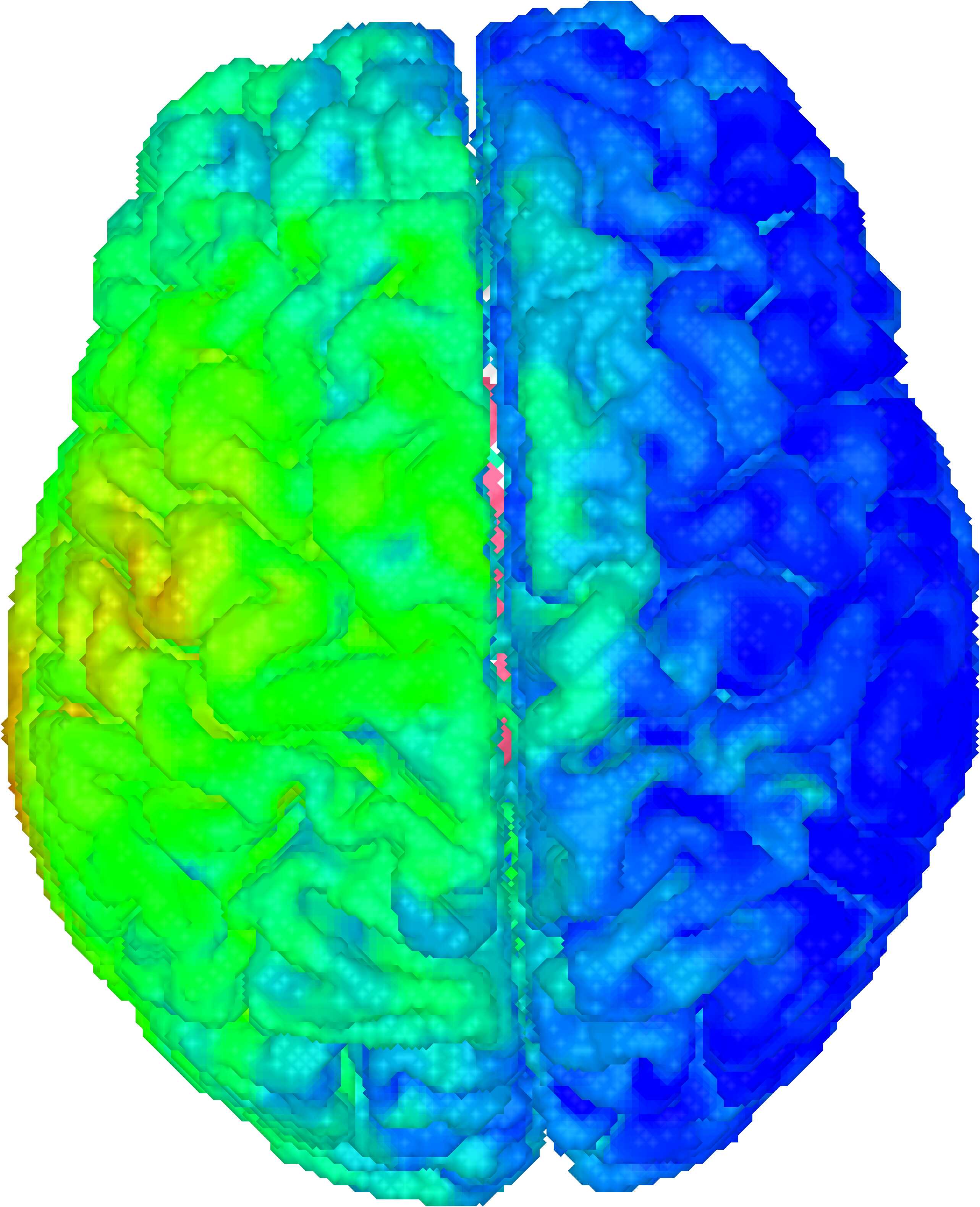}
        \caption{$\azimuth=\qty{0}{\degree}$, $\elevation=\qty{90}{\degree}$.}
        \end{subfigure}
        \hfill
        \begin{subfigure}[b]{\subFigWidth}
        \includegraphics[width=\linewidth]{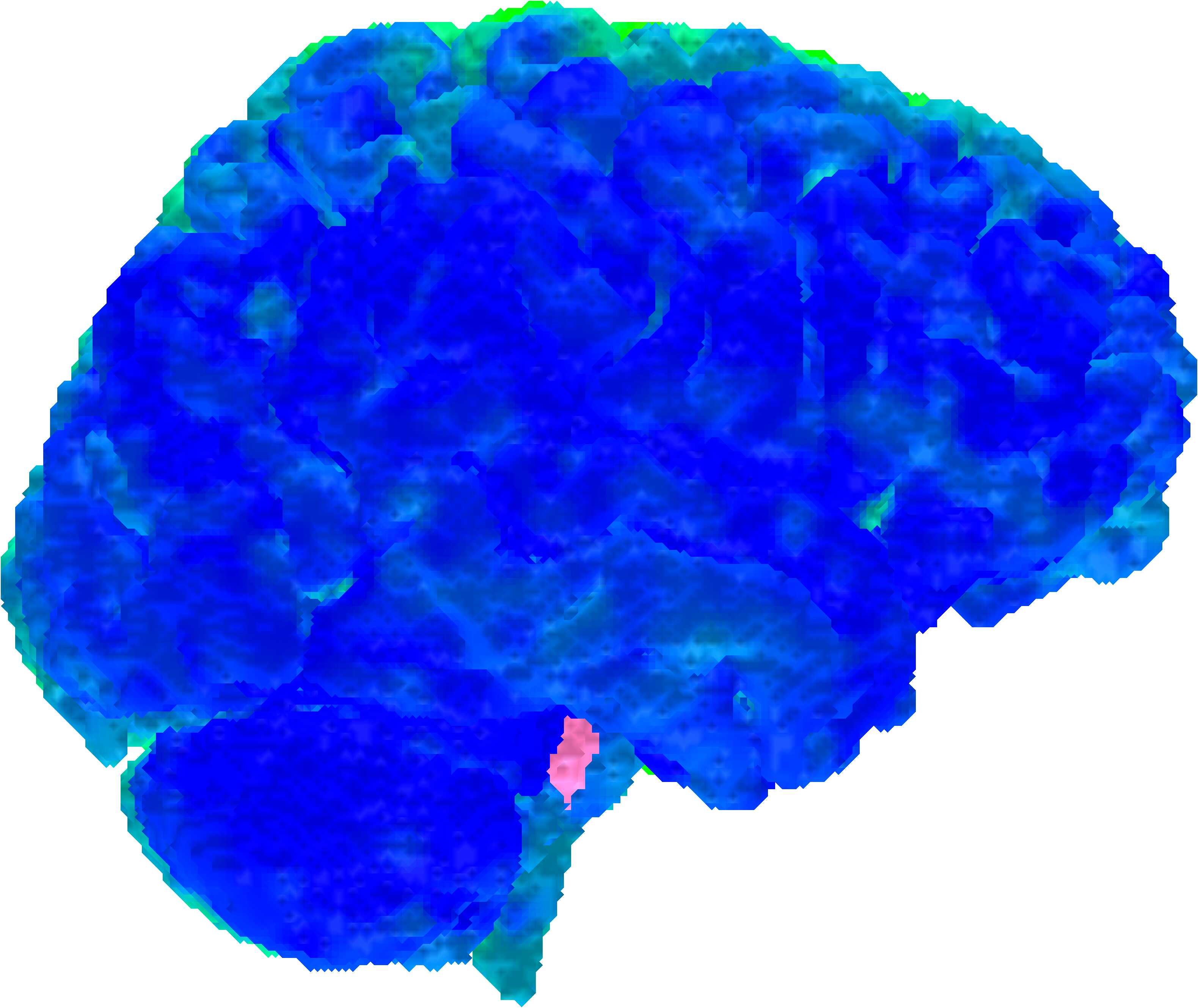}
            \caption{$\azimuth=\qty{90}{\degree}$, $\elevation=\qty{0}{\degree}$.}
        \end{subfigure}

        \def\subFigWidth{0.45\linewidth}

        \begin{subfigure}[b]{\subFigWidth}
            \includegraphics[width=\linewidth]{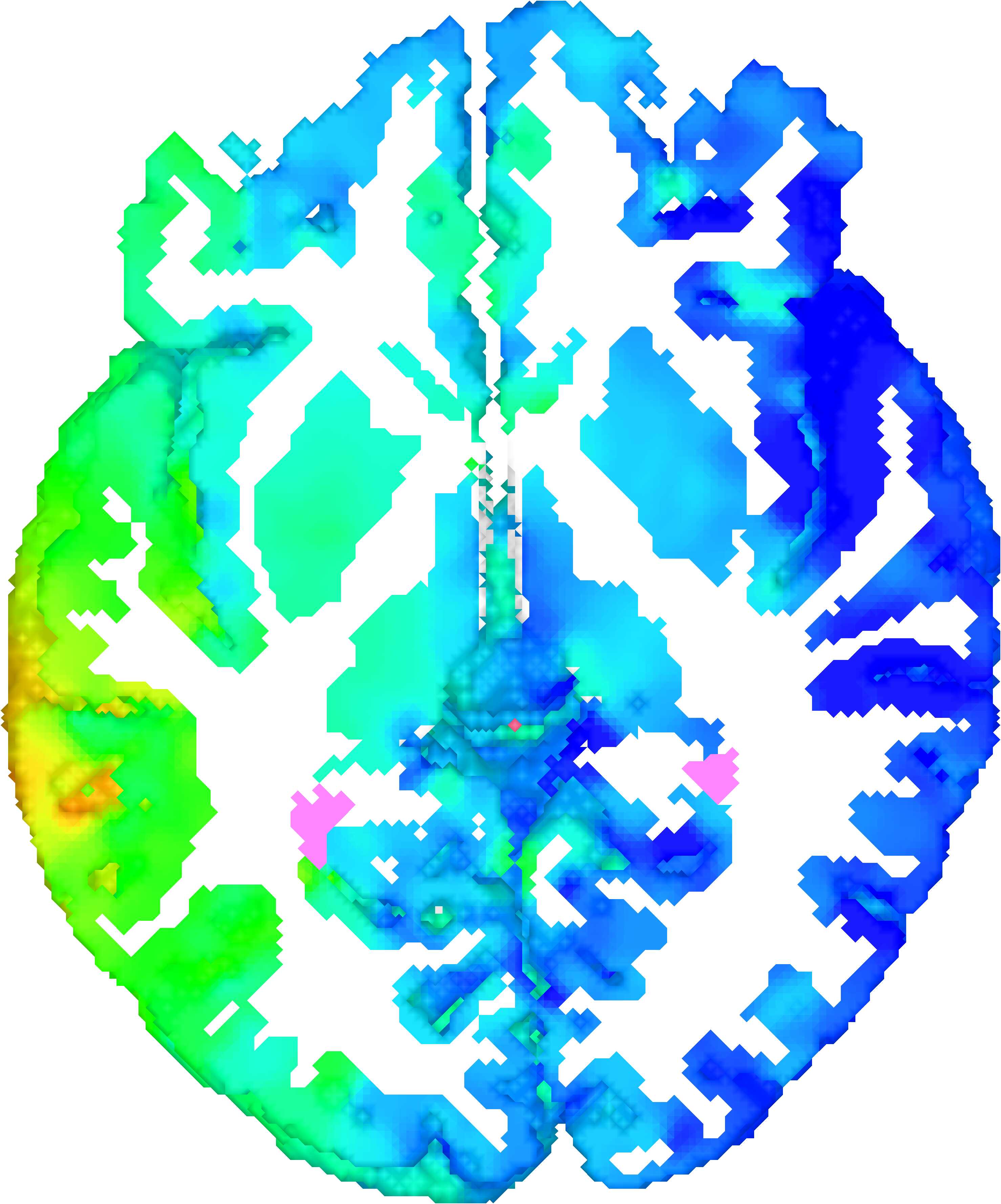}
            \caption{$\azimuth=\qty{0}{\degree}$, $\elevation=\qty{90}{\degree}$, $z=\qty{30}{\milli\meter}$.}
        \end{subfigure}
       \hfill
        \begin{subfigure}[b]{\subFigWidth}
            \includegraphics[width=\linewidth]{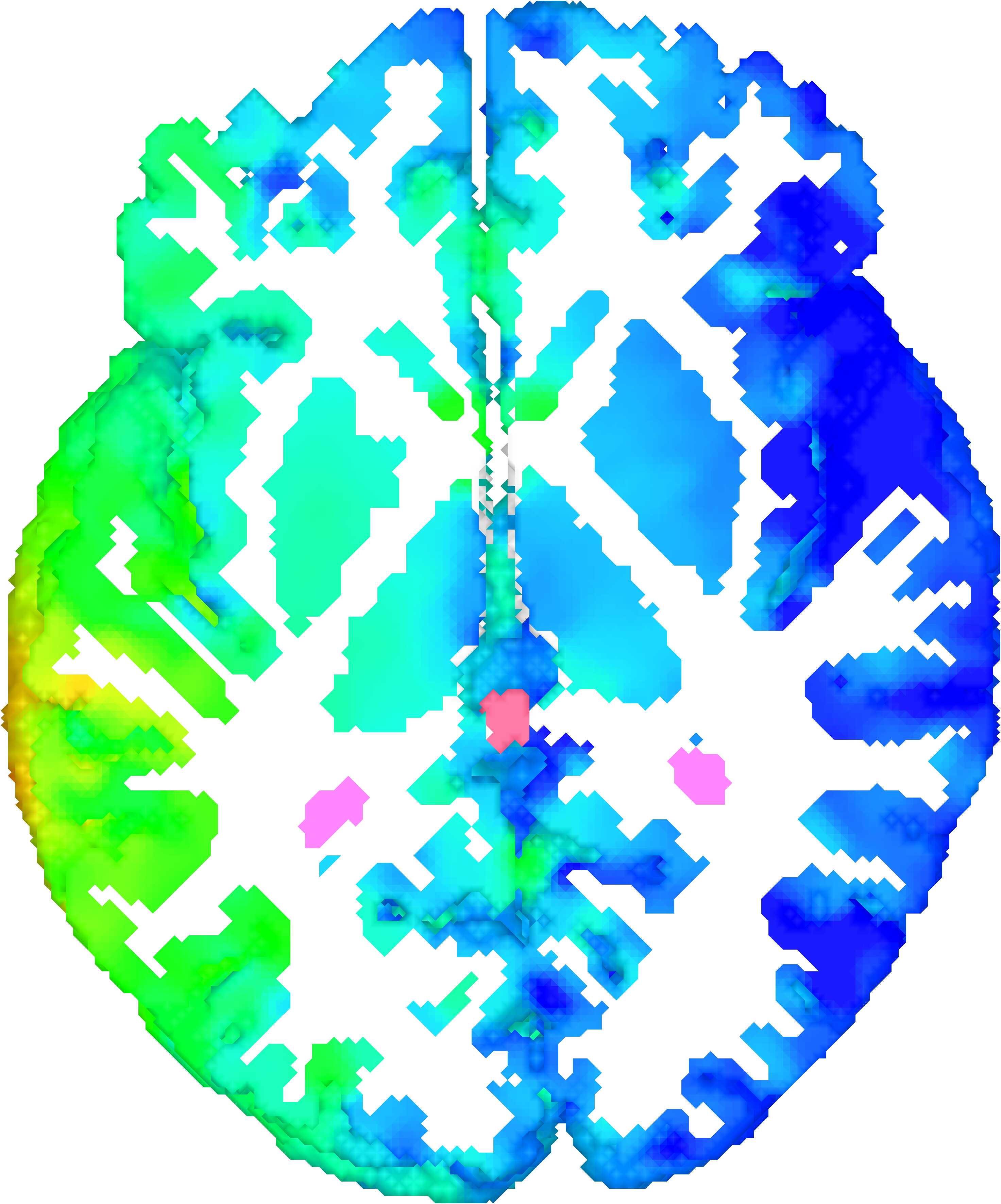}
            \caption{$\azimuth=\qty{0}{\degree}$, $\elevation=\qty{90}{\degree}$, $z=\qty{33}{\milli\meter}$.}
        \end{subfigure}

        \begin{subfigure}[b]{\subFigWidth}
            \includegraphics[width=\linewidth]{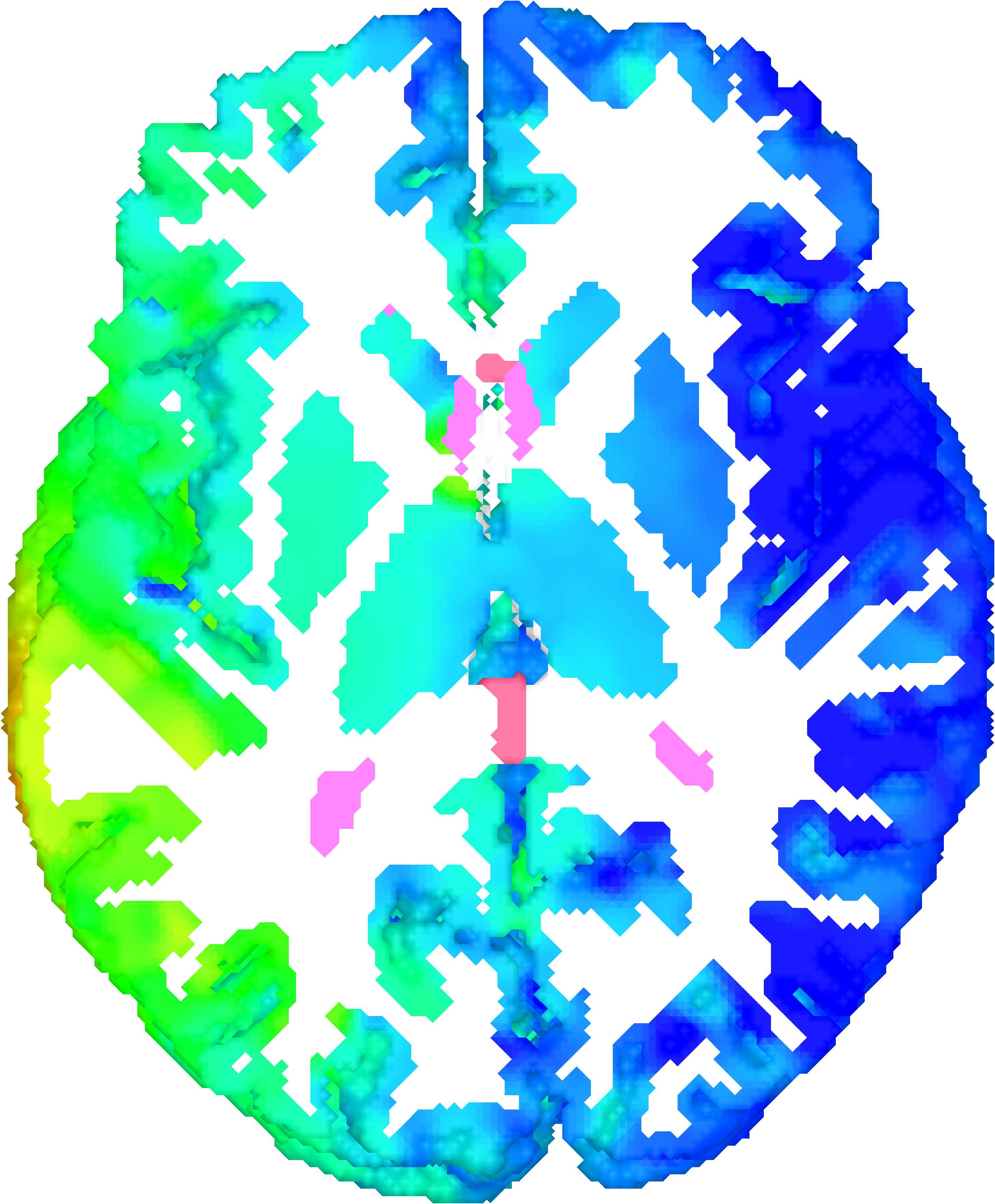}
            \caption{$\azimuth=\qty{0}{\degree}$, $\elevation=\qty{90}{\degree}$, $z=\qty{37}{\milli\meter}$.}
        \end{subfigure}
        \hfill
        \begin{subfigure}[b]{\subFigWidth}
            \includegraphics[width=\linewidth]{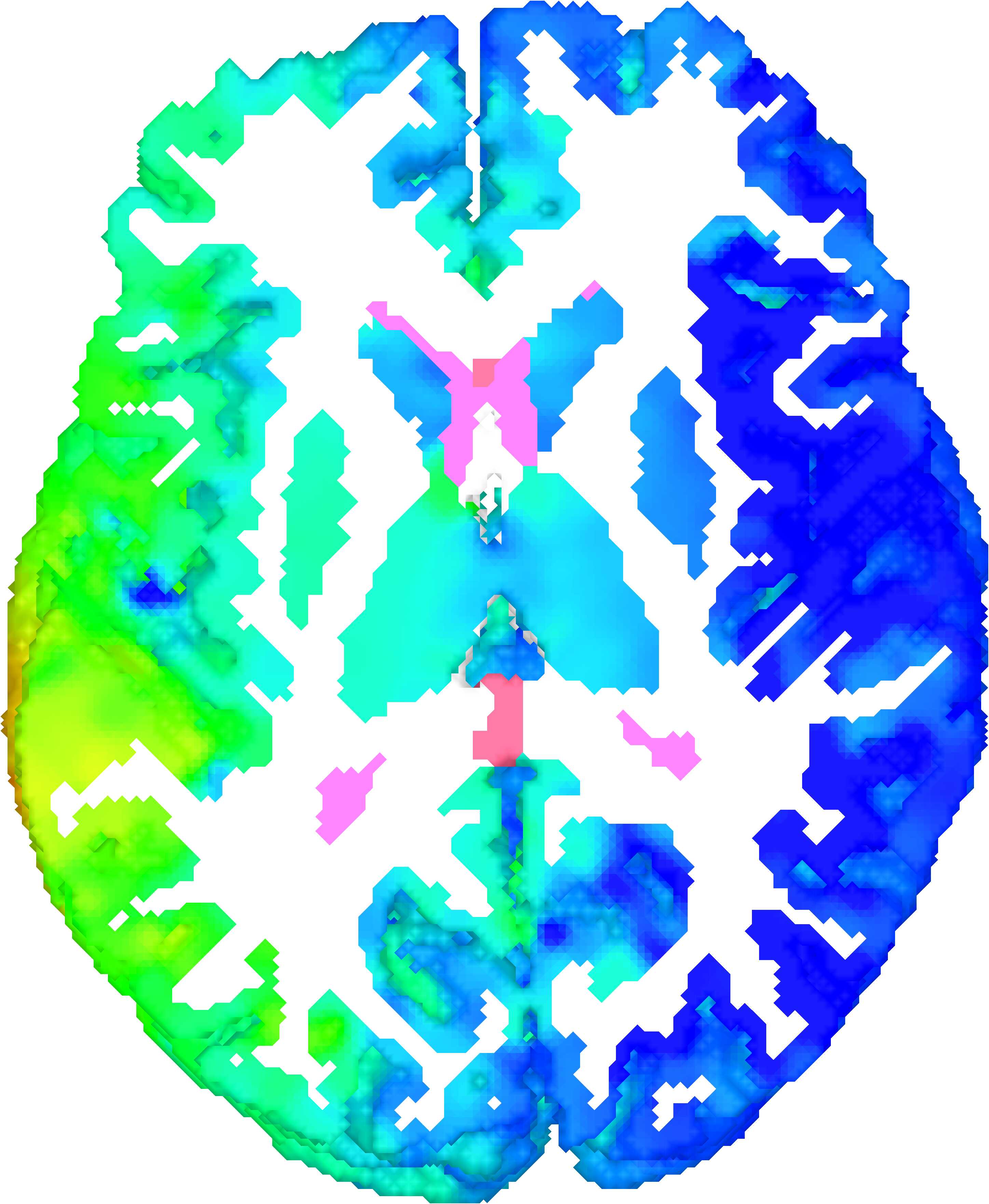}
            \caption{$\azimuth=\qty{0}{\degree}$, $\elevation=\qty{90}{\degree}$, $z=\qty{40}{\milli\meter}$.}
        \end{subfigure}
    \end{minipage}
    \hfill
    \begin{minipage}{0.09\linewidth}
    \includegraphics[width=\linewidth]{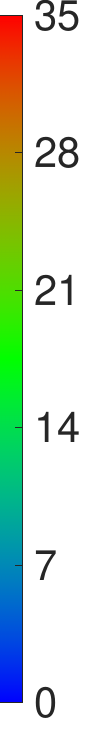}
    \end{minipage}
    \caption{Smoothed current density distribution $\dBT\args{\volumeCurrentDensity/\max\volumeCurrentDensity,\qty{35}{\decibel}}$ of equation \eqref{eq:visual-log-transform}, produced by the left pair of electrodes at different azimuths $\azimuth$ and elevations $\elevation$, with base stimulation frequency $\frequency=\qty{1000}{\hertz}$. The current pattern was $\currentPattern=[-1.0,1.0]\,\si{\milli\ampere}$ and the maximal value of $\volumeCurrentDensity$ roughly equal to \qty{2.04e-4}{\milli\ampere\per\milli\meter\squared}.}
    \label{fig:left-current-field}
\end{figure}

\begin{figure}[!h]
    \centering
    \begin{minipage}{0.85\linewidth}
        \def\subFigWidth{0.31\linewidth}
        \begin{subfigure}[b]{\subFigWidth}
            \includegraphics[width=\linewidth]{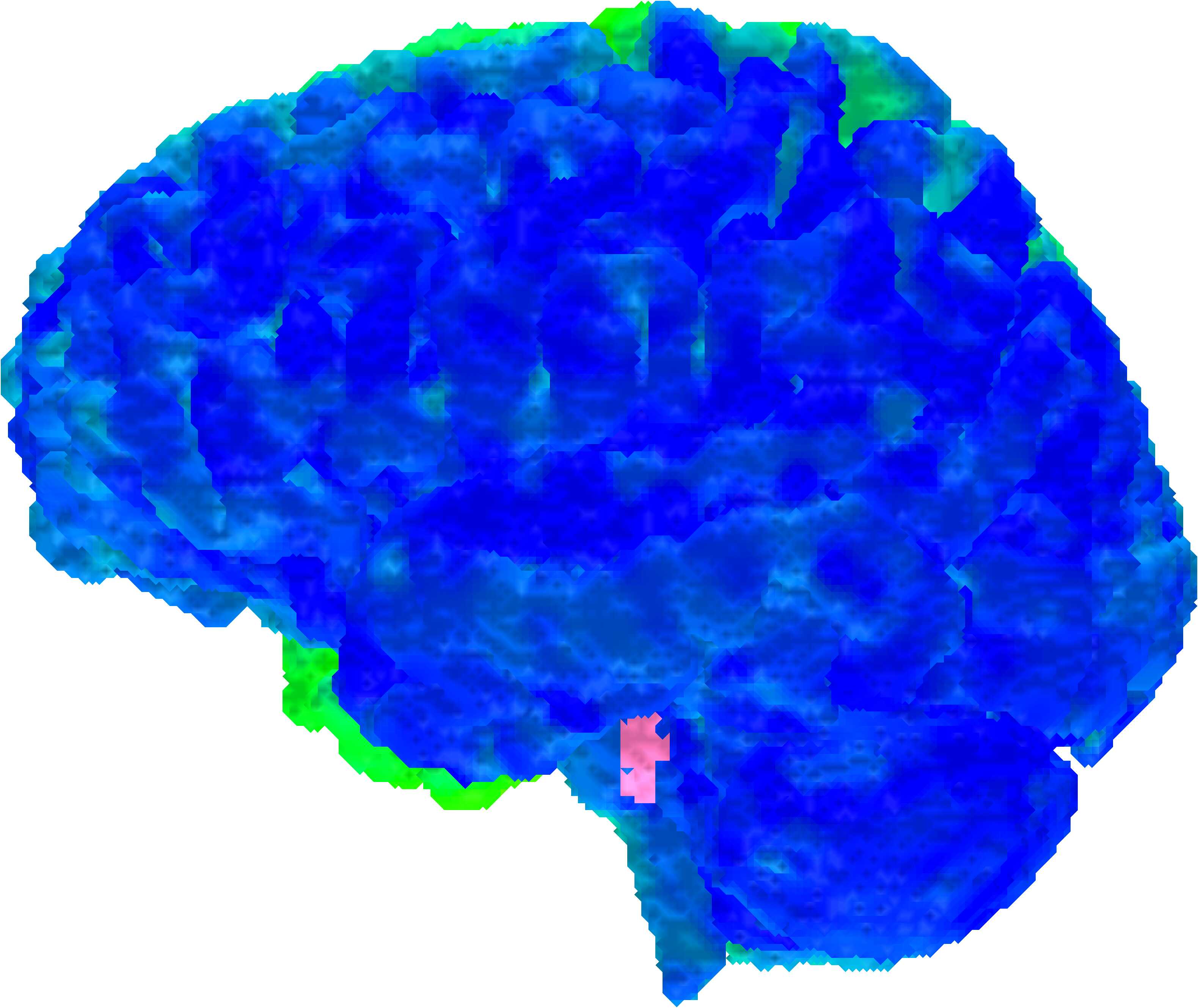}
            \caption{$\azimuth=\qty{-90}{\degree}$, $\elevation=\qty{0}{\degree}$.}
        \end{subfigure}
        \hfill
        \begin{subfigure}[b]{\subFigWidth}
        \includegraphics[width=\linewidth]{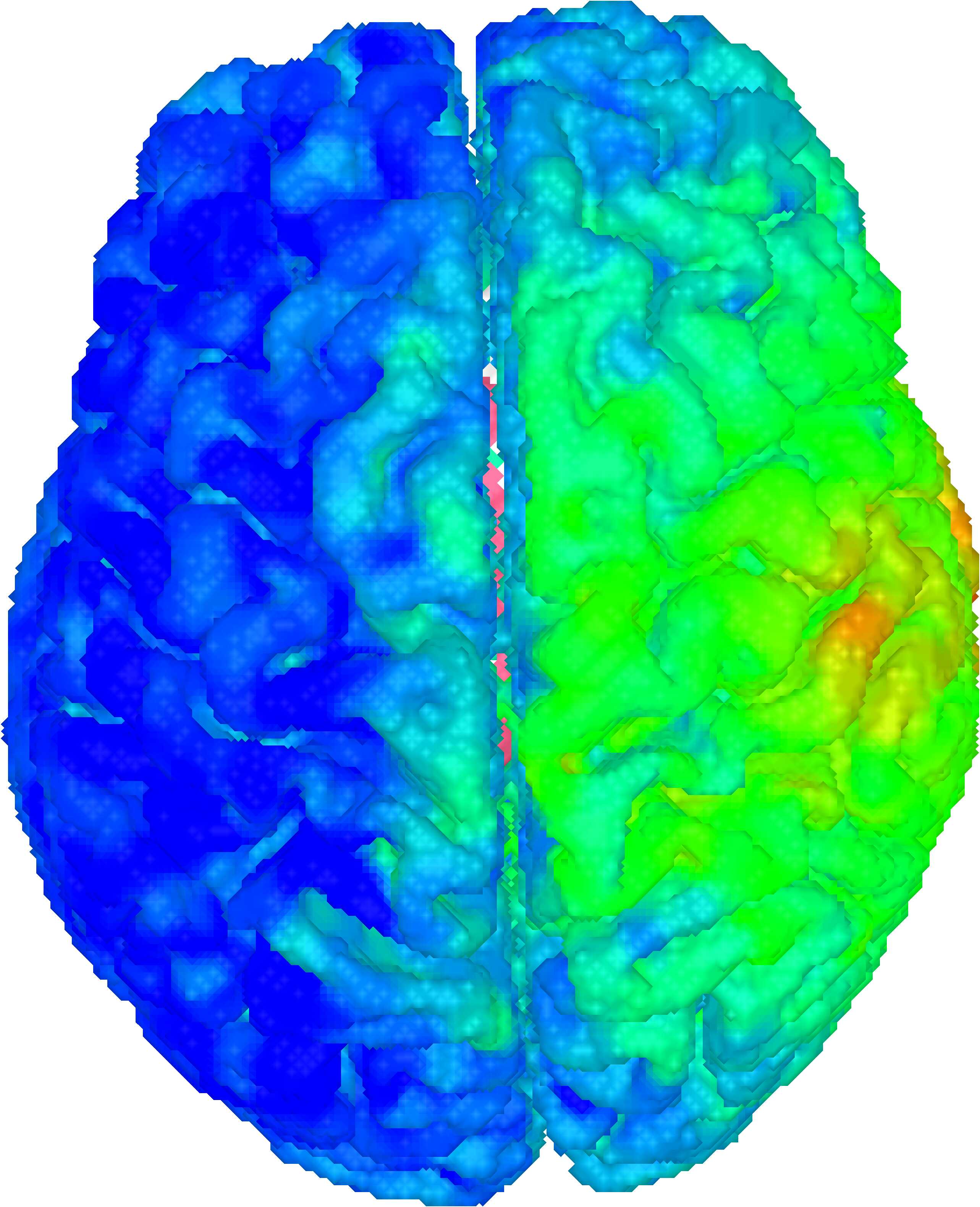}
        \caption{$\azimuth=\qty{0}{\degree}$, $\elevation=\qty{90}{\degree}$.}
        \end{subfigure}
        \hfill
        \begin{subfigure}[b]{\subFigWidth}
        \includegraphics[width=\linewidth]{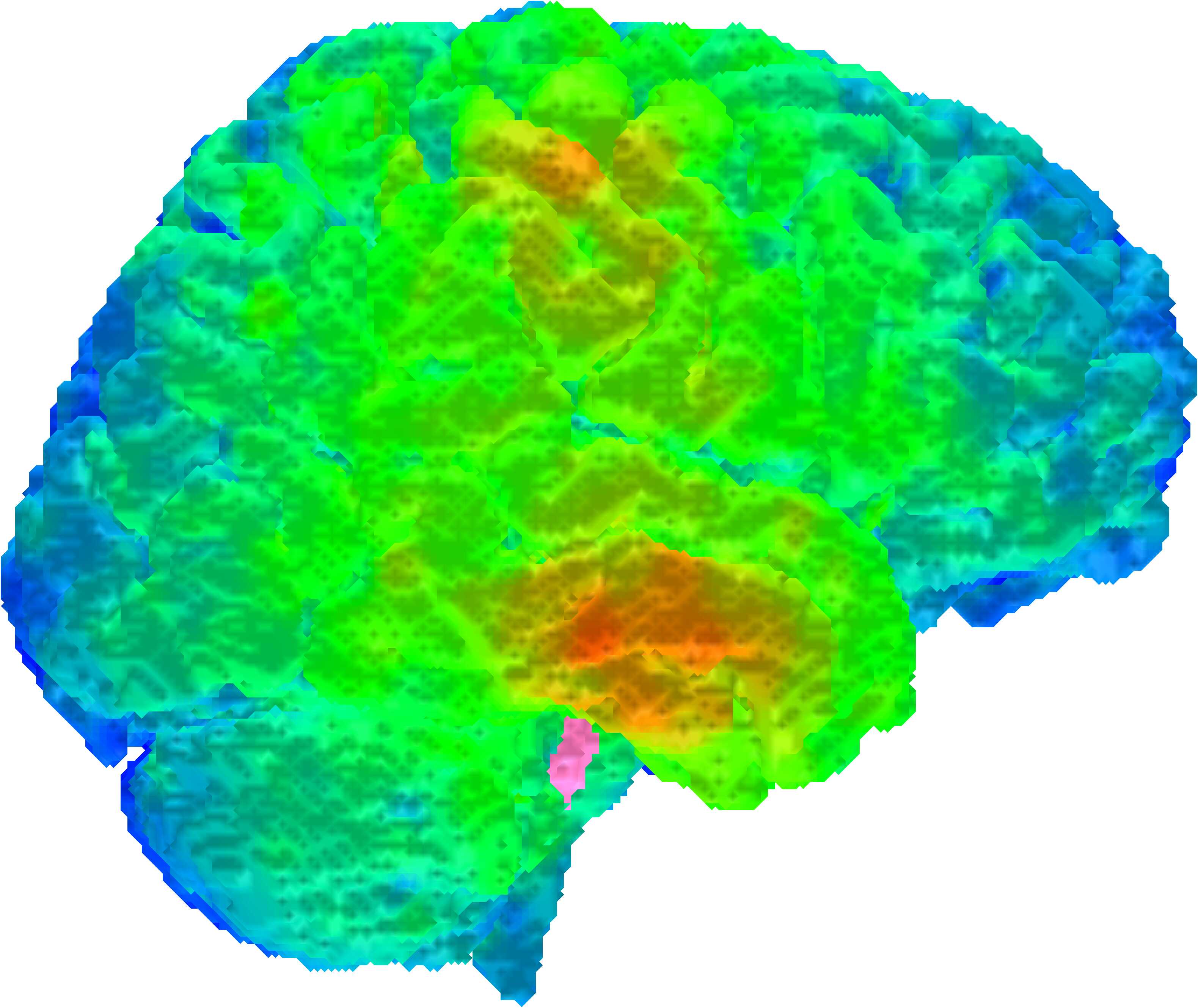}
            \caption{$\azimuth=\qty{90}{\degree}$, $\elevation=\qty{0}{\degree}$.}
        \end{subfigure}

        \def\subFigWidth{0.45\linewidth}

        \begin{subfigure}[b]{\subFigWidth}
            \includegraphics[width=\linewidth]{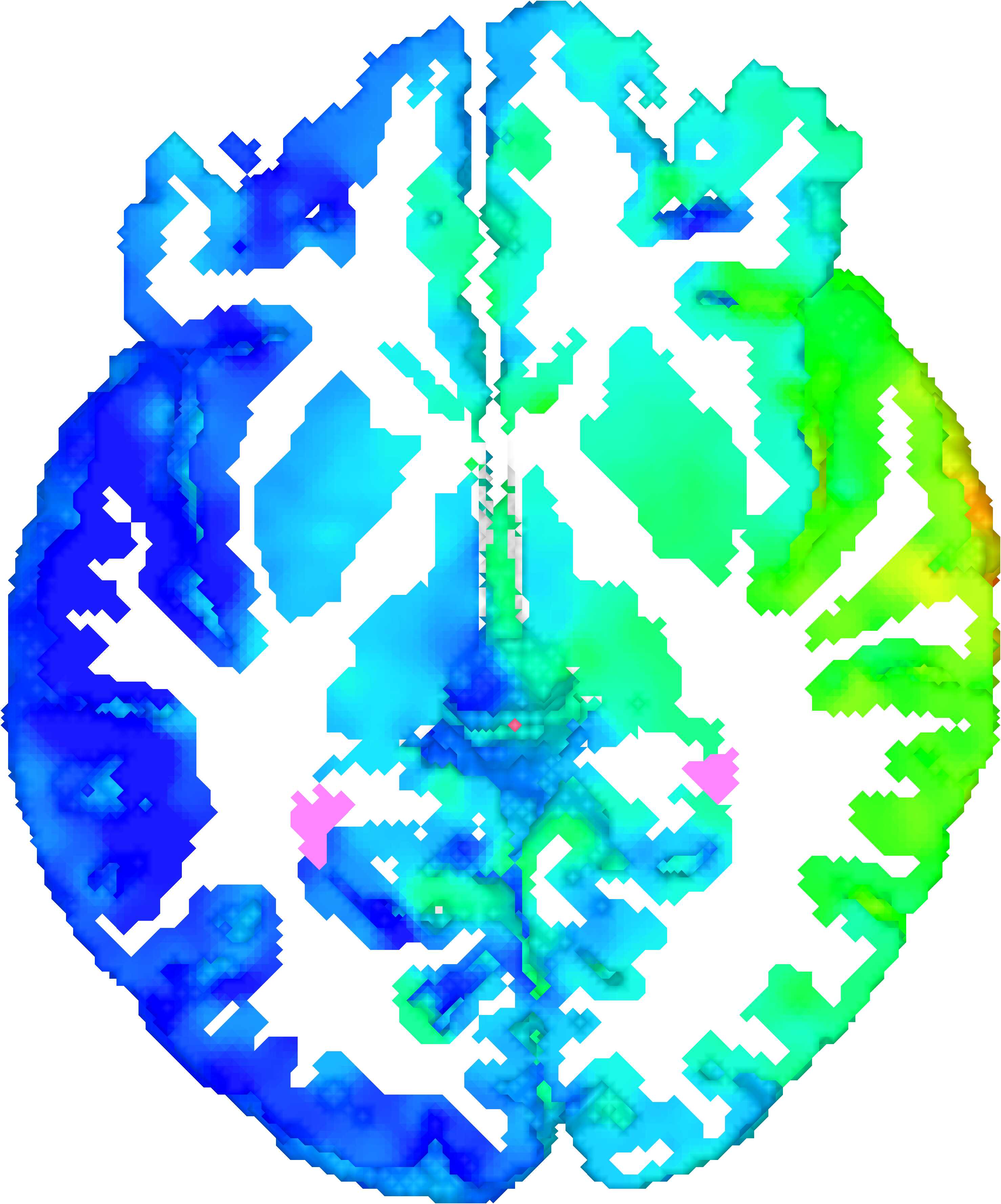}
            \caption{$\azimuth=\qty{0}{\degree}$, $\elevation=\qty{90}{\degree}$, $z=\qty{30}{\milli\meter}$.}
        \end{subfigure}
       \hfill
        \begin{subfigure}[b]{\subFigWidth}
            \includegraphics[width=\linewidth]{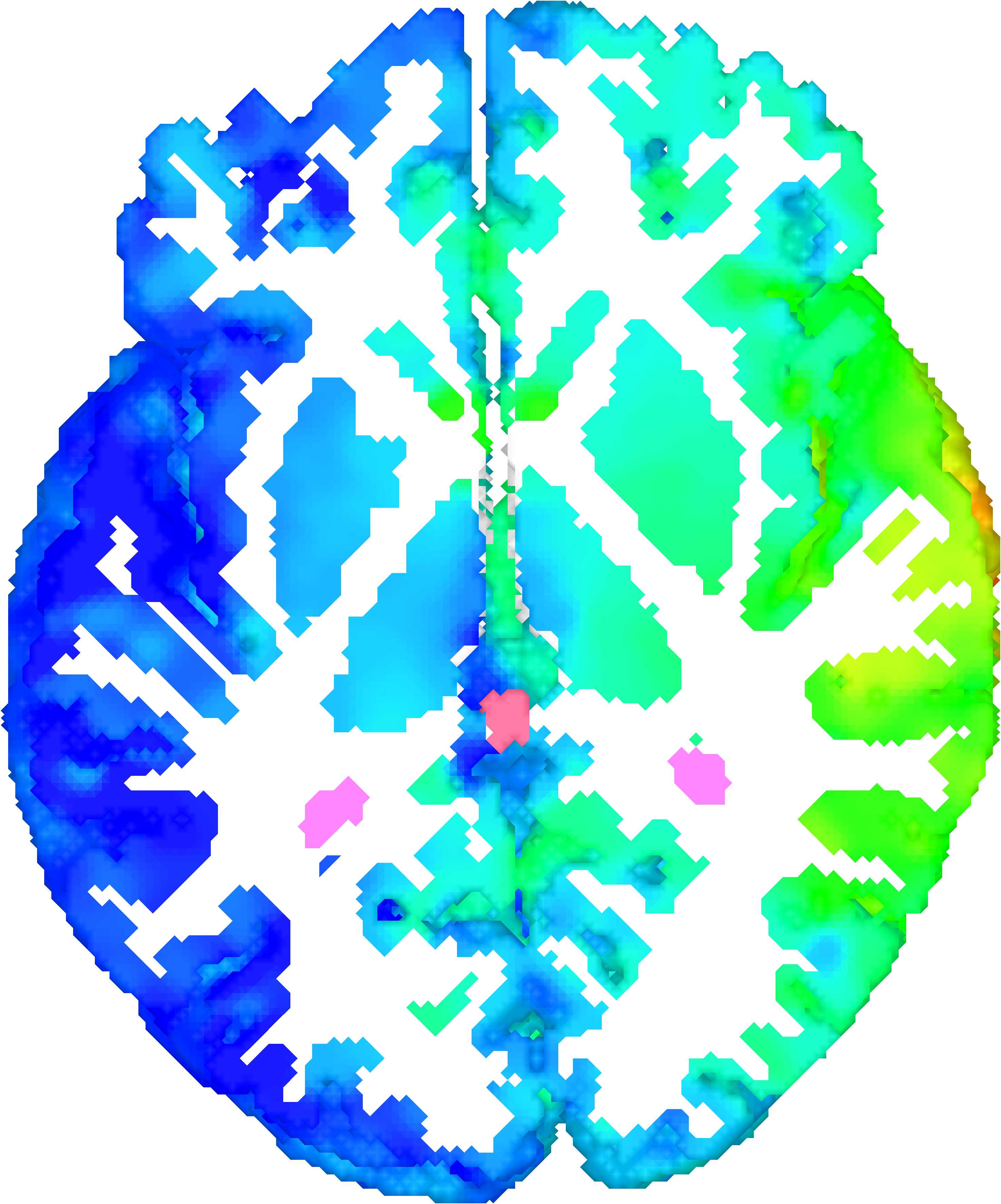}
            \caption{$\azimuth=\qty{0}{\degree}$, $\elevation=\qty{90}{\degree}$, $z=\qty{33}{\milli\meter}$.}
        \end{subfigure}

        \begin{subfigure}[b]{\subFigWidth}
            \includegraphics[width=\linewidth]{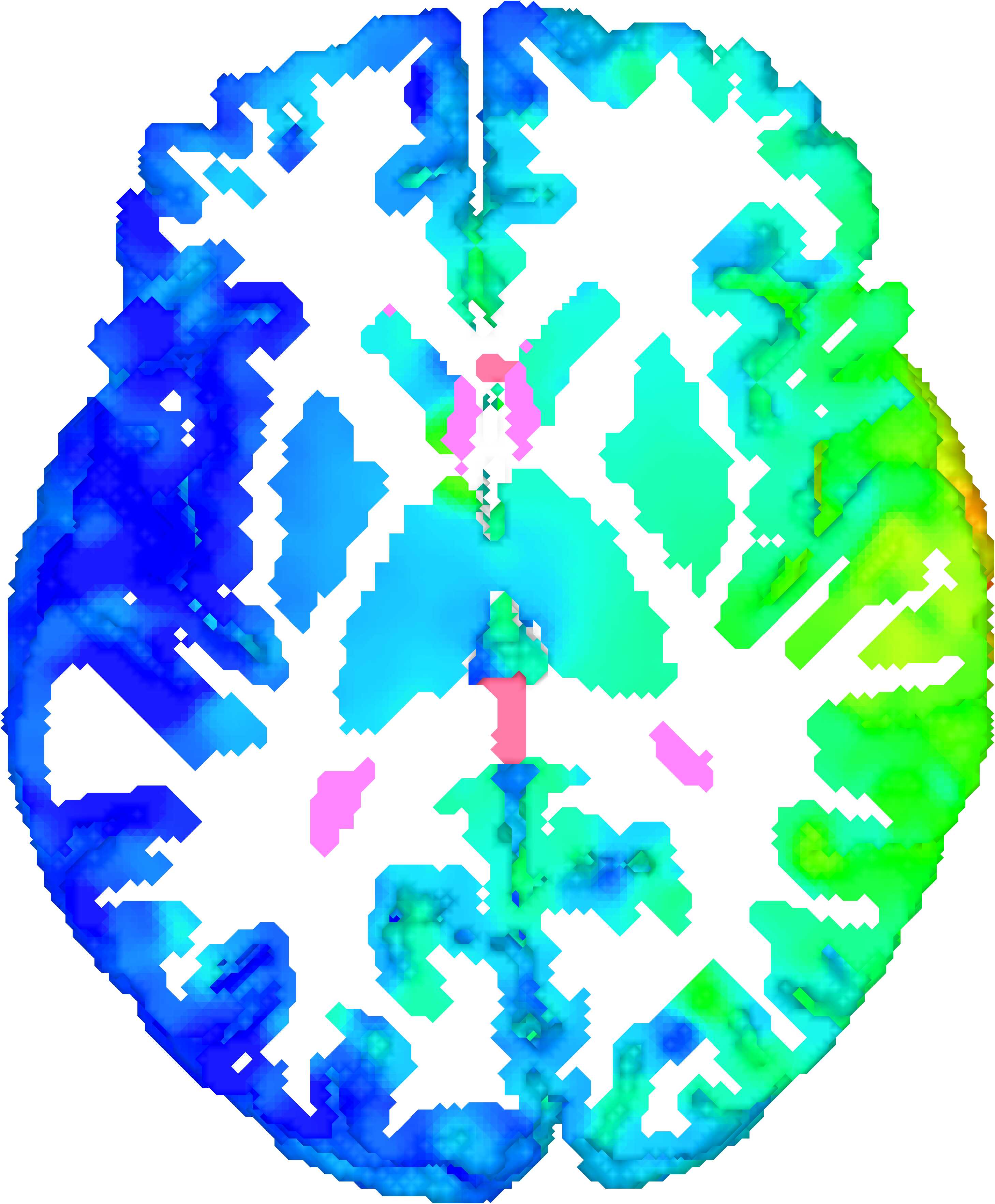}
            \caption{$\azimuth=\qty{0}{\degree}$, $\elevation=\qty{90}{\degree}$, $z=\qty{37}{\milli\meter}$.}
        \end{subfigure}
        \hfill
        \begin{subfigure}[b]{\subFigWidth}
            \includegraphics[width=\linewidth]{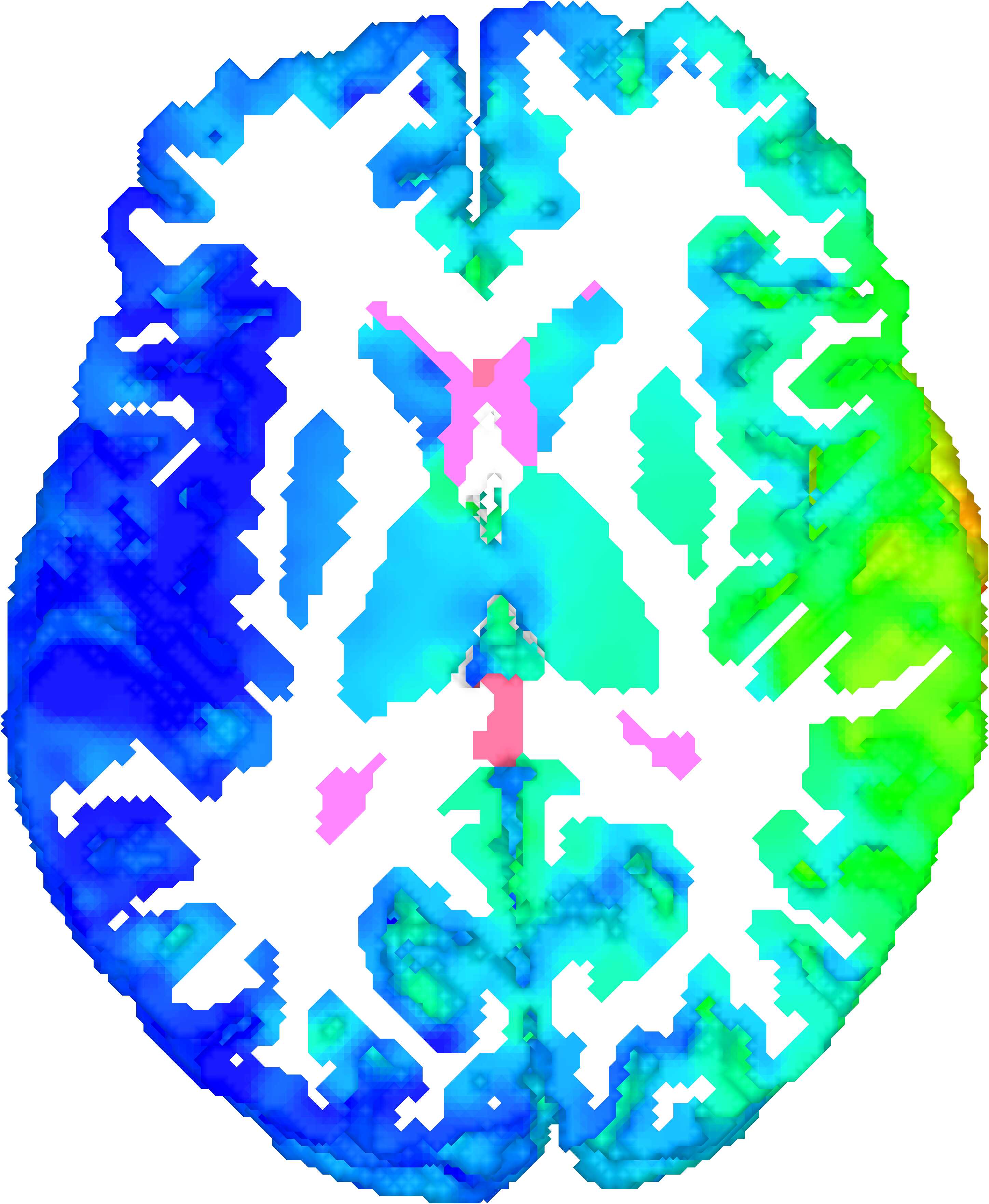}
            \caption{$\azimuth=\qty{0}{\degree}$, $\elevation=\qty{90}{\degree}$, $z=\qty{40}{\milli\meter}$.}
        \end{subfigure}
    \end{minipage}
    \hfill
    \begin{minipage}{0.09\linewidth}
    \includegraphics[width=\linewidth]{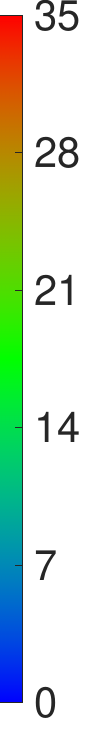}
    \end{minipage}
    \caption{Smoothed current density distribution $\dBT\args{\volumeCurrentDensity/\max\volumeCurrentDensity,\qty{35}{\decibel}}$ of equation \eqref{eq:visual-log-transform} produced by the right pair of electrodes C4 and FT10 at different azimuths $\azimuth$ and elevations $\elevation$, with base stimulation frequency $\frequency=\qty{1000}{\hertz}$. The current pattern was $\currentPattern=[-1.0,1.0]$ and the maximal value of $\volumeCurrentDensity$ roughly equal to \qty{2.06e-4}{\milli\ampere\per\milli\meter\squared}.}
    \label{fig:right-current-field}
\end{figure}

Before drawing the set of interference fields, we started
off by observing how the individual current fields
$\volumeCurrentDensity$, with $\unitOf\volumeCurrentDensity
= \si{\milli\ampere\per\milli\meter\squared}$ emanating from the left and right
electrode pairs C3--TP9 and C4--TF10, behaved individually.
This is seen in Figures~\ref{fig:left-current-field} and
\ref{fig:right-current-field}, where we see the left and
right unit current fields diminishing as they progress
from the electrodes to the opposite side of the volume
conductor. The choice of maximum contrast value was chosen
as $\errorDiscrepancy = \qty{-35}{\decibel}$, such that
the visualized current field $\volumeCurrentDensity$ drops
to $0$ on the other side of the volume conductor. From
the same figures, we see that the current fields have
diminished in roughly equal amounts near the center of the
volume conductor to a value of roughly \qty{-18}{\decibel}
after having moved past the thalamus from either
side. A maximal value for $\volumeCurrentDensity$
near the electrodes was approximately
\num{2.04}--\qty{2.05e-4}{\milli\ampere\per\milli\meter\squared}, in both cases.

\begin{figure}[!h]
    \centering
    \begin{minipage}{0.85\linewidth}
        \def\subFigWidth{0.31\linewidth}
        \begin{subfigure}[b]{\subFigWidth}
            \includegraphics[width=\linewidth]{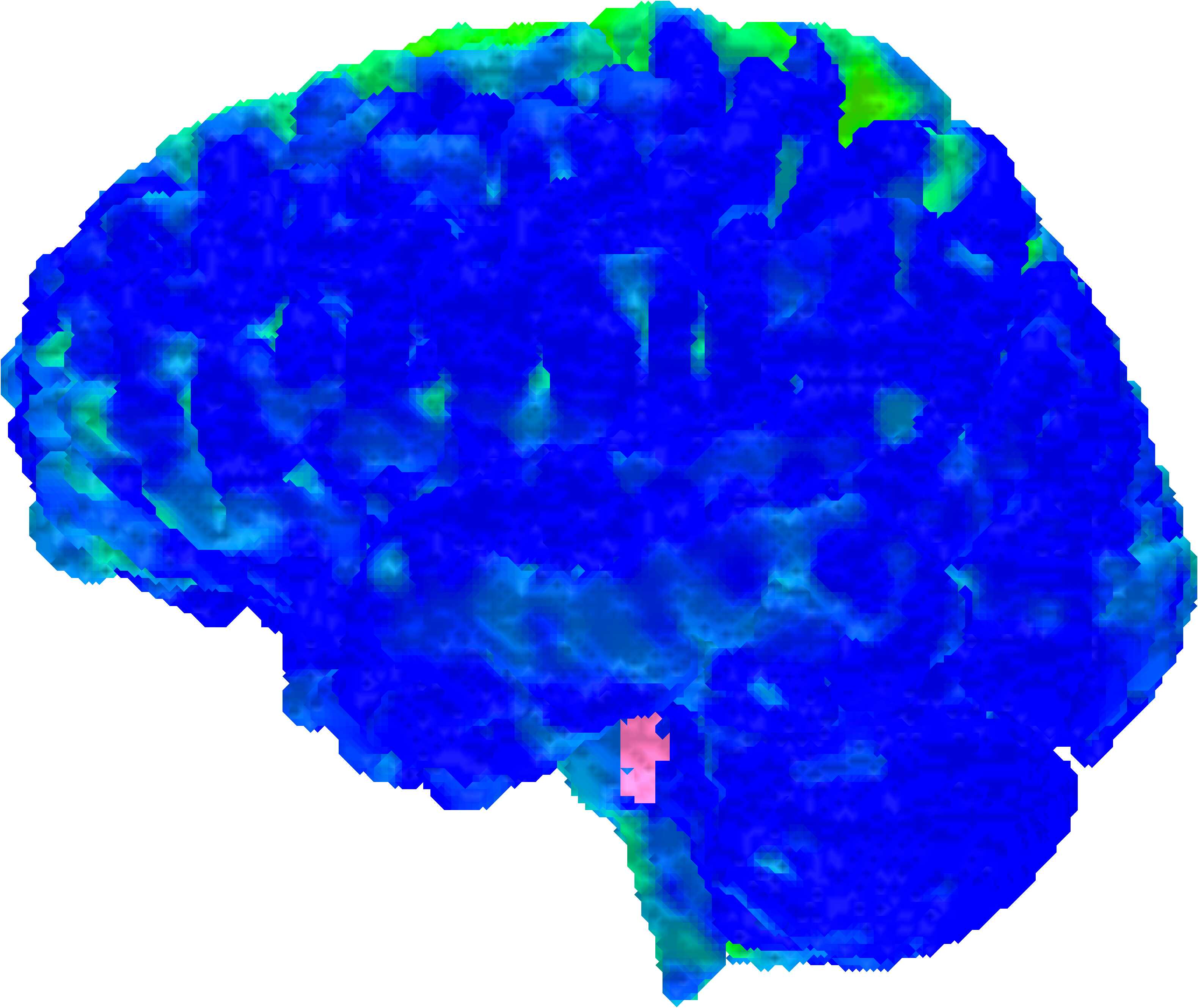}
            \caption{$\azimuth=\qty{-90}{\degree}$, $\elevation=\qty{0}{\degree}$.}
        \end{subfigure}
        \hfill
        \begin{subfigure}[b]{\subFigWidth}
        \includegraphics[width=\linewidth]{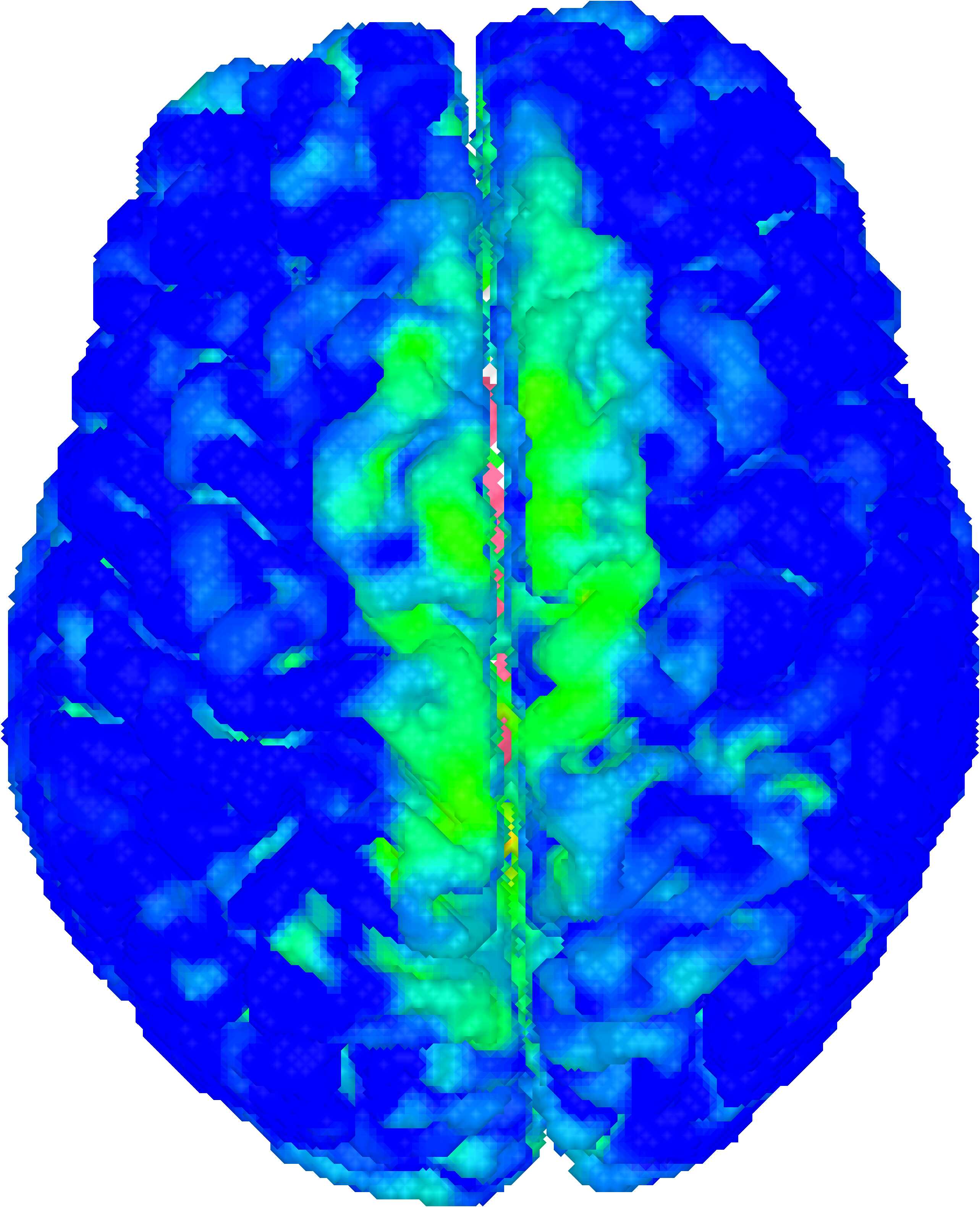}
        \caption{$\azimuth=\qty{0}{\degree}$, $\elevation=\qty{90}{\degree}$.}
        \end{subfigure}
        \hfill
        \begin{subfigure}[b]{\subFigWidth}
        \includegraphics[width=\linewidth]{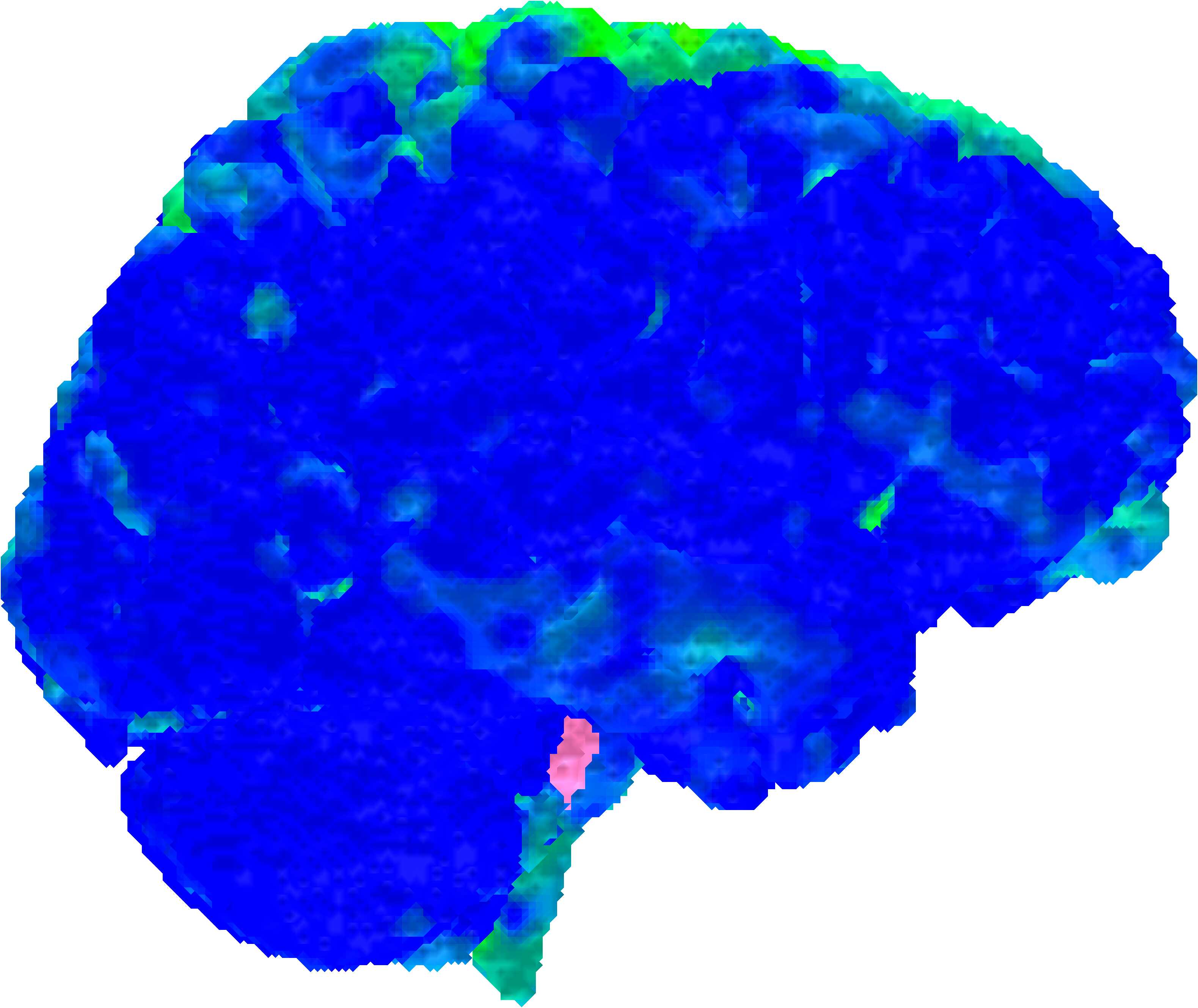}
            \caption{$\azimuth=\qty{90}{\degree}$, $\elevation=\qty{0}{\degree}$.}
        \end{subfigure}

        \def\subFigWidth{0.45\linewidth}

        \begin{subfigure}[b]{\subFigWidth}
            \includegraphics[width=\linewidth]{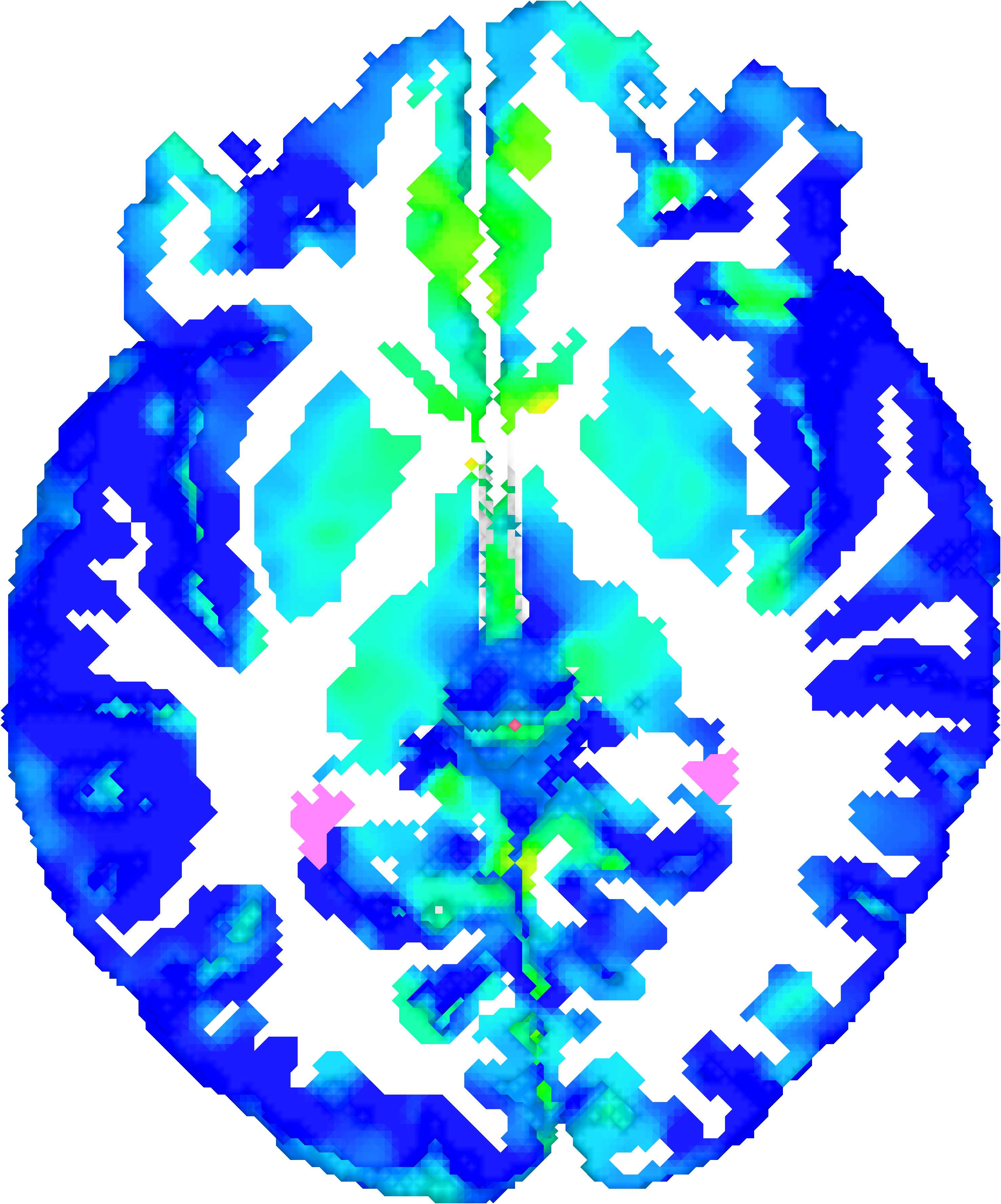}
            \caption{$\azimuth=\qty{0}{\degree}$, $\elevation=\qty{90}{\degree}$, $z=\qty{30}{\milli\meter}$.}
        \end{subfigure}
       \hfill
        \begin{subfigure}[b]{\subFigWidth}
            \includegraphics[width=\linewidth]{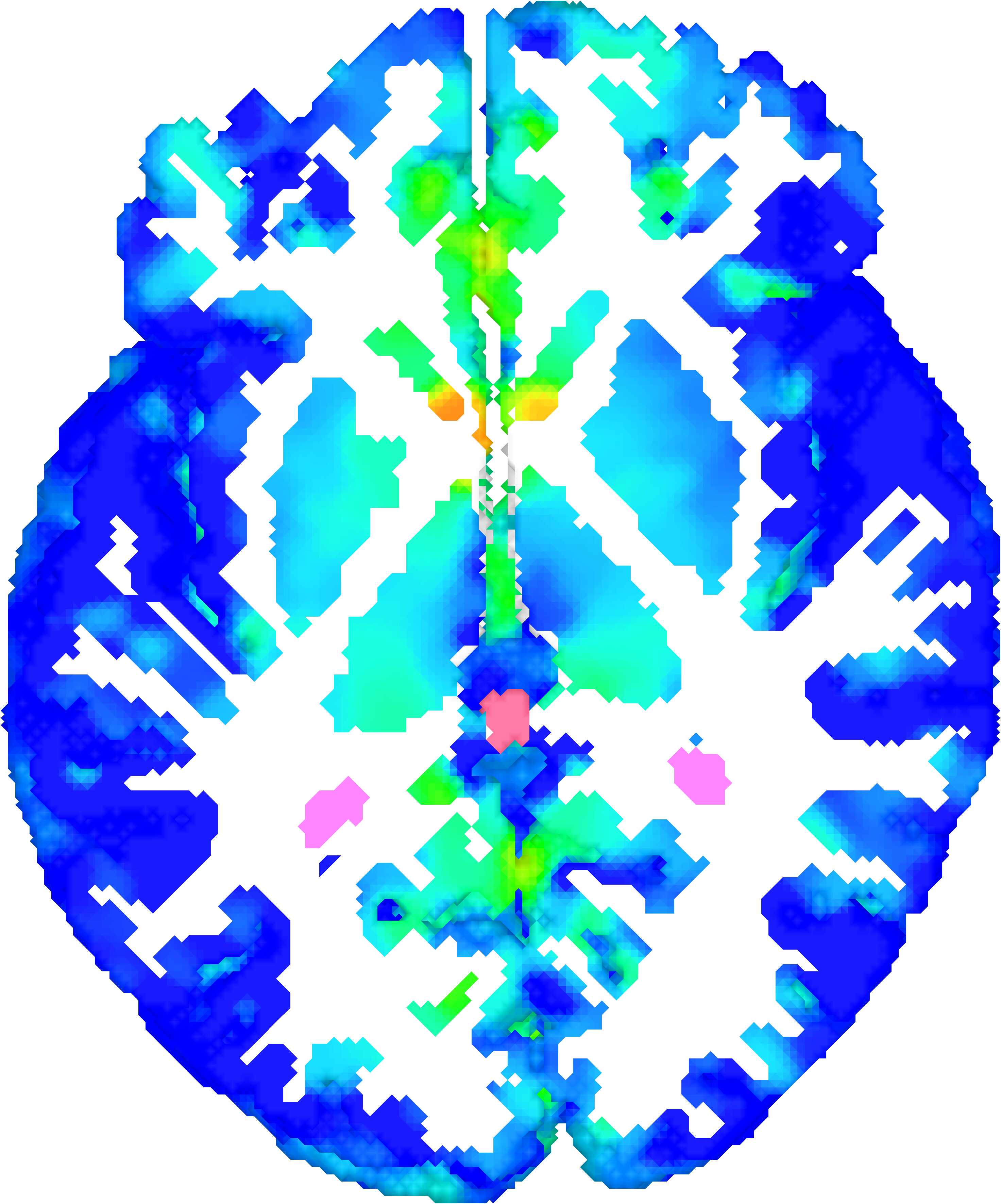}
            \caption{$\azimuth=\qty{0}{\degree}$, $\elevation=\qty{90}{\degree}$, $z=\qty{33}{\milli\meter}$.}
        \end{subfigure}

        \begin{subfigure}[b]{\subFigWidth}
            \includegraphics[width=\linewidth]{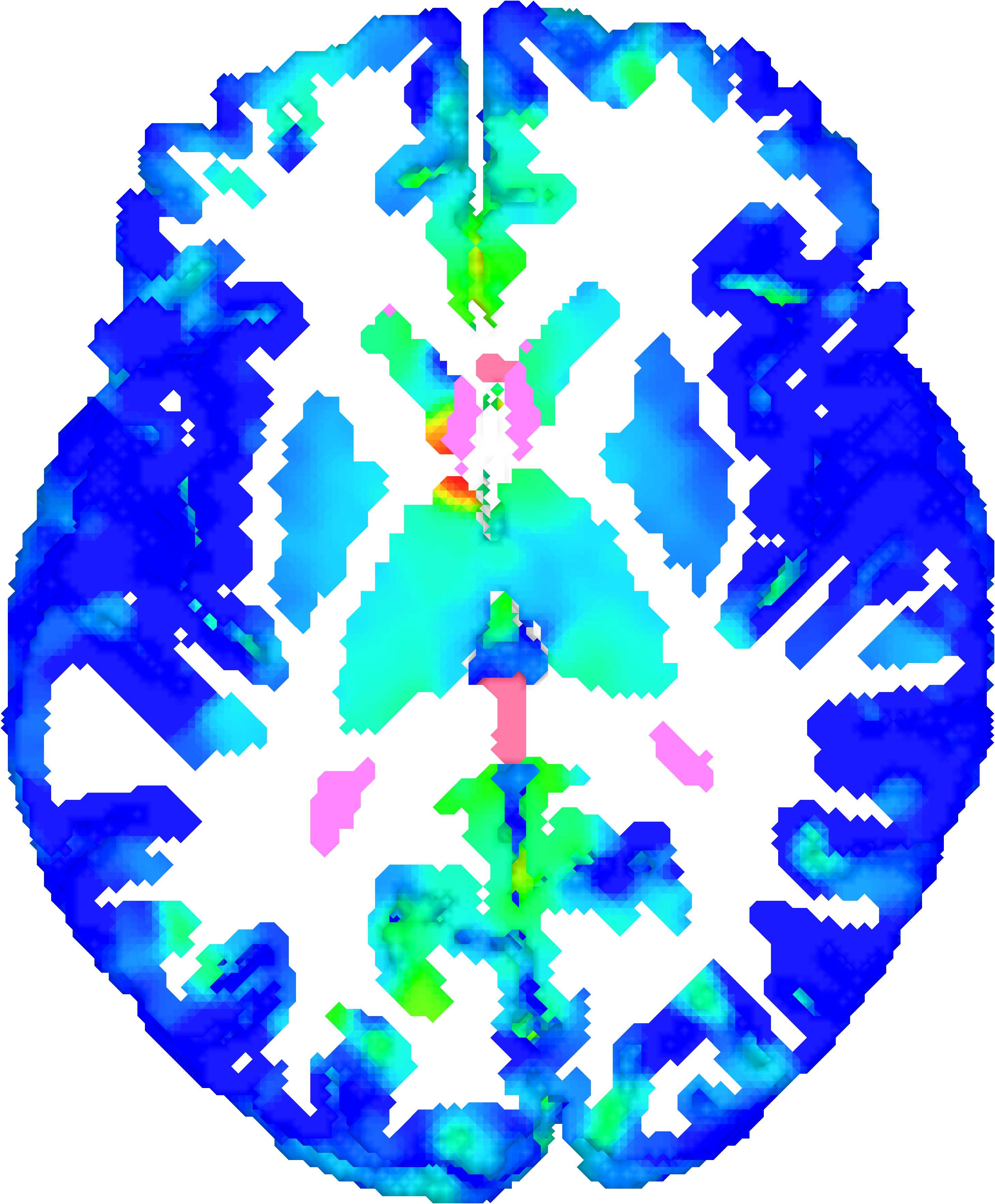}
            \caption{$\azimuth=\qty{0}{\degree}$, $\elevation=\qty{90}{\degree}$, $z=\qty{37}{\milli\meter}$.}
        \end{subfigure}
        \hfill
        \begin{subfigure}[b]{\subFigWidth}
            \includegraphics[width=\linewidth]{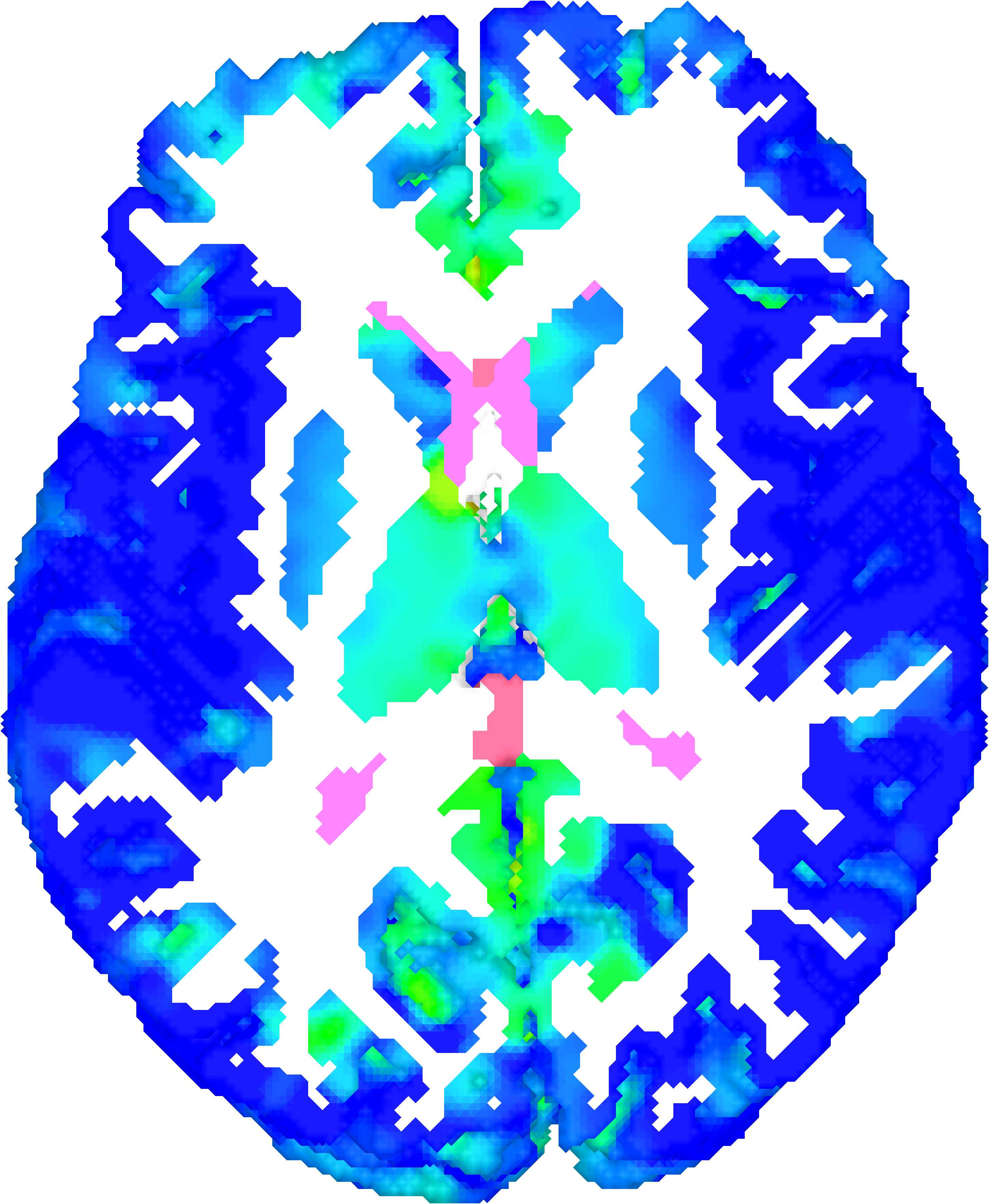}
            \caption{$\azimuth=\qty{0}{\degree}$, $\elevation=\qty{90}{\degree}$, $z=\qty{40}{\milli\meter}$.}
        \end{subfigure}
    \end{minipage}
    \hfill
    \begin{minipage}{0.09\linewidth}
    \includegraphics[width=\linewidth]{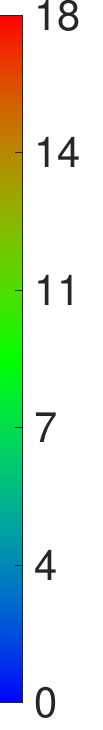}
    \end{minipage}
    \caption{Volume current interference field \eqref{eq:interference-field} of $\volumeCurrentDensity$, smoothed via the function $\dBT\args{\someField/\max\someField,\qty{18}{\decibel}}$ of equation \eqref{eq:visual-log-transform} at different azimuths $\azimuth$ and elevations $\elevation$, with base stimulation frequency $\frequency=\qty{1000}{\hertz}$. The current pattern was $\currentPattern=[-1.0,1.0,-1.0,1.0]\,\si{\milli\ampere}$ and the maximal magnitude of the output $\volumeCurrentDensity$ was approximately \qty{6.4e-5}{\milli\ampere\per\milli\meter\squared}.}
    \label{fig:interference-field-equal-currents}
\end{figure}

\begin{figure}[!h]
    \centering
    \begin{minipage}{0.85\linewidth}
        \def\subFigWidth{0.31\linewidth}
        \begin{subfigure}[b]{\subFigWidth}
            \includegraphics[width=\linewidth]{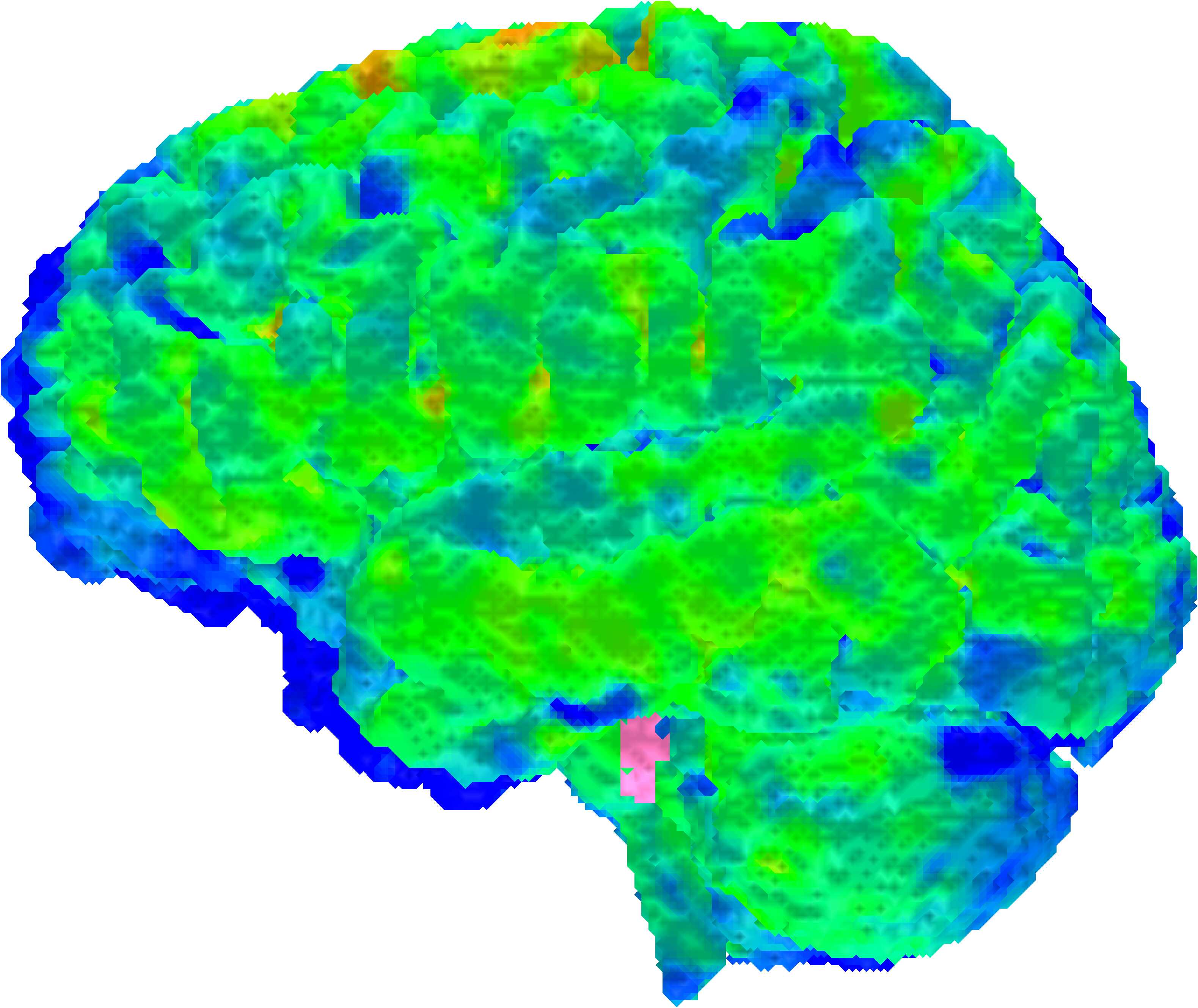}
            \caption{$\azimuth=\qty{-90}{\degree}$, $\elevation=\qty{0}{\degree}$.}
        \end{subfigure}
        \hfill
        \begin{subfigure}[b]{\subFigWidth}
        \includegraphics[width=\linewidth]{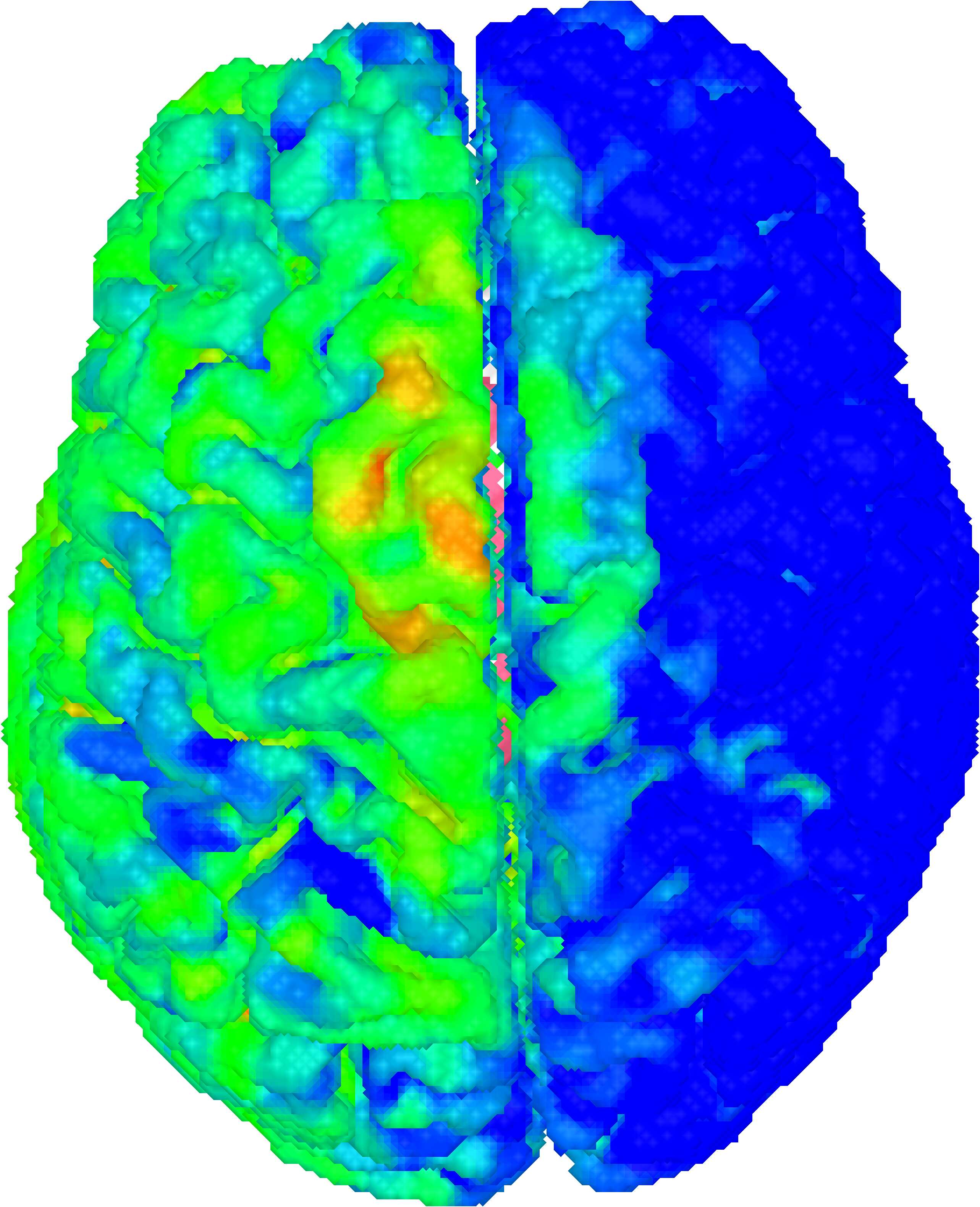}
        \caption{$\azimuth=\qty{0}{\degree}$, $\elevation=\qty{90}{\degree}$.}
        \end{subfigure}
        \hfill
        \begin{subfigure}[b]{\subFigWidth}
        \includegraphics[width=\linewidth]{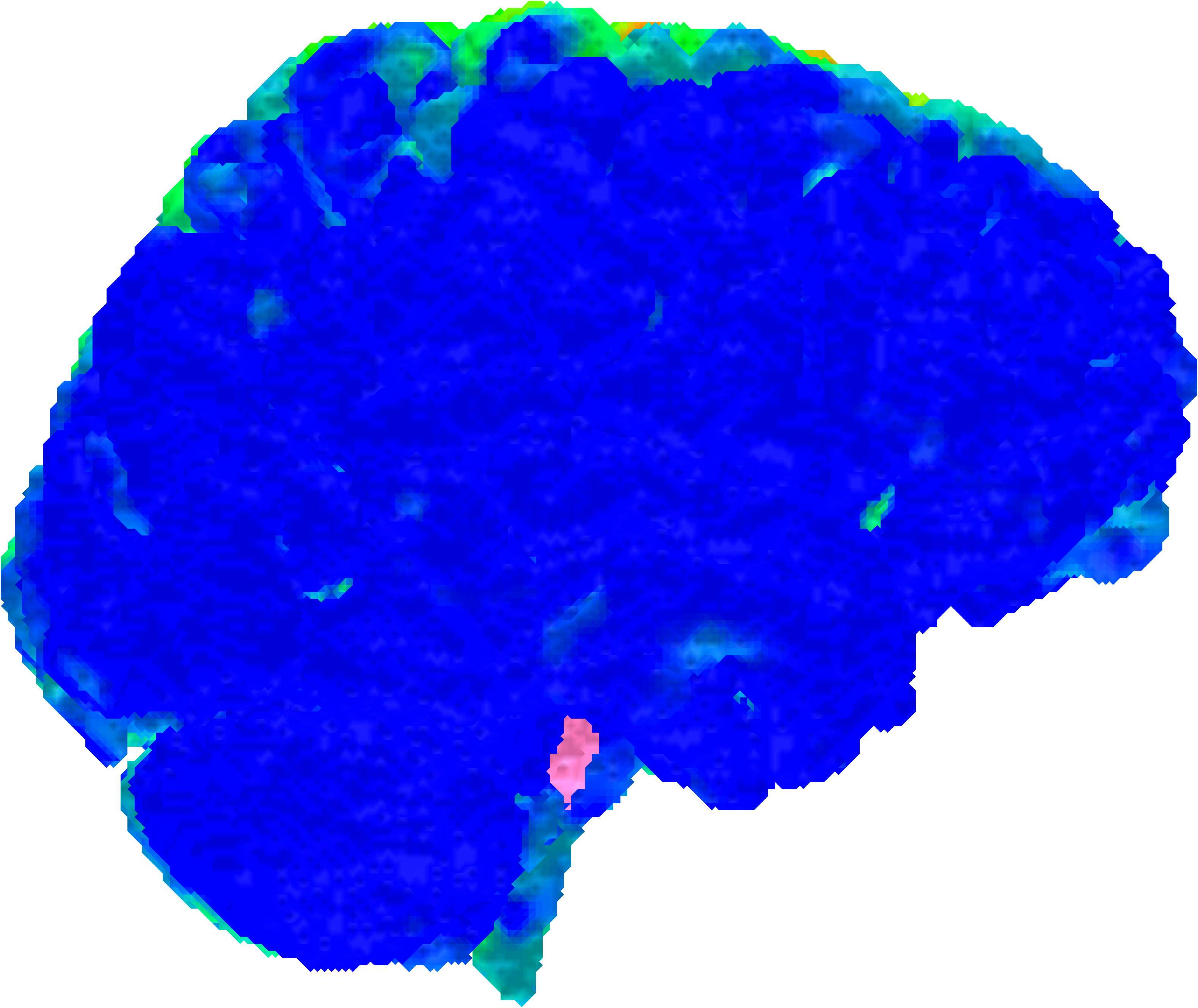}
            \caption{$\azimuth=\qty{90}{\degree}$, $\elevation=\qty{0}{\degree}$.}
        \end{subfigure}

        \def\subFigWidth{0.45\linewidth}

        \begin{subfigure}[b]{\subFigWidth}
            \includegraphics[width=\linewidth]{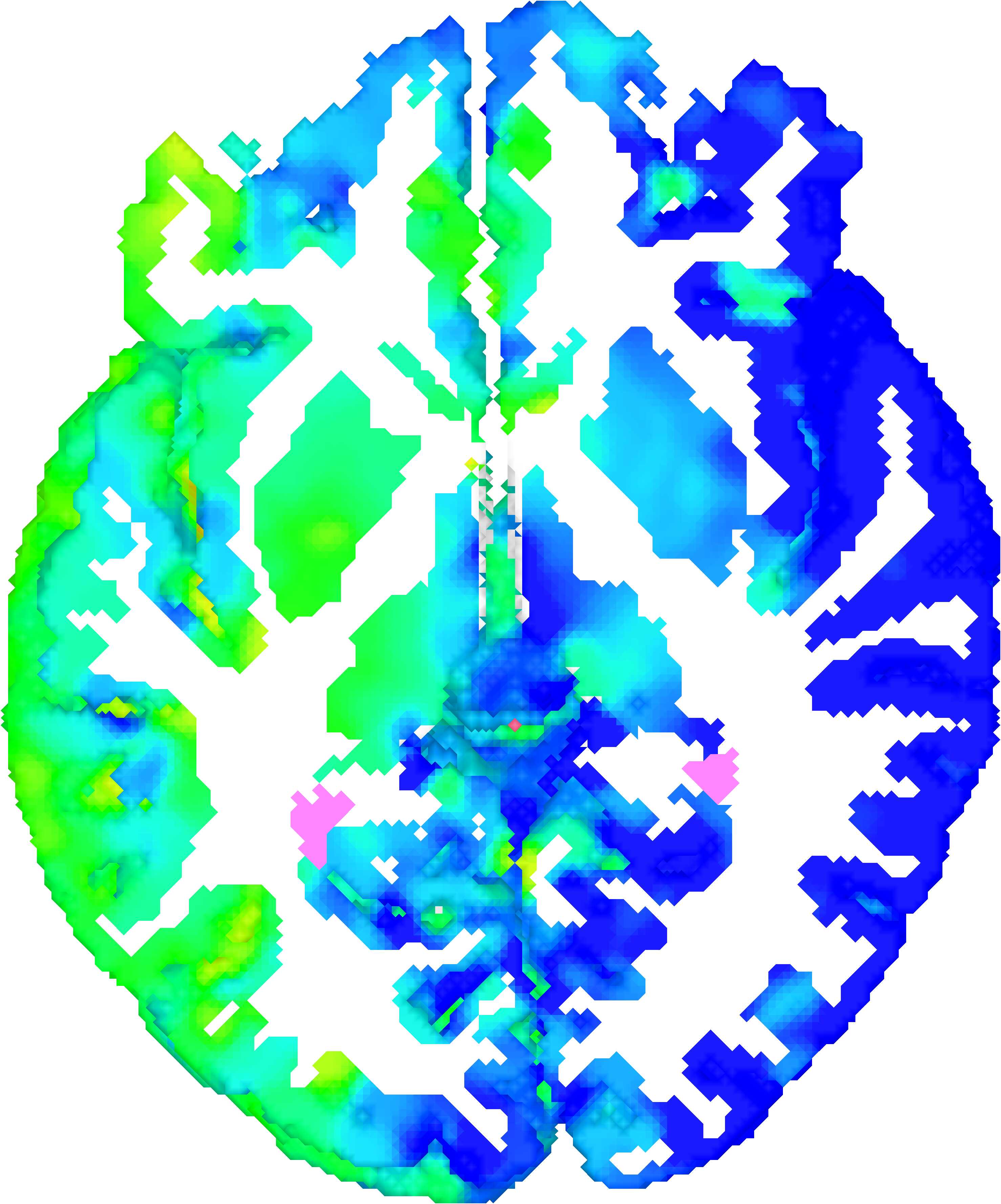}
            \caption{$\azimuth=\qty{0}{\degree}$, $\elevation=\qty{90}{\degree}$, $z=\qty{30}{\milli\meter}$.}
        \end{subfigure}
       \hfill
        \begin{subfigure}[b]{\subFigWidth}
            \includegraphics[width=\linewidth]{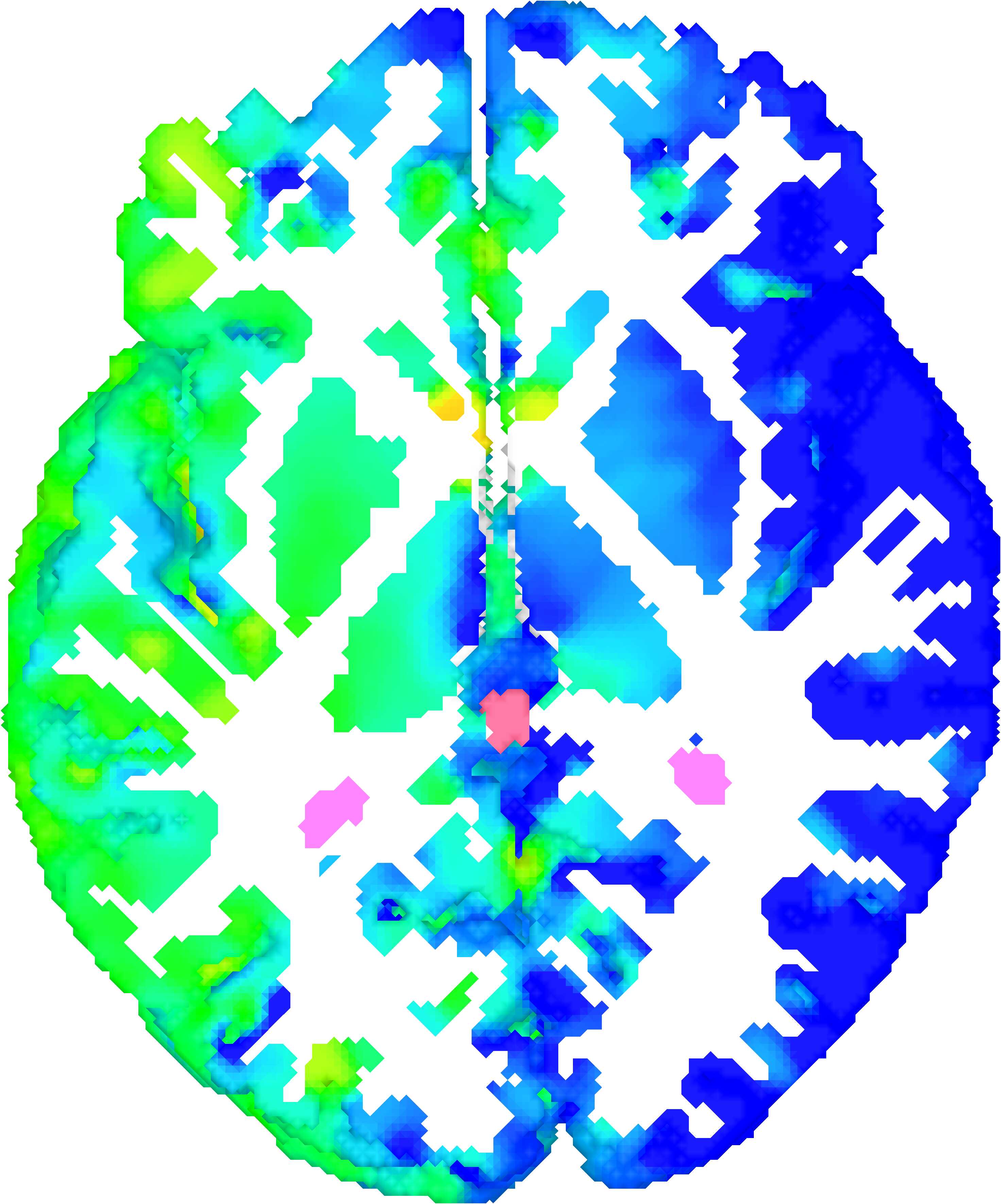}
            \caption{$\azimuth=\qty{0}{\degree}$, $\elevation=\qty{90}{\degree}$, $z=\qty{33}{\milli\meter}$.}
        \end{subfigure}

        \begin{subfigure}[b]{\subFigWidth}
            \includegraphics[width=\linewidth]{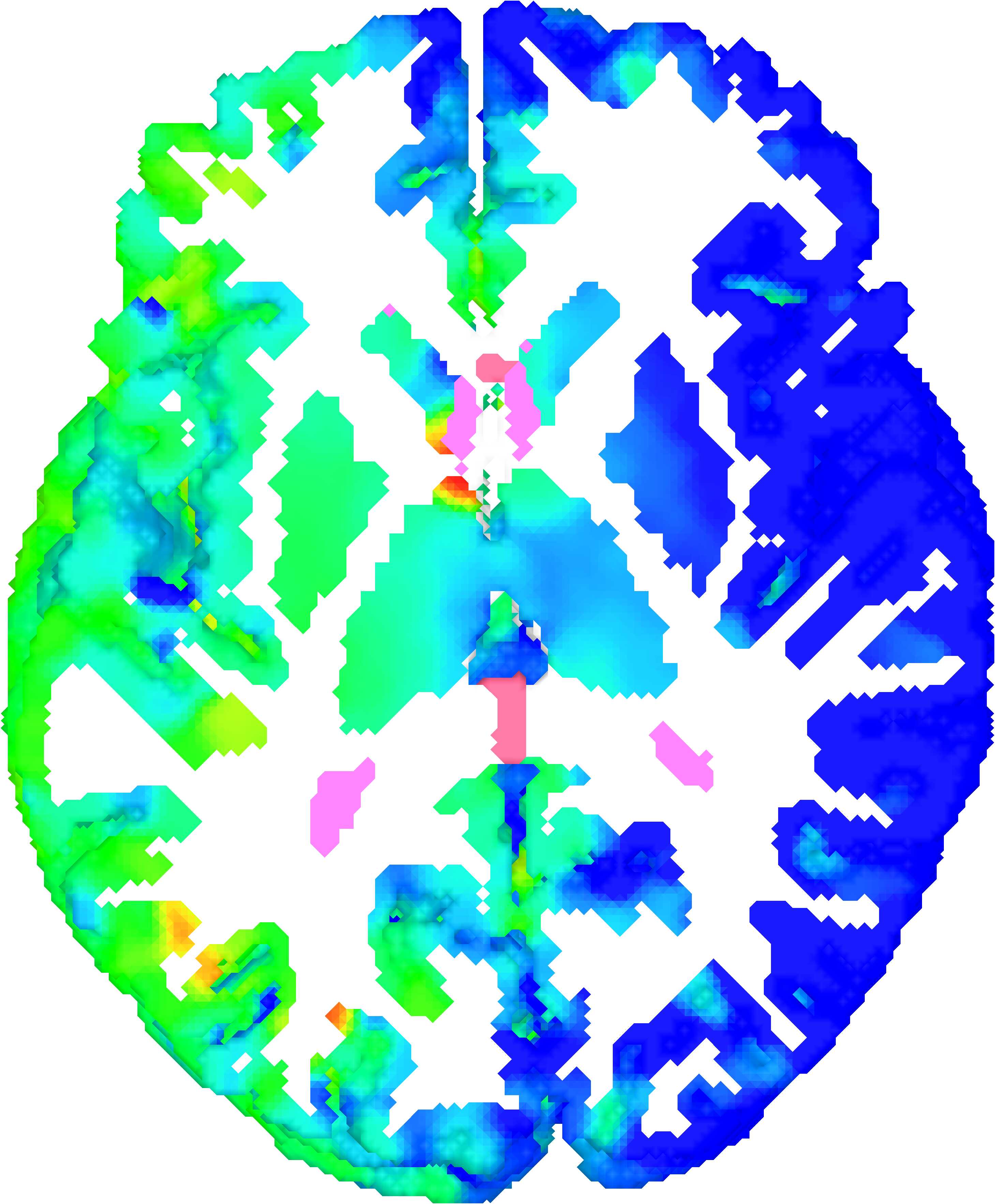}
            \caption{$\azimuth=\qty{0}{\degree}$, $\elevation=\qty{90}{\degree}$, $z=\qty{37}{\milli\meter}$.}
        \end{subfigure}
        \hfill
        \begin{subfigure}[b]{\subFigWidth}
            \includegraphics[width=\linewidth]{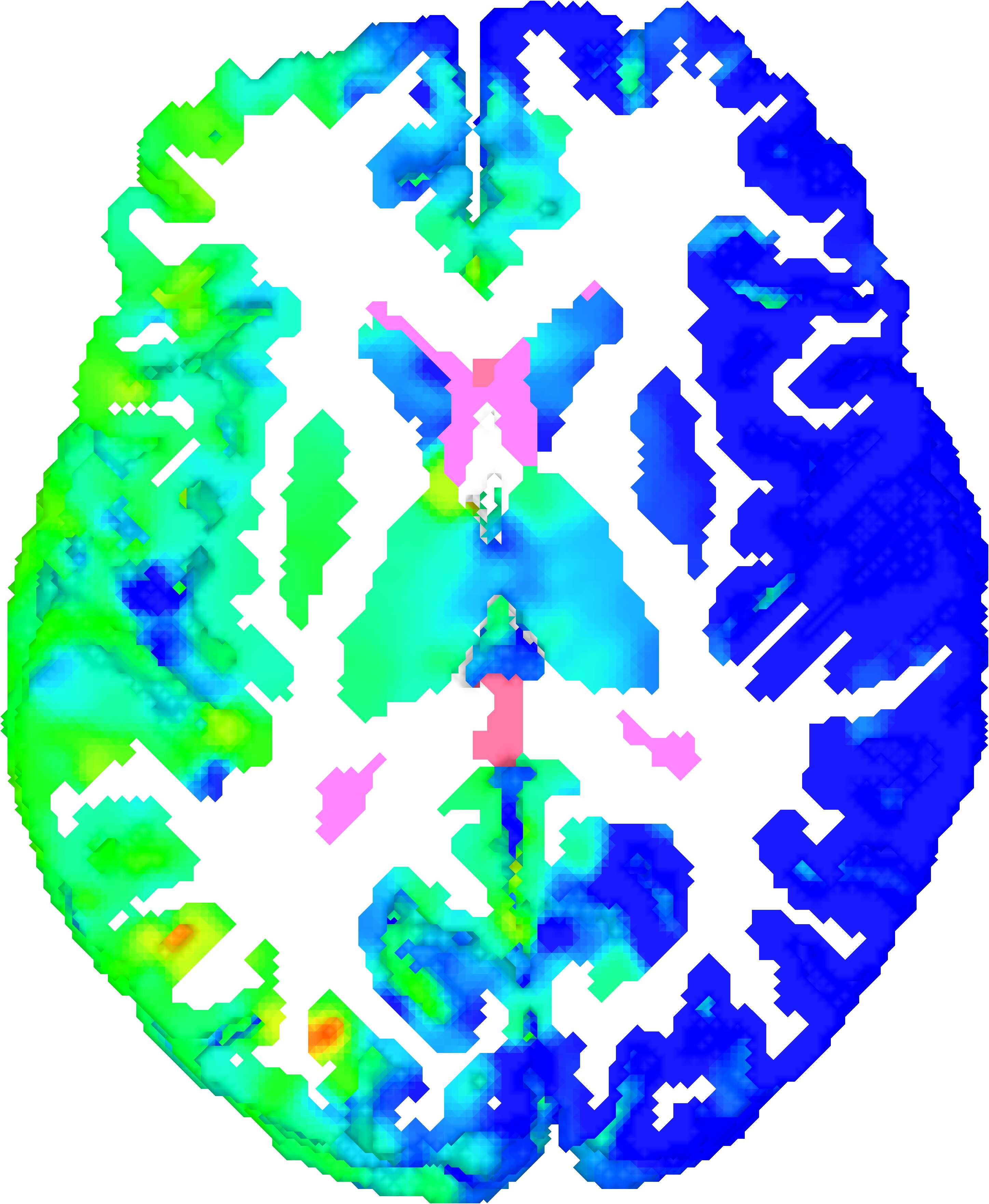}
            \caption{$\azimuth=\qty{0}{\degree}$, $\elevation=\qty{90}{\degree}$, $z=\qty{40}{\milli\meter}$.}
        \end{subfigure}
    \end{minipage}
    \hfill
    \begin{minipage}{0.09\linewidth}
    \includegraphics[width=\linewidth]{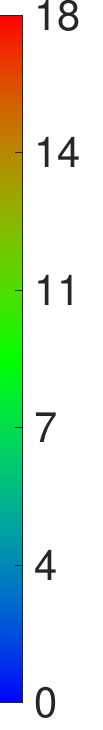}
    \end{minipage}
    \caption{Volume current interference field \eqref{eq:interference-field} of $\volumeCurrentDensity$, smoothed via the function $\dBT\args{\someField/\max\someField,\qty{18}{\decibel}}$ of equation \eqref{eq:visual-log-transform} at different azimuths $\azimuth$ and elevations $\elevation$, with base stimulation frequency $\frequency=\qty{1000}{\hertz}$. The current pattern was $\currentPattern=[-0.5,0.5,-1.5,1.5]\,\si{\milli\ampere}$ and the maximal value of $\volumeCurrentDensity$ roughly equal to \qty{4.0e-5}{\milli\ampere\per\milli\meter\squared}.}
    \label{fig:interference-field-weaker-left}
\end{figure}

\begin{figure}[!h]
    \centering
    \begin{minipage}{0.85\linewidth}
        \def\subFigWidth{0.31\linewidth}
        \begin{subfigure}[b]{\subFigWidth}
            \includegraphics[width=\linewidth]{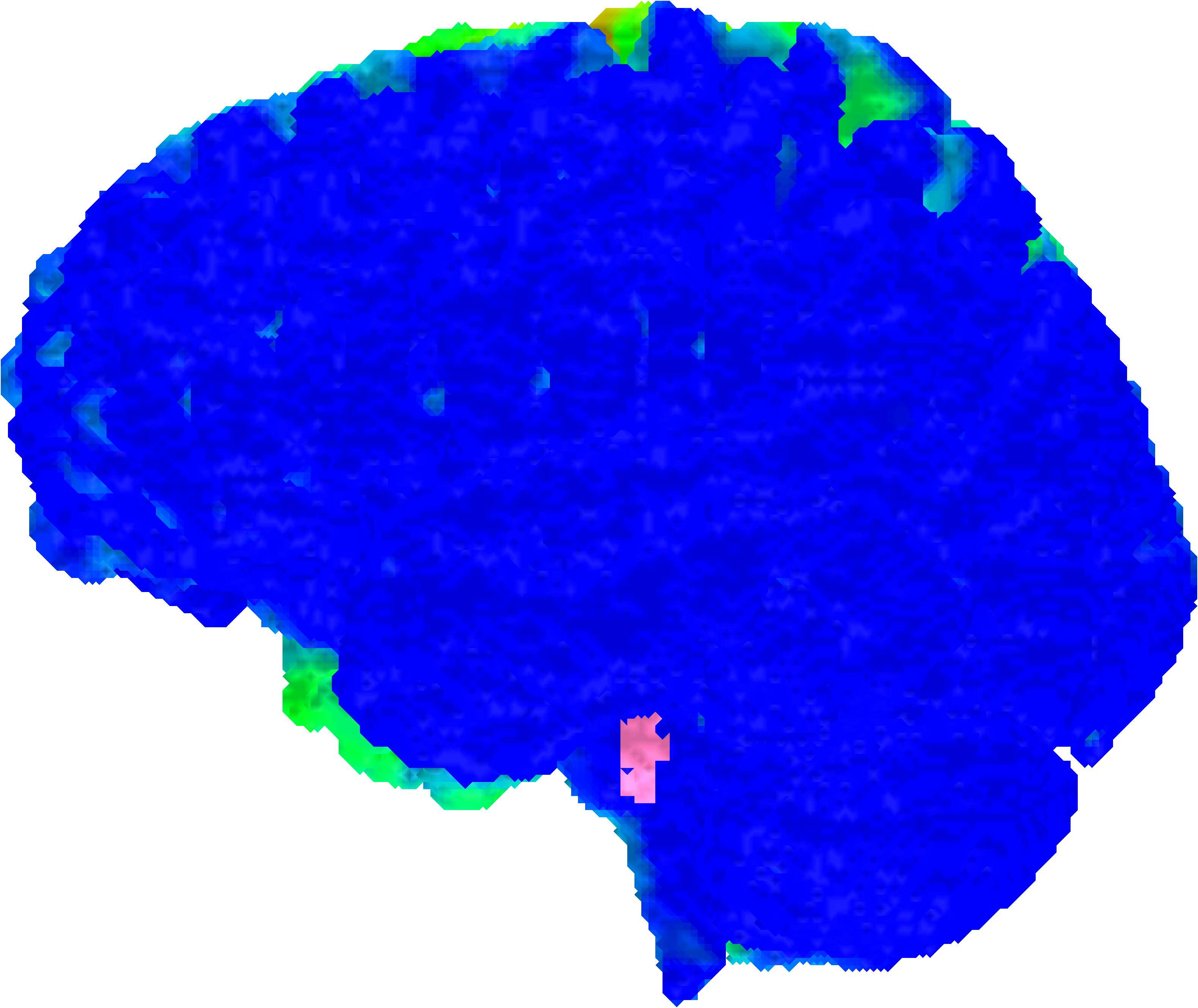}
            \caption{$\azimuth=\qty{-90}{\degree}$, $\elevation=\qty{0}{\degree}$.}
        \end{subfigure}
        \hfill
        \begin{subfigure}[b]{\subFigWidth}
        \includegraphics[width=\linewidth]{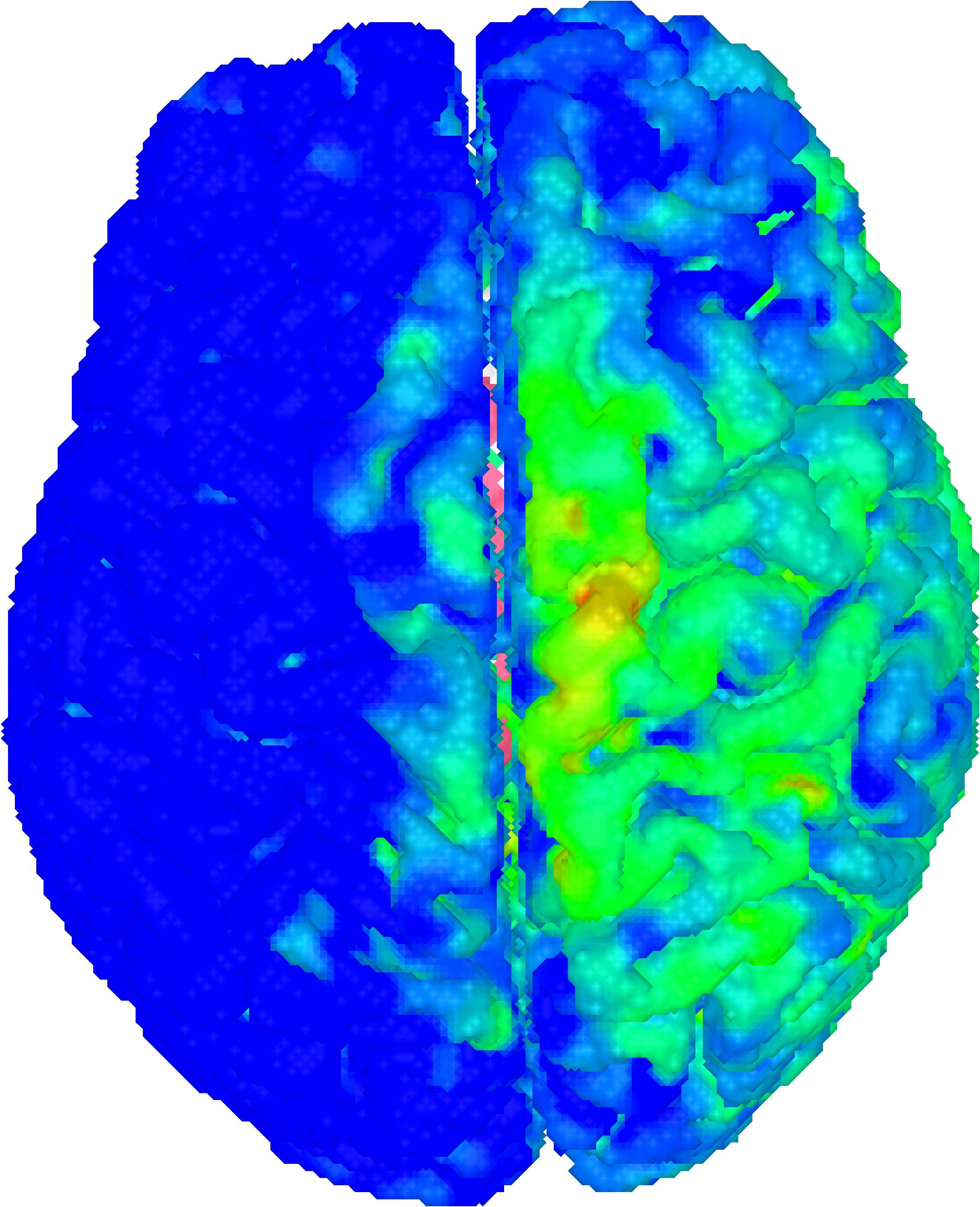}
        \caption{$\azimuth=\qty{0}{\degree}$, $\elevation=\qty{90}{\degree}$.}
        \end{subfigure}
        \hfill
        \begin{subfigure}[b]{\subFigWidth}
        \includegraphics[width=\linewidth]{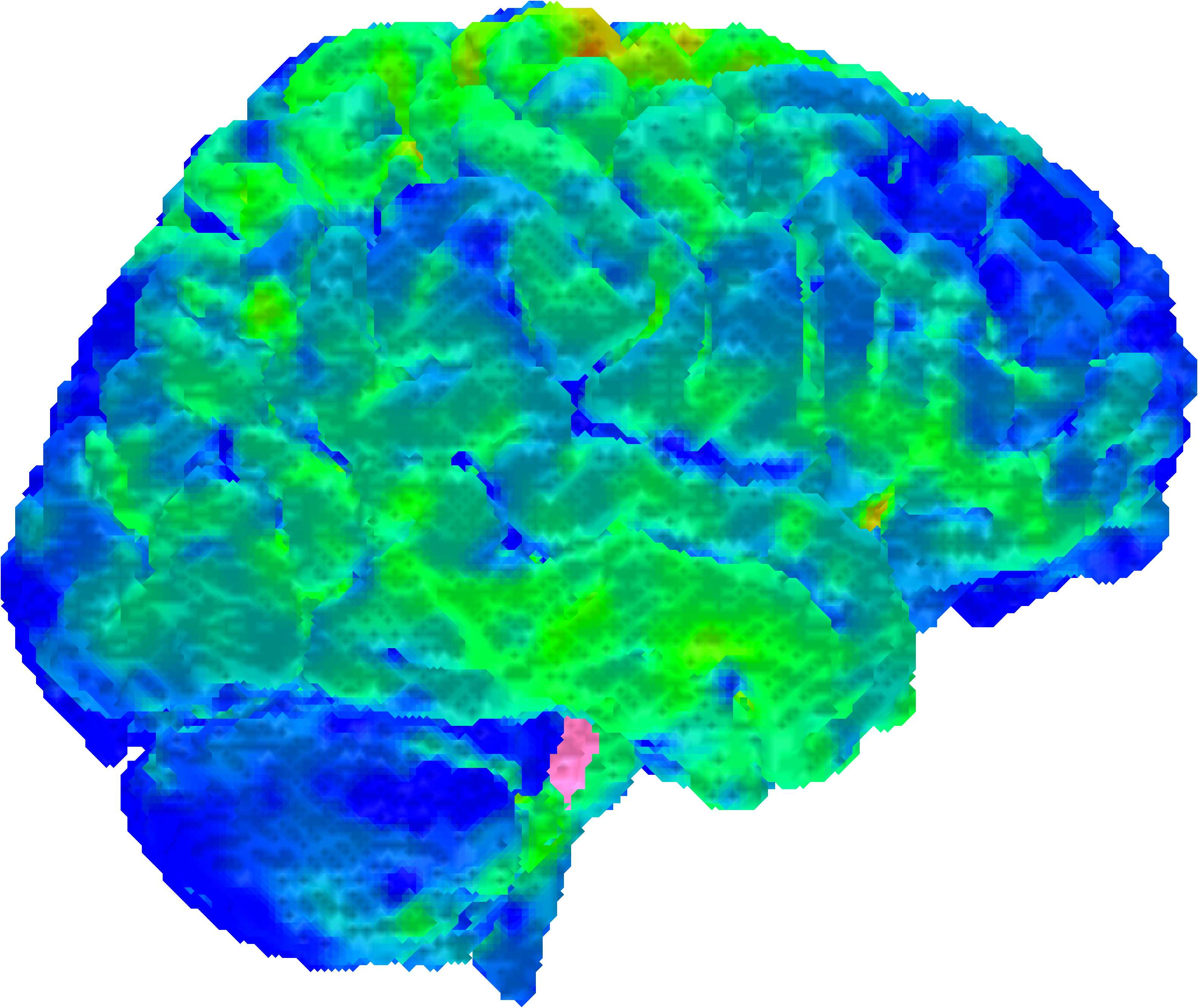}
            \caption{$\azimuth=\qty{90}{\degree}$, $\elevation=\qty{0}{\degree}$.}
        \end{subfigure}

        \def\subFigWidth{0.45\linewidth}

        \begin{subfigure}[b]{\subFigWidth}
            \includegraphics[width=\linewidth]{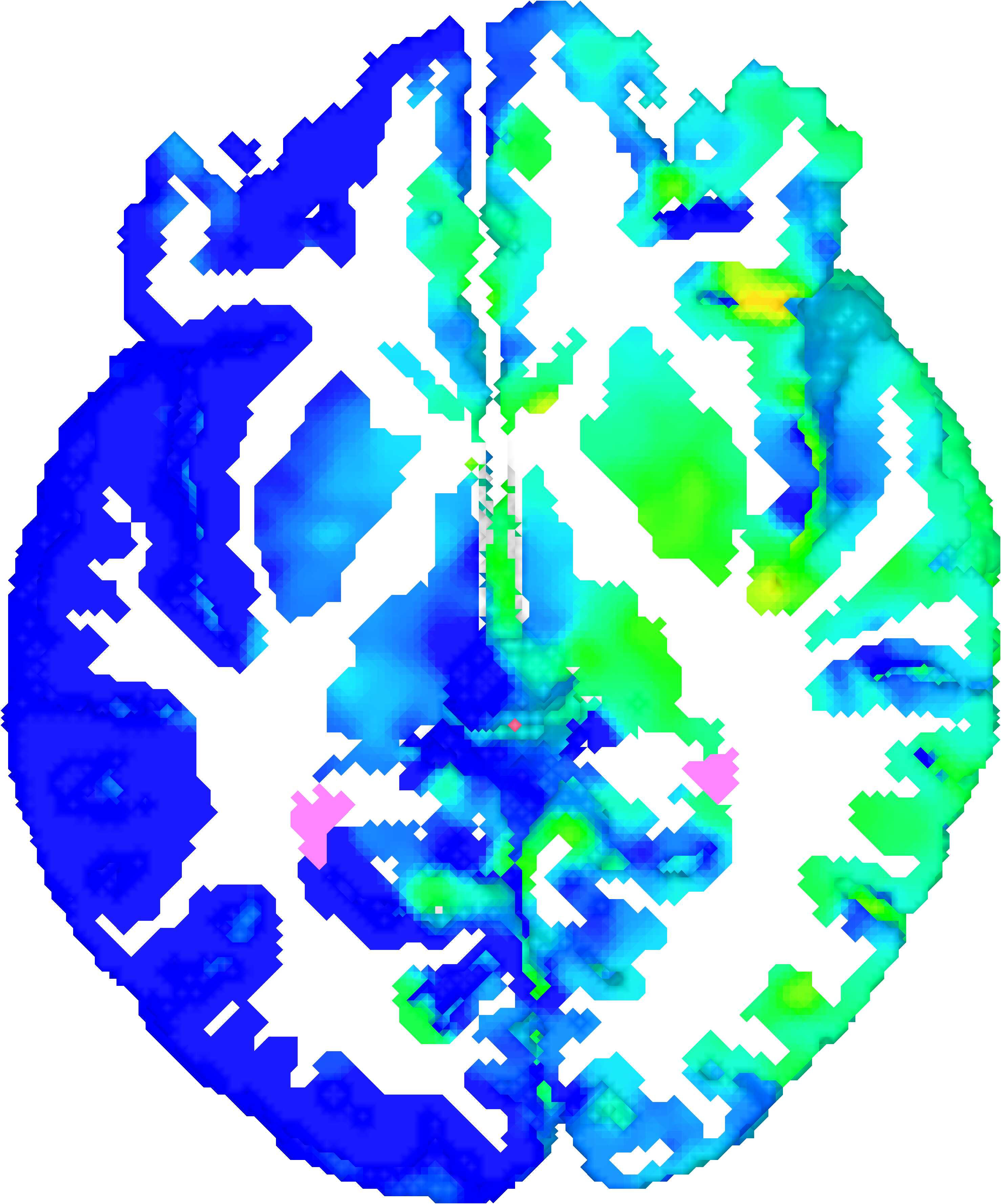}
            \caption{$\azimuth=\qty{0}{\degree}$, $\elevation=\qty{90}{\degree}$, $z=\qty{30}{\milli\meter}$.}
        \end{subfigure}
       \hfill
        \begin{subfigure}[b]{\subFigWidth}
            \includegraphics[width=\linewidth]{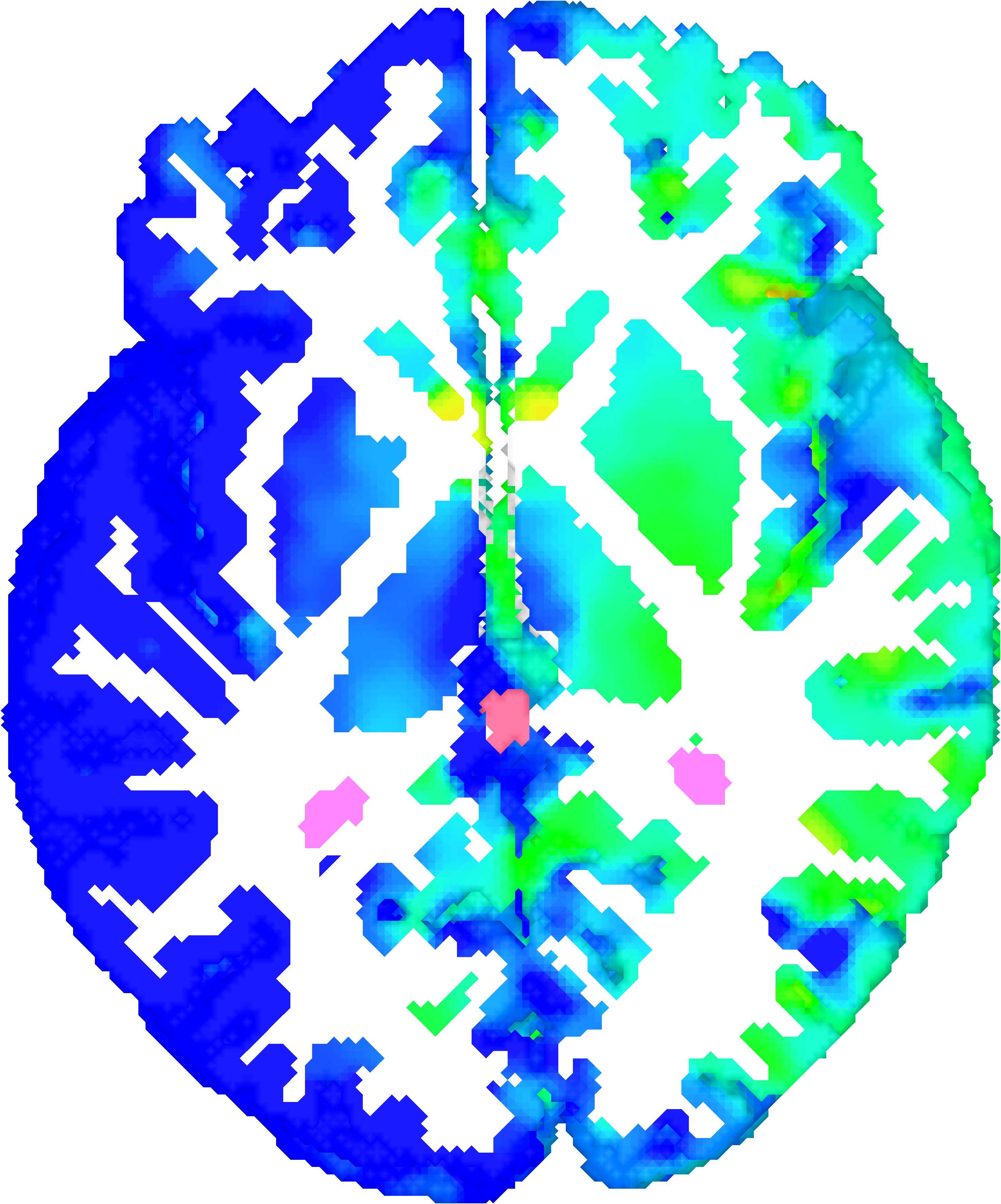}
            \caption{$\azimuth=\qty{0}{\degree}$, $\elevation=\qty{90}{\degree}$, $z=\qty{33}{\milli\meter}$.}
        \end{subfigure}

        \begin{subfigure}[b]{\subFigWidth}
            \includegraphics[width=\linewidth]{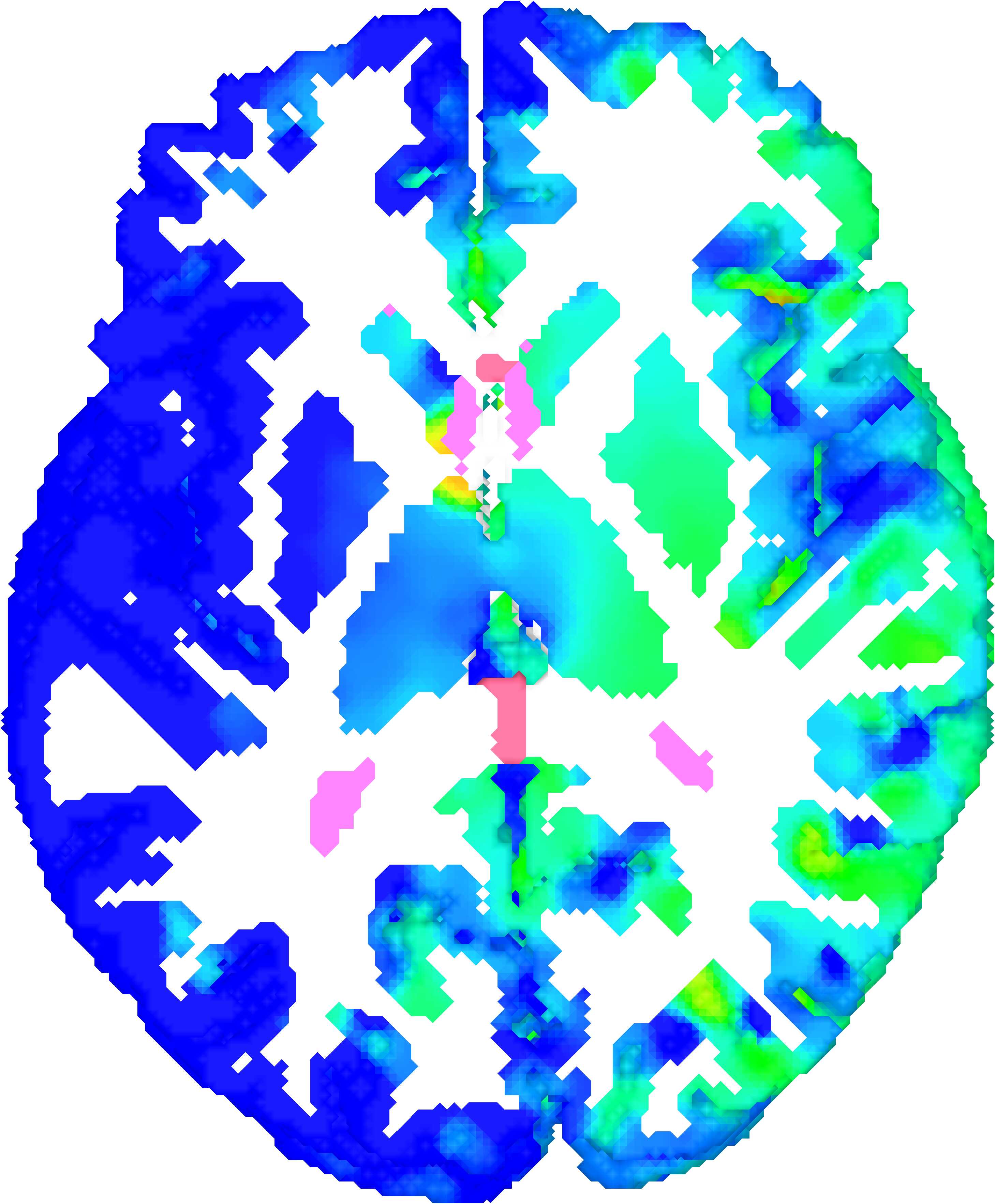}
            \caption{$\azimuth=\qty{0}{\degree}$, $\elevation=\qty{90}{\degree}$, $z=\qty{37}{\milli\meter}$.}
        \end{subfigure}
        \hfill
        \begin{subfigure}[b]{\subFigWidth}
            \includegraphics[width=\linewidth]{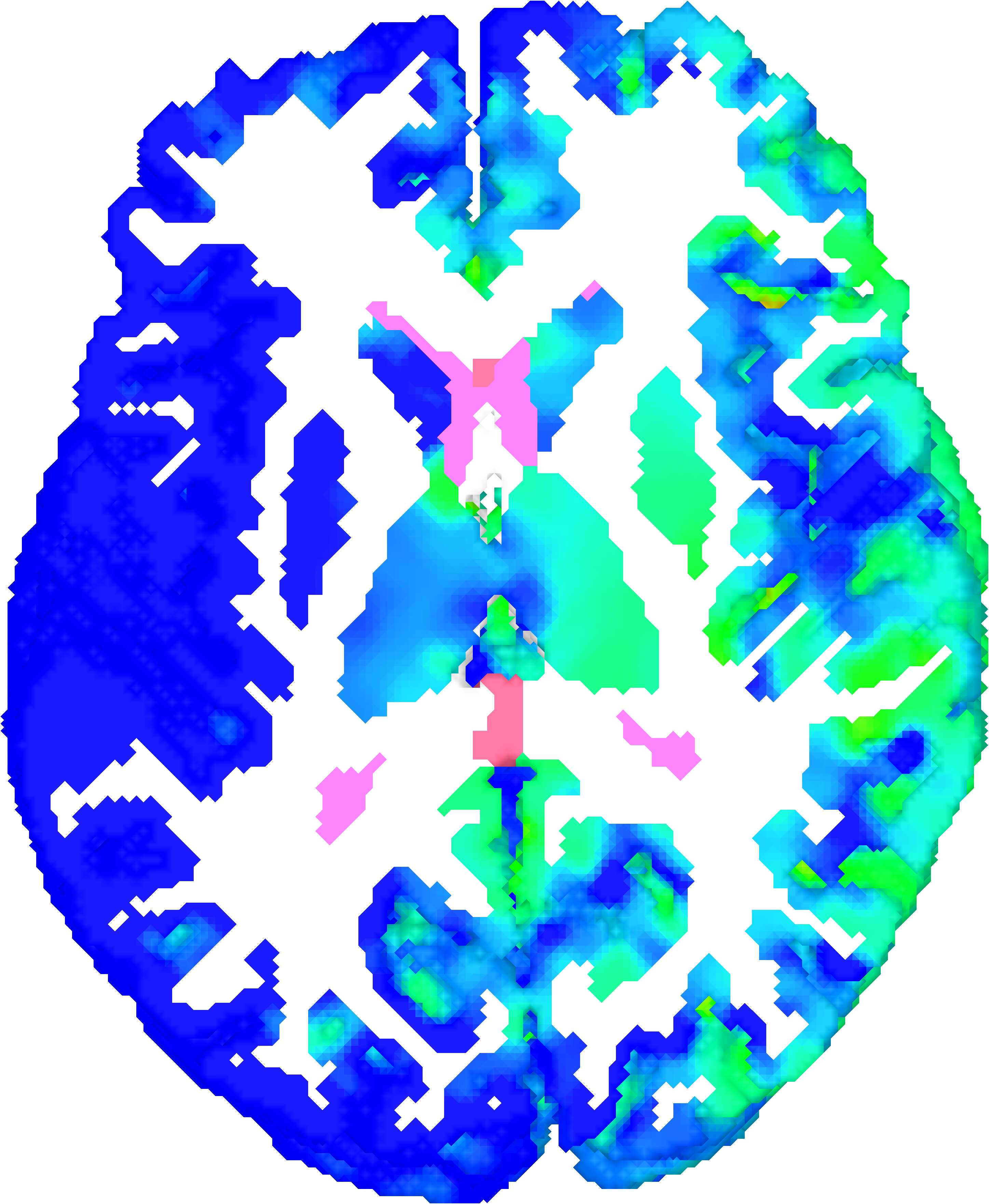}
            \caption{$\azimuth=\qty{0}{\degree}$, $\elevation=\qty{90}{\degree}$, $z=\qty{40}{\milli\meter}$.}
        \end{subfigure}
    \end{minipage}
    \hfill
    \begin{minipage}{0.09\linewidth}
    \includegraphics[width=\linewidth]{c,f=1000,li=0.5,ri=1.5,dB=18,az=0,el=90.colorbar.pdf}
    \end{minipage}
    \caption{Volume current interference field \eqref{eq:interference-field} of $\volumeCurrentDensity$, smoothed via the function $\dBT\args{\someField/\max\someField,\qty{18}{\decibel}}$ of equation \eqref{eq:visual-log-transform} at different azimuths $\azimuth$ and elevations $\elevation$, with base stimulation frequency $\frequency=\qty{1000}{\hertz}$. The current pattern was $\currentPattern=[-1.5,1.5,-0.5,0.5]\,\si{\milli\ampere}$  and the maximal value of $\volumeCurrentDensity$ roughly equal to \qty{5.0e-5}{\milli\ampere\per\milli\meter\squared}.}
    \label{fig:interference-field-weaker-right}
\end{figure}

Knowing a sensible maximal contrast value
$\errorDiscrepancy$ near the stimulation area of
interest, the thalamus, we set out to depict the
interference fields \eqref{eq:interference-field}
with $\volumeCurrentDensity_1$ as the field
emanating from the left electrode pair C3--TP9
and $\volumeCurrentDensity_2$ as the field
emanating from the right electrode pair C4--TF10.
Figures~\ref{fig:interference-field-equal-currents}--\ref{fig:interference-field-weaker-right}
show how the interference field behaves, when the sum of
total currents following Kirchhoffs first law is kept at a
constant $\qty{2}{\milli\ampere}$, but is varied such that
the currents are either equal in strength, the left current
source is weaker than the right, or vice versa. With equal
currents, the maximum of the distribution is seen to occur
near the front of the thalamus, and when the input currents
are changed, the peak always moves towards the electrodes
with a weaker current throughput. The maximal value for
the inteference field in the volume for equal currents was
$\volumeCurrentDensity\approx\qty{6.4e-5}{\milli\ampere\per\milli\meter\squared}$.
With the right current source outputting smaller
current, the maximum is slightly lowered to
\qty{5.04e-5}{\milli\ampere\per\milli\meter\squared},
and when the left current source is upholding
a smaller amplitude, we have an interference
current field with a magnitude of approximately
\qty{3.99e-5}{\milli\ampere\per\milli\meter\squared} at
most.

\begin{figure}[!h]
    \centering
    \begin{minipage}{0.85\linewidth}
        \def\subFigWidth{0.31\linewidth}
        \begin{subfigure}[b]{\subFigWidth}
            \includegraphics[width=\linewidth]{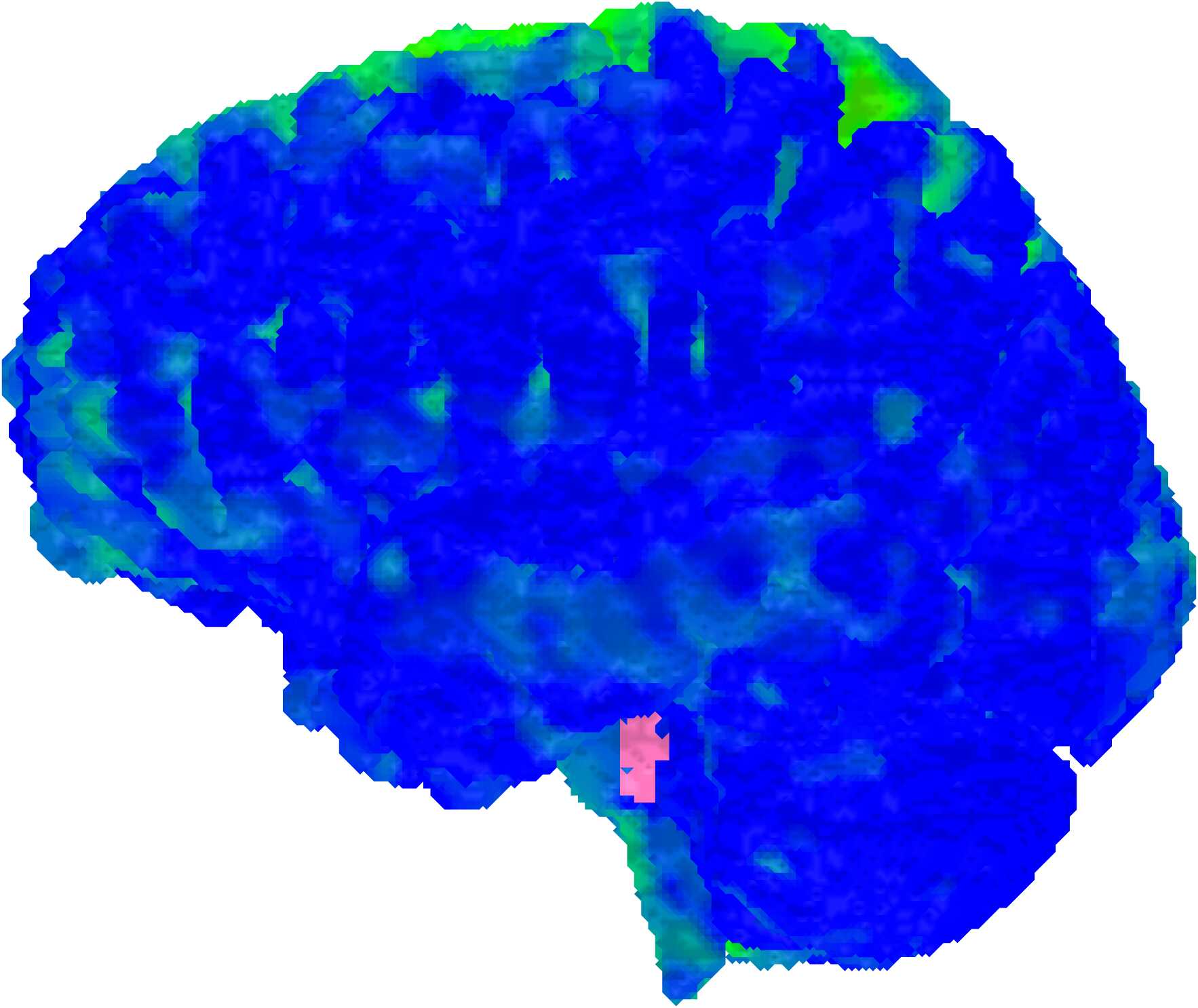}
            \caption{$\azimuth=\qty{-90}{\degree}$, $\elevation=\qty{0}{\degree}$.}
        \end{subfigure}
        \hfill
        \begin{subfigure}[b]{\subFigWidth}
        \includegraphics[width=\linewidth]{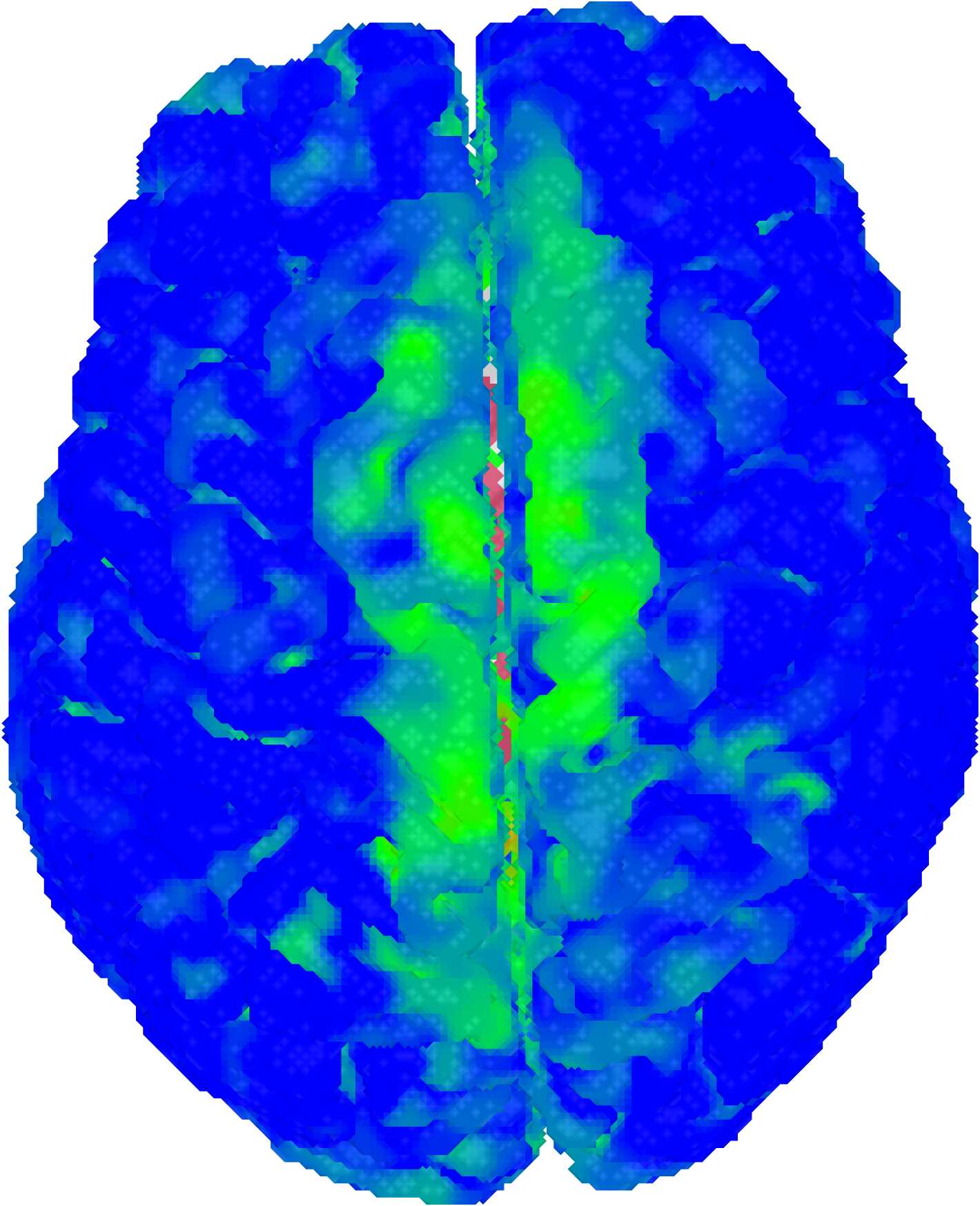}
        \caption{$\azimuth=\qty{0}{\degree}$, $\elevation=\qty{90}{\degree}$.}
        \end{subfigure}
        \hfill
        \begin{subfigure}[b]{\subFigWidth}
        \includegraphics[width=\linewidth]{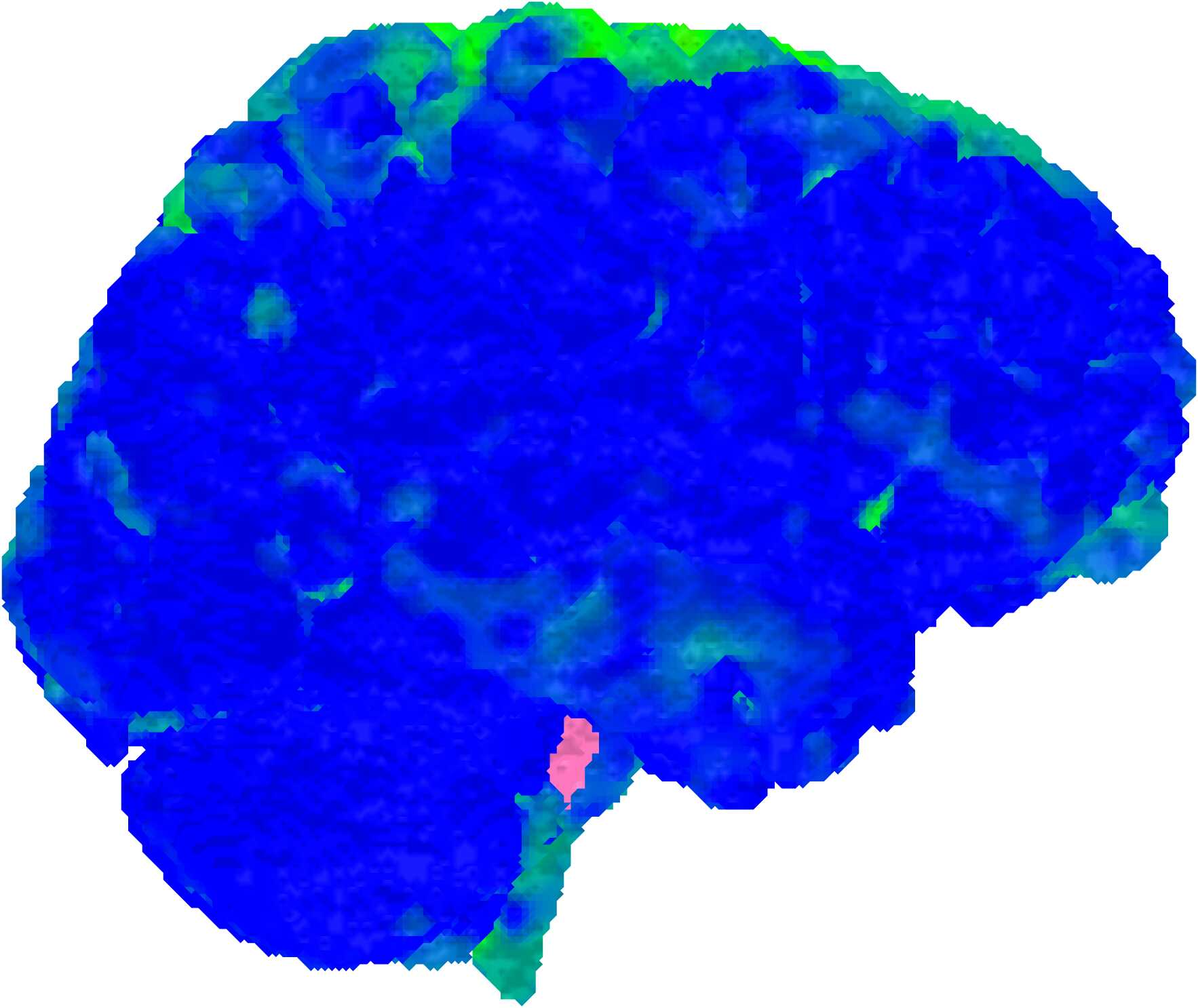}
            \caption{$\azimuth=\qty{90}{\degree}$, $\elevation=\qty{0}{\degree}$.}
        \end{subfigure}

        \def\subFigWidth{0.45\linewidth}

        \begin{subfigure}[b]{\subFigWidth}
            \includegraphics[width=\linewidth]{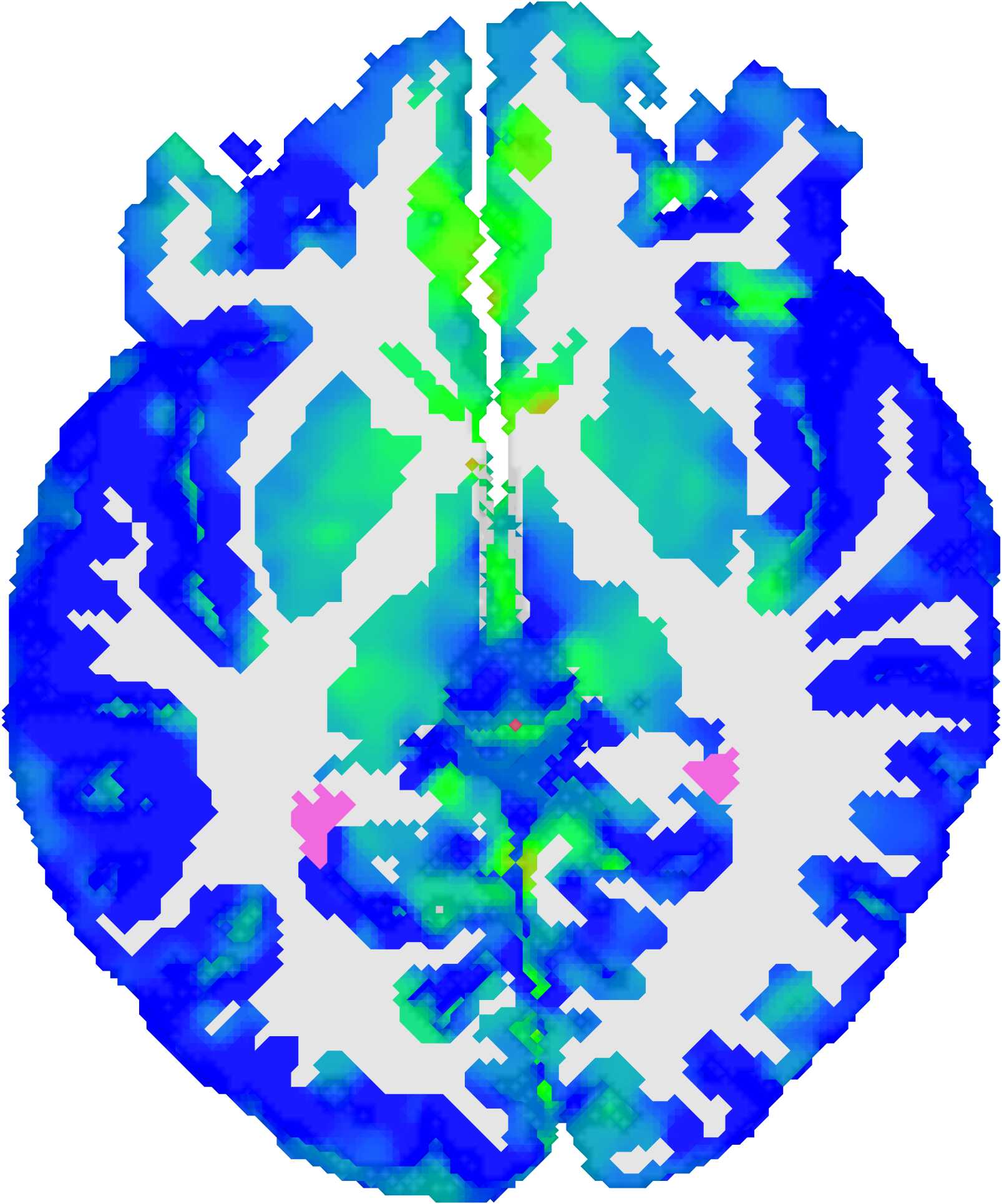}
            \caption{$\azimuth=\qty{0}{\degree}$, $\elevation=\qty{90}{\degree}$, $z=\qty{30}{\milli\meter}$.}
        \end{subfigure}
       \hfill
        \begin{subfigure}[b]{\subFigWidth}
            \includegraphics[width=\linewidth]{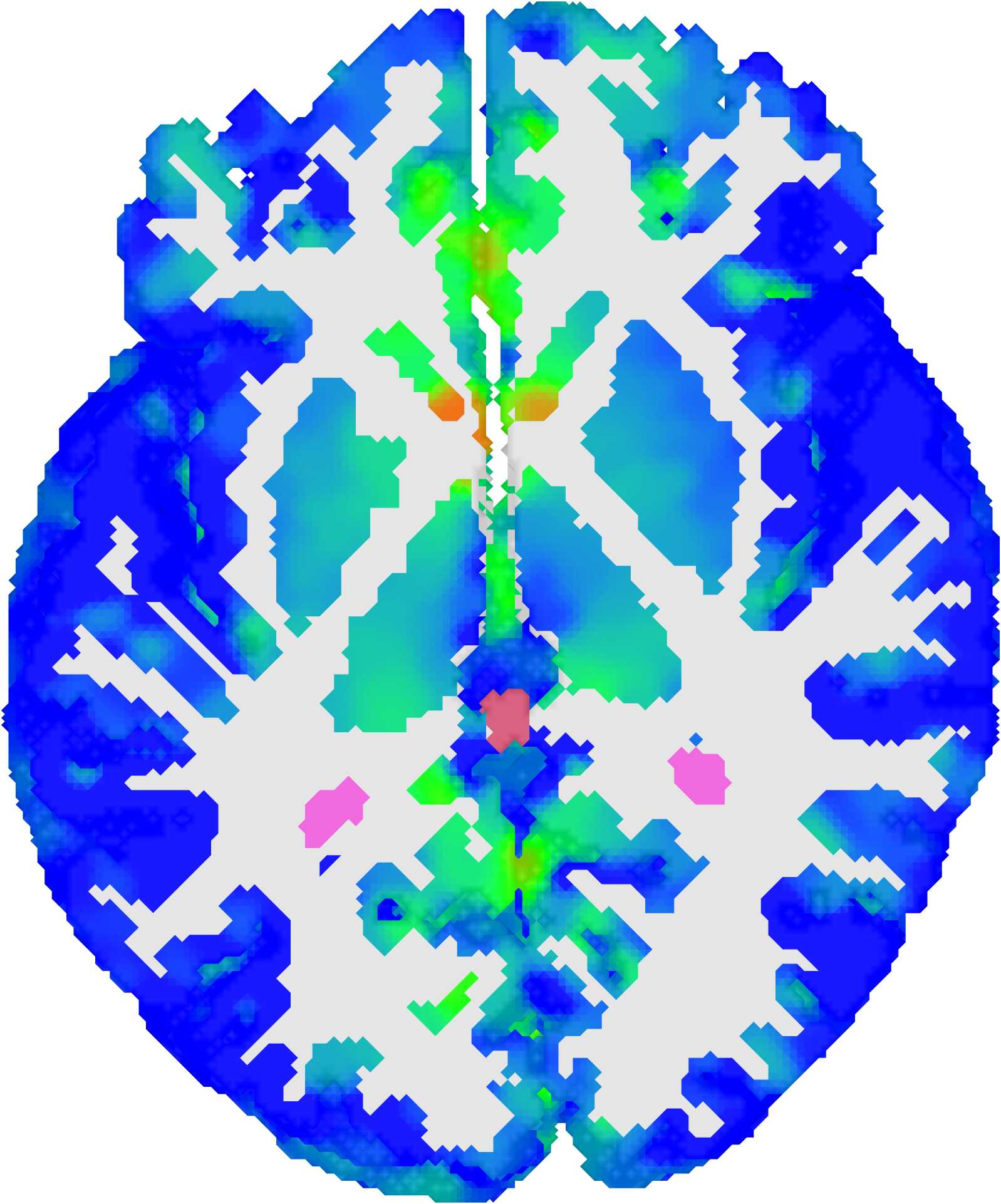}
            \caption{$\azimuth=\qty{0}{\degree}$, $\elevation=\qty{90}{\degree}$, $z=\qty{33}{\milli\meter}$.}
        \end{subfigure}

        \begin{subfigure}[b]{\subFigWidth}
            \includegraphics[width=\linewidth]{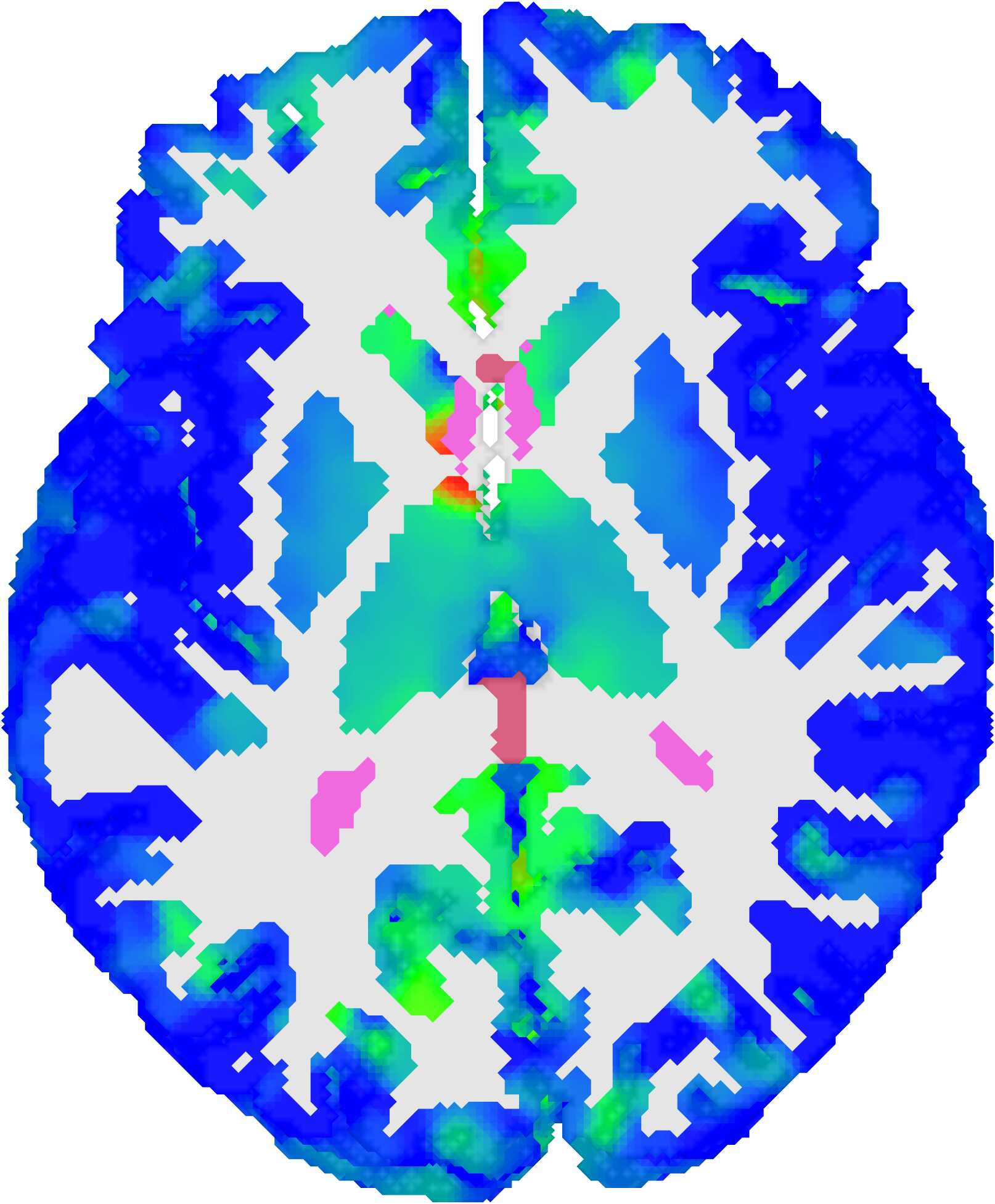}
            \caption{$\azimuth=\qty{0}{\degree}$, $\elevation=\qty{90}{\degree}$, $z=\qty{37}{\milli\meter}$.}
        \end{subfigure}
        \hfill
        \begin{subfigure}[b]{\subFigWidth}
            \includegraphics[width=\linewidth]{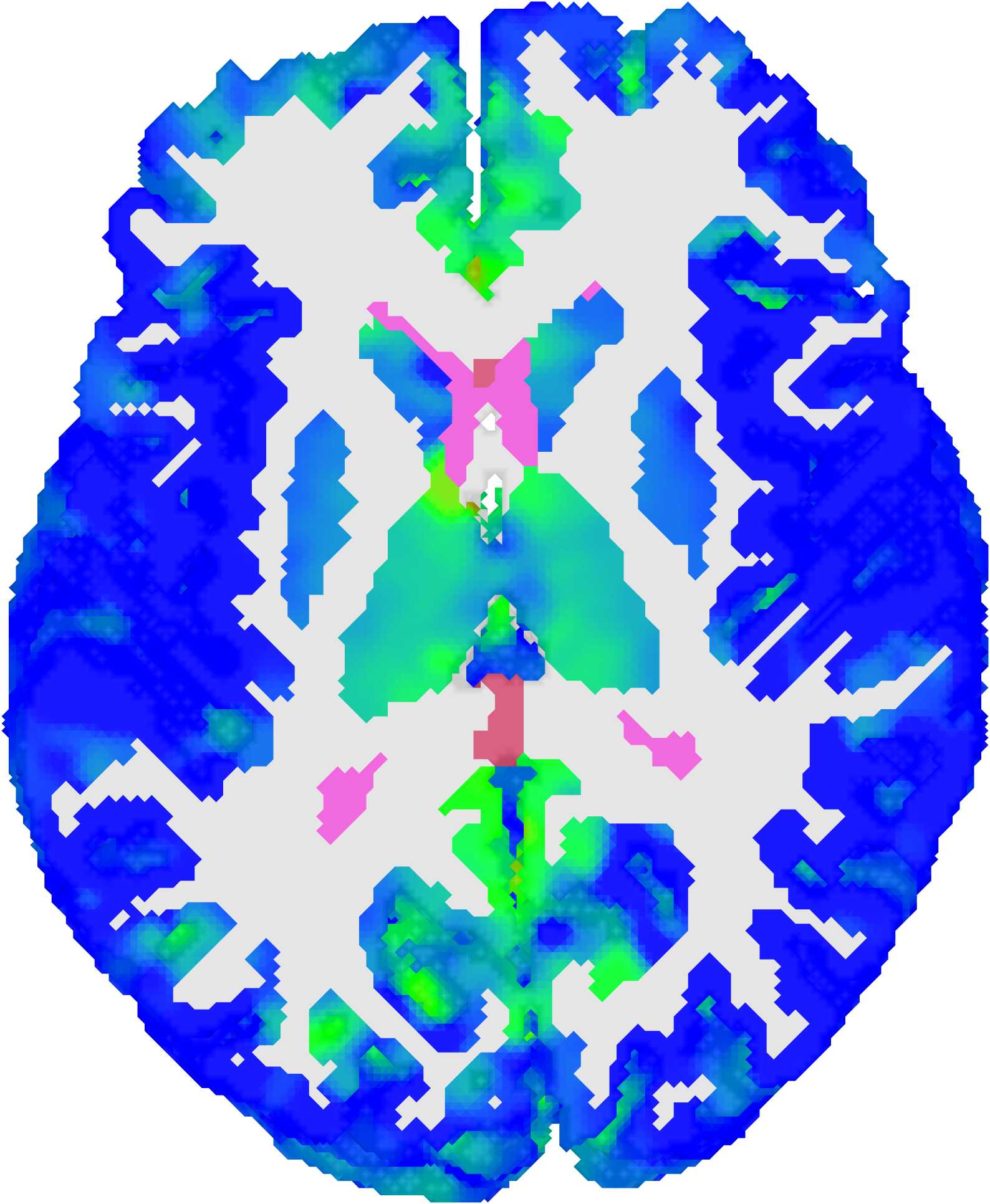}
            \caption{$\azimuth=\qty{0}{\degree}$, $\elevation=\qty{90}{\degree}$, $z=\qty{40}{\milli\meter}$.}
        \end{subfigure}
    \end{minipage}
    \hfill
    \begin{minipage}{0.09\linewidth}
    \includegraphics[width=\linewidth]{c,f=1000,li=0.5,ri=1.5,dB=18,az=0,el=90.colorbar.pdf}
    \end{minipage}
    \caption{Volume current interference field \eqref{eq:interference-field} of $\refJv$, smoothed via the function $\dBT\args{\someField/\max\someField,\qty{18}{\decibel}}$ of equation \eqref{eq:visual-log-transform} at different azimuths $\azimuth$ and elevations $\elevation$, with base stimulation frequency $\frequency=\qty{1000}{\hertz}$ and an updated $\contactResistance=\qty{1270}{\ohm}$. The current pattern was $\currentPattern=[-1.0,1.0,-1.0,1.0]\,\si{\milli\ampere}$ and the maximal magnitude of the output $\refJv$ was approximately \qty{6.5109e-05}{\milli\ampere\per\milli\meter\squared}.}
    \label{fig:interference-reffield-equal-currents-1270}
\end{figure}

\begin{figure}[!h]
    \centering
    \begin{minipage}{0.85\linewidth}
        \def\subFigWidth{0.31\linewidth}
        \begin{subfigure}[b]{\subFigWidth}
            \includegraphics[width=\linewidth]{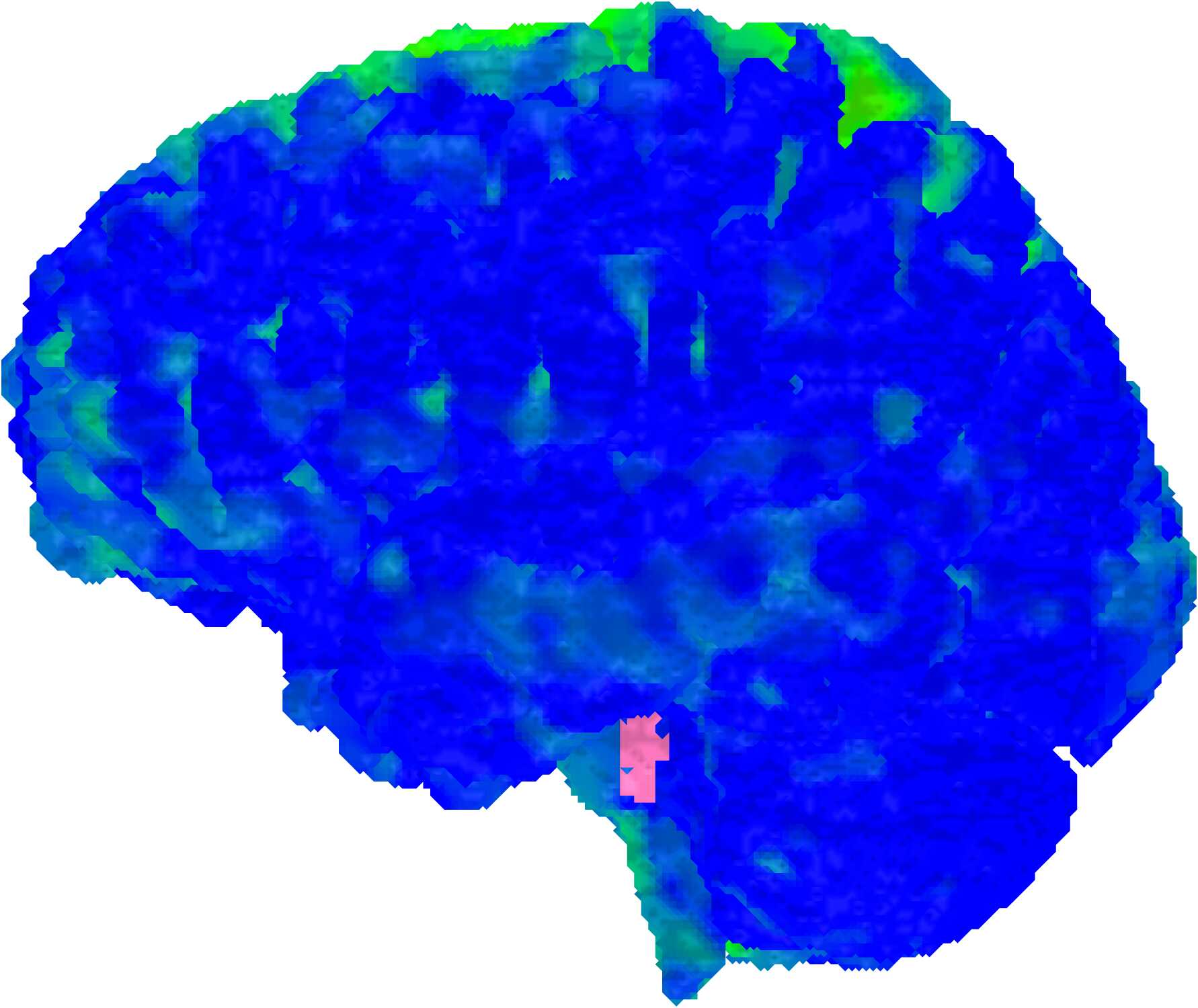}
            \caption{$\azimuth=\qty{-90}{\degree}$, $\elevation=\qty{0}{\degree}$.}
        \end{subfigure}
        \hfill
        \begin{subfigure}[b]{\subFigWidth}
        \includegraphics[width=\linewidth]{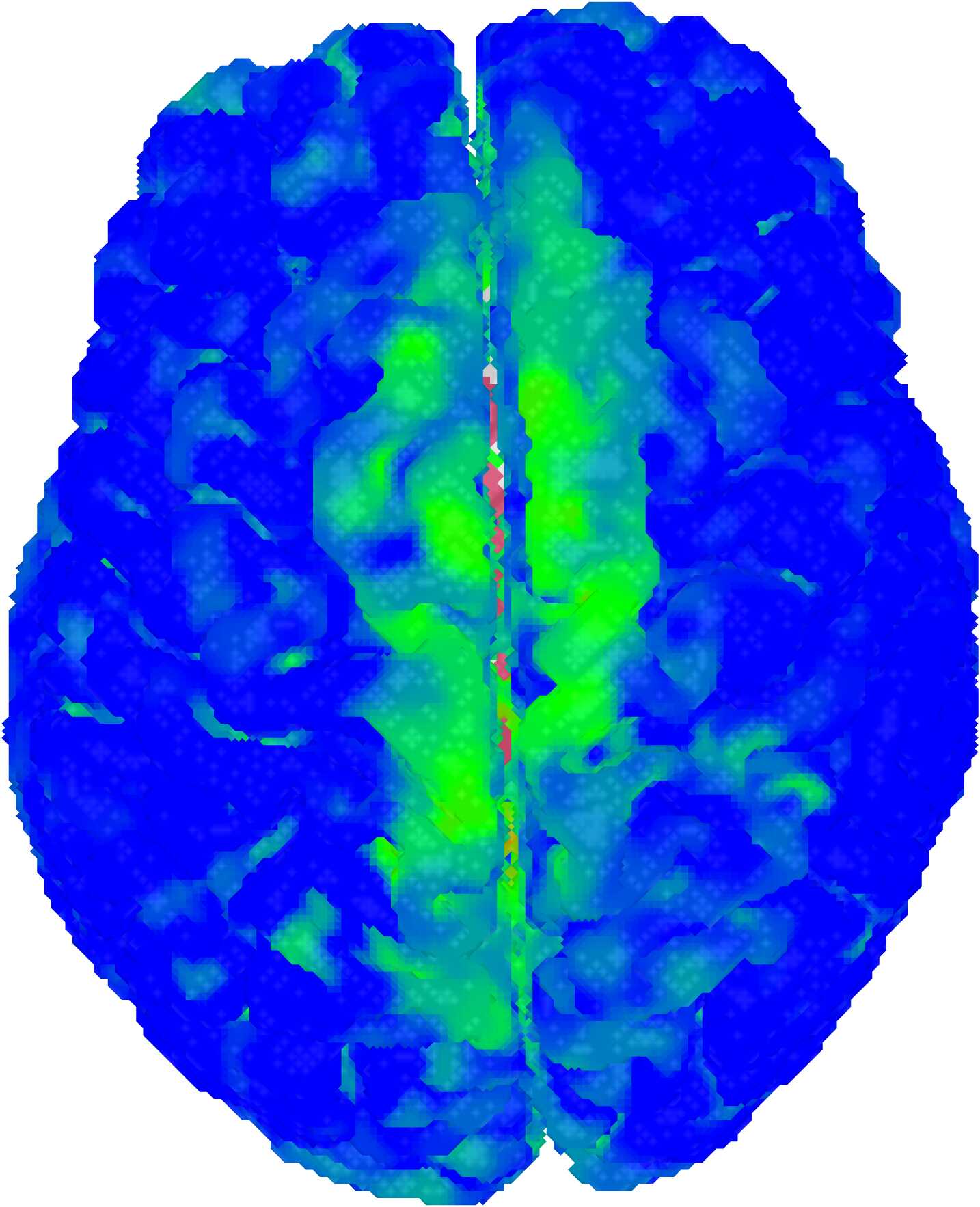}
        \caption{$\azimuth=\qty{0}{\degree}$, $\elevation=\qty{90}{\degree}$.}
        \end{subfigure}
        \hfill
        \begin{subfigure}[b]{\subFigWidth}
        \includegraphics[width=\linewidth]{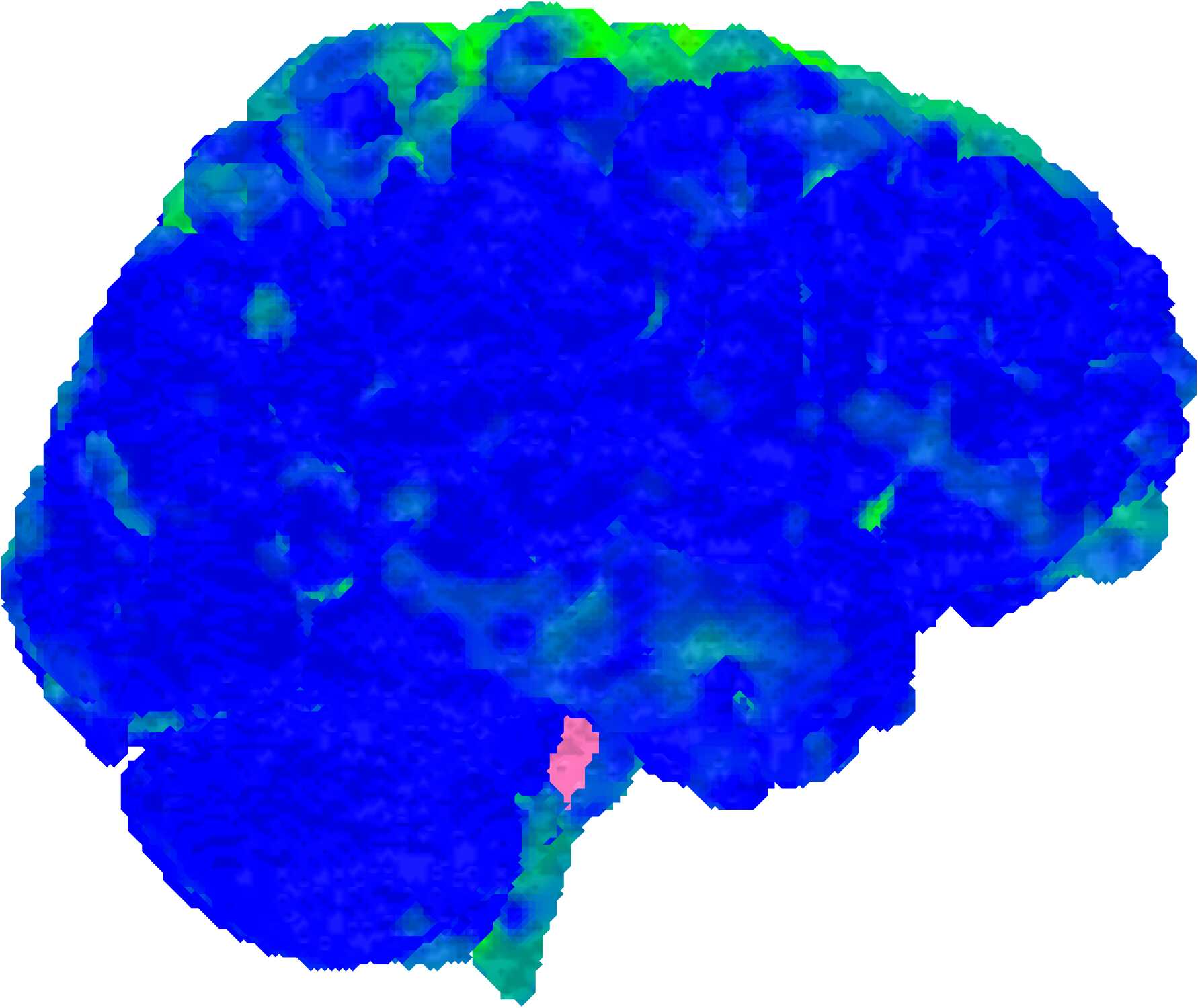}
            \caption{$\azimuth=\qty{90}{\degree}$, $\elevation=\qty{0}{\degree}$.}
        \end{subfigure}

        \def\subFigWidth{0.45\linewidth}

        \begin{subfigure}[b]{\subFigWidth}
            \includegraphics[width=\linewidth]{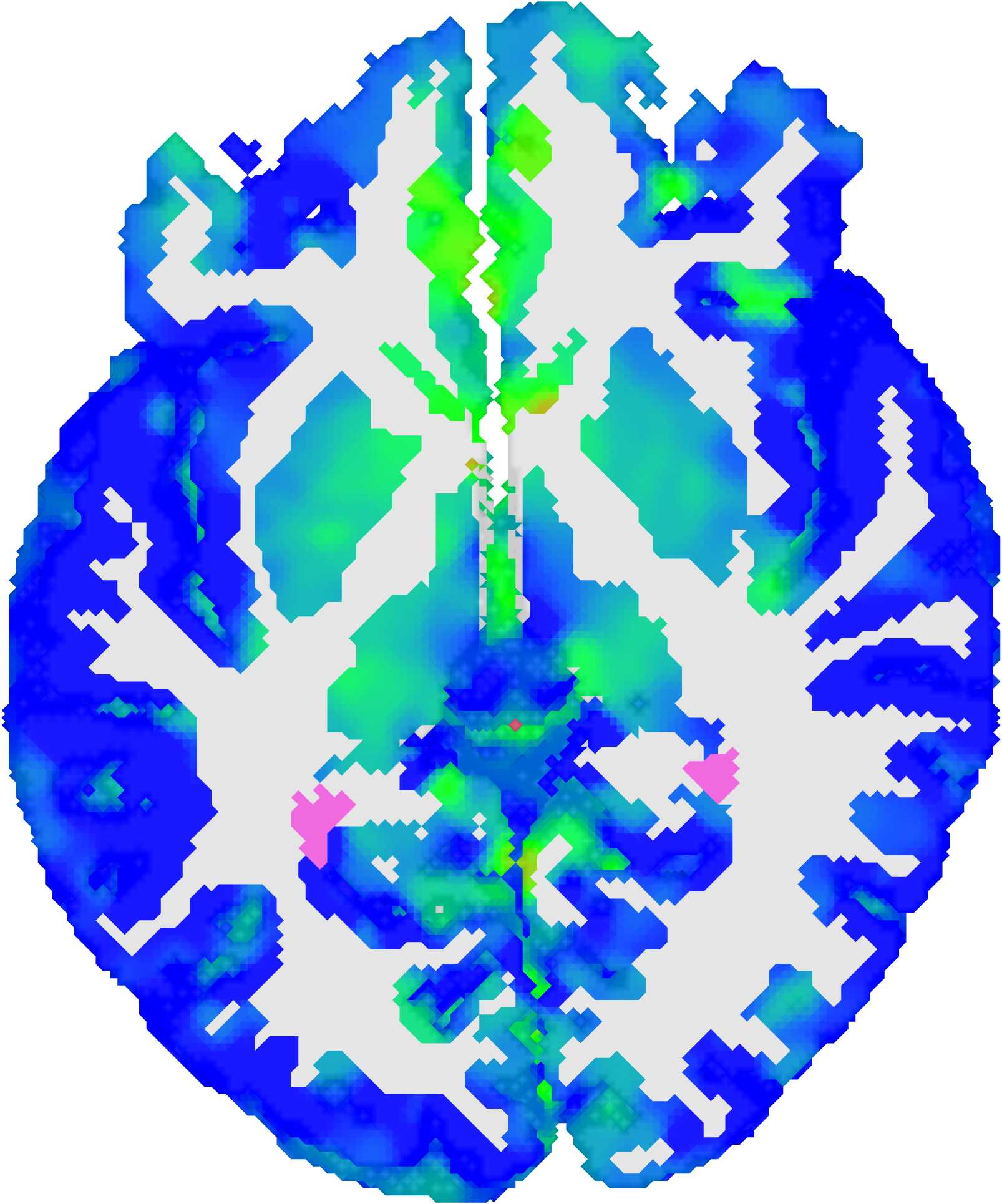}
            \caption{$\azimuth=\qty{0}{\degree}$, $\elevation=\qty{90}{\degree}$, $z=\qty{30}{\milli\meter}$.}
        \end{subfigure}
       \hfill
        \begin{subfigure}[b]{\subFigWidth}
            \includegraphics[width=\linewidth]{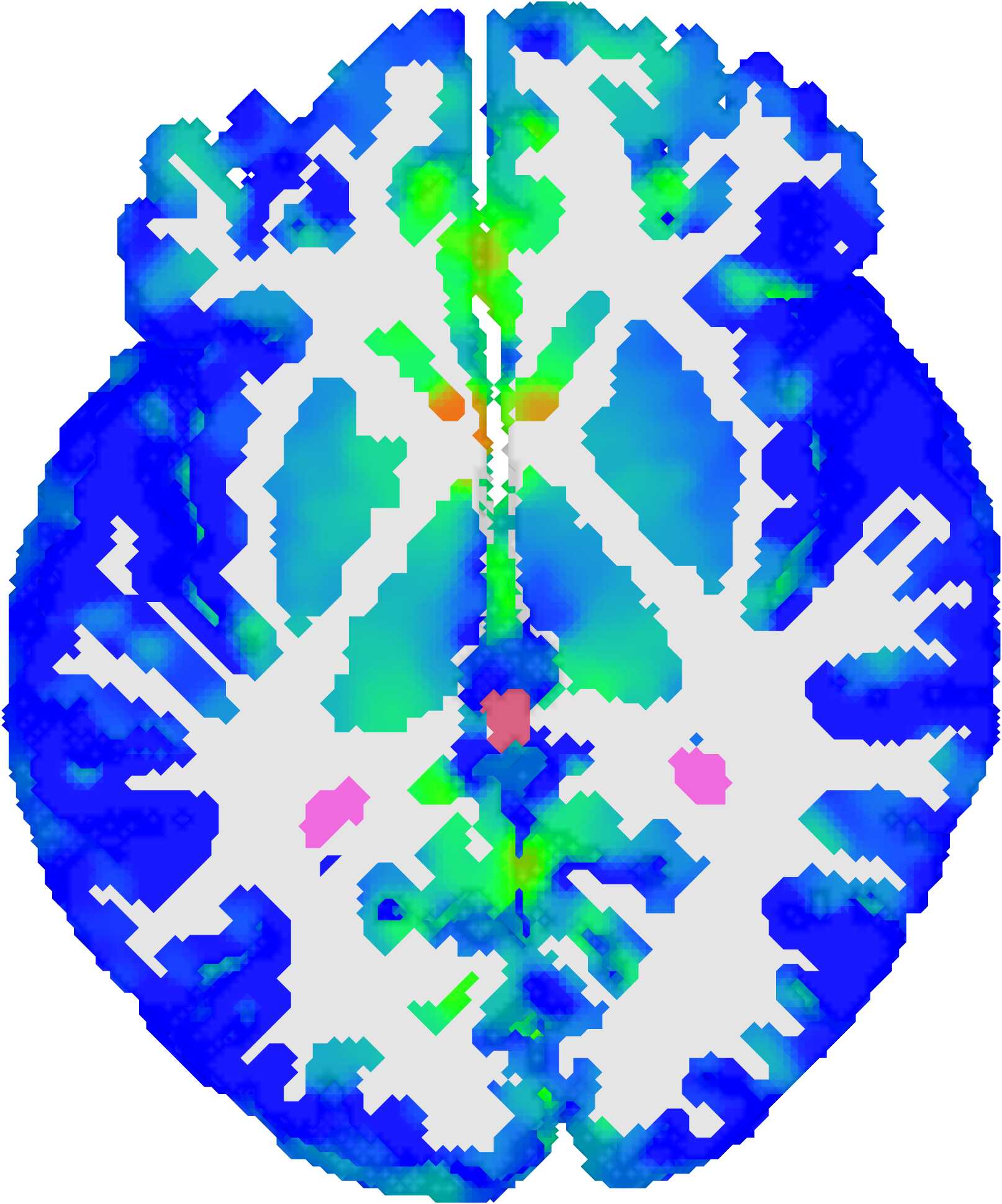}
            \caption{$\azimuth=\qty{0}{\degree}$, $\elevation=\qty{90}{\degree}$, $z=\qty{33}{\milli\meter}$.}
        \end{subfigure}

        \begin{subfigure}[b]{\subFigWidth}
            \includegraphics[width=\linewidth]{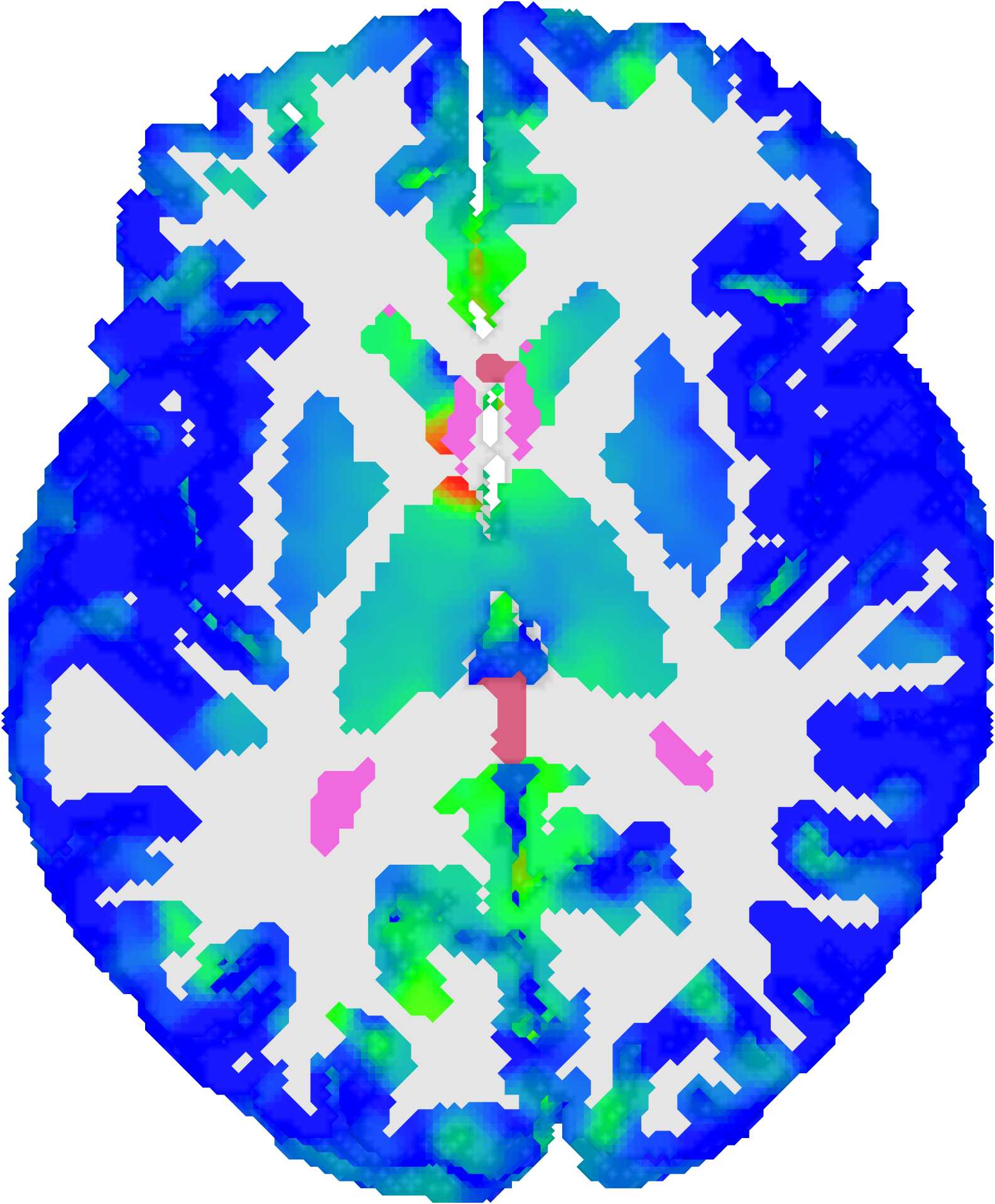}
            \caption{$\azimuth=\qty{0}{\degree}$, $\elevation=\qty{90}{\degree}$, $z=\qty{37}{\milli\meter}$.}
        \end{subfigure}
        \hfill
        \begin{subfigure}[b]{\subFigWidth}
            \includegraphics[width=\linewidth]{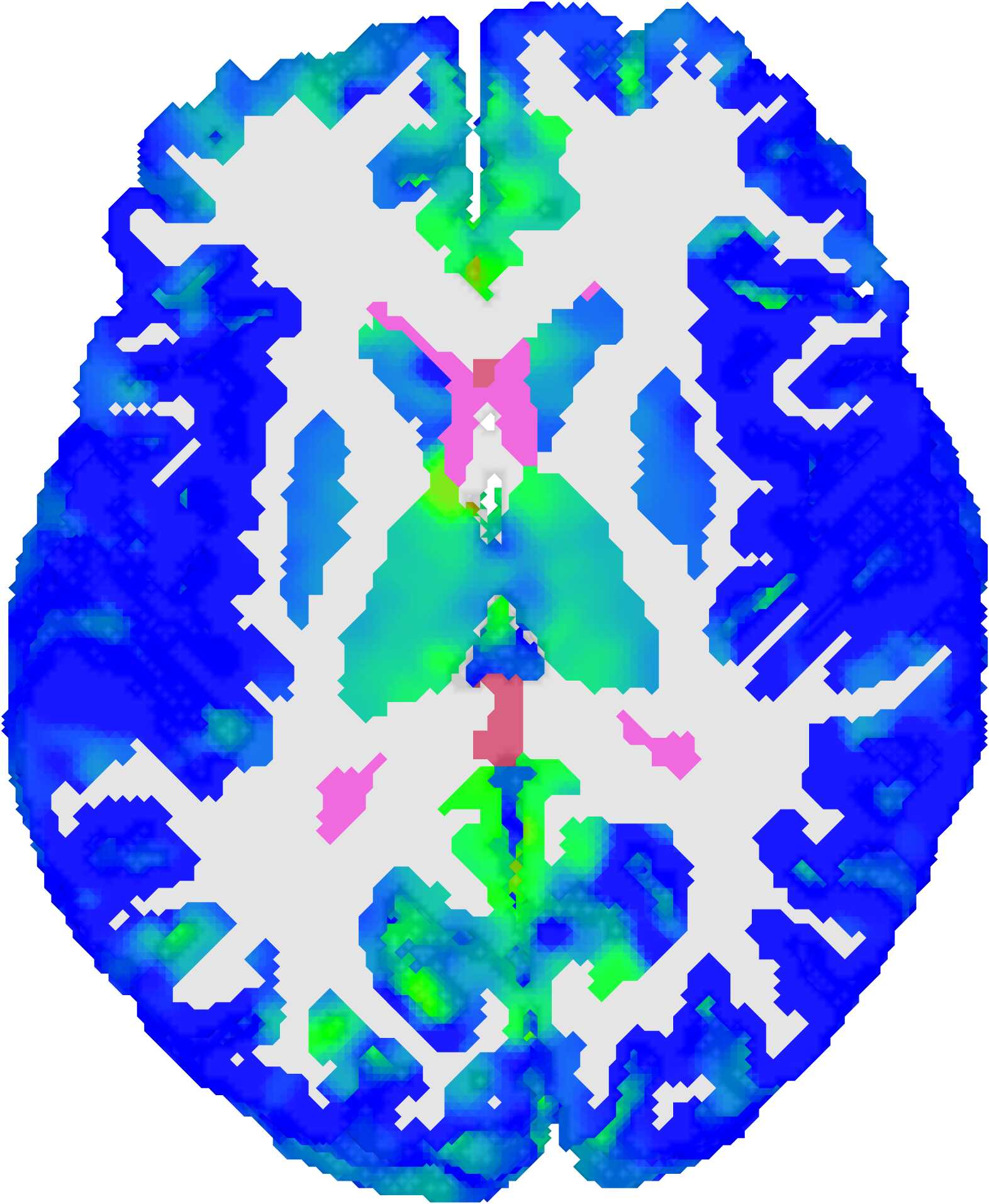}
            \caption{$\azimuth=\qty{0}{\degree}$, $\elevation=\qty{90}{\degree}$, $z=\qty{40}{\milli\meter}$.}
        \end{subfigure}
    \end{minipage}
    \hfill
    \begin{minipage}{0.09\linewidth}
    \includegraphics[width=\linewidth]{c,f=1000,li=0.5,ri=1.5,dB=18,az=0,el=90.colorbar.pdf}
    \end{minipage}
    \caption{Volume current interference field \eqref{eq:interference-field} of $\refJv$, smoothed via the function $\dBT\args{\someField/\max\someField,\qty{18}{\decibel}}$ of equation \eqref{eq:visual-log-transform} at different azimuths $\azimuth$ and elevations $\elevation$, with base stimulation frequency $\frequency=\qty{1000}{\hertz}$ and an updated $\contactResistance=\qty{5270}{\ohm}$. The current pattern was $\currentPattern=[-1.0,1.0,-1.0,1.0]\,\si{\milli\ampere}$ and the maximal magnitude of the output $\refJv$ was approximately \qty{6.4419e-05}{\milli\ampere\per\milli\meter\squared}.}
    \label{fig:interference-reffield-equal-currents-5270}
\end{figure}

\begin{figure}[!h]
    \centering
    \begin{minipage}{0.85\linewidth}
        \def\subFigWidth{0.31\linewidth}
        \begin{subfigure}[b]{\subFigWidth}
            \includegraphics[width=\linewidth]{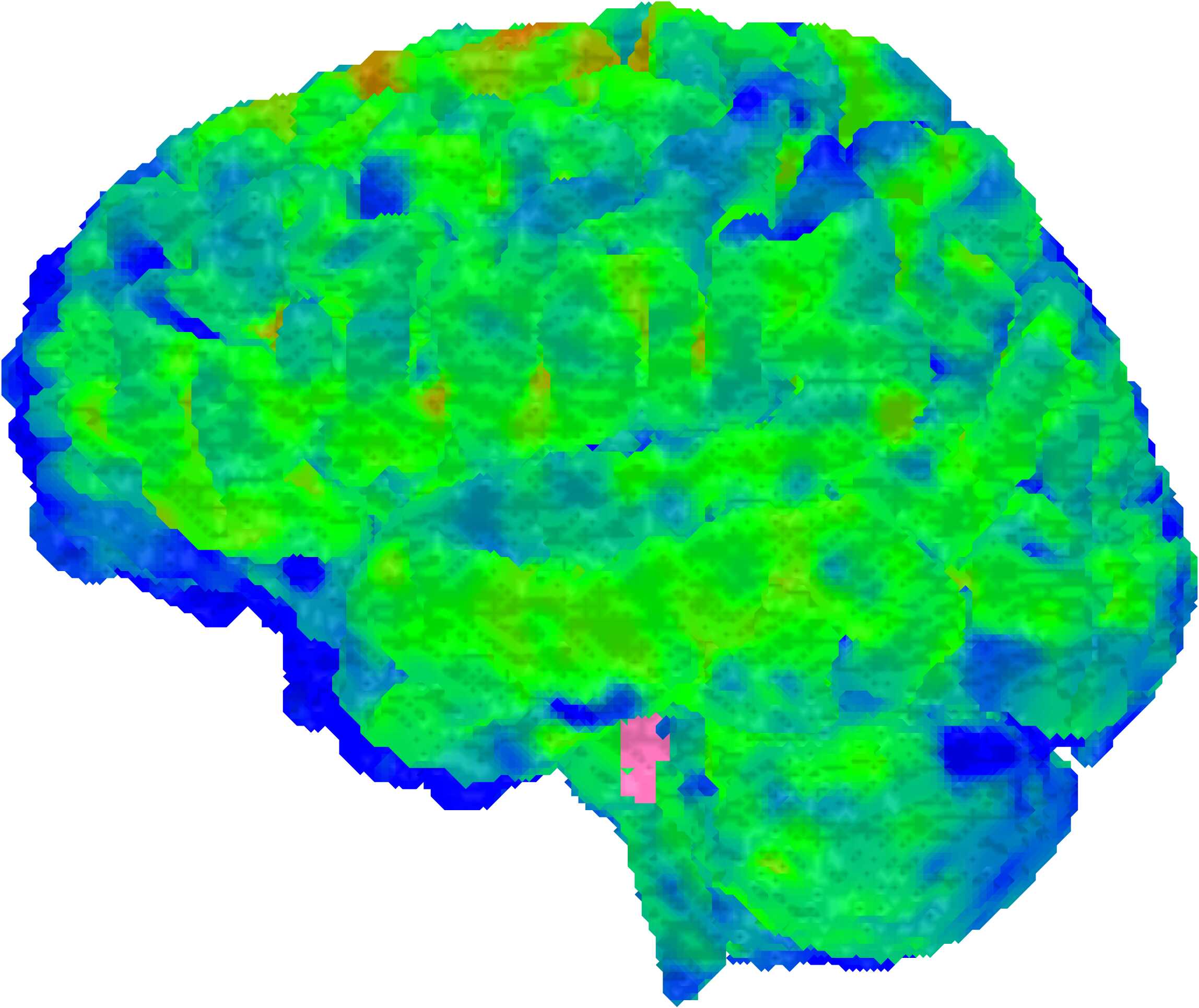}
            \caption{$\azimuth=\qty{-90}{\degree}$, $\elevation=\qty{0}{\degree}$.}
        \end{subfigure}
        \hfill
        \begin{subfigure}[b]{\subFigWidth}
        \includegraphics[width=\linewidth]{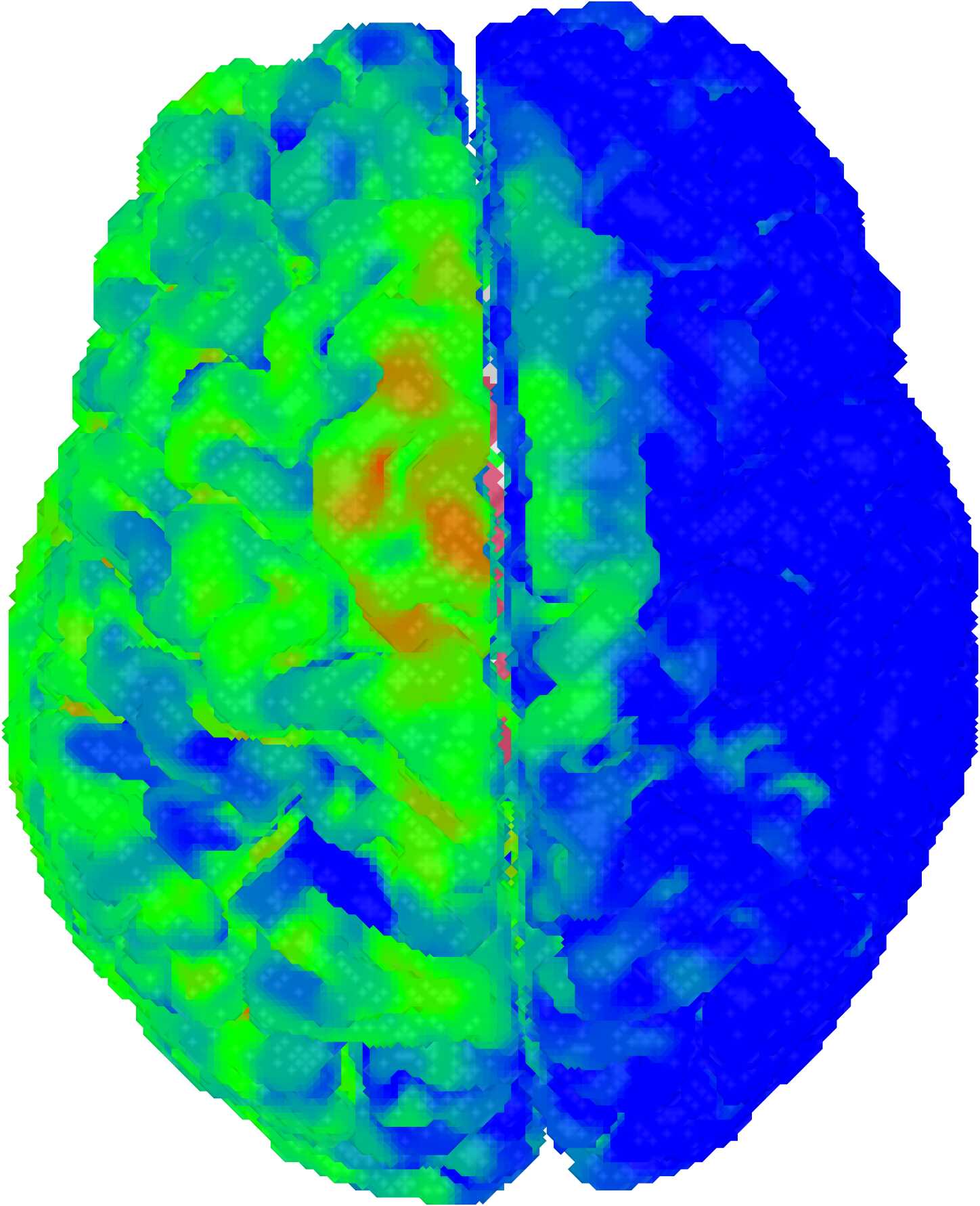}
        \caption{$\azimuth=\qty{0}{\degree}$, $\elevation=\qty{90}{\degree}$.}
        \end{subfigure}
        \hfill
        \begin{subfigure}[b]{\subFigWidth}
        \includegraphics[width=\linewidth]{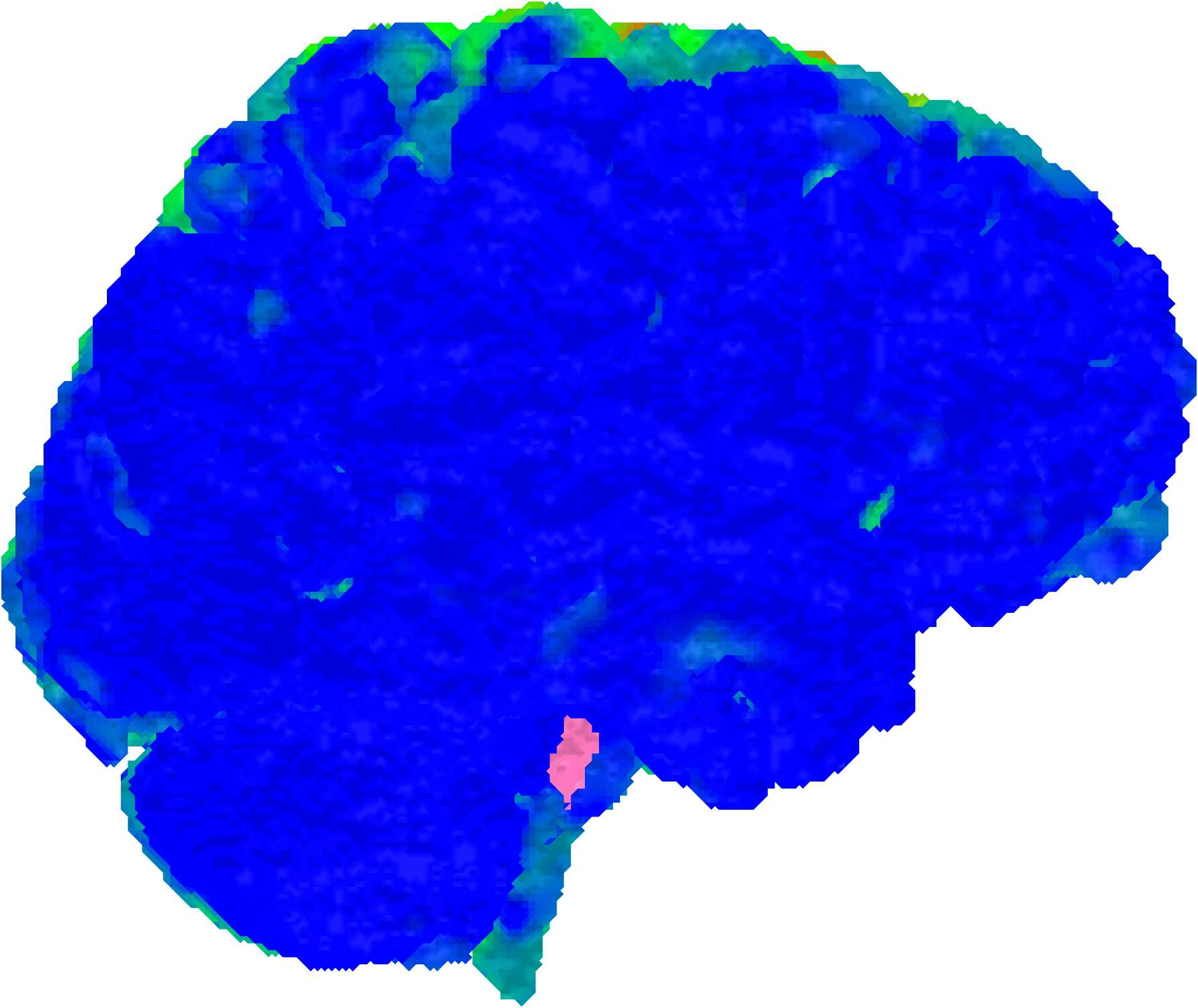}
            \caption{$\azimuth=\qty{90}{\degree}$, $\elevation=\qty{0}{\degree}$.}
        \end{subfigure}

        \def\subFigWidth{0.45\linewidth}

        \begin{subfigure}[b]{\subFigWidth}
            \includegraphics[width=\linewidth]{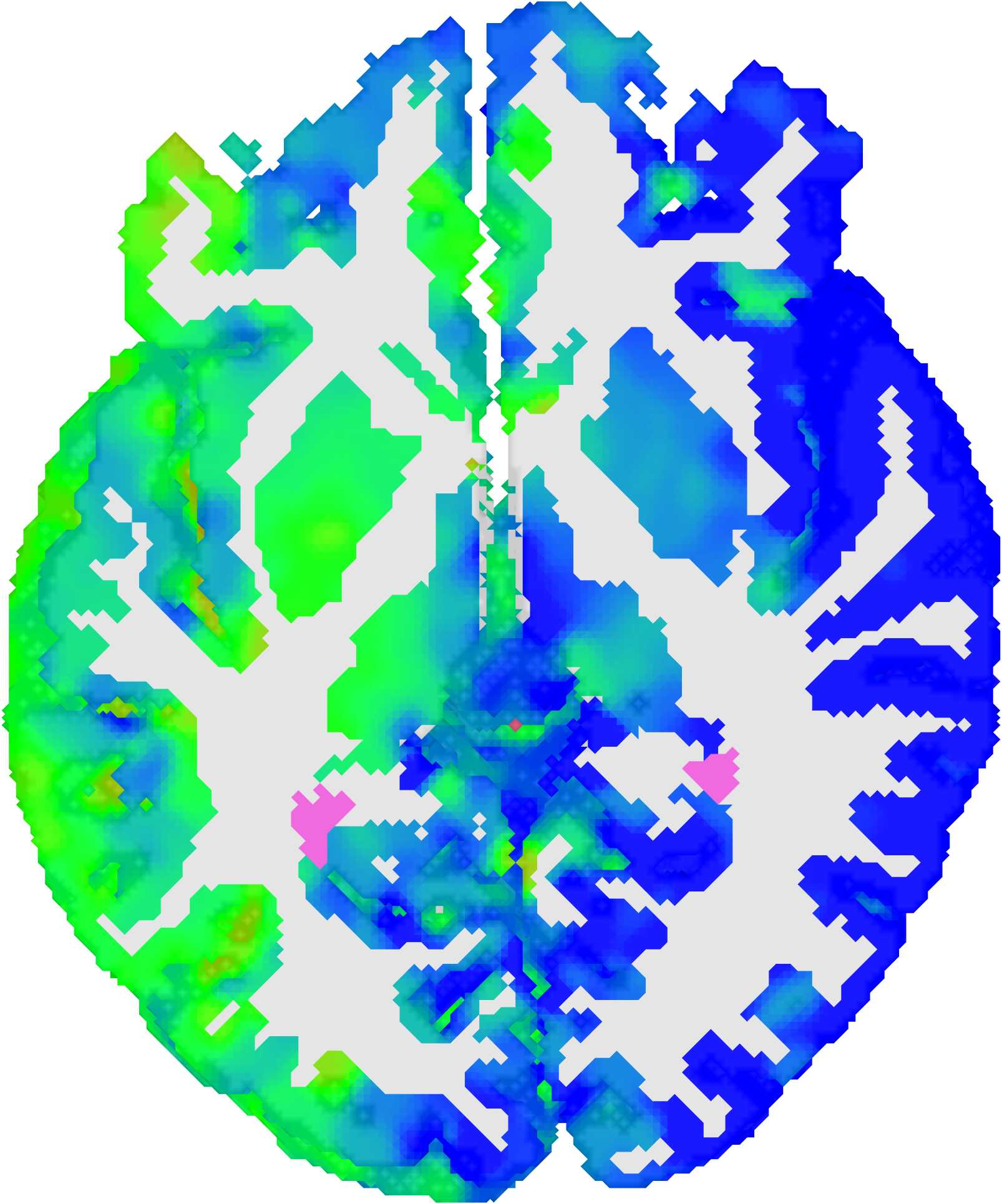}
            \caption{$\azimuth=\qty{0}{\degree}$, $\elevation=\qty{90}{\degree}$, $z=\qty{30}{\milli\meter}$.}
        \end{subfigure}
       \hfill
        \begin{subfigure}[b]{\subFigWidth}
            \includegraphics[width=\linewidth]{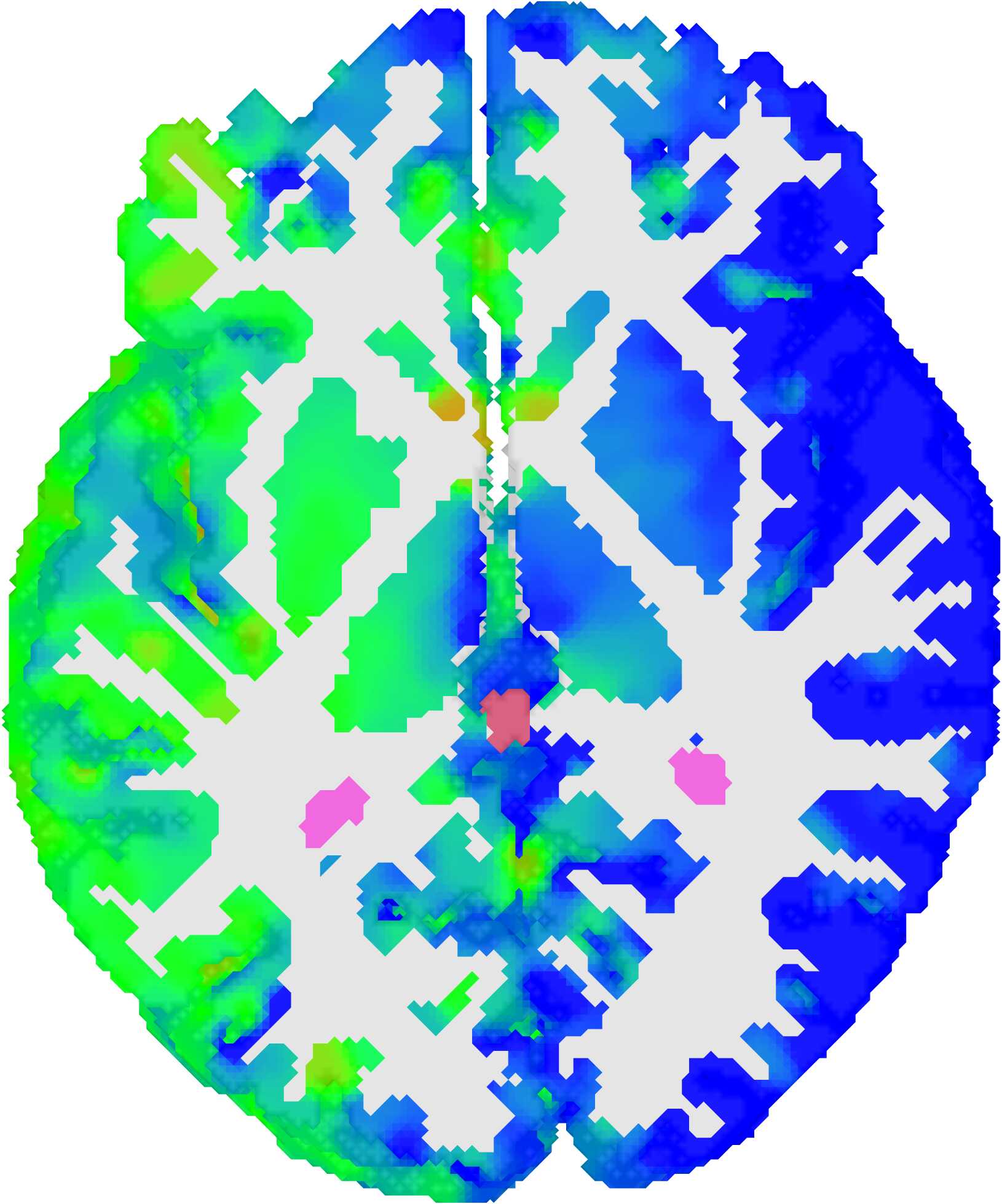}
            \caption{$\azimuth=\qty{0}{\degree}$, $\elevation=\qty{90}{\degree}$, $z=\qty{33}{\milli\meter}$.}
        \end{subfigure}

        \begin{subfigure}[b]{\subFigWidth}
            \includegraphics[width=\linewidth]{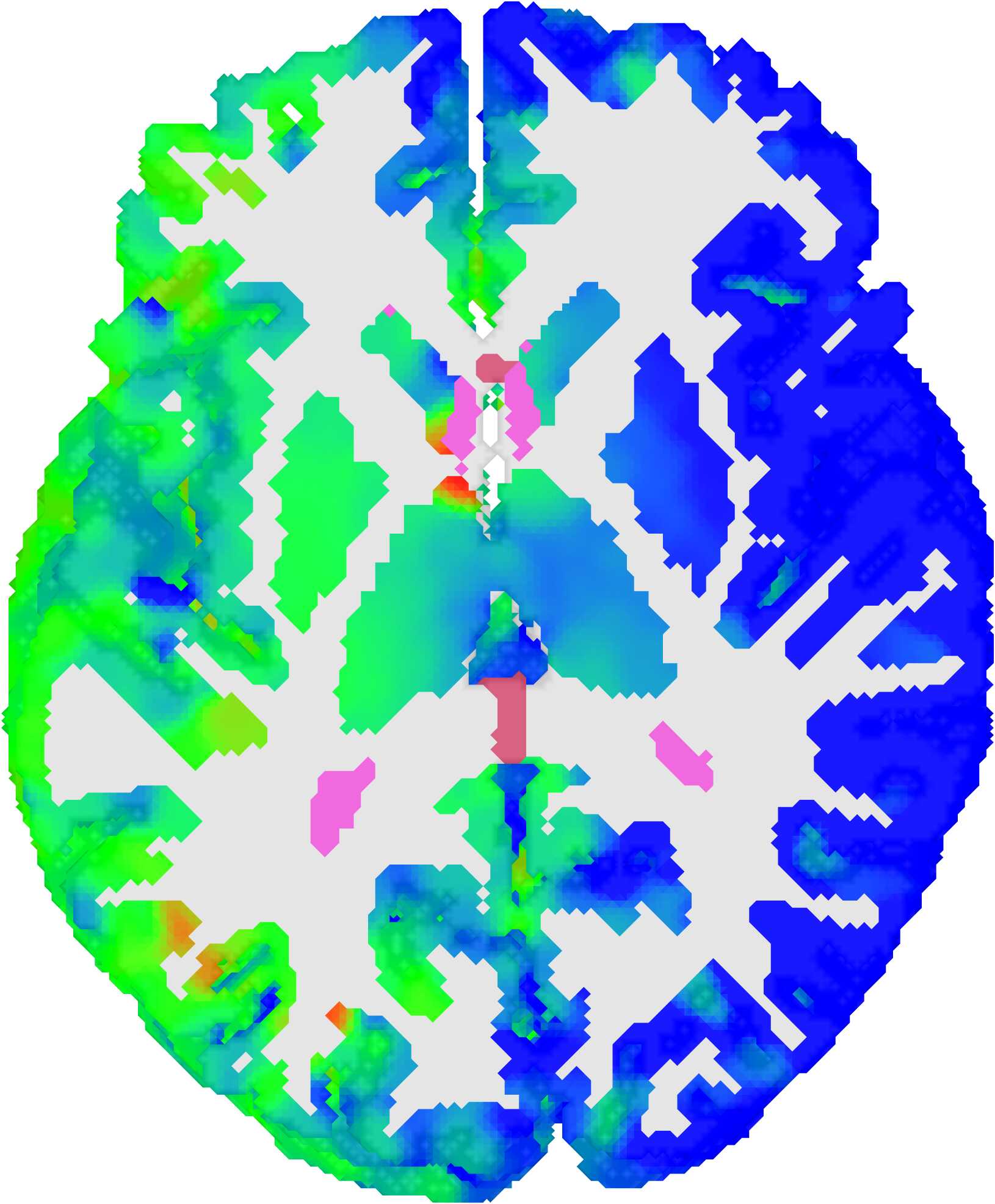}
            \caption{$\azimuth=\qty{0}{\degree}$, $\elevation=\qty{90}{\degree}$, $z=\qty{37}{\milli\meter}$.}
        \end{subfigure}
        \hfill
        \begin{subfigure}[b]{\subFigWidth}
            \includegraphics[width=\linewidth]{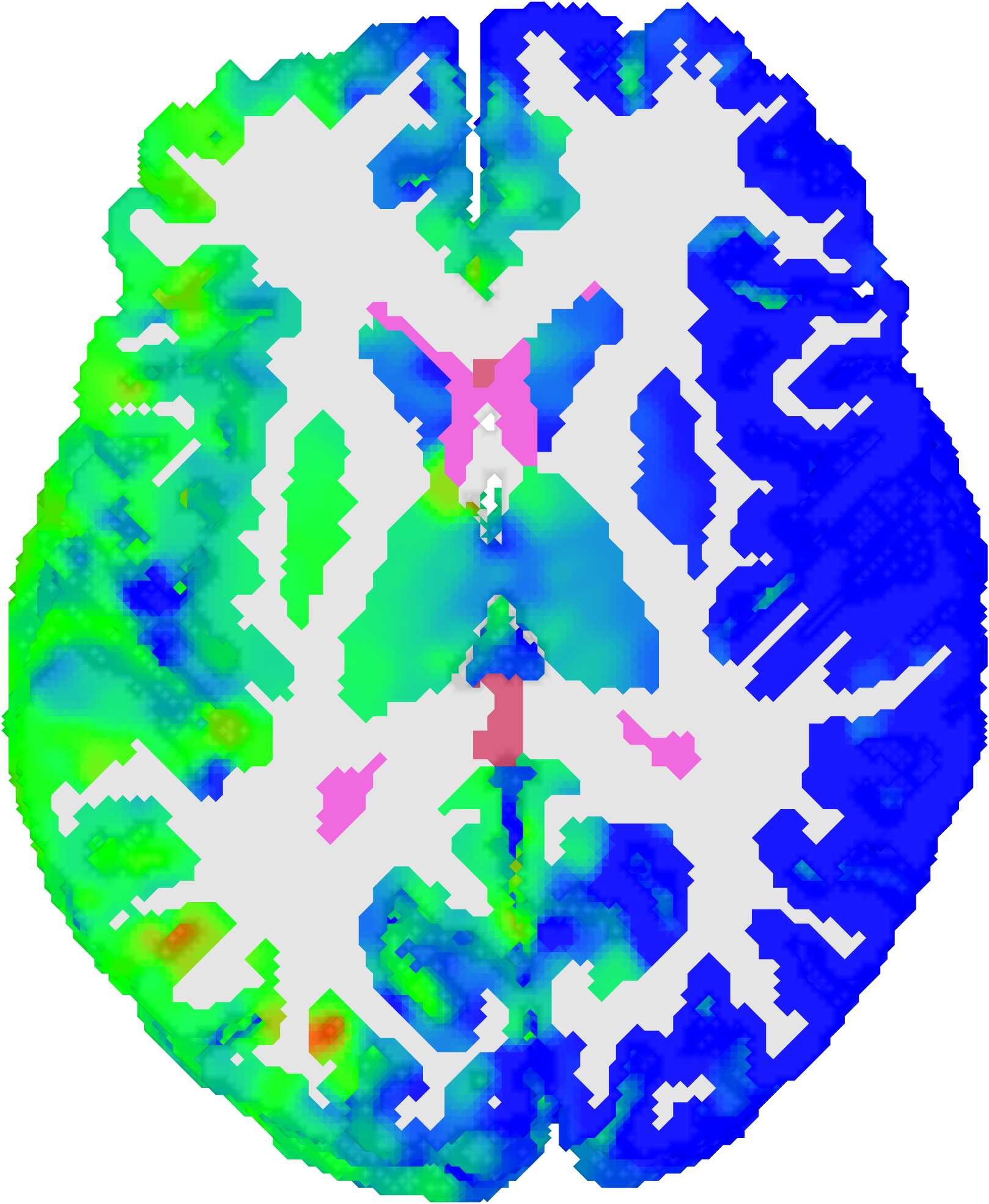}
            \caption{$\azimuth=\qty{0}{\degree}$, $\elevation=\qty{90}{\degree}$, $z=\qty{40}{\milli\meter}$.}
        \end{subfigure}
    \end{minipage}
    \hfill
    \begin{minipage}{0.09\linewidth}
    \includegraphics[width=\linewidth]{c,f=1000,li=0.5,ri=1.5,dB=18,az=0,el=90.colorbar.pdf}
    \end{minipage}
    \caption{Volume current interference field \eqref{eq:interference-field} of $\refJv$, smoothed via the function $\dBT\args{\someField/\max\someField,\qty{18}{\decibel}}$ of equation \eqref{eq:visual-log-transform} at different azimuths $\azimuth$ and elevations $\elevation$, with base stimulation frequency $\frequency=\qty{1000}{\hertz}$ and an updated $\contactResistance=\qty{1270}{\ohm}$. The current pattern was $\currentPattern=[-0.5,0.5,-1.5,1.5]\,\si{\milli\ampere}$ and the maximal value of $\refJv$ roughly equal to \qty{4.0239e-05}{\milli\ampere\per\milli\meter\squared}.}
    \label{fig:interference-reffield-weaker-left-1270}
\end{figure}

\begin{figure}[!h]
    \centering
    \begin{minipage}{0.85\linewidth}
        \def\subFigWidth{0.31\linewidth}
        \begin{subfigure}[b]{\subFigWidth}
            \includegraphics[width=\linewidth]{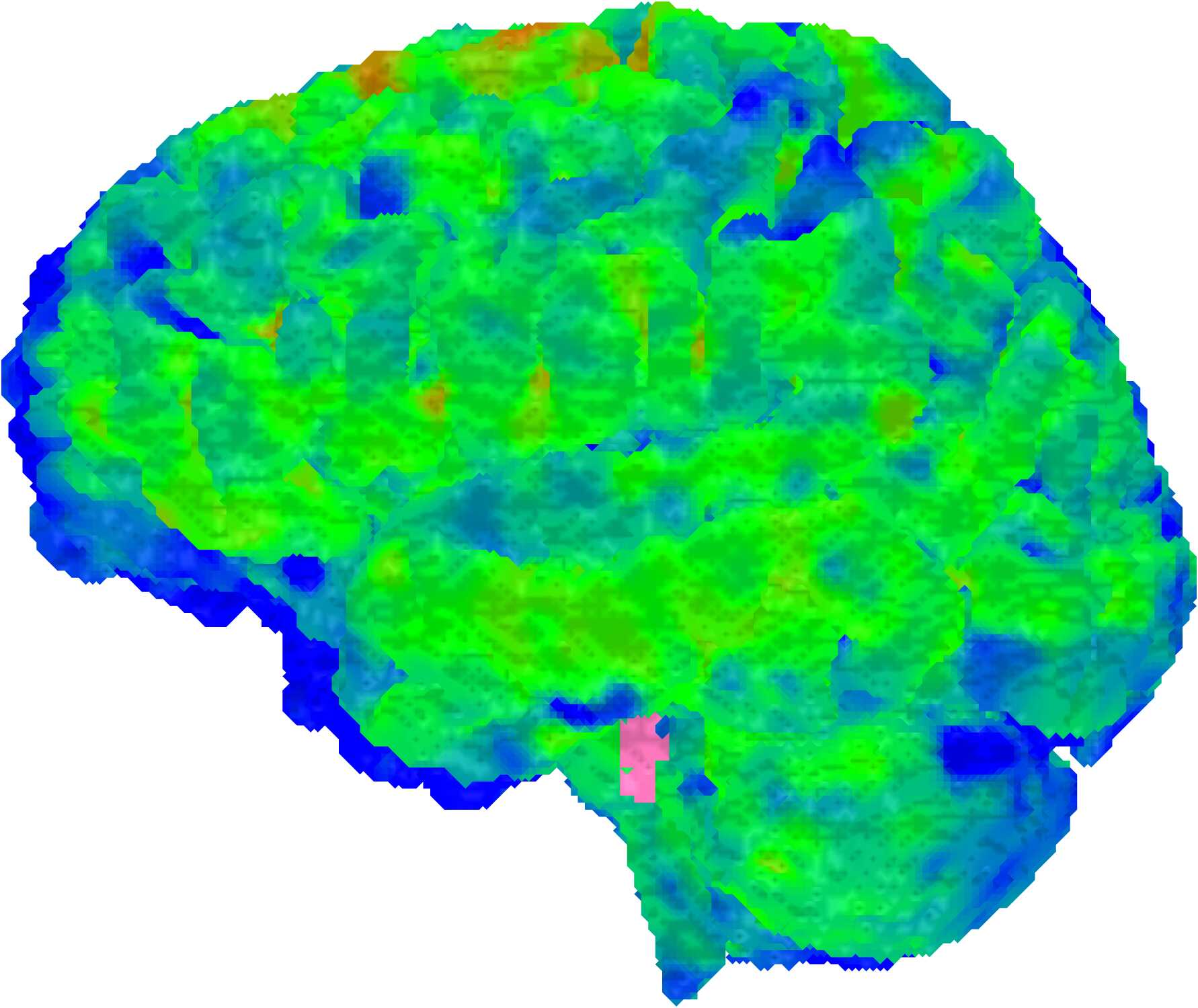}
            \caption{$\azimuth=\qty{-90}{\degree}$, $\elevation=\qty{0}{\degree}$.}
        \end{subfigure}
        \hfill
        \begin{subfigure}[b]{\subFigWidth}
        \includegraphics[width=\linewidth]{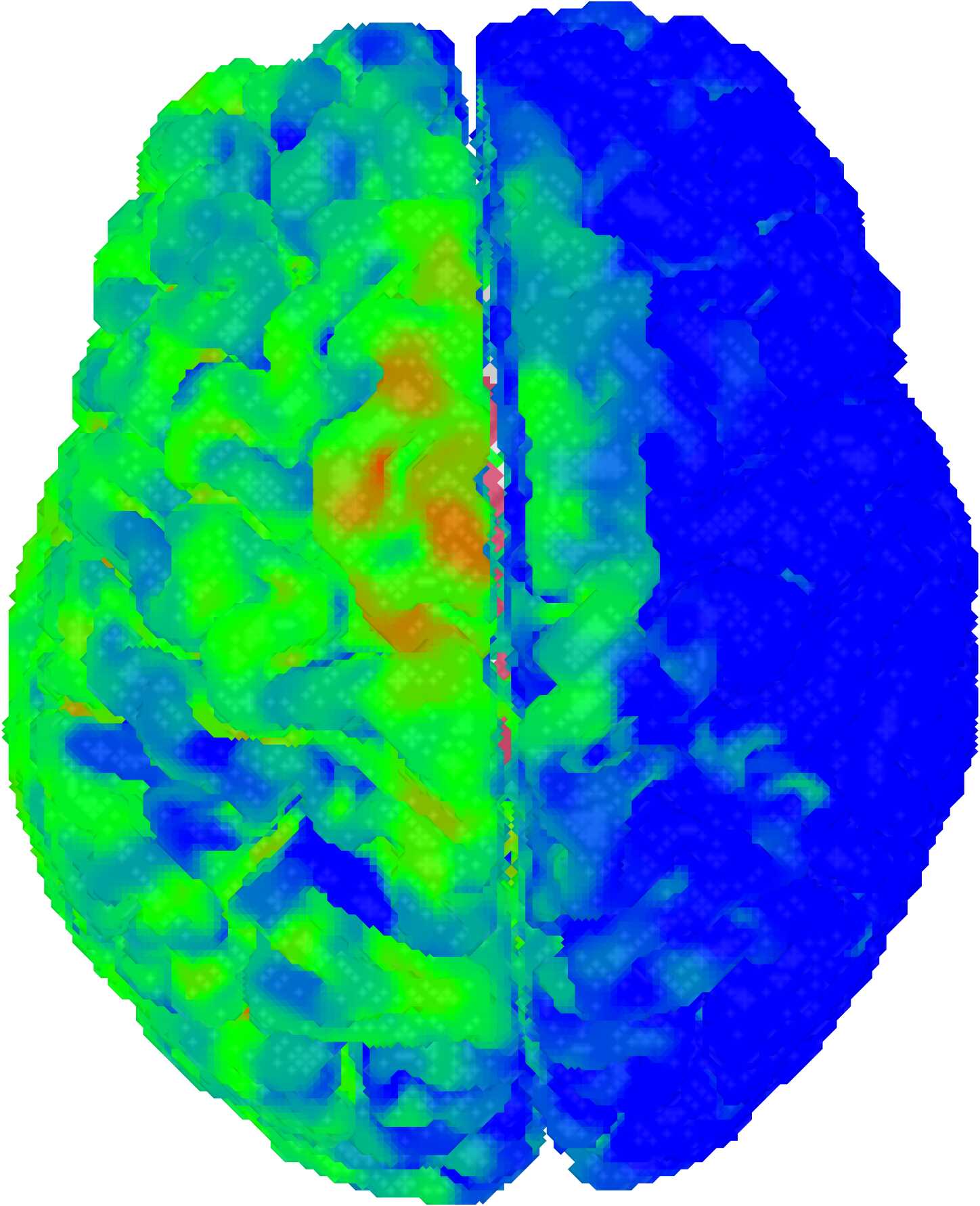}
        \caption{$\azimuth=\qty{0}{\degree}$, $\elevation=\qty{90}{\degree}$.}
        \end{subfigure}
        \hfill
        \begin{subfigure}[b]{\subFigWidth}
        \includegraphics[width=\linewidth]{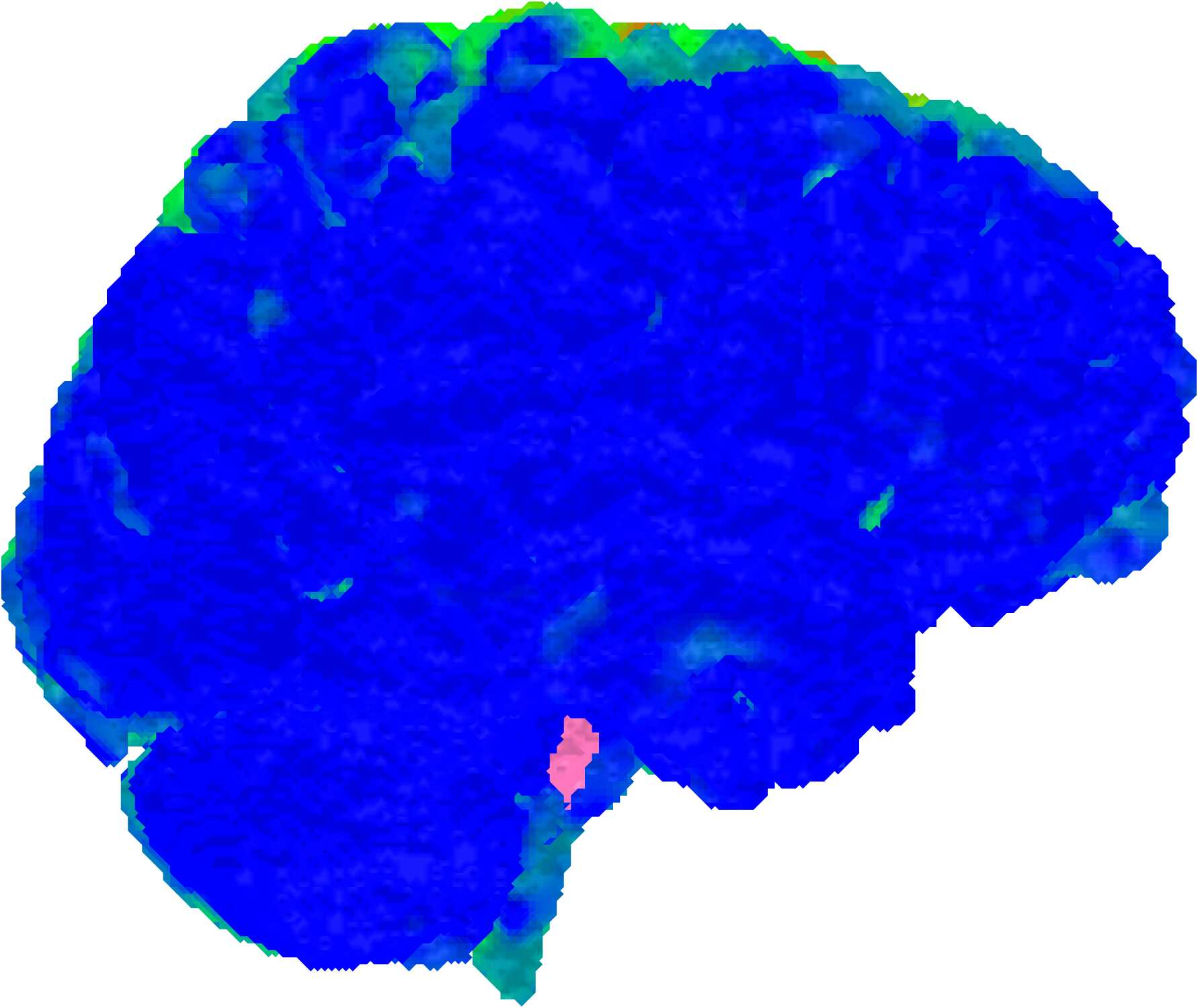}
            \caption{$\azimuth=\qty{90}{\degree}$, $\elevation=\qty{0}{\degree}$.}
        \end{subfigure}

        \def\subFigWidth{0.45\linewidth}

        \begin{subfigure}[b]{\subFigWidth}
            \includegraphics[width=\linewidth]{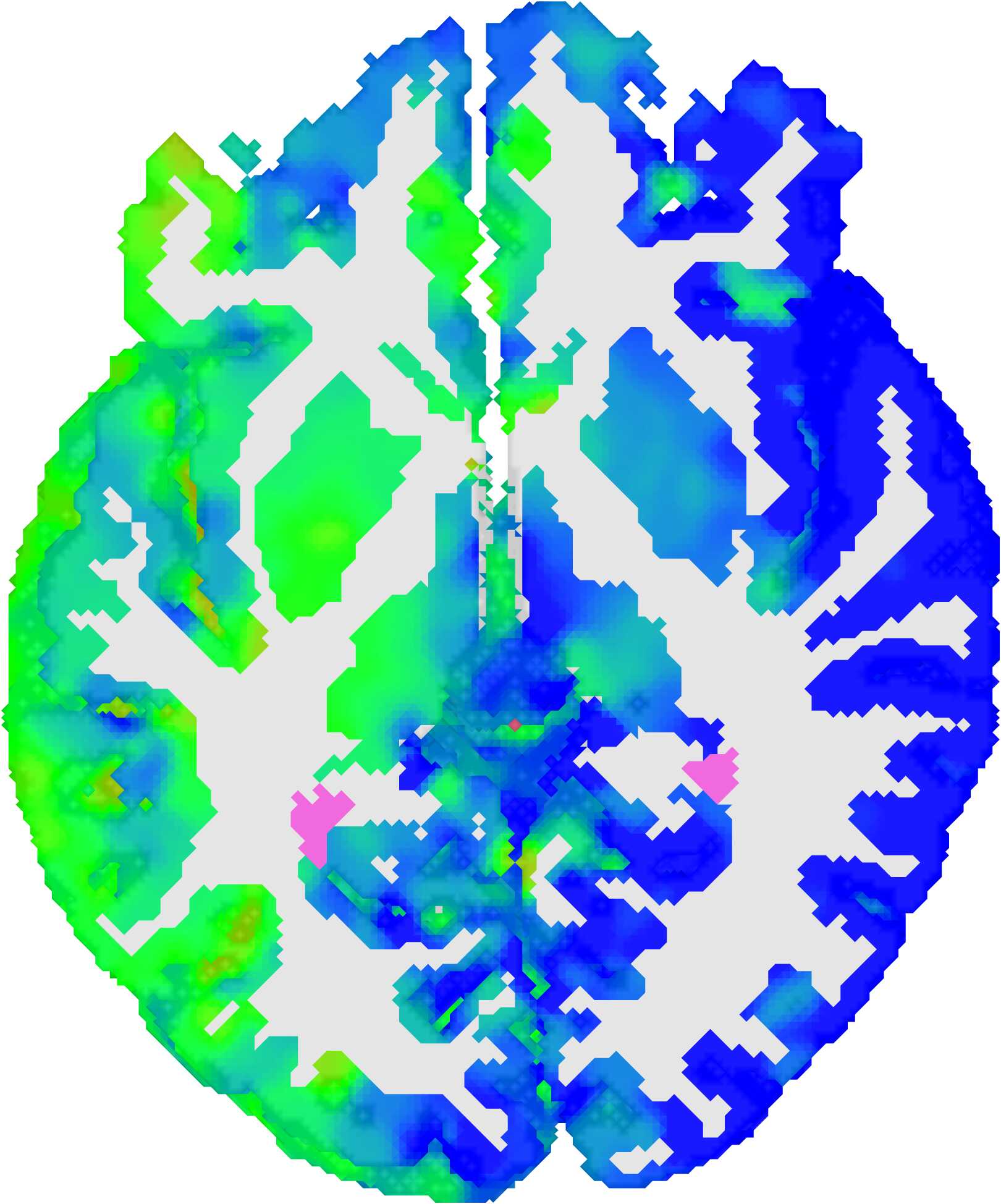}
            \caption{$\azimuth=\qty{0}{\degree}$, $\elevation=\qty{90}{\degree}$, $z=\qty{30}{\milli\meter}$.}
        \end{subfigure}
       \hfill
        \begin{subfigure}[b]{\subFigWidth}
            \includegraphics[width=\linewidth]{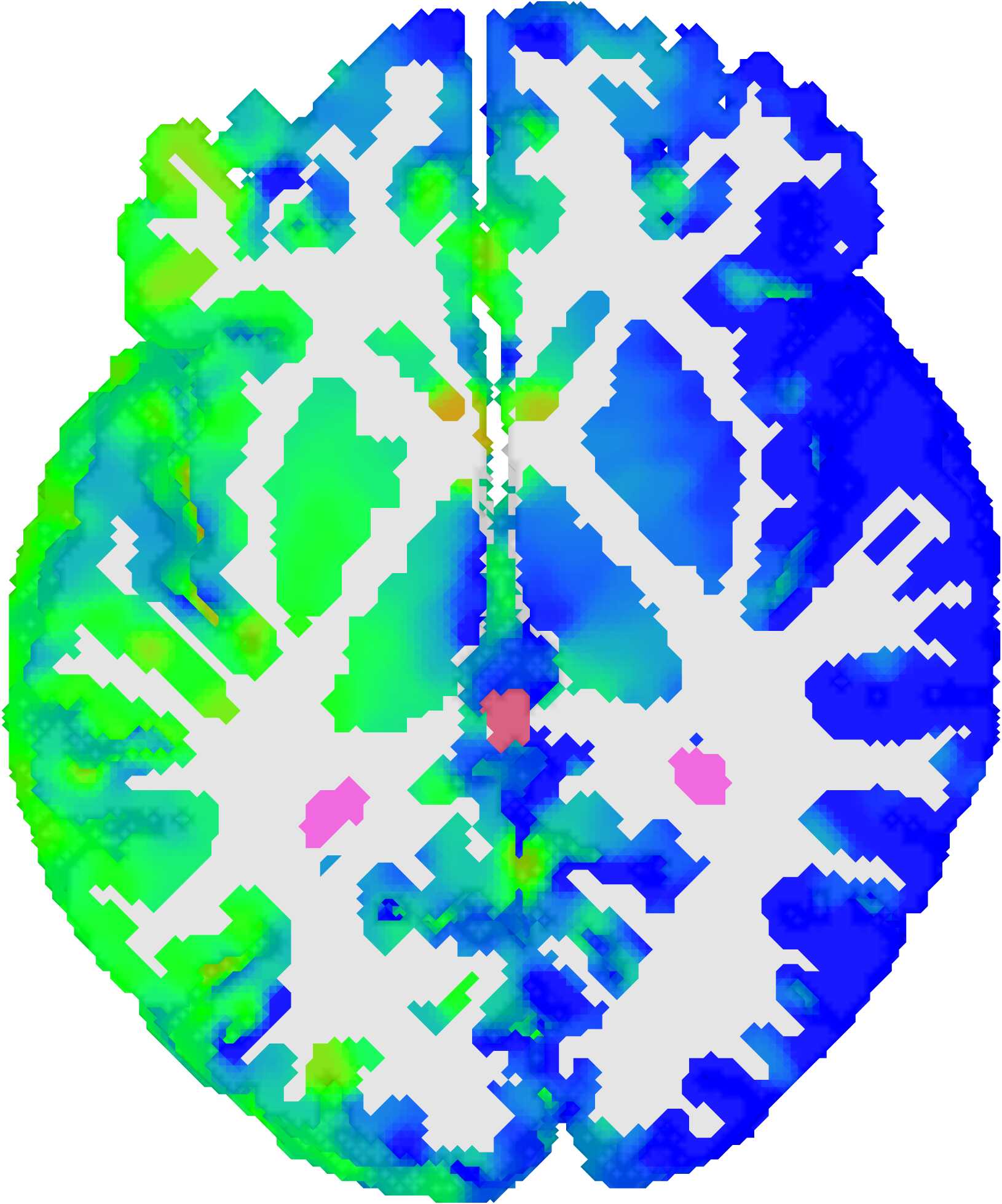}
            \caption{$\azimuth=\qty{0}{\degree}$, $\elevation=\qty{90}{\degree}$, $z=\qty{33}{\milli\meter}$.}
        \end{subfigure}

        \begin{subfigure}[b]{\subFigWidth}
            \includegraphics[width=\linewidth]{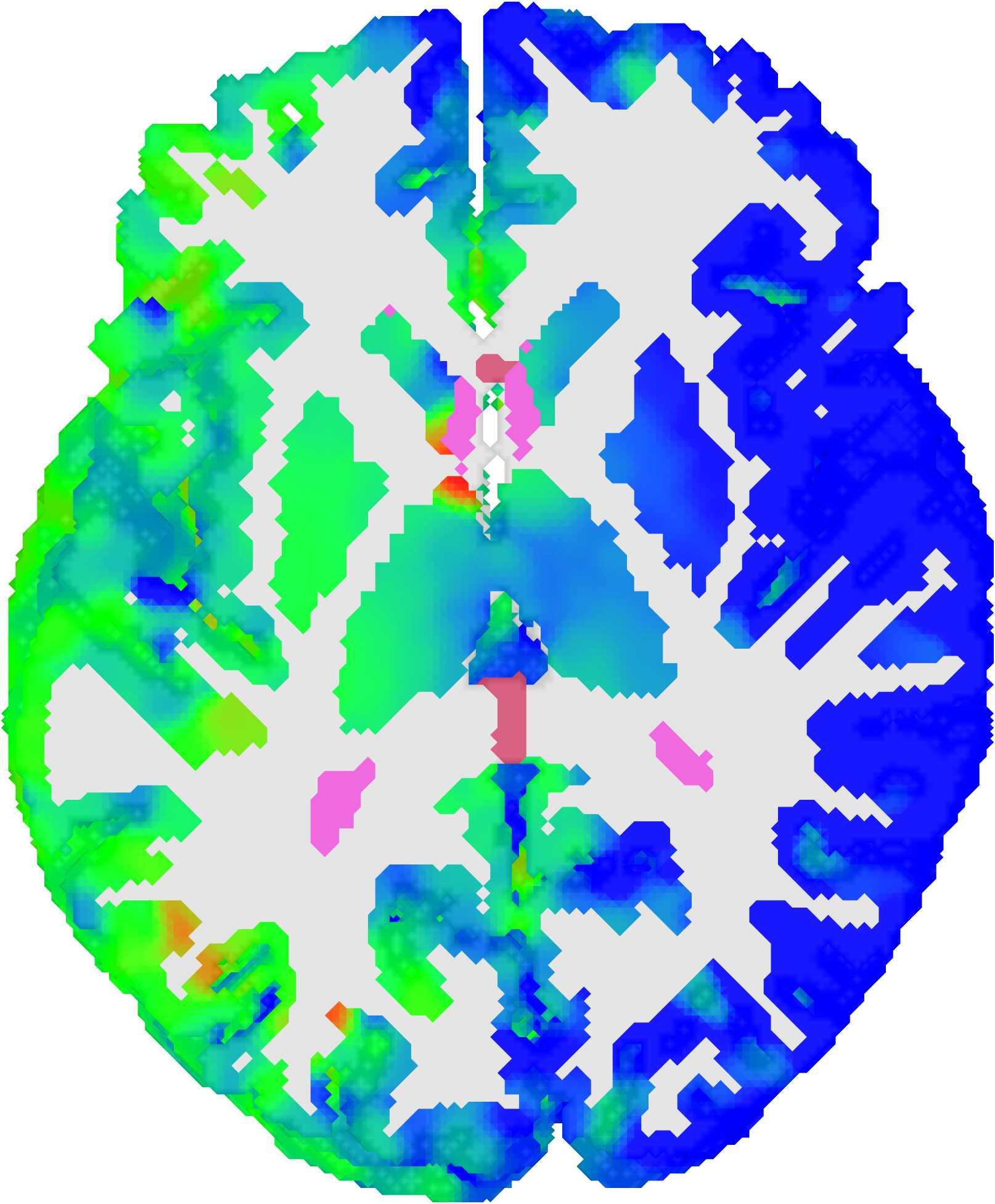}
            \caption{$\azimuth=\qty{0}{\degree}$, $\elevation=\qty{90}{\degree}$, $z=\qty{37}{\milli\meter}$.}
        \end{subfigure}
        \hfill
        \begin{subfigure}[b]{\subFigWidth}
            \includegraphics[width=\linewidth]{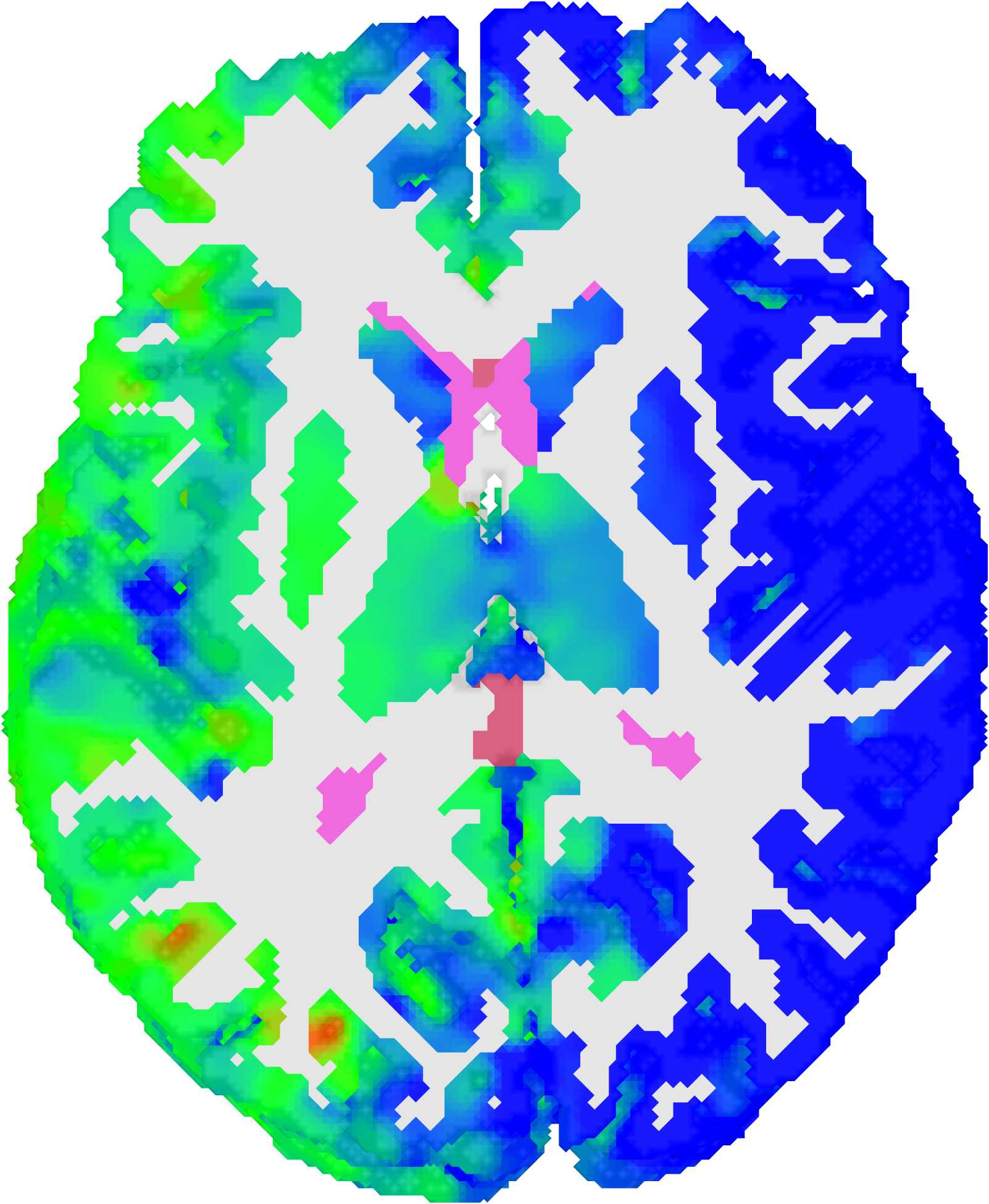}
            \caption{$\azimuth=\qty{0}{\degree}$, $\elevation=\qty{90}{\degree}$, $z=\qty{40}{\milli\meter}$.}
        \end{subfigure}
    \end{minipage}
    \hfill
    \begin{minipage}{0.09\linewidth}
    \includegraphics[width=\linewidth]{c,f=1000,li=0.5,ri=1.5,dB=18,az=0,el=90.colorbar.pdf}
    \end{minipage}
    \caption{Volume current interference field \eqref{eq:interference-field} of $\refJv$, smoothed via the function $\dBT\args{\someField/\max\someField,\qty{18}{\decibel}}$ of equation \eqref{eq:visual-log-transform} at different azimuths $\azimuth$ and elevations $\elevation$, with base stimulation frequency $\frequency=\qty{1000}{\hertz}$ and an updated $\contactResistance=\qty{5270}{\ohm}$. The current pattern was $\currentPattern=[-0.5,0.5,-1.5,1.5]\,\si{\milli\ampere}$ and the maximal value of $\refJv$ roughly equal to \qty{3.9774e-05}{\milli\ampere\per\milli\meter\squared}.}
    \label{fig:interference-reffield-weaker-left-5270}
\end{figure}

\begin{figure}[!h]
    \centering
    \begin{minipage}{0.85\linewidth}
        \def\subFigWidth{0.31\linewidth}
        \begin{subfigure}[b]{\subFigWidth}
            \includegraphics[width=\linewidth]{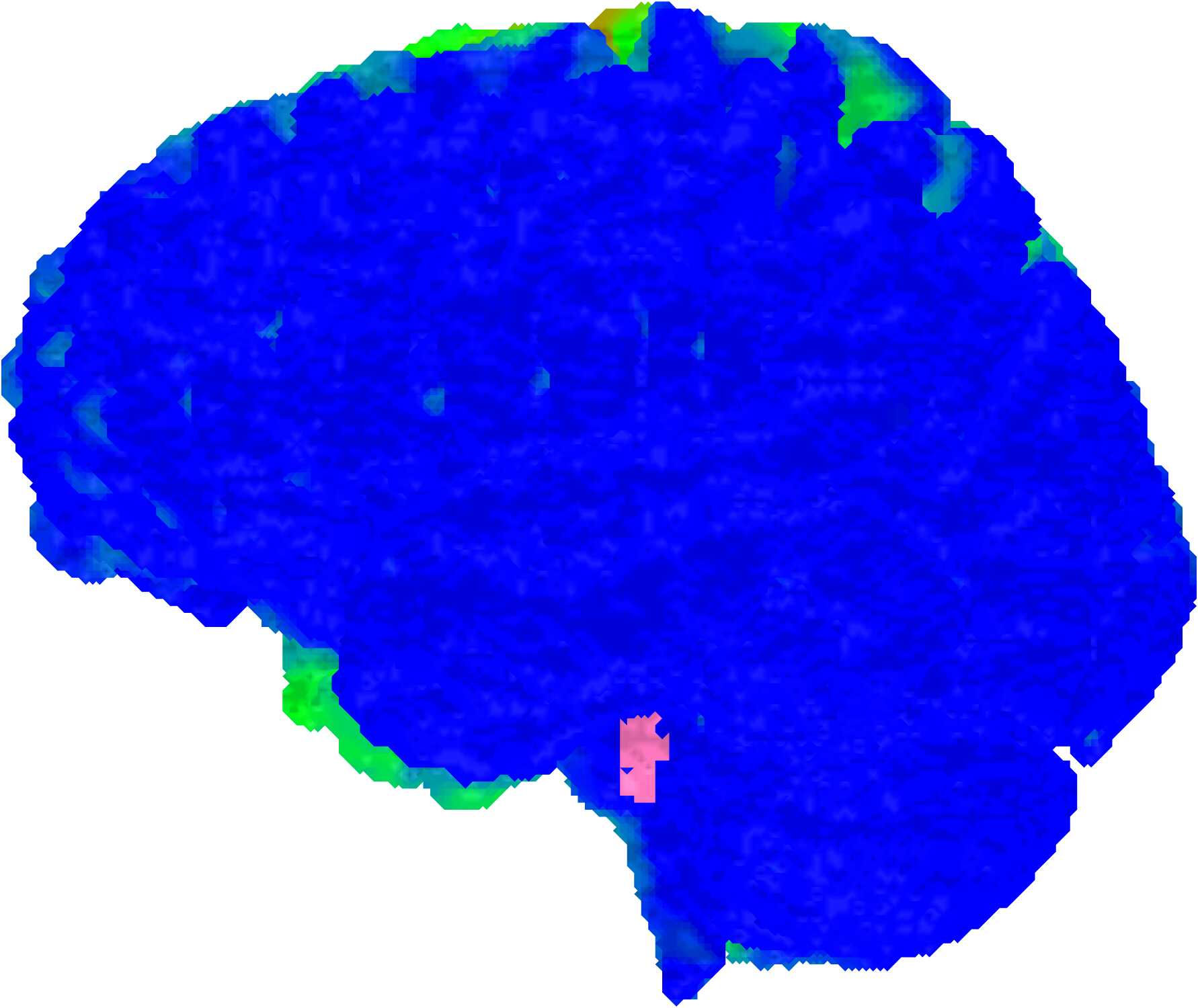}
            \caption{$\azimuth=\qty{-90}{\degree}$, $\elevation=\qty{0}{\degree}$.}
        \end{subfigure}
        \hfill
        \begin{subfigure}[b]{\subFigWidth}
        \includegraphics[width=\linewidth]{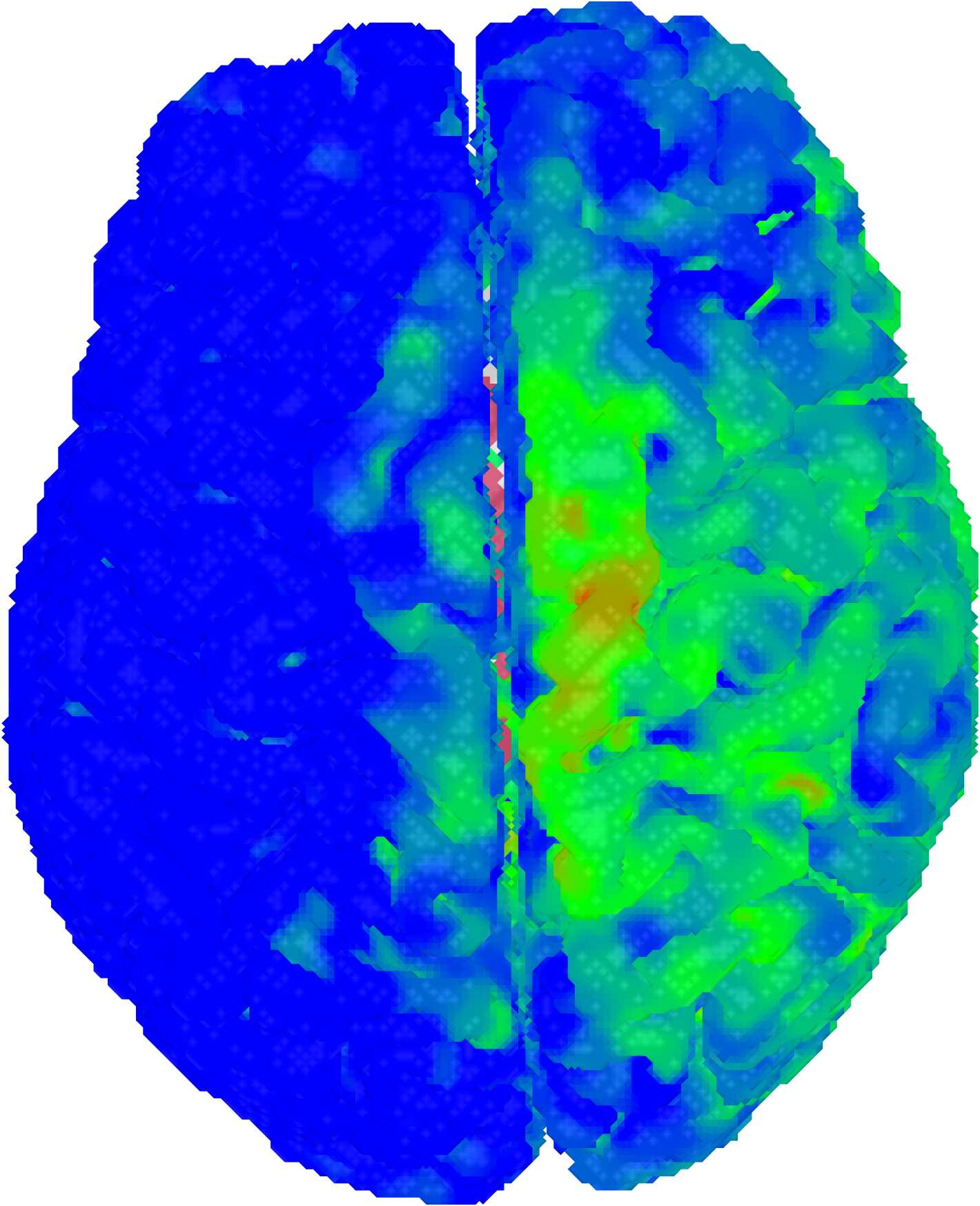}
        \caption{$\azimuth=\qty{0}{\degree}$, $\elevation=\qty{90}{\degree}$.}
        \end{subfigure}
        \hfill
        \begin{subfigure}[b]{\subFigWidth}
        \includegraphics[width=\linewidth]{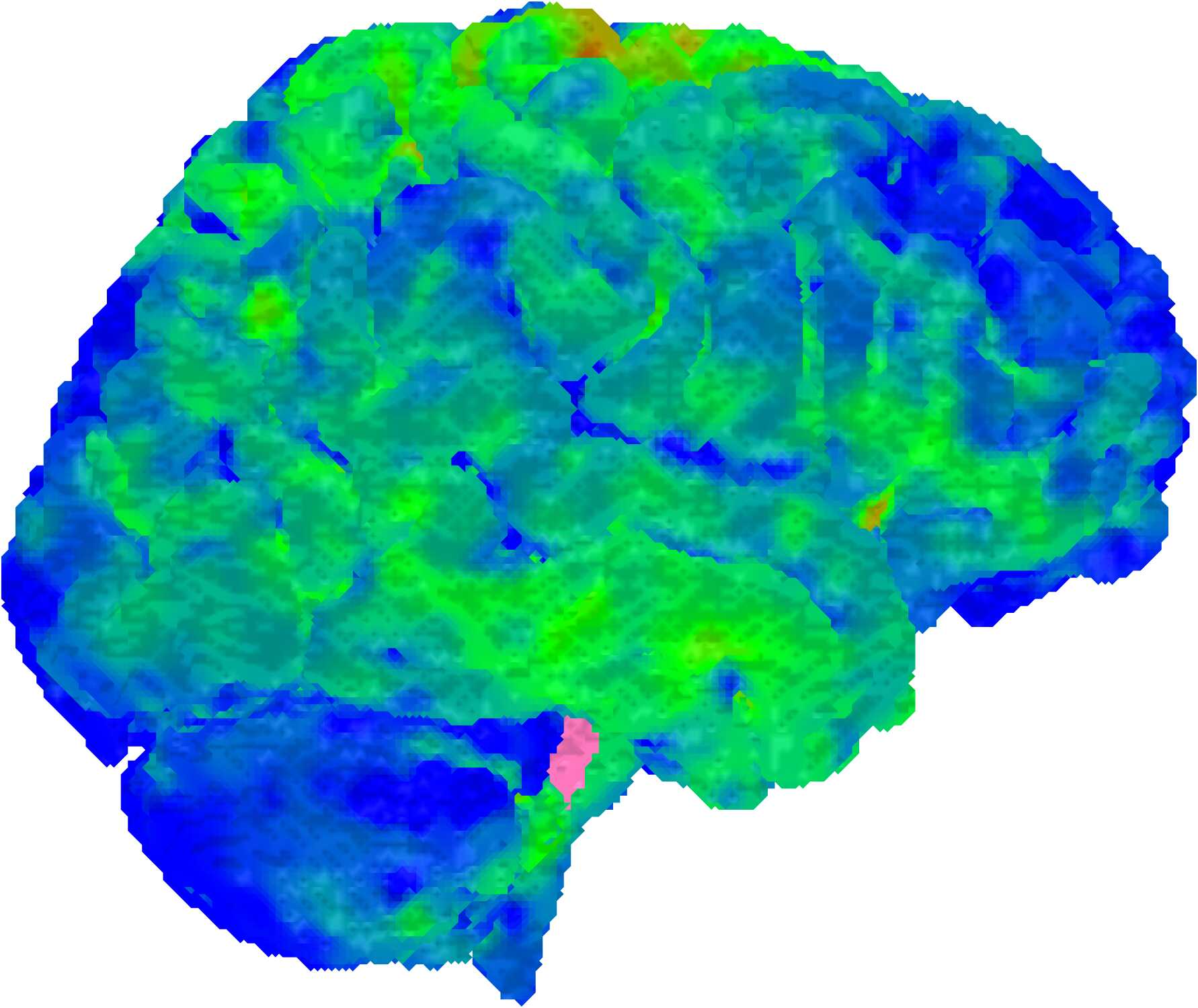}
            \caption{$\azimuth=\qty{90}{\degree}$, $\elevation=\qty{0}{\degree}$.}
        \end{subfigure}

        \def\subFigWidth{0.45\linewidth}

        \begin{subfigure}[b]{\subFigWidth}
            \includegraphics[width=\linewidth]{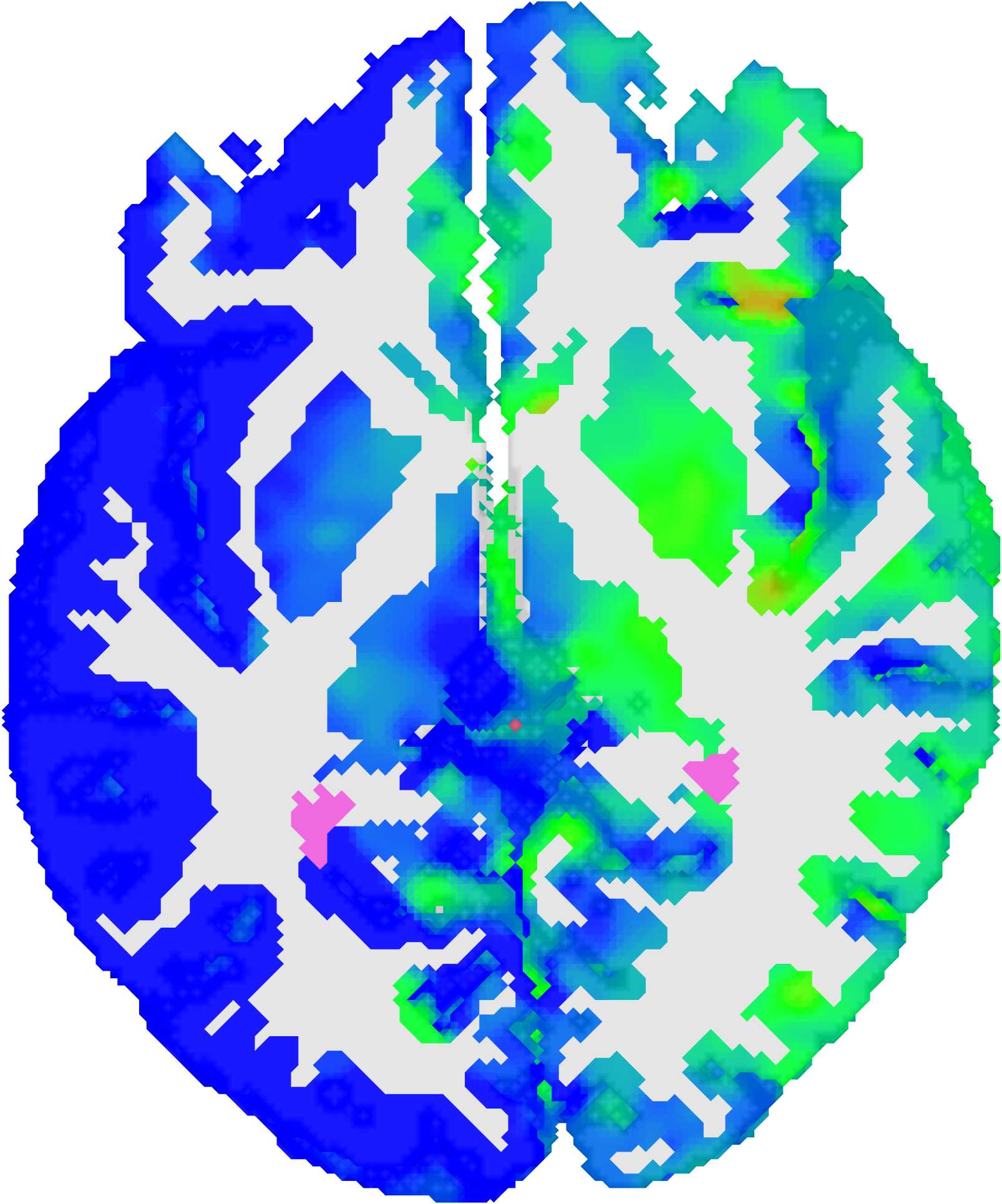}
            \caption{$\azimuth=\qty{0}{\degree}$, $\elevation=\qty{90}{\degree}$, $z=\qty{30}{\milli\meter}$.}
        \end{subfigure}
       \hfill
        \begin{subfigure}[b]{\subFigWidth}
            \includegraphics[width=\linewidth]{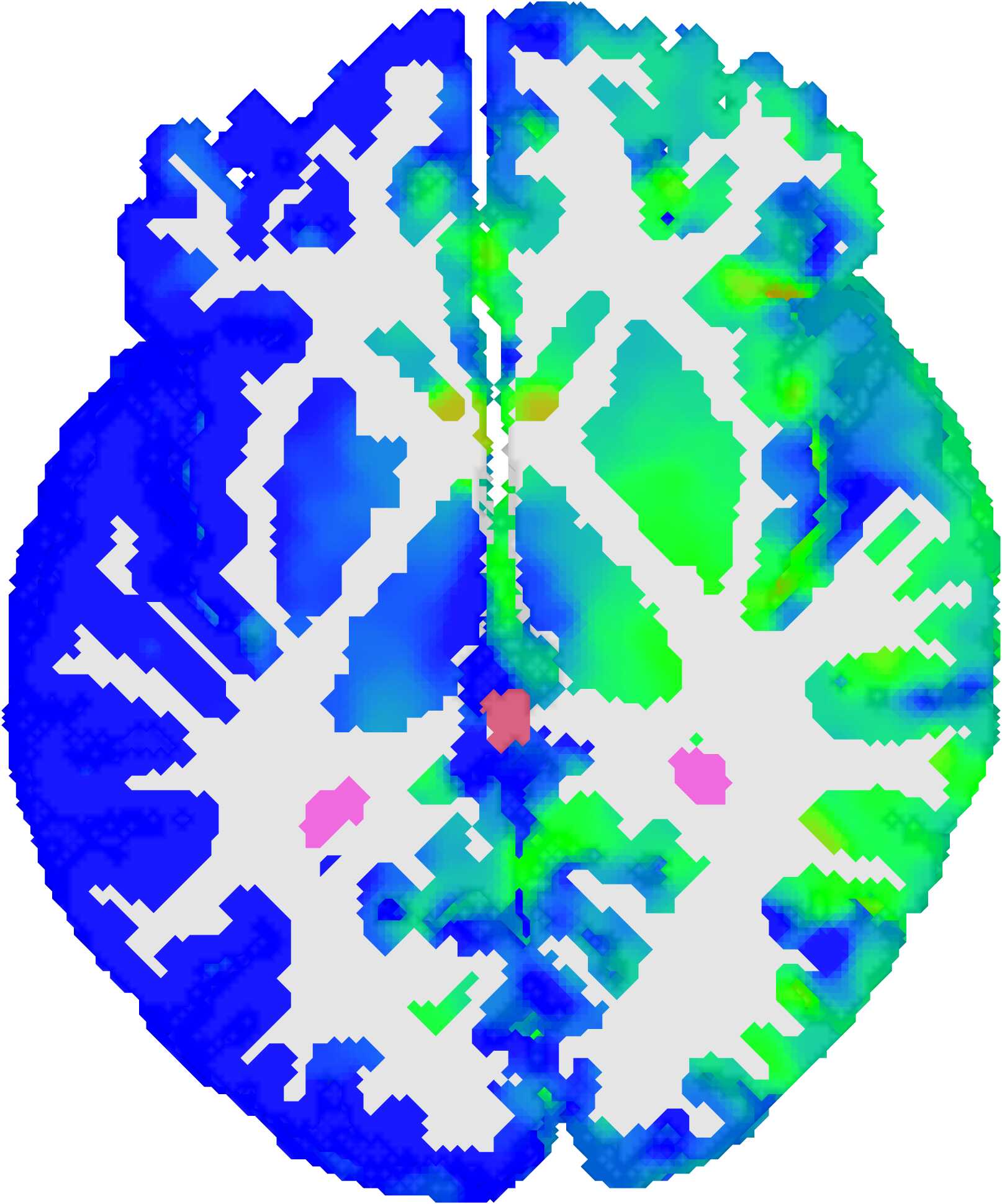}
            \caption{$\azimuth=\qty{0}{\degree}$, $\elevation=\qty{90}{\degree}$, $z=\qty{33}{\milli\meter}$.}
        \end{subfigure}

        \begin{subfigure}[b]{\subFigWidth}
            \includegraphics[width=\linewidth]{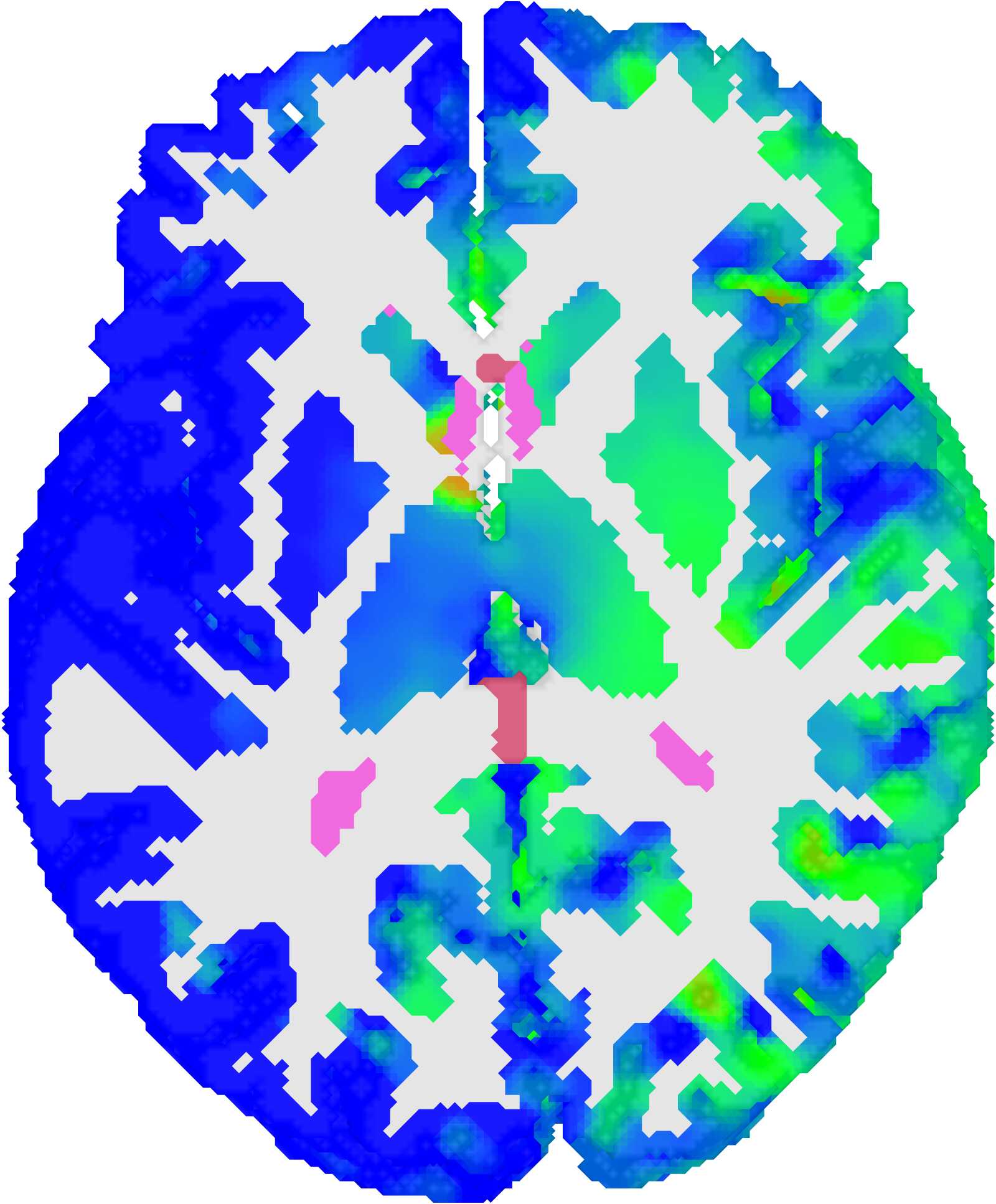}
            \caption{$\azimuth=\qty{0}{\degree}$, $\elevation=\qty{90}{\degree}$, $z=\qty{37}{\milli\meter}$.}
        \end{subfigure}
        \hfill
        \begin{subfigure}[b]{\subFigWidth}
            \includegraphics[width=\linewidth]{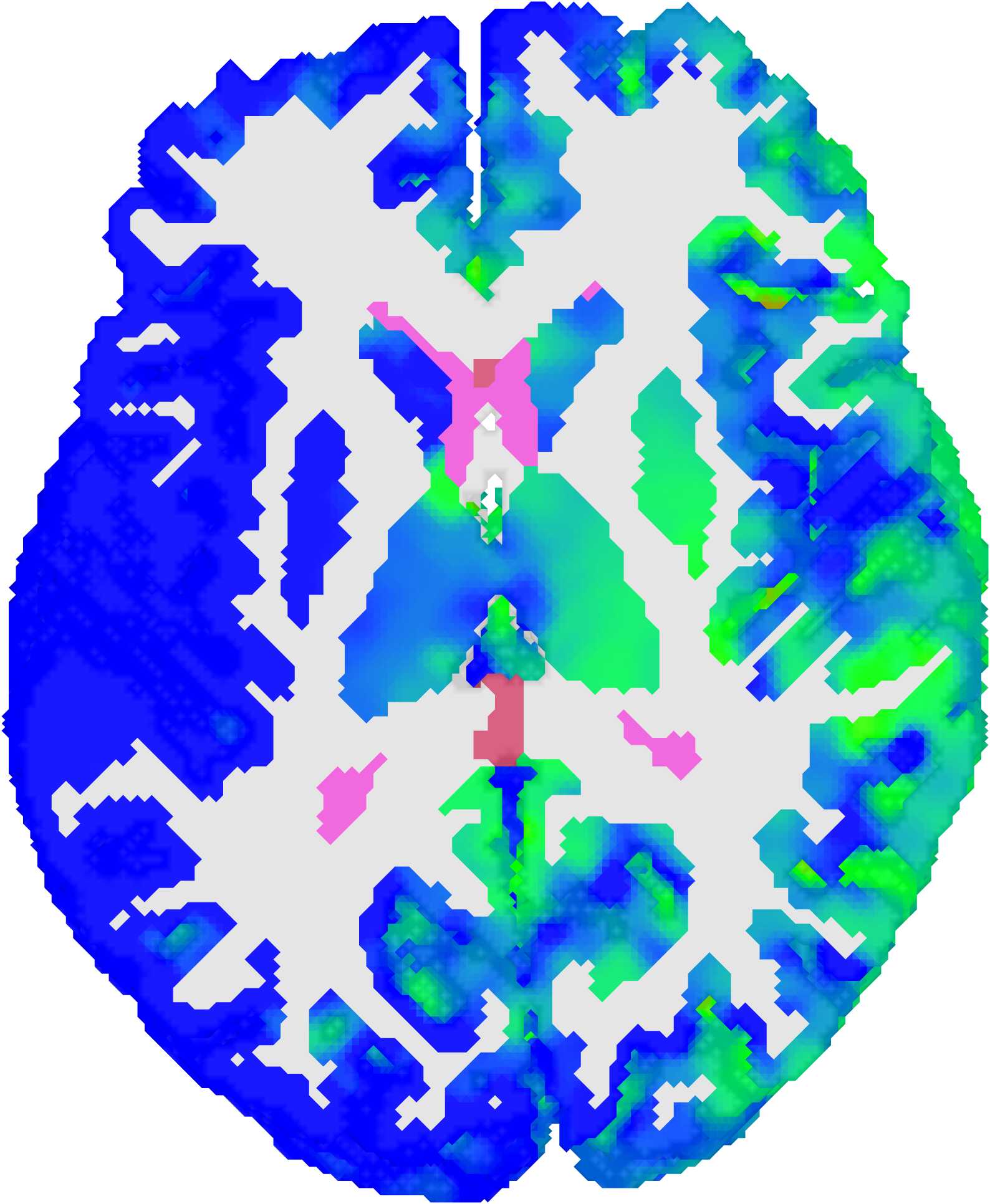}
            \caption{$\azimuth=\qty{0}{\degree}$, $\elevation=\qty{90}{\degree}$, $z=\qty{40}{\milli\meter}$.}
        \end{subfigure}
    \end{minipage}
    \hfill
    \begin{minipage}{0.09\linewidth}
    \includegraphics[width=\linewidth]{c,f=1000,li=0.5,ri=1.5,dB=18,az=0,el=90.colorbar.pdf}
    \end{minipage}
    \caption{Volume current interference field \eqref{eq:interference-field} of $\refJv$, smoothed via the function $\dBT\args{\someField/\max\someField,\qty{18}{\decibel}}$ of equation \eqref{eq:visual-log-transform} at different azimuths $\azimuth$ and elevations $\elevation$, with base stimulation frequency $\frequency=\qty{1000}{\hertz}$ and an updated $\contactResistance=\qty{1270}{\ohm}$. The current pattern was $\currentPattern=[-1.5,1.5,-0.5,0.5]\,\si{\milli\ampere}$  and the maximal value of $\refJv$ roughly equal to \qty{5.0823e-05}{\milli\ampere\per\milli\meter\squared}.}
    \label{fig:interference-reffield-weaker-right-1270}
\end{figure}

\begin{figure}[!h]
    \centering
    \begin{minipage}{0.85\linewidth}
        \def\subFigWidth{0.31\linewidth}
        \begin{subfigure}[b]{\subFigWidth}
            \includegraphics[width=\linewidth]{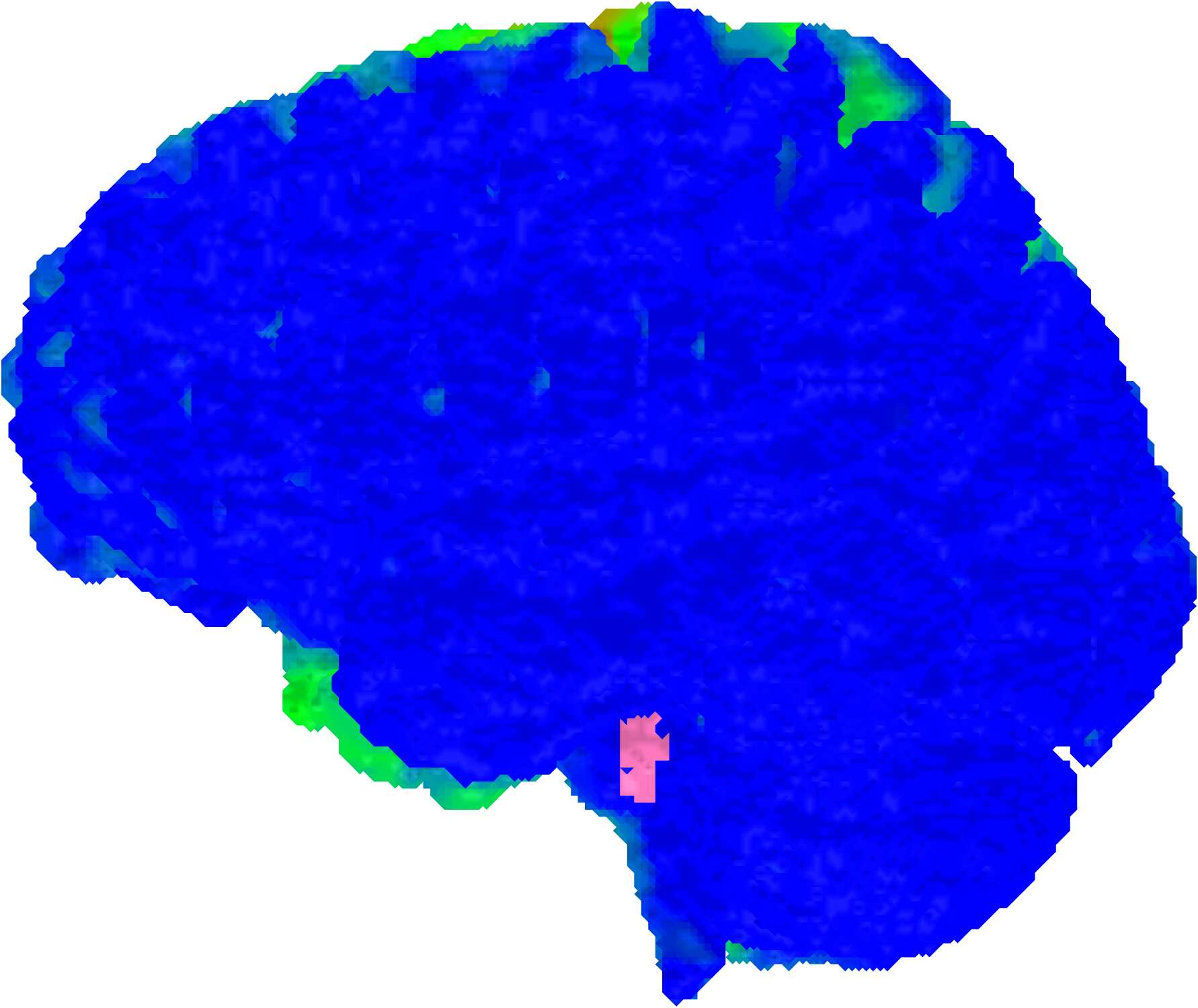}
            \caption{$\azimuth=\qty{-90}{\degree}$, $\elevation=\qty{0}{\degree}$.}
        \end{subfigure}
        \hfill
        \begin{subfigure}[b]{\subFigWidth}
        \includegraphics[width=\linewidth]{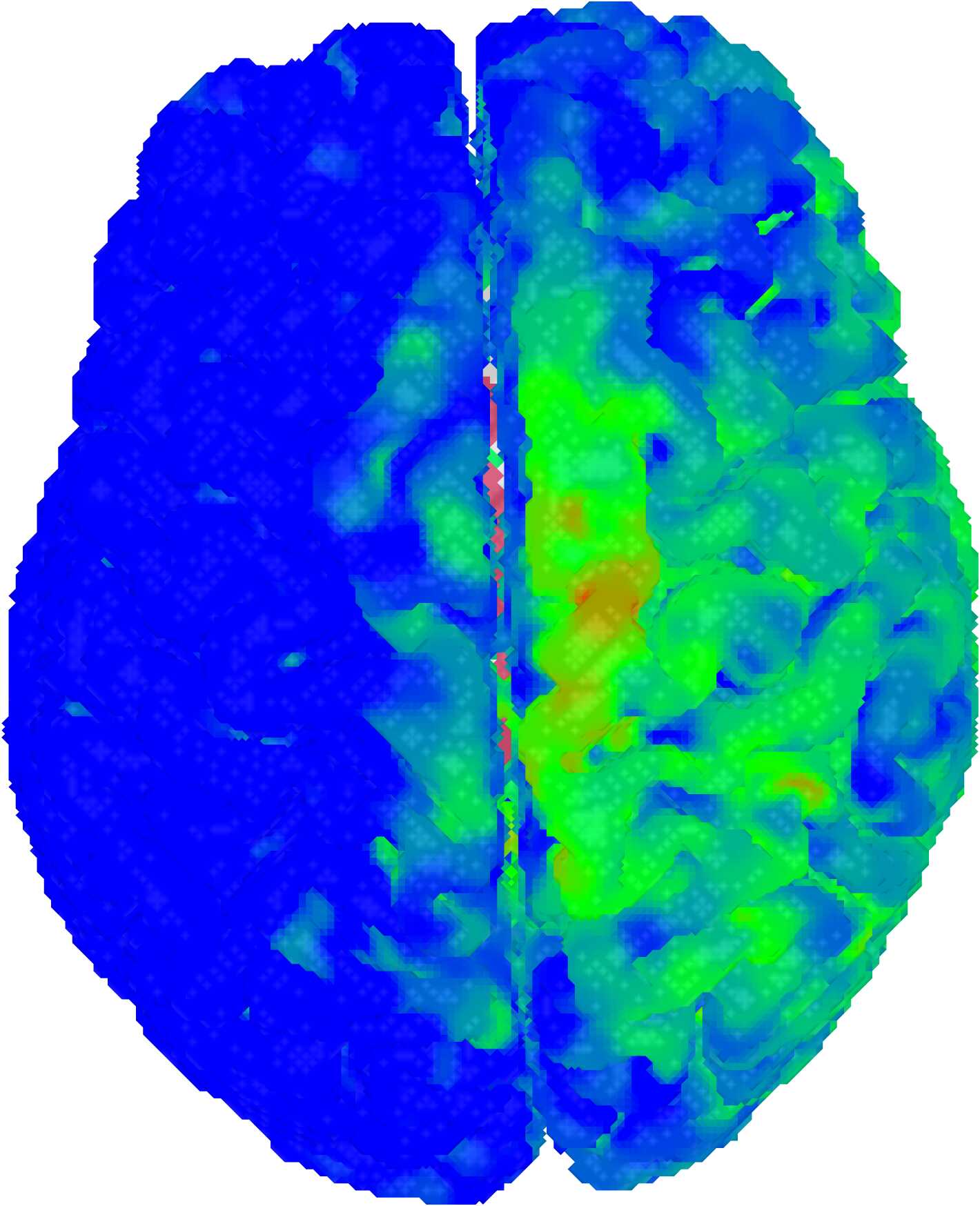}
        \caption{$\azimuth=\qty{0}{\degree}$, $\elevation=\qty{90}{\degree}$.}
        \end{subfigure}
        \hfill
        \begin{subfigure}[b]{\subFigWidth}
        \includegraphics[width=\linewidth]{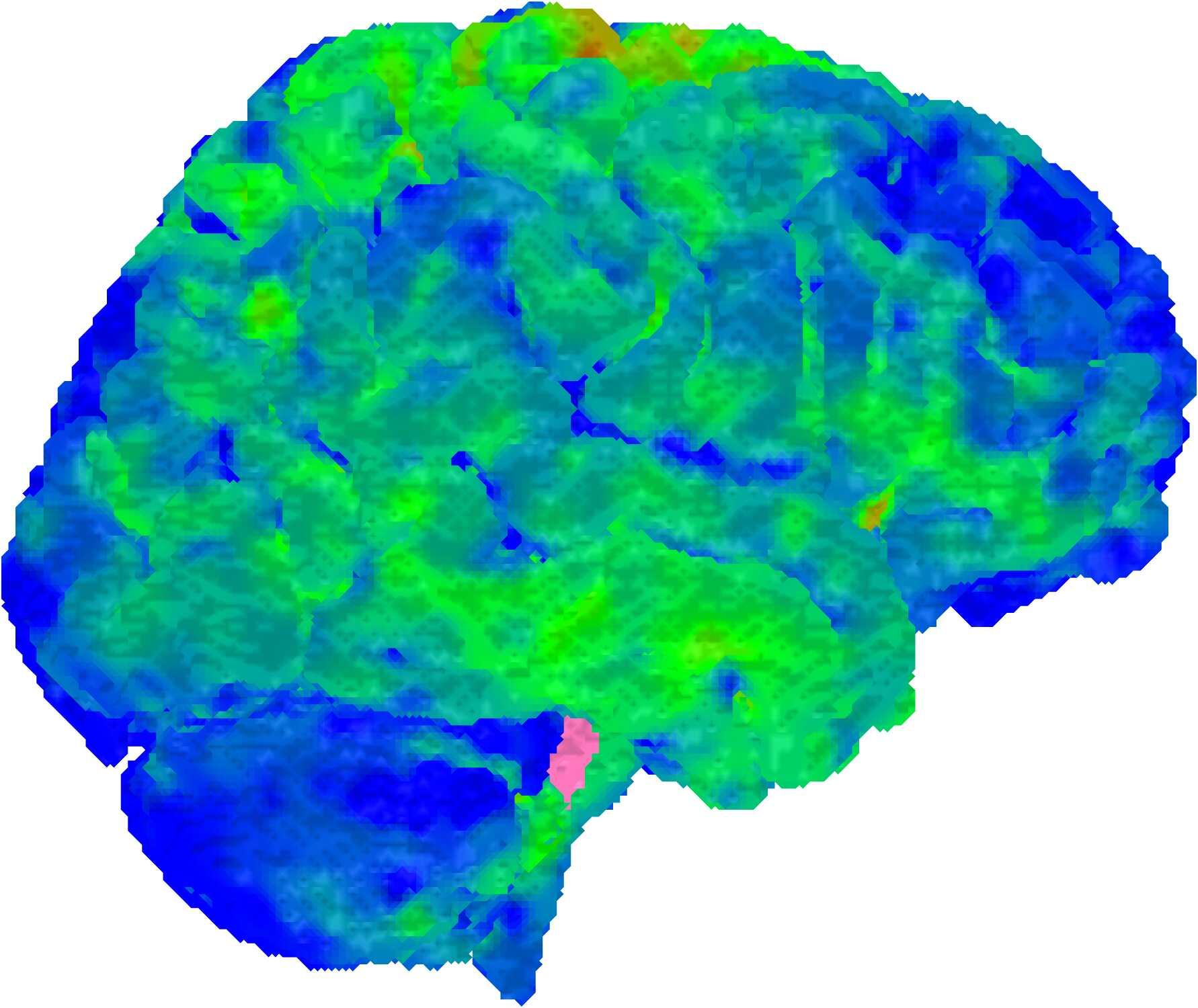}
            \caption{$\azimuth=\qty{90}{\degree}$, $\elevation=\qty{0}{\degree}$.}
        \end{subfigure}

        \def\subFigWidth{0.45\linewidth}

        \begin{subfigure}[b]{\subFigWidth}
            \includegraphics[width=\linewidth]{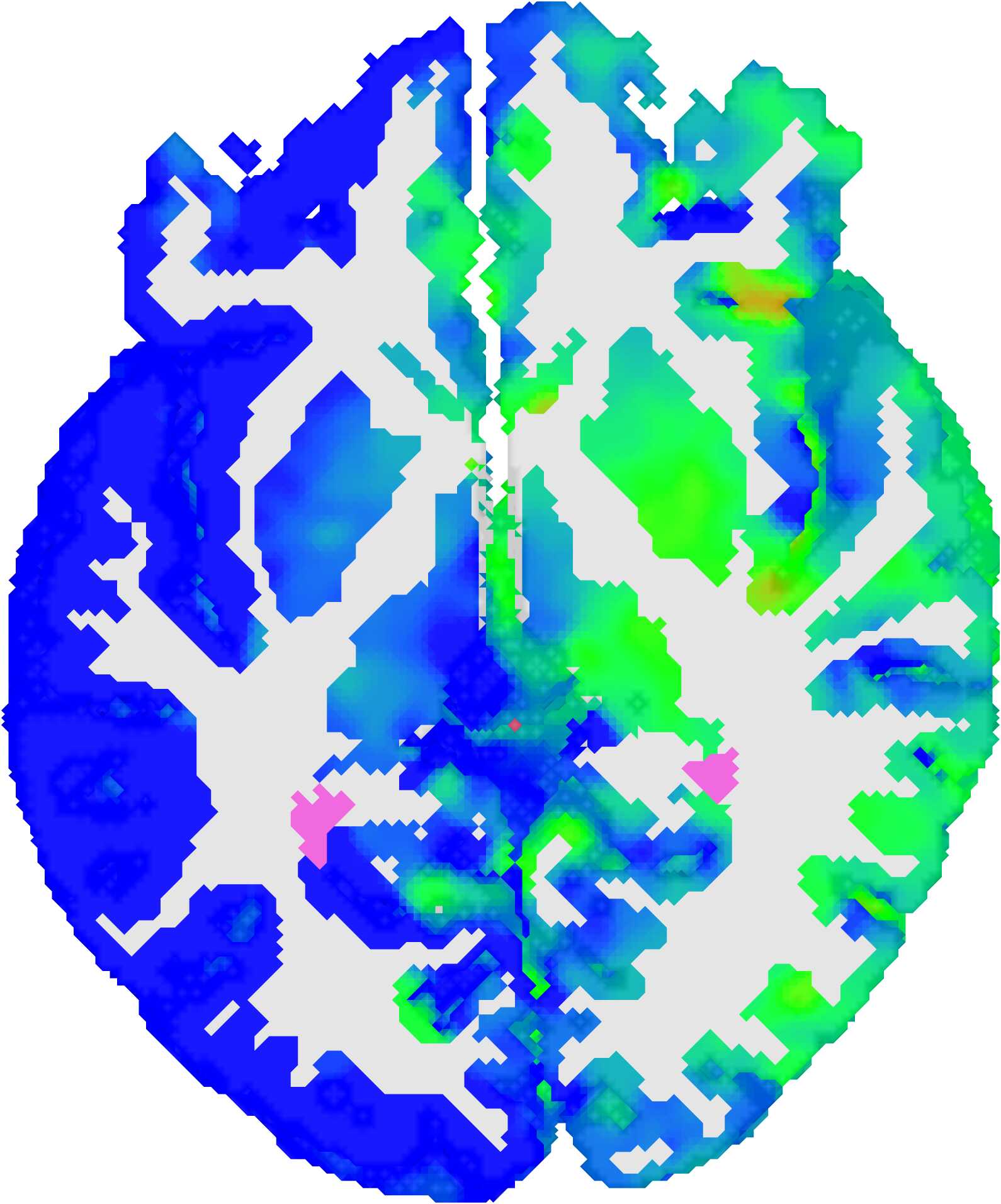}
            \caption{$\azimuth=\qty{0}{\degree}$, $\elevation=\qty{90}{\degree}$, $z=\qty{30}{\milli\meter}$.}
        \end{subfigure}
       \hfill
        \begin{subfigure}[b]{\subFigWidth}
            \includegraphics[width=\linewidth]{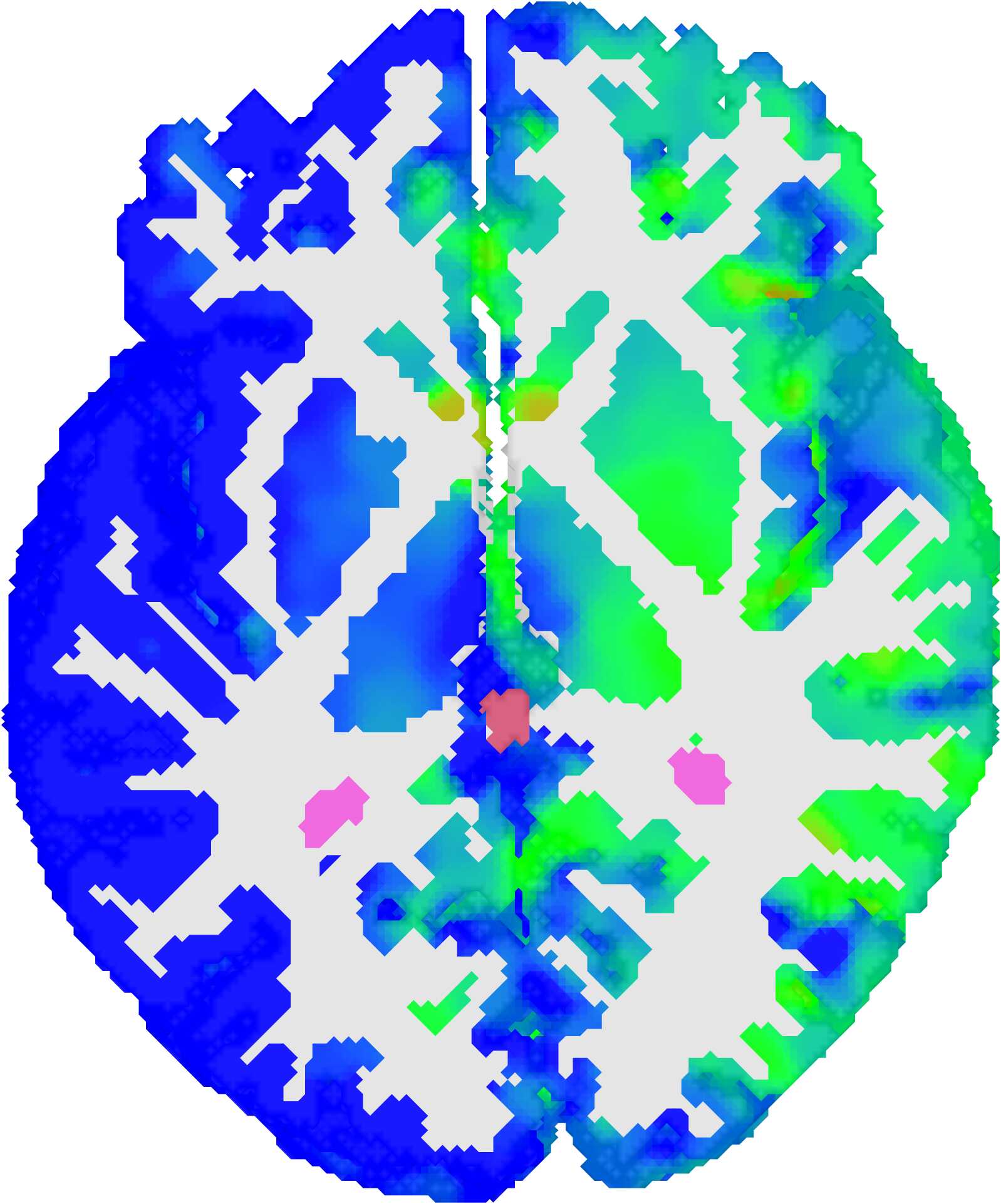}
            \caption{$\azimuth=\qty{0}{\degree}$, $\elevation=\qty{90}{\degree}$, $z=\qty{33}{\milli\meter}$.}
        \end{subfigure}

        \begin{subfigure}[b]{\subFigWidth}
            \includegraphics[width=\linewidth]{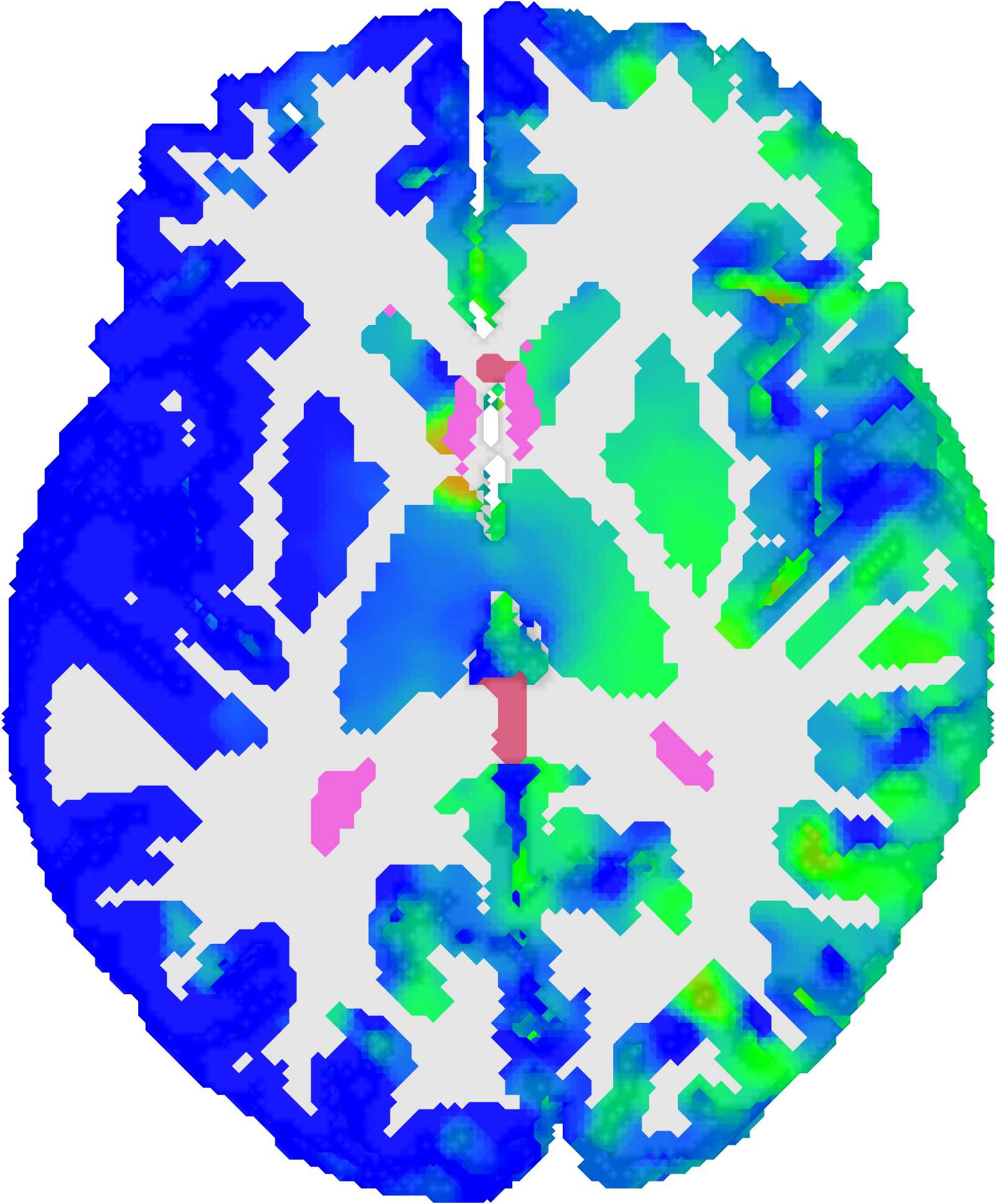}
            \caption{$\azimuth=\qty{0}{\degree}$, $\elevation=\qty{90}{\degree}$, $z=\qty{37}{\milli\meter}$.}
        \end{subfigure}
        \hfill
        \begin{subfigure}[b]{\subFigWidth}
            \includegraphics[width=\linewidth]{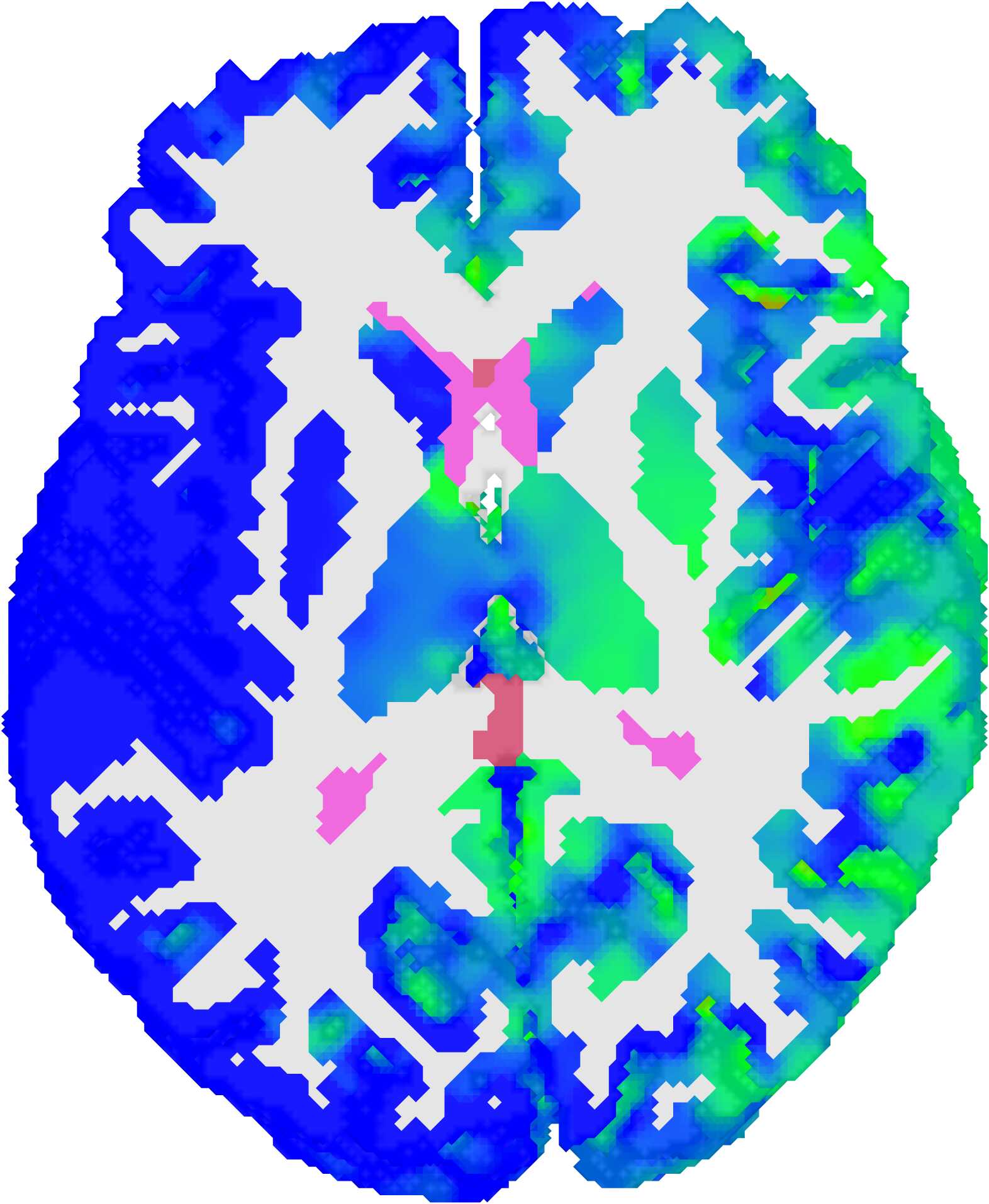}
            \caption{$\azimuth=\qty{0}{\degree}$, $\elevation=\qty{90}{\degree}$, $z=\qty{40}{\milli\meter}$.}
        \end{subfigure}
    \end{minipage}
    \hfill
    \begin{minipage}{0.09\linewidth}
    \includegraphics[width=\linewidth]{c,f=1000,li=0.5,ri=1.5,dB=18,az=0,el=90.colorbar.pdf}
    \end{minipage}
    \caption{Volume current interference field \eqref{eq:interference-field} of $\refJv$, smoothed via the function $\dBT\args{\someField/\max\someField,\qty{18}{\decibel}}$ of equation \eqref{eq:visual-log-transform} at different azimuths $\azimuth$ and elevations $\elevation$, with base stimulation frequency $\frequency=\qty{1000}{\hertz}$ and an updated $\contactResistance=\qty{5270}{\ohm}$. The current pattern was $\currentPattern=[-1.5,1.5,-0.5,0.5]\,\si{\milli\ampere}$  and the maximal value of $\refJv$ roughly equal to \qty{5.0097e-05}{\milli\ampere\per\milli\meter\squared}.}
    \label{fig:interference-reffield-weaker-right-5270}
\end{figure}

We also observe how the interference field behaves
after the contact resistance $\contactResistance$
of electrode TP9 has worsened to both
\qty{1270}{\ohm} and \qty{5270}{\ohm}. The effects
of the worsening $\contactResistance$ are seen in
Figures~\ref{fig:interference-reffield-equal-currents-1270}--\ref{fig:interference-reffield-weaker-right-5270}.
The effect of the increasing resistance in the electrode
seems to slightly modify the resulting inteference field
magnitude, with the lower resistance of the two producing
an ever-so-slightly stronger stimulation. The steering
effect seen in the original interference field also
remains, and the actual shape of the field regarding
its magnitudes across the colume conductor also does not
visibly change.

\begin{figure}[!h]
    \centering
    \begin{minipage}{0.85\linewidth}
        \def\subFigWidth{0.31\linewidth}
        \begin{subfigure}[b]{\subFigWidth}
            \includegraphics[width=\linewidth]{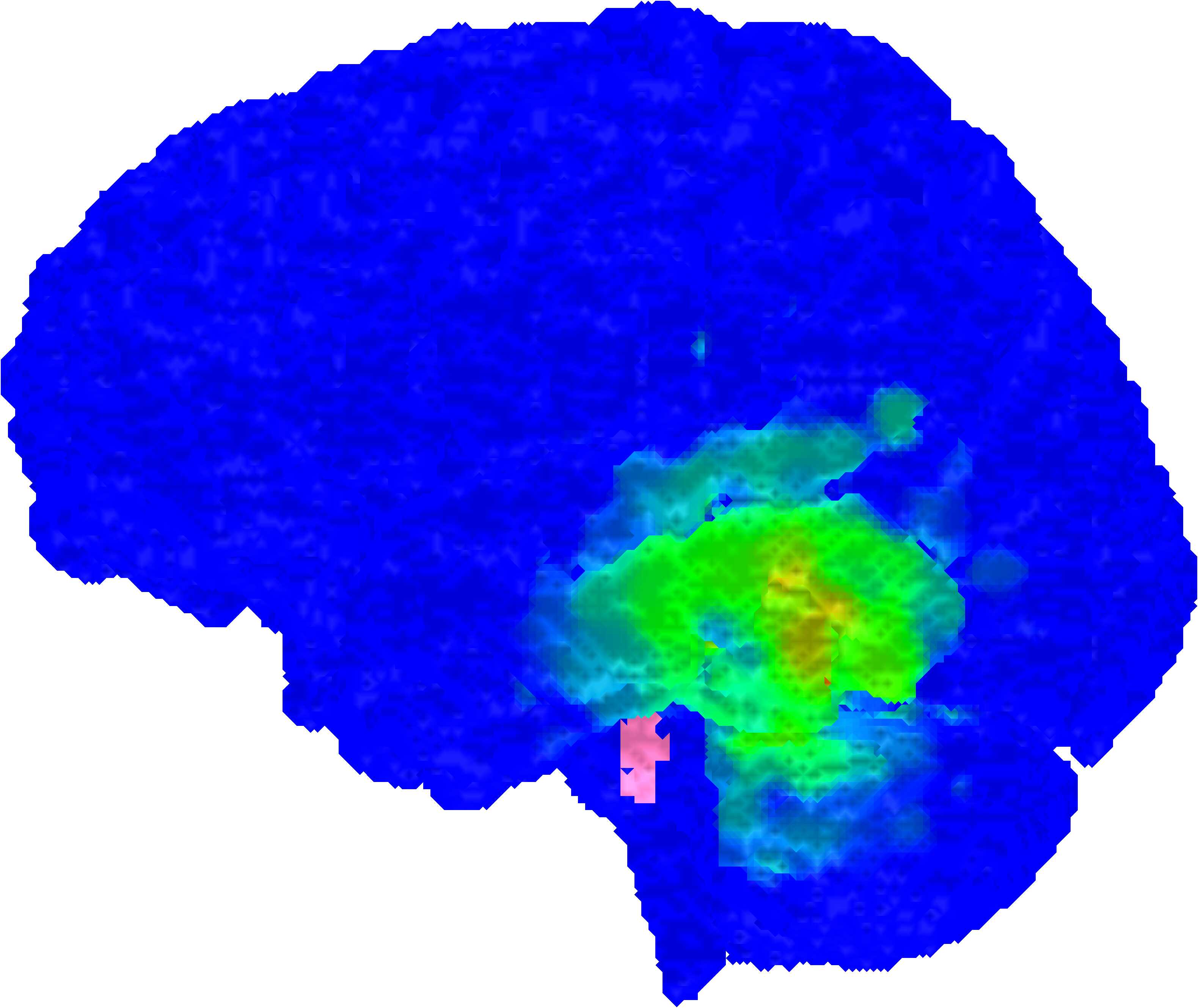}
            \caption{$\azimuth=\qty{-90}{\degree}$, $\elevation=\qty{0}{\degree}$.}
        \end{subfigure}
        \hfill
        \begin{subfigure}[b]{\subFigWidth}
        \includegraphics[width=\linewidth]{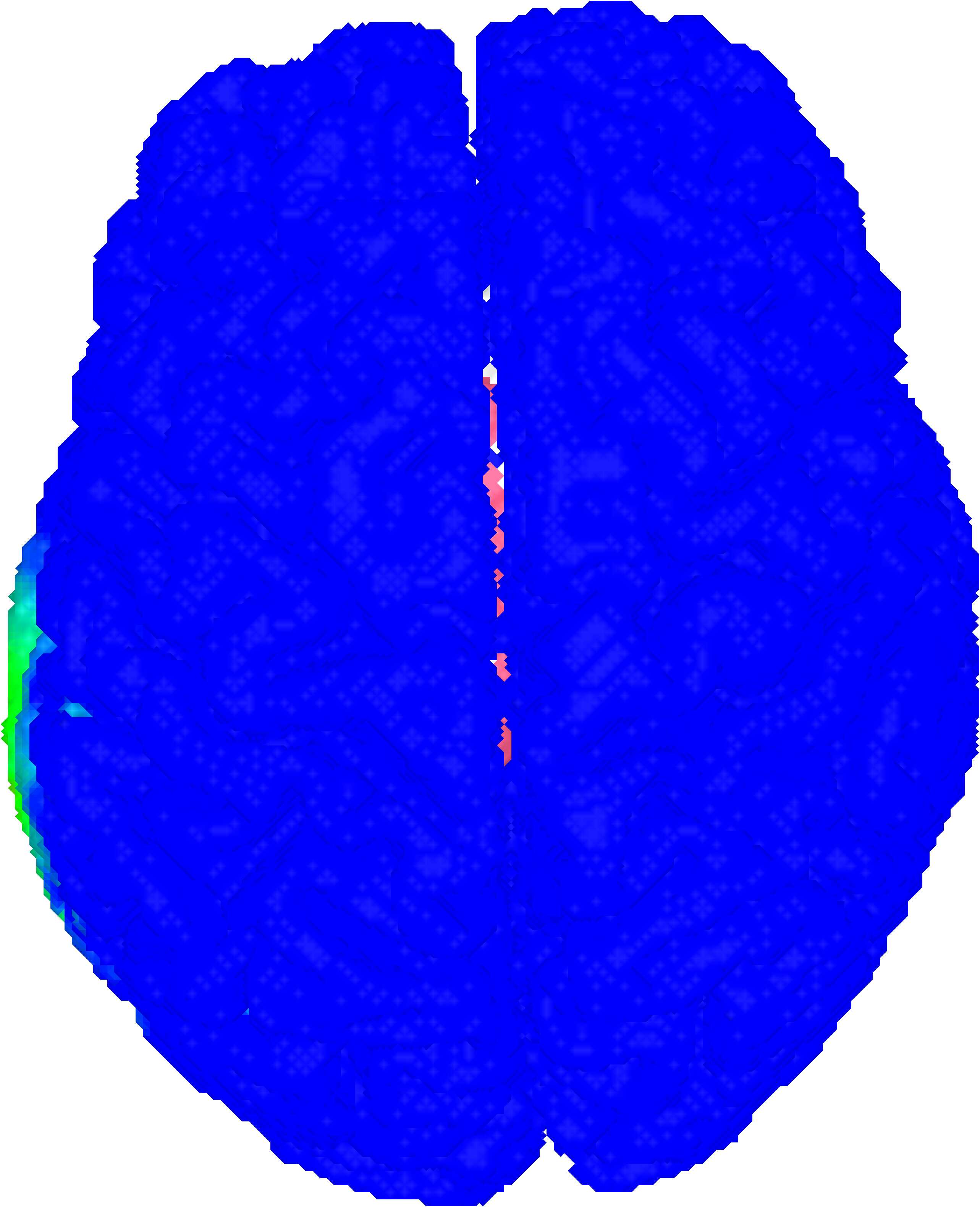}
        \caption{$\azimuth=\qty{0}{\degree}$, $\elevation=\qty{90}{\degree}$.}
        \end{subfigure}
        \hfill
        \begin{subfigure}[b]{\subFigWidth}
        \includegraphics[width=\linewidth]{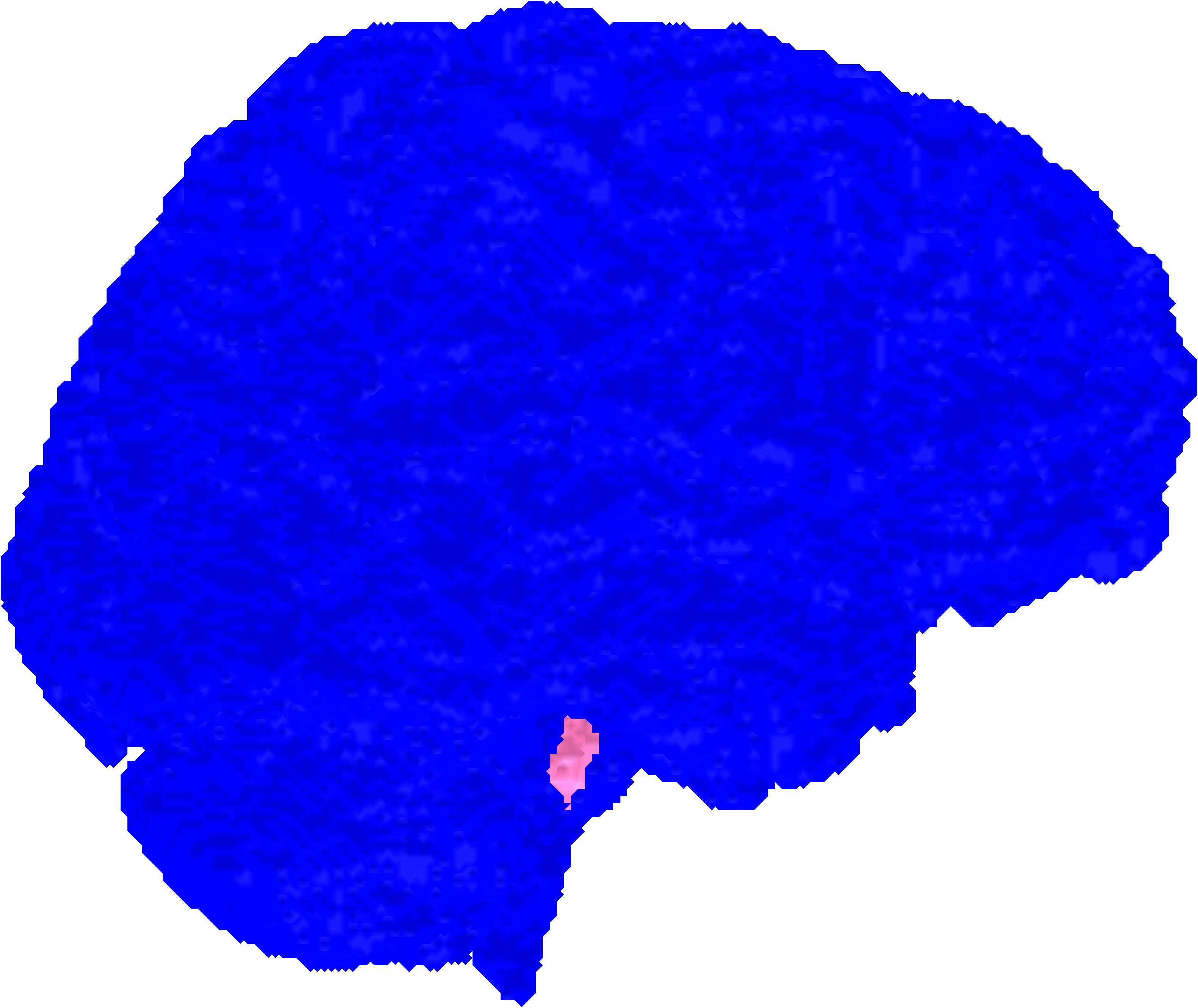}
            \caption{$\azimuth=\qty{90}{\degree}$, $\elevation=\qty{0}{\degree}$.}
        \end{subfigure}

        \def\subFigWidth{0.45\linewidth}

        \begin{subfigure}[b]{\subFigWidth}
            \includegraphics[width=\linewidth]{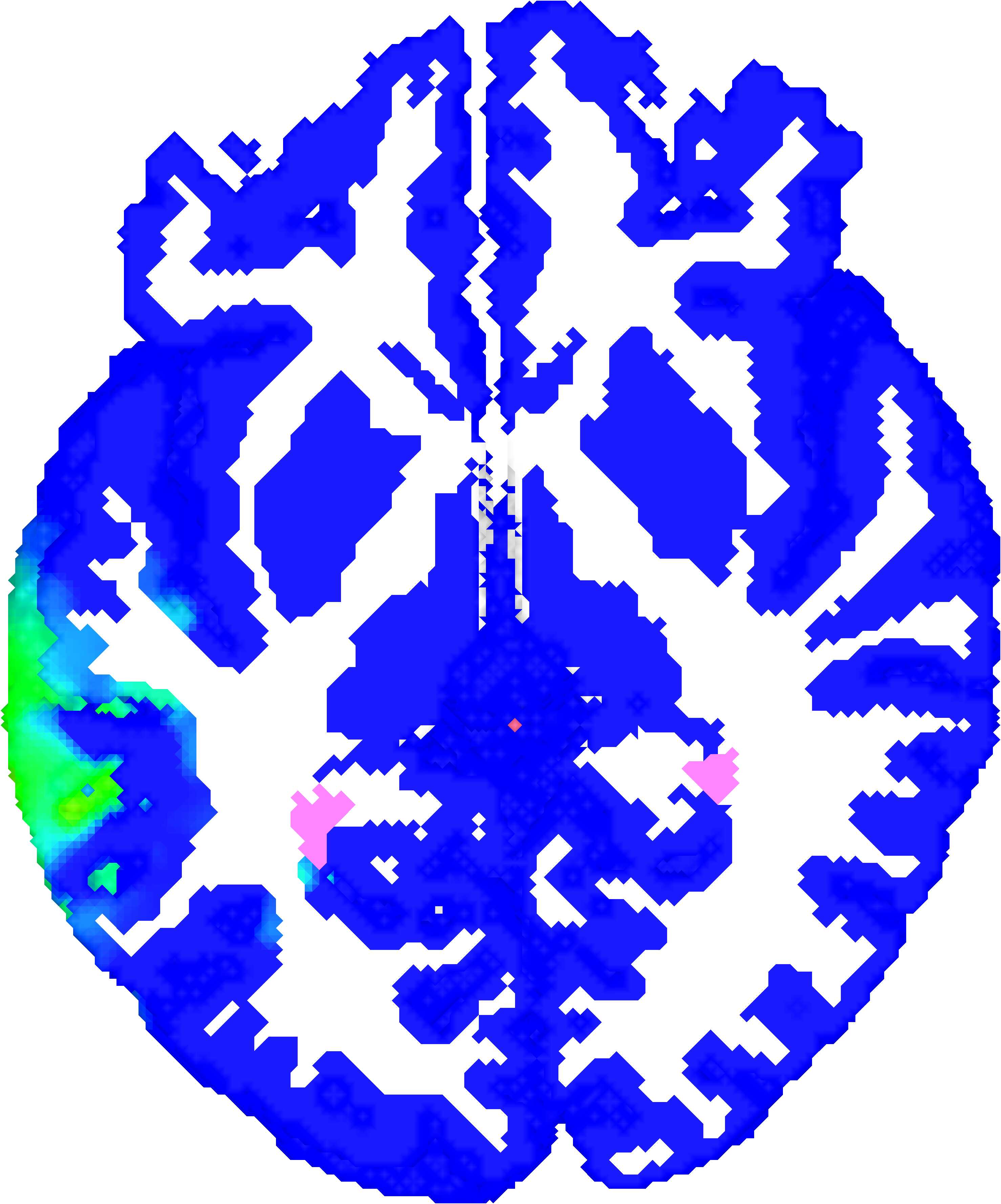}
            \caption{$\azimuth=\qty{0}{\degree}$, $\elevation=\qty{90}{\degree}$, $z=\qty{30}{\milli\meter}$.}
        \end{subfigure}
       \hfill
        \begin{subfigure}[b]{\subFigWidth}
            \includegraphics[width=\linewidth]{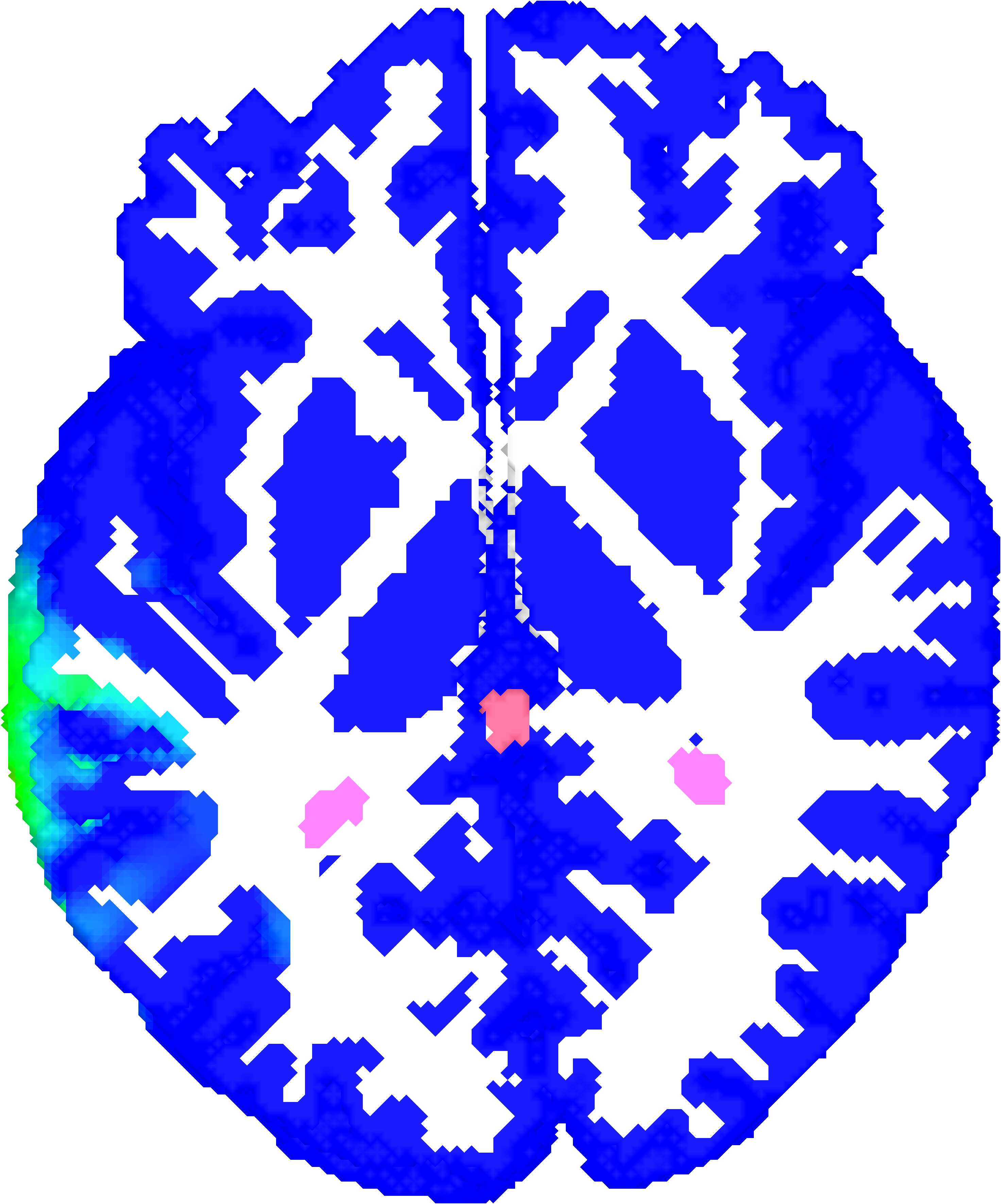}
            \caption{$\azimuth=\qty{0}{\degree}$, $\elevation=\qty{90}{\degree}$, $z=\qty{33}{\milli\meter}$.}
        \end{subfigure}

        \begin{subfigure}[b]{\subFigWidth}
            \includegraphics[width=\linewidth]{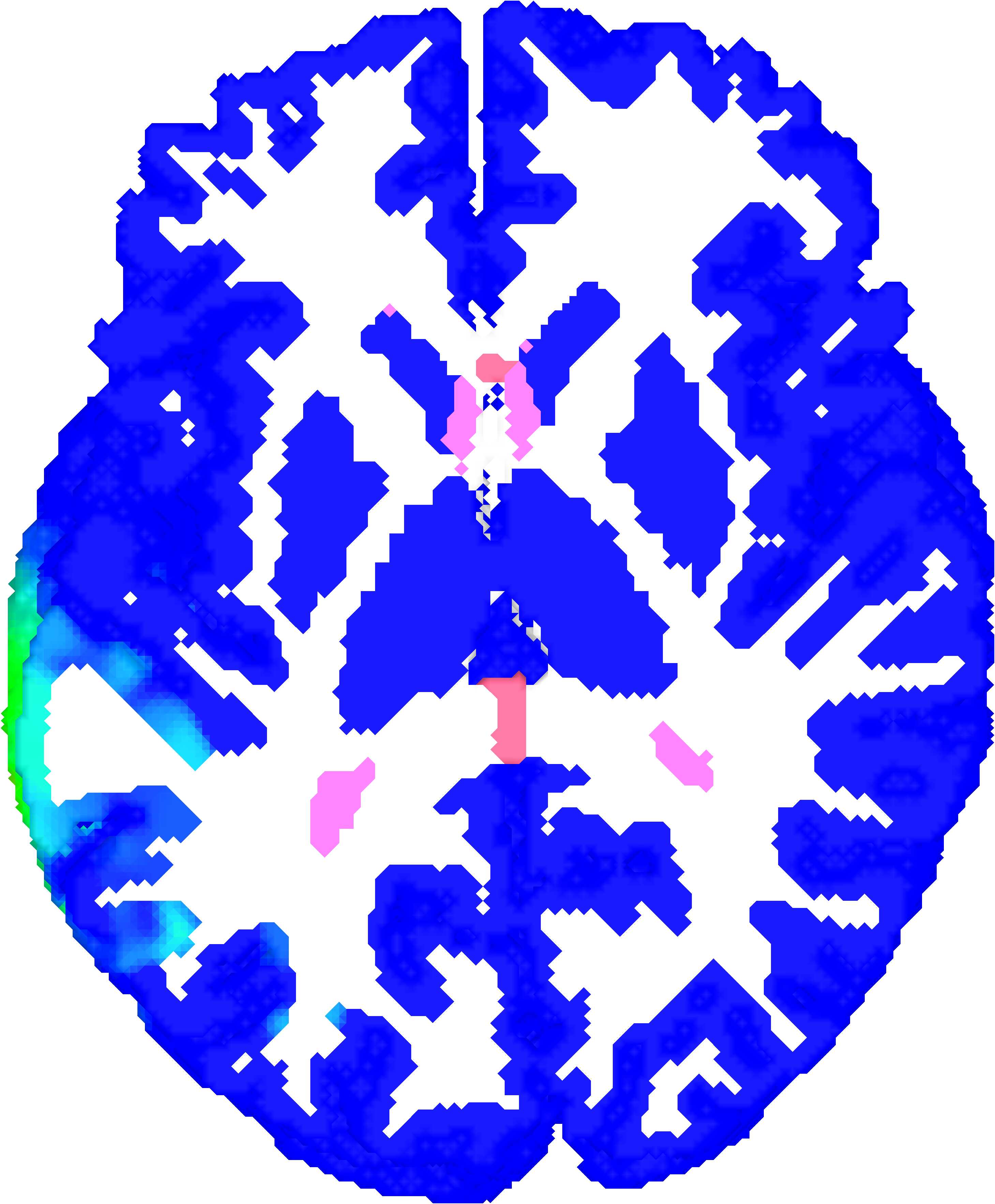}
            \caption{$\azimuth=\qty{0}{\degree}$, $\elevation=\qty{90}{\degree}$, $z=\qty{37}{\milli\meter}$.}
        \end{subfigure}
        \hfill
        \begin{subfigure}[b]{\subFigWidth}
            \includegraphics[width=\linewidth]{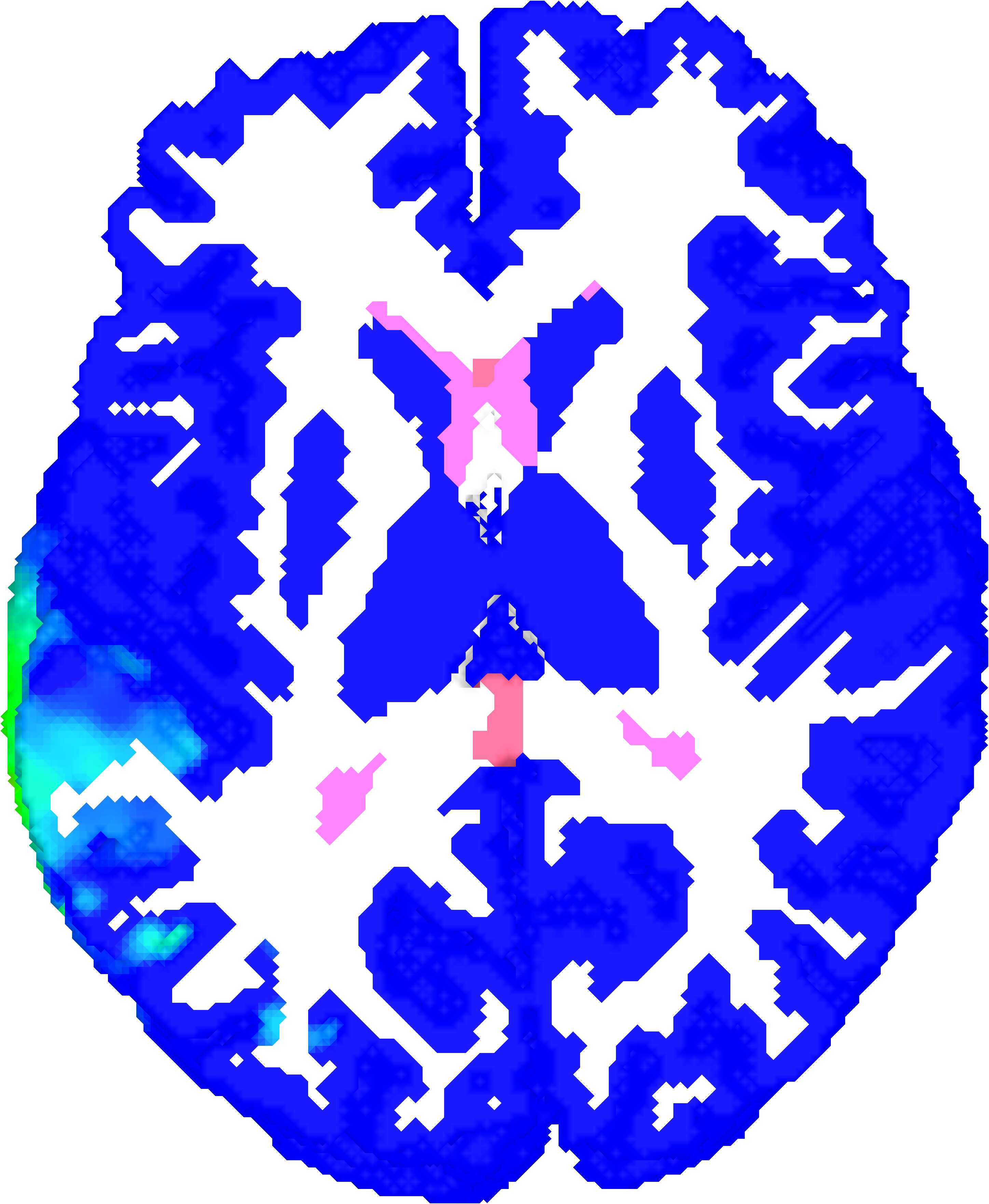}
            \caption{$\azimuth=\qty{0}{\degree}$, $\elevation=\qty{90}{\degree}$, $z=\qty{40}{\milli\meter}$.}
        \end{subfigure}
    \end{minipage}
    \hfill
    \begin{minipage}{0.09\linewidth}
    \includegraphics[width=\linewidth]{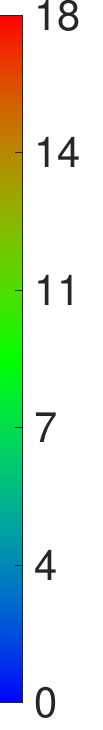}
    \end{minipage}
    \caption{Relative differences \eqref{eq:rel-diff} between $\linL$ and $\refL$, smoothed via the function $\dBT\args{\someField/\max\someField,\qty{18}{\decibel}}$ of equation \eqref{eq:visual-log-transform} at different azimuths $\azimuth$ and elevations $\elevation$, with base stimulation frequency $\frequency=\qty{1000}{\hertz}$ and an updated $\contactResistance=\qty{1270}{\ohm}$. The maximal value for the absolute difference of the lead fields was approximately \qty{0.13}{\per\meter\squared}, when  $\norm\leadFieldMat\approx\qty{0.19e-3}{}$, $\norm\refL\approx\qty{19e-3}{}$ and $\norm\linL\approx\qty{0.21e-3}{\per\meter\squared}$.}
    \label{fig:linrefLreldiff-Rc=1270}
\end{figure}

\begin{figure}[!h]
    \centering
    \begin{minipage}{0.85\linewidth}
        \def\subFigWidth{0.31\linewidth}
        \begin{subfigure}[b]{\subFigWidth}
            \includegraphics[width=\linewidth]{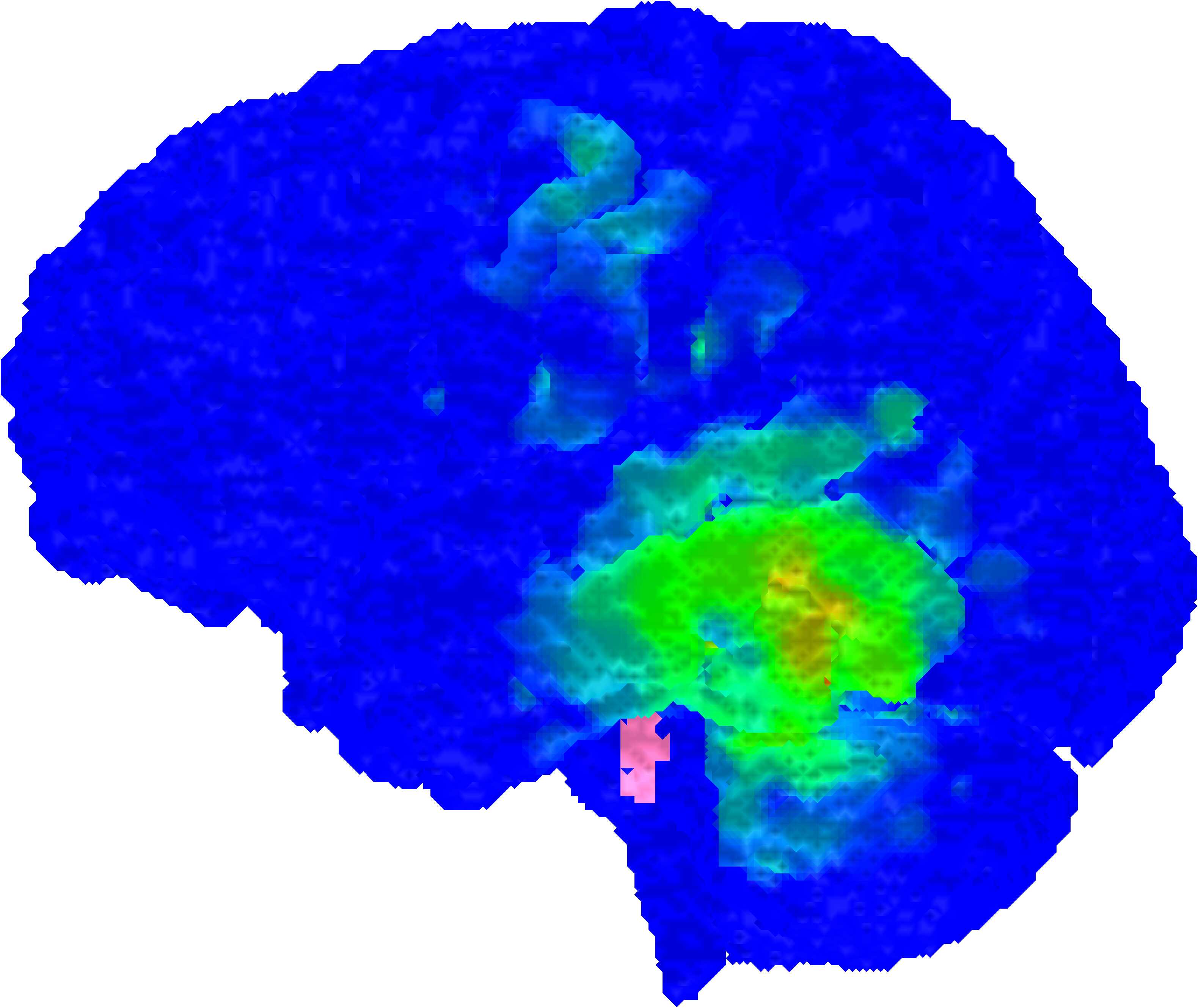}
            \caption{$\azimuth=\qty{-90}{\degree}$, $\elevation=\qty{0}{\degree}$.}
        \end{subfigure}
        \hfill
        \begin{subfigure}[b]{\subFigWidth}
        \includegraphics[width=\linewidth]{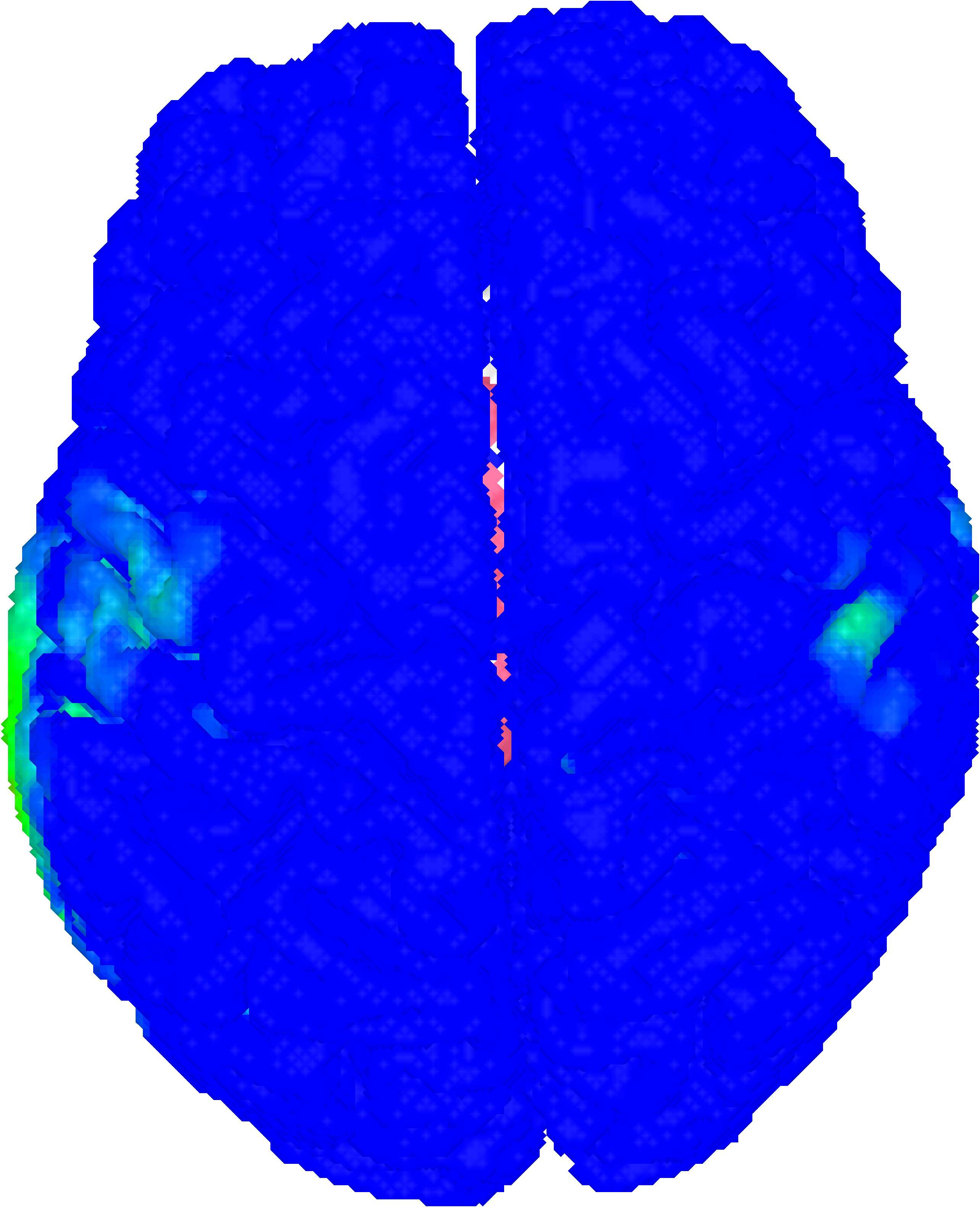}
        \caption{$\azimuth=\qty{0}{\degree}$, $\elevation=\qty{90}{\degree}$.}
        \end{subfigure}
        \hfill
        \begin{subfigure}[b]{\subFigWidth}
        \includegraphics[width=\linewidth]{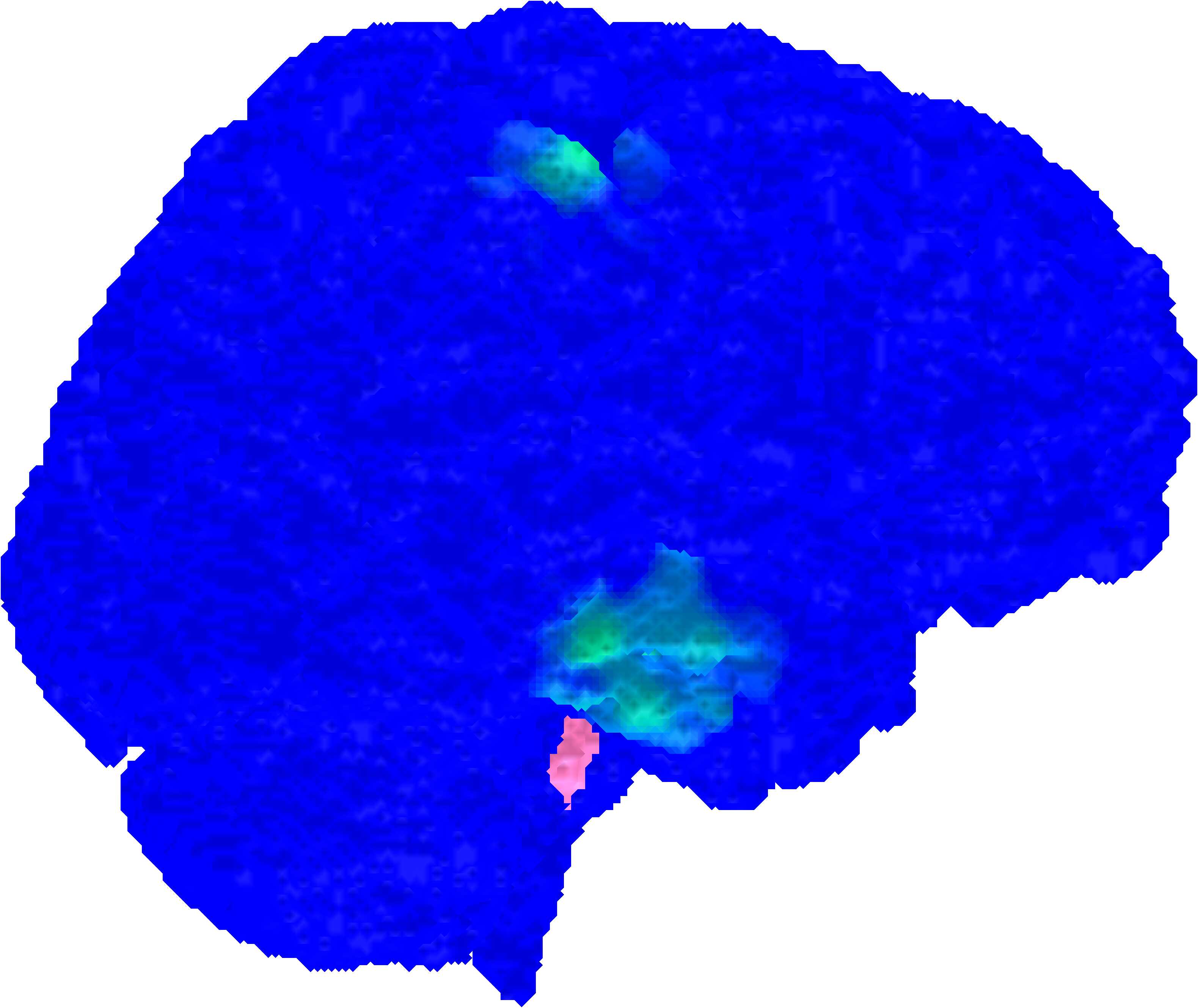}
            \caption{$\azimuth=\qty{90}{\degree}$, $\elevation=\qty{0}{\degree}$.}
        \end{subfigure}

        \def\subFigWidth{0.45\linewidth}

        \begin{subfigure}[b]{\subFigWidth}
            \includegraphics[width=\linewidth]{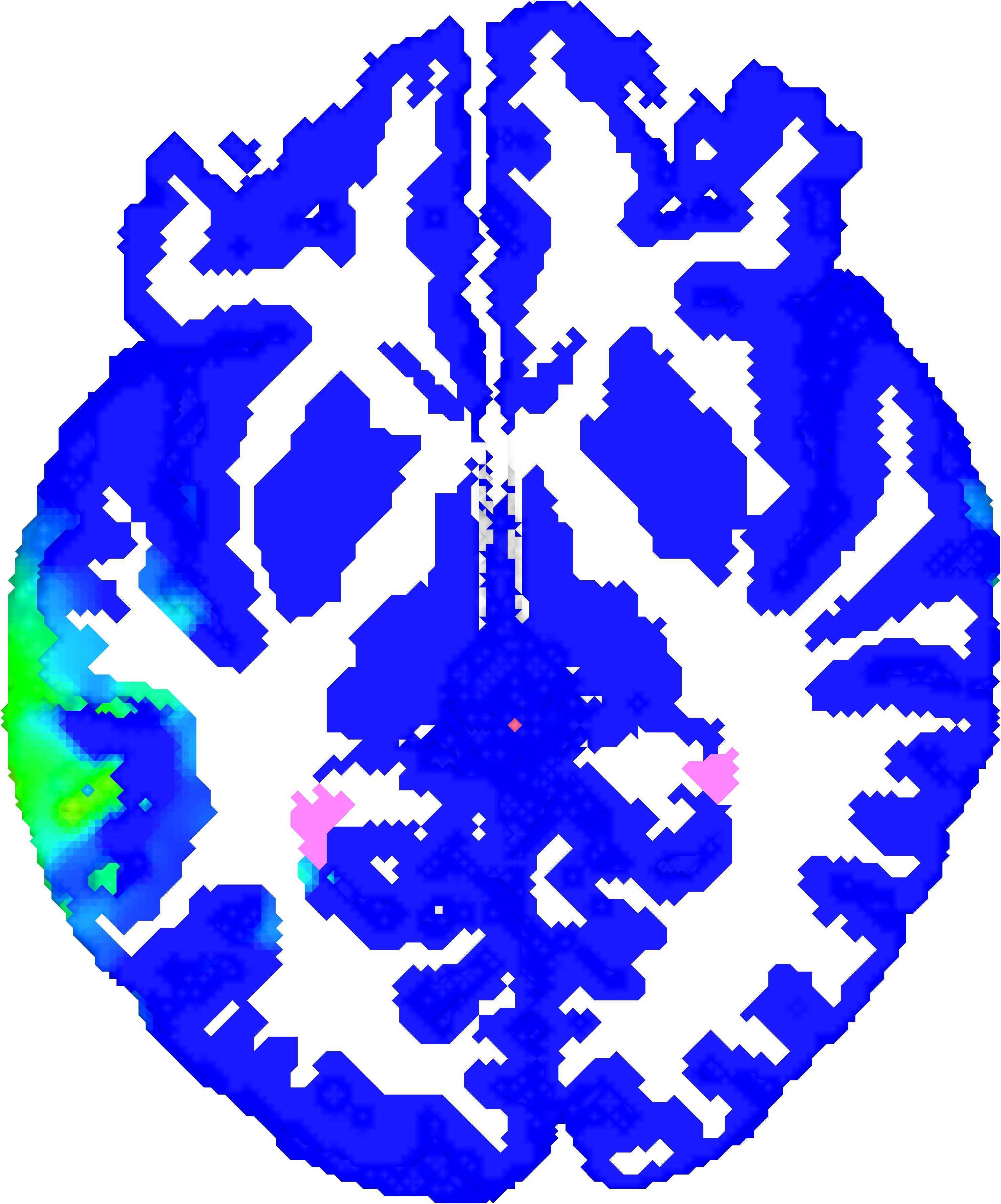}
            \caption{$\azimuth=\qty{0}{\degree}$, $\elevation=\qty{90}{\degree}$, $z=\qty{30}{\milli\meter}$.}
        \end{subfigure}
       \hfill
        \begin{subfigure}[b]{\subFigWidth}
            \includegraphics[width=\linewidth]{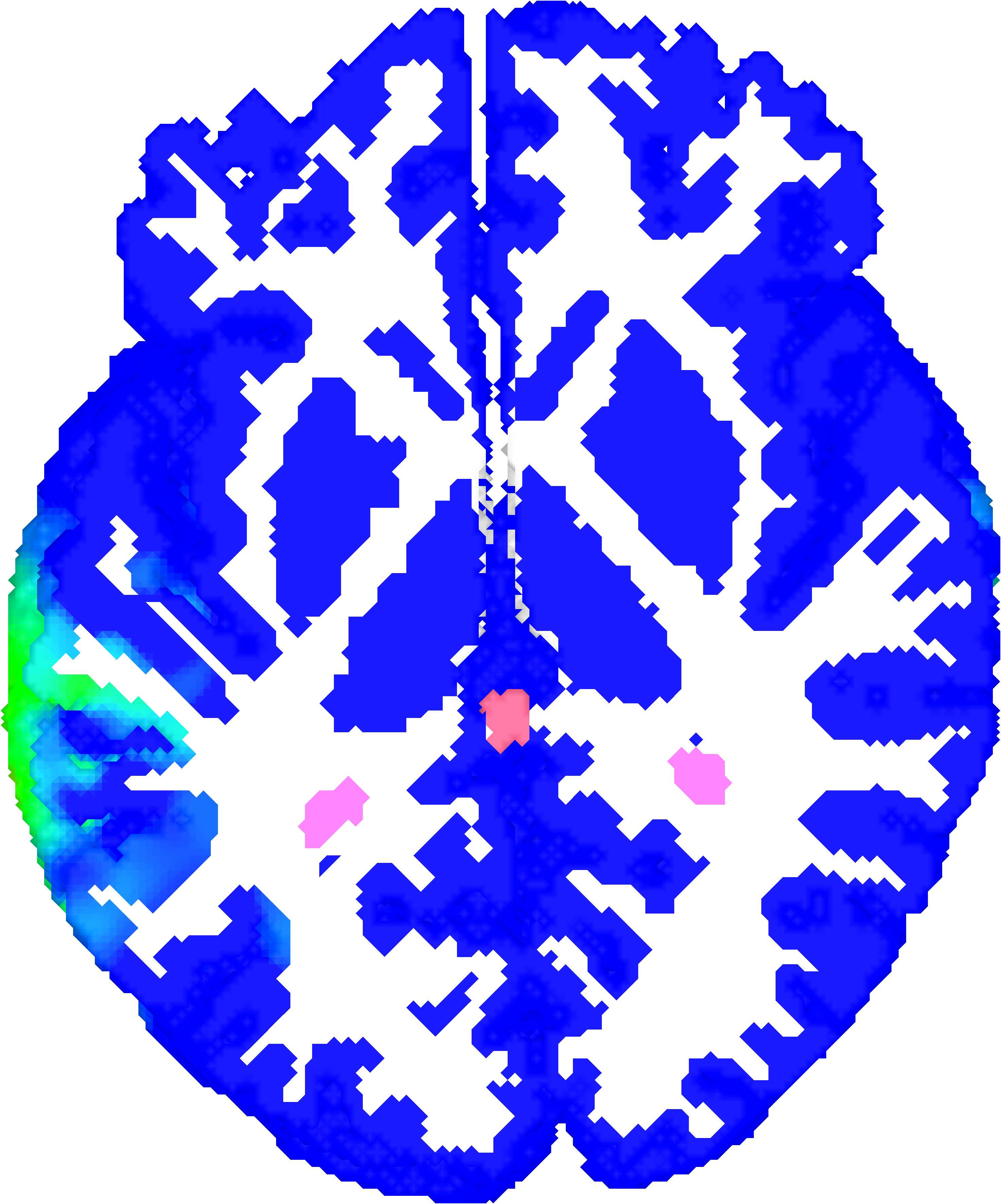}
            \caption{$\azimuth=\qty{0}{\degree}$, $\elevation=\qty{90}{\degree}$, $z=\qty{33}{\milli\meter}$.}
        \end{subfigure}

        \begin{subfigure}[b]{\subFigWidth}
            \includegraphics[width=\linewidth]{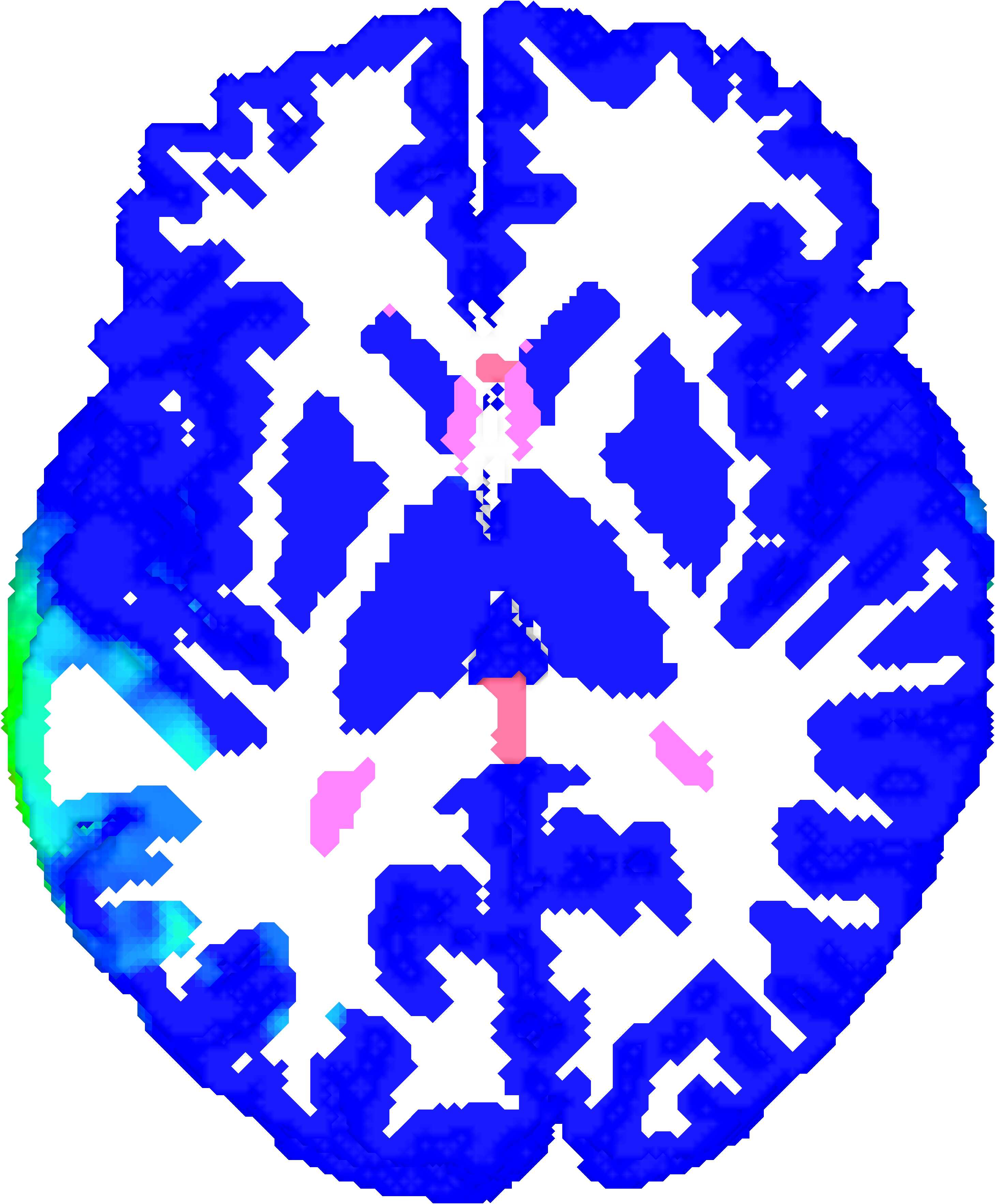}
            \caption{$\azimuth=\qty{0}{\degree}$, $\elevation=\qty{90}{\degree}$, $z=\qty{37}{\milli\meter}$.}
        \end{subfigure}
        \hfill
        \begin{subfigure}[b]{\subFigWidth}
            \includegraphics[width=\linewidth]{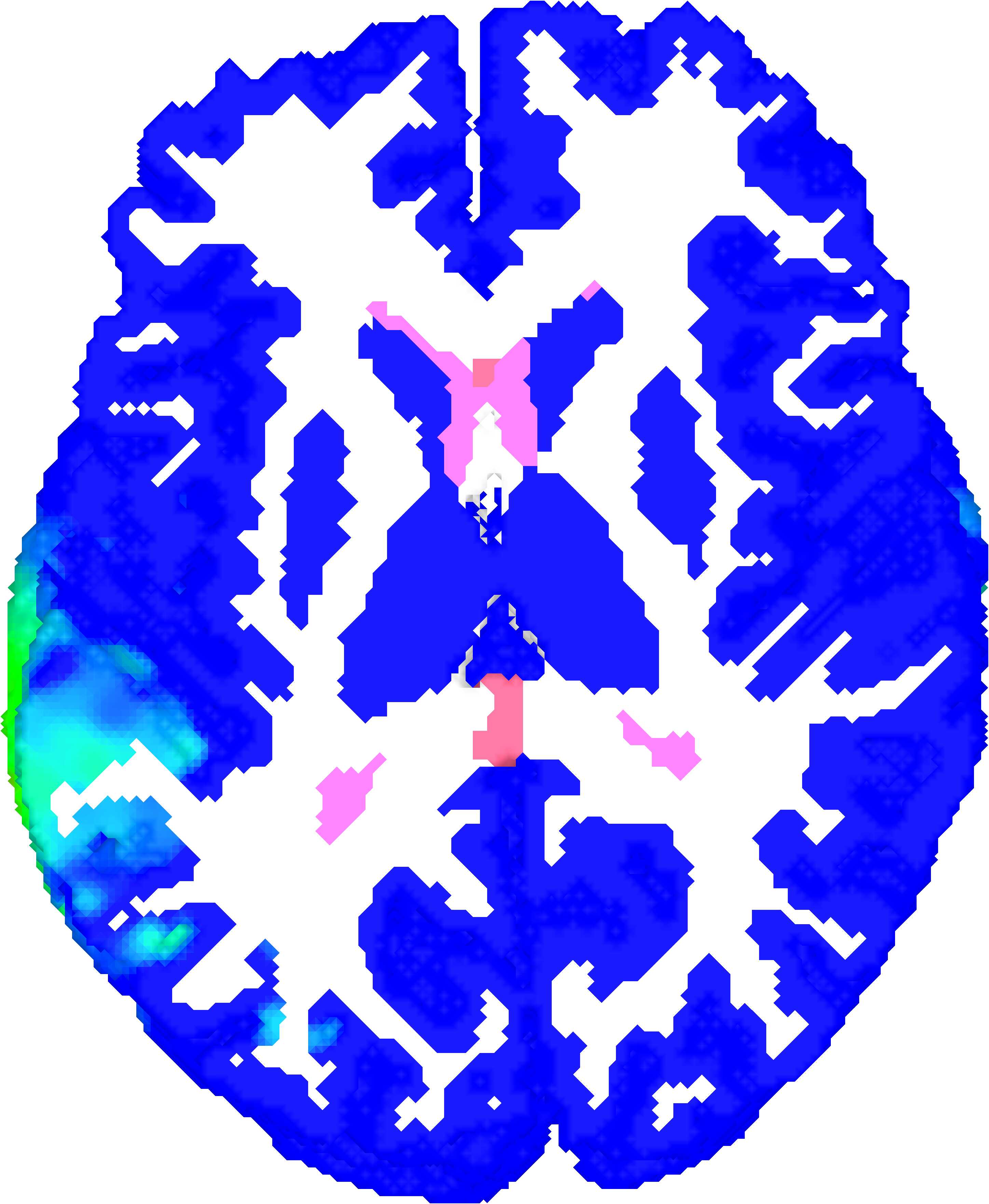}
            \caption{$\azimuth=\qty{0}{\degree}$, $\elevation=\qty{90}{\degree}$, $z=\qty{40}{\milli\meter}$.}
        \end{subfigure}
    \end{minipage}
    \hfill
    \begin{minipage}{0.09\linewidth}
    \includegraphics[width=\linewidth]{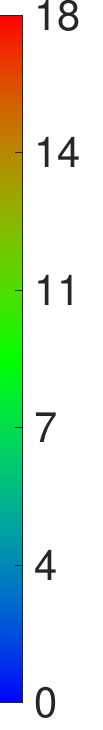}
    \end{minipage}
    \caption{Relative differences \eqref{eq:rel-diff} between $\linL$ and $\refL$, smoothed via the function $\dBT\args{\someField/\max\someField,\qty{18}{\decibel}}$ of equation \eqref{eq:visual-log-transform} at different azimuths $\azimuth$ and elevations $\elevation$, with base stimulation frequency $\frequency=\qty{1000}{\hertz}$ and an updated $\contactResistance=\qty{5270}{\ohm}$. The maximal value for the absolute difference of the lead fields was approximately \qty{0.16e-3}{\per\meter\squared}, when  $\norm\iniL\approx\qty{0.17e-3}{\per\meter\squared}$, $\norm\refL\approx\qty{21e-3}{\per\meter\squared}$ and $\norm\linL\approx\qty{0.21e-3}{\per\meter\squared}$.}
    \label{fig:linrefLreldiff-Rc=5270}
\end{figure}

Having observed how the initial interference
fields with contact resistances at a steady
$\contactResistance=\qty{270}{\ohm}$ behave,
Figures~\ref{fig:linrefLreldiff-Rc=1270} and
\ref{fig:linrefLreldiff-Rc=5270} display relative
differences~\eqref{eq:rel-diff} between the column
norms of the lead fields $\linL$ and $\refL$
with the contact resistance having increased to
$\contactResistance\in\set{\num{1270},\num{5270}}\,\si\ohm$
at the modified electrode TP9. We see the greatest changes
near the modified electrode itself, whereas the effects
on the rest of the volume $\domain$ seem relatively
unchanged at $\contactResistance=\qty{1270}{\ohm}$. At
$\contactResistance=\qty{5270}{\ohm}$, visible differences
are also seen near the electrodes that were modified.
However, at this sensitivity scale there is no visible
relative difference between $\linL$ and $\refL$ in the
deeper brain compartments in either case.

\begin{figure}[!h]
    \centering
    \begin{minipage}{0.85\linewidth}
        \def\subFigWidth{0.31\linewidth}
        \begin{subfigure}[b]{\subFigWidth}
            \includegraphics[width=\linewidth]{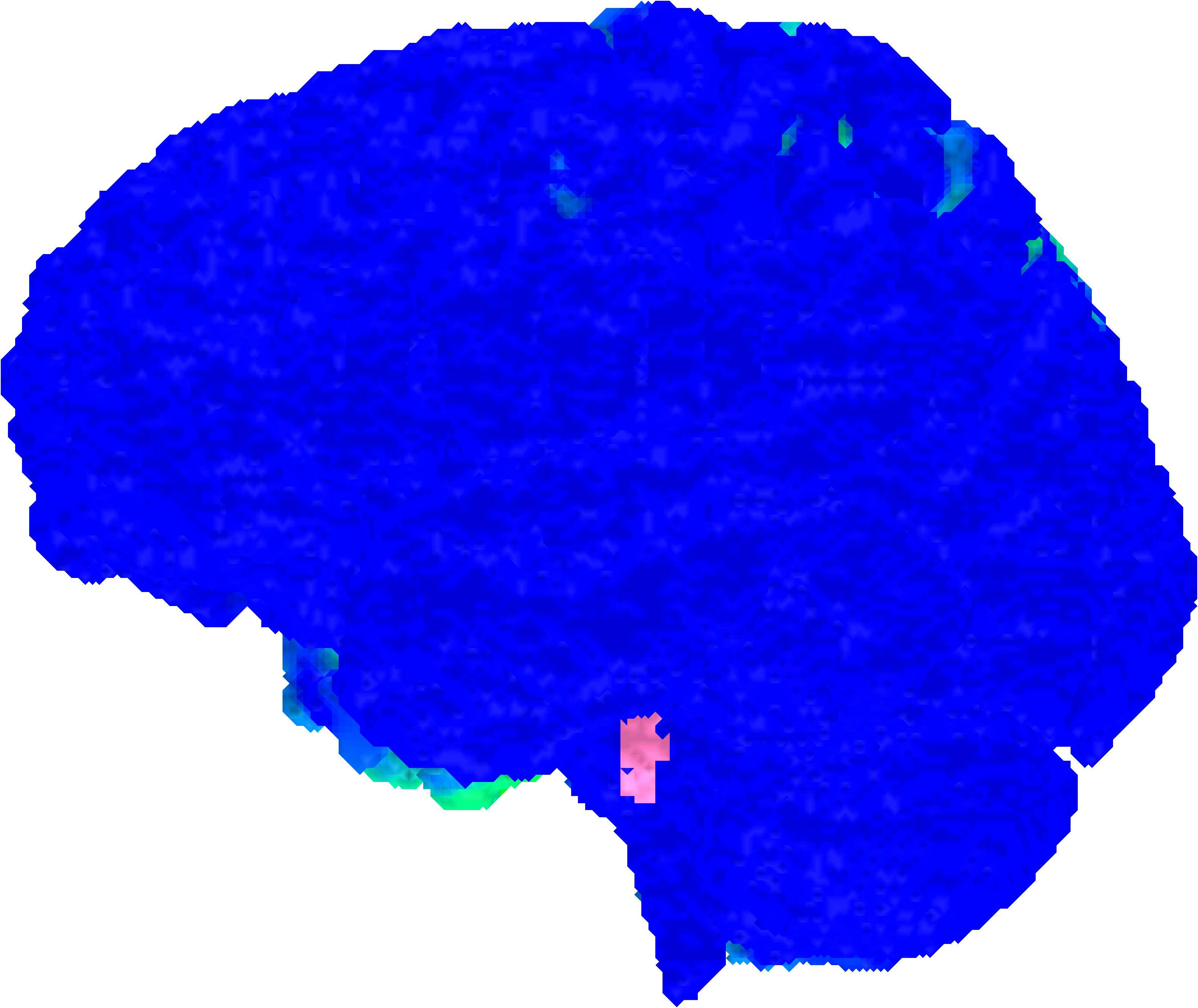}
            \caption{$\azimuth=\qty{-90}{\degree}$, $\elevation=\qty{0}{\degree}$.}
        \end{subfigure}
        \hfill
        \begin{subfigure}[b]{\subFigWidth}
        \includegraphics[width=\linewidth]{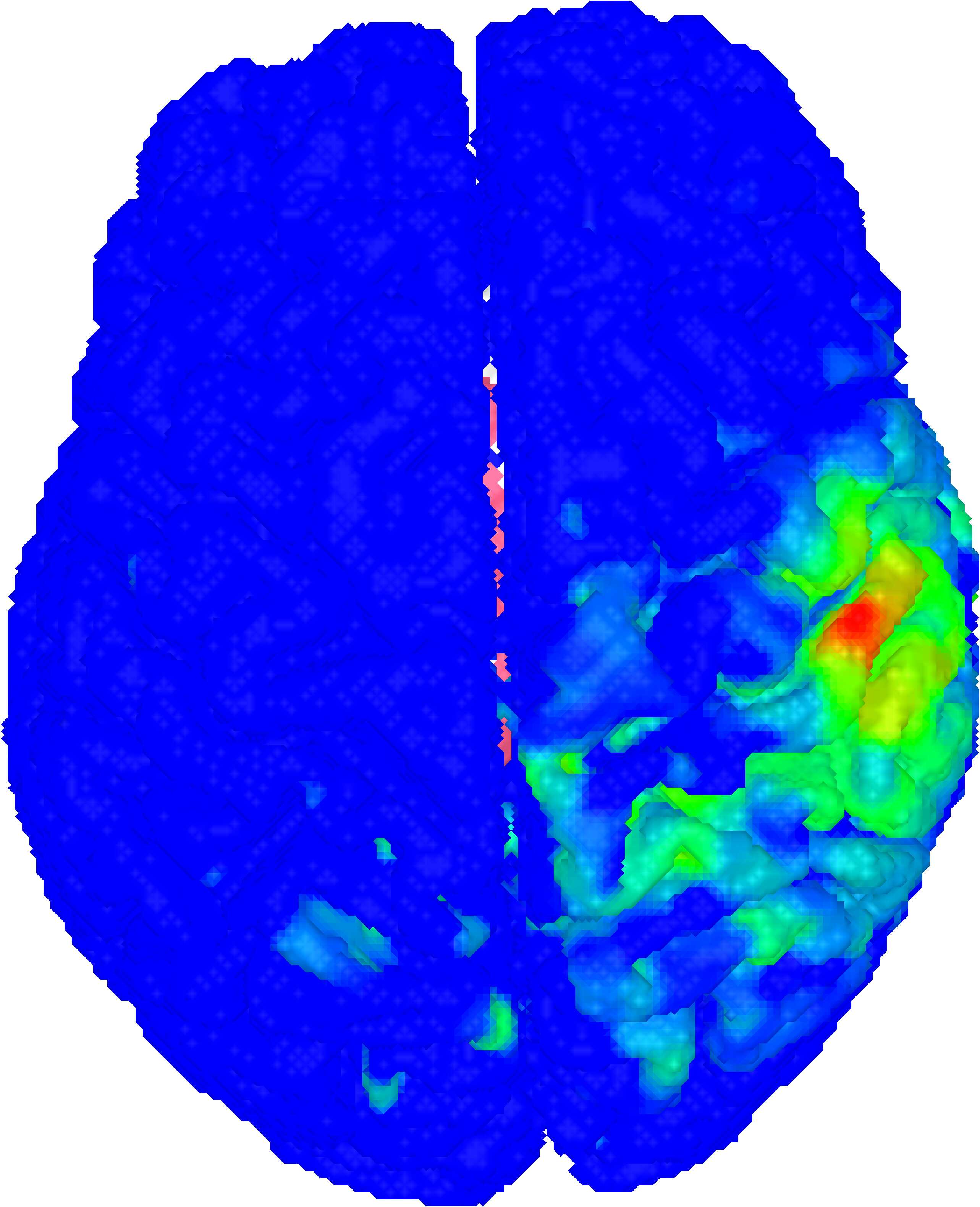}
        \caption{$\azimuth=\qty{0}{\degree}$, $\elevation=\qty{90}{\degree}$.}
        \end{subfigure}
        \hfill
        \begin{subfigure}[b]{\subFigWidth}
        \includegraphics[width=\linewidth]{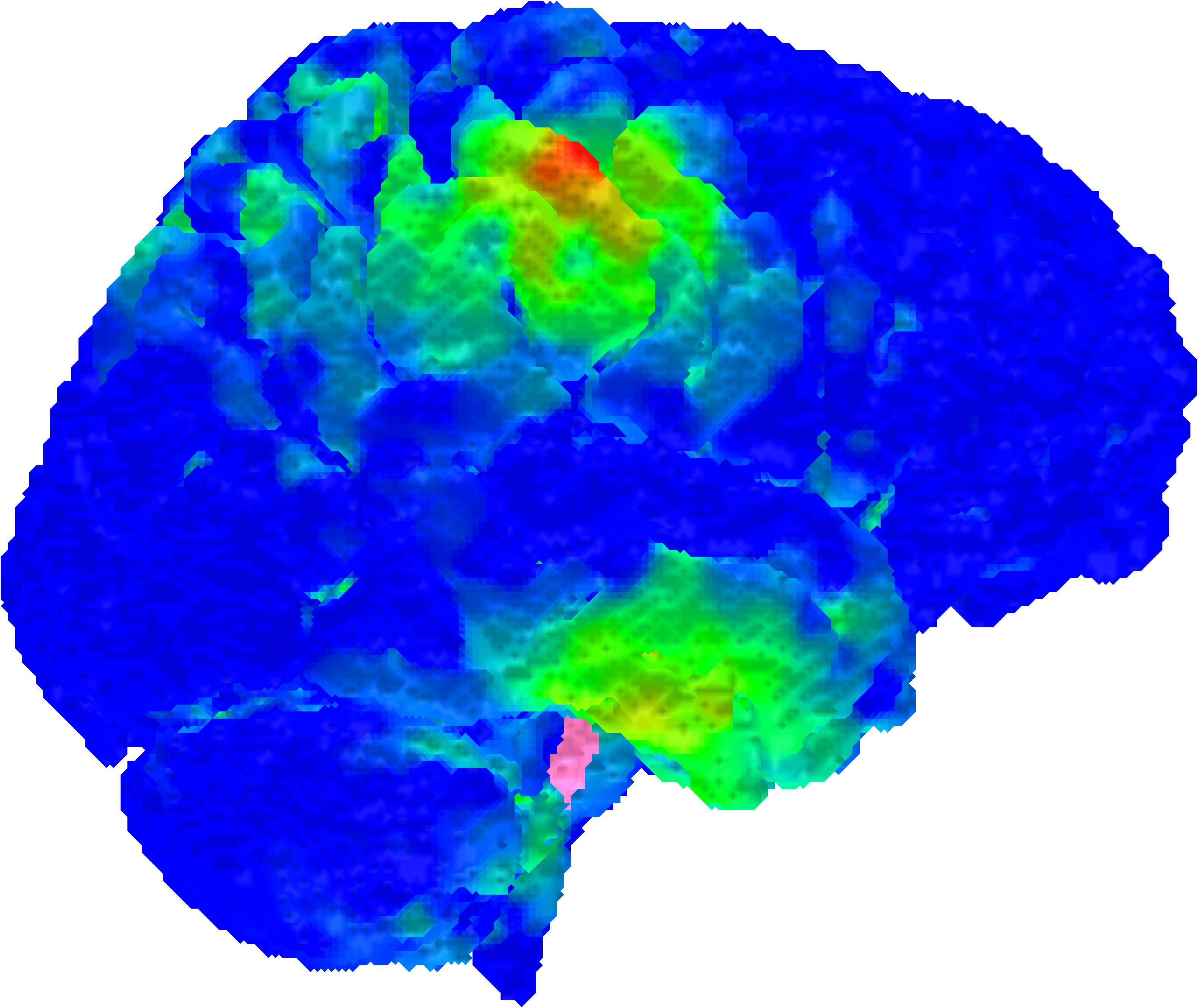}
            \caption{$\azimuth=\qty{90}{\degree}$, $\elevation=\qty{0}{\degree}$.}
        \end{subfigure}

        \def\subFigWidth{0.45\linewidth}

        \begin{subfigure}[b]{\subFigWidth}
            \includegraphics[width=\linewidth]{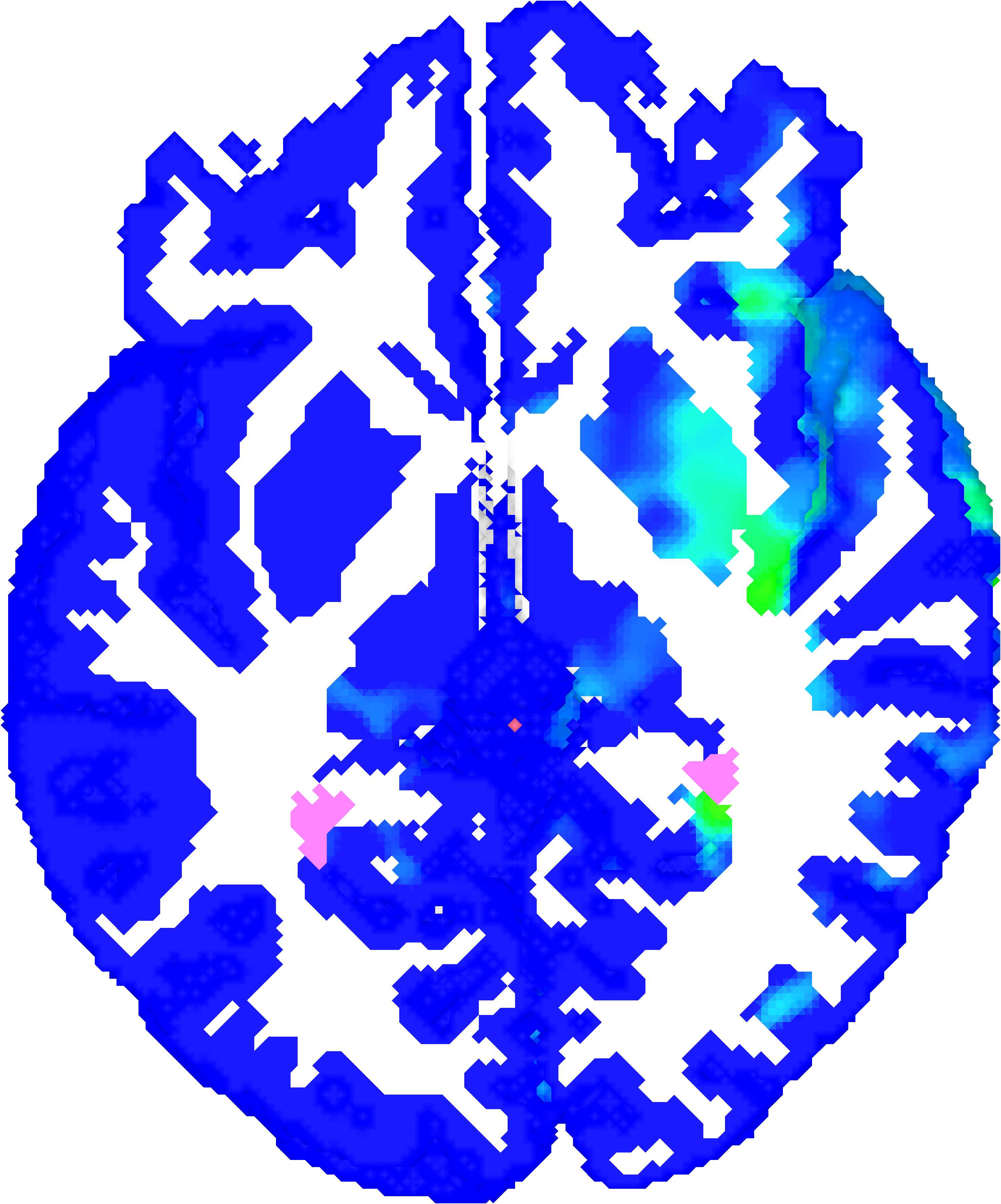}
            \caption{$\azimuth=\qty{0}{\degree}$, $\elevation=\qty{90}{\degree}$, $z=\qty{30}{\milli\meter}$.}
        \end{subfigure}
       \hfill
        \begin{subfigure}[b]{\subFigWidth}
            \includegraphics[width=\linewidth]{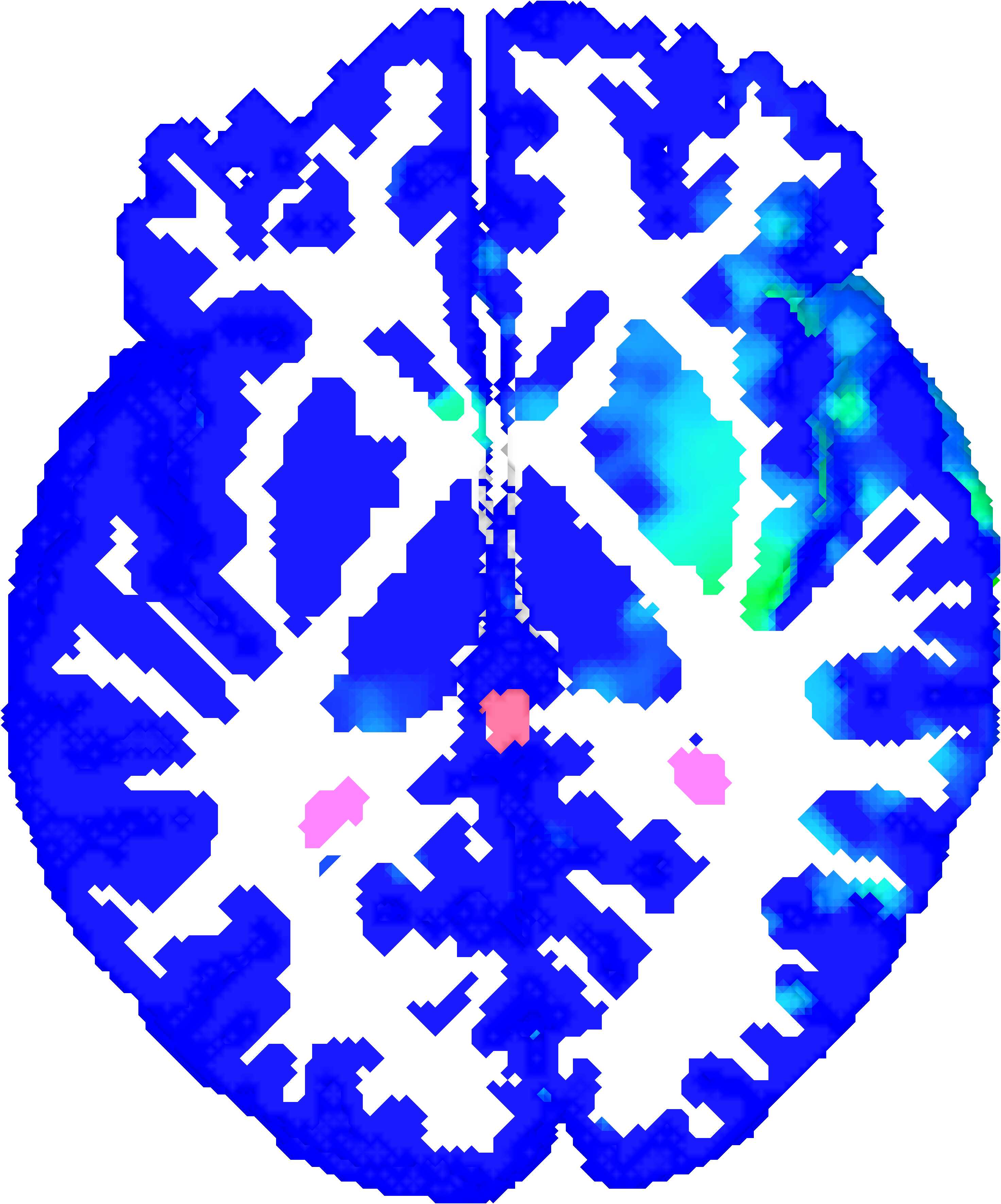}
            \caption{$\azimuth=\qty{0}{\degree}$, $\elevation=\qty{90}{\degree}$, $z=\qty{33}{\milli\meter}$.}
        \end{subfigure}

        \begin{subfigure}[b]{\subFigWidth}
            \includegraphics[width=\linewidth]{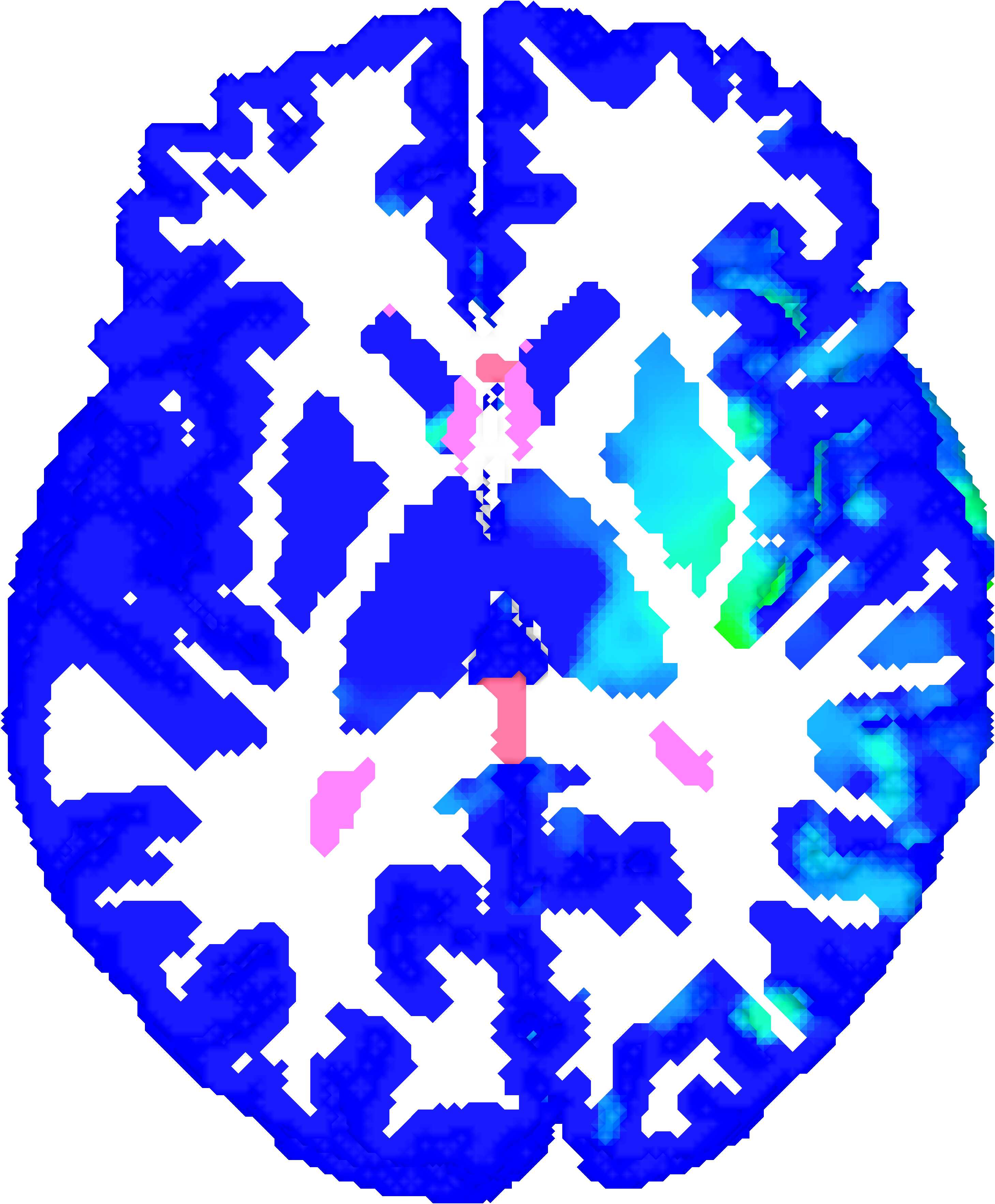}
            \caption{$\azimuth=\qty{0}{\degree}$, $\elevation=\qty{90}{\degree}$, $z=\qty{37}{\milli\meter}$.}
        \end{subfigure}
        \hfill
        \begin{subfigure}[b]{\subFigWidth}
            \includegraphics[width=\linewidth]{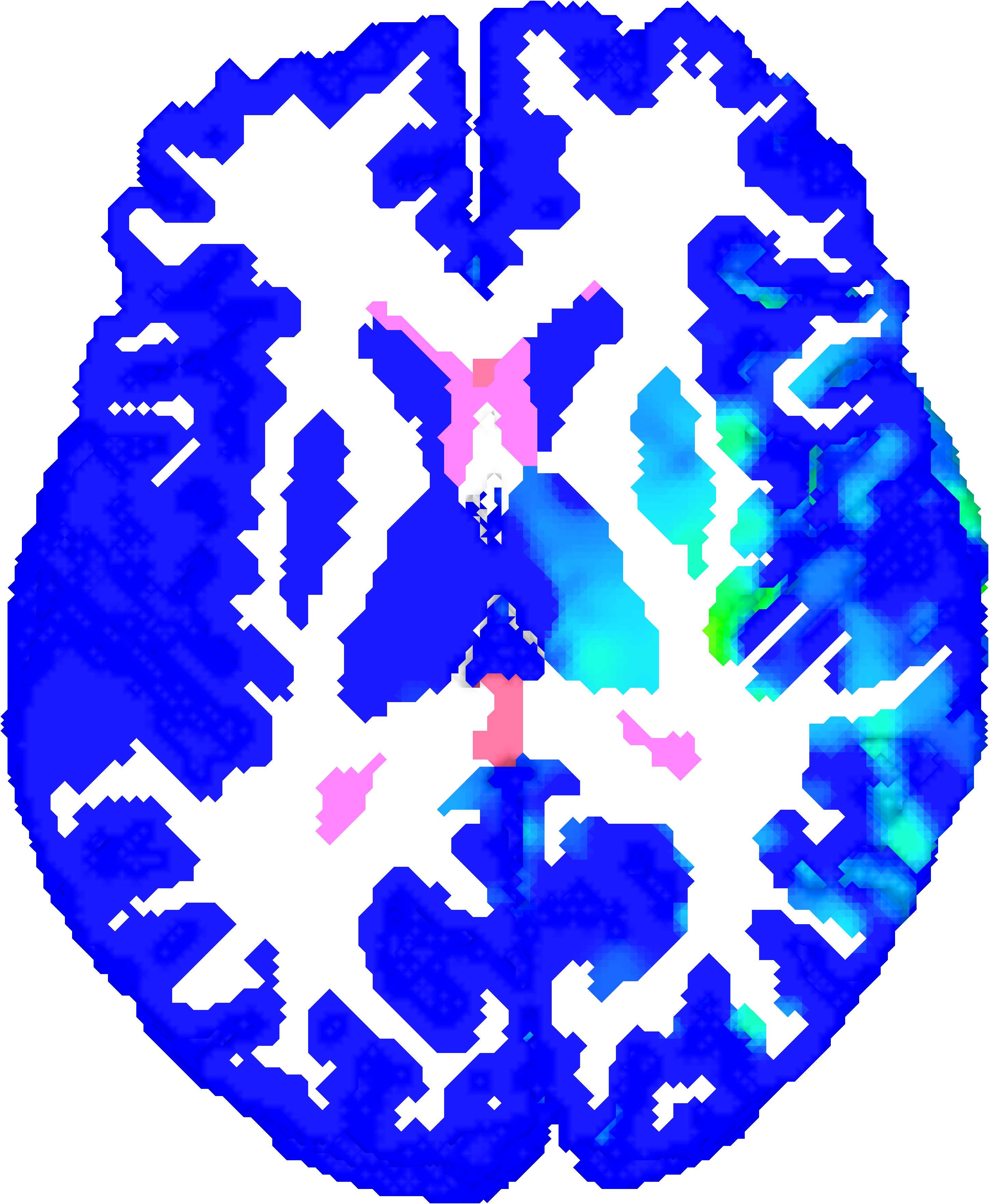}
            \caption{$\azimuth=\qty{0}{\degree}$, $\elevation=\qty{90}{\degree}$, $z=\qty{40}{\milli\meter}$.}
        \end{subfigure}
    \end{minipage}
    \hfill
    \begin{minipage}{0.09\linewidth}
    \includegraphics[width=\linewidth]{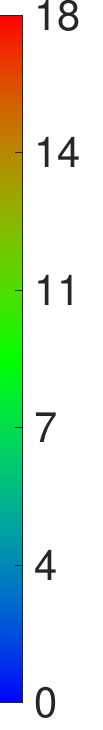}
    \end{minipage}
    \caption{Relative differences \eqref{eq:rel-diff} between the interference fields \eqref{eq:interference-field} of $\linJv$ and $\refJv$, smoothed via the function $\dBT\args{\someField/\max\someField,\qty{18}{\decibel}}$ of equation \eqref{eq:visual-log-transform} at different azimuths $\azimuth$ and elevations $\elevation$, with base stimulation frequency $\frequency=\qty{1000}{\hertz}$ and an updated $\contactResistance=\qty{1270}{\ohm}$. The current pattern was $\currentPattern=[-1.0,1.0,-1.0,1.0]\,\si{\milli\ampere}$. The maximal value of the absolute differences was roughly \qty{2.08e-6}{\milli\ampere\per\milli\meter\squared}.}
    \label{fig:linrefreldiff-Rc=1270-i=1.0-1.0}
\end{figure}

\begin{figure}[!h]
    \centering
    \begin{minipage}{0.85\linewidth}
        \def\subFigWidth{0.31\linewidth}
        \begin{subfigure}[b]{\subFigWidth}
            \includegraphics[width=\linewidth]{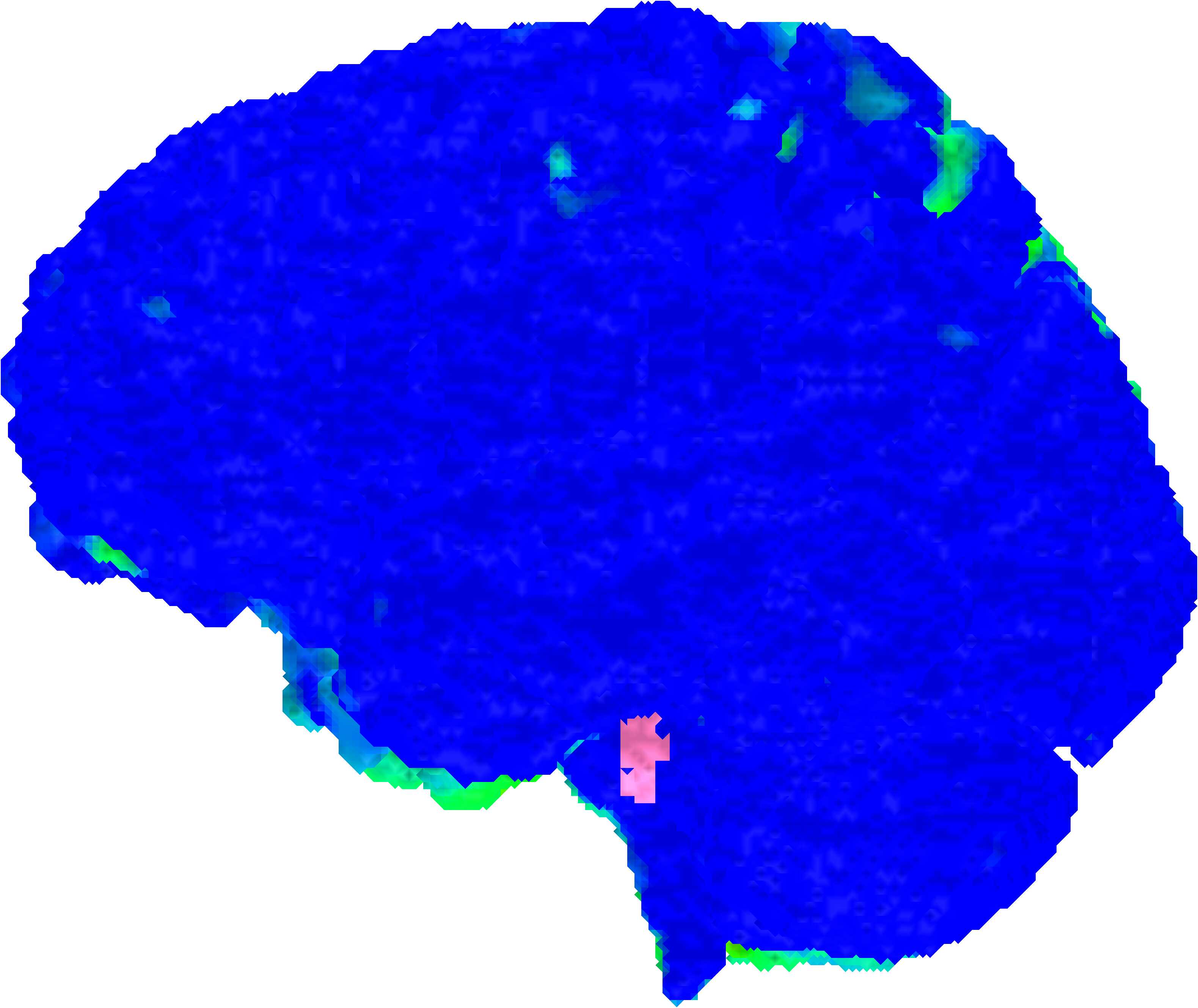}
            \caption{$\azimuth=\qty{-90}{\degree}$, $\elevation=\qty{0}{\degree}$.}
        \end{subfigure}
        \hfill
        \begin{subfigure}[b]{\subFigWidth}
        \includegraphics[width=\linewidth]{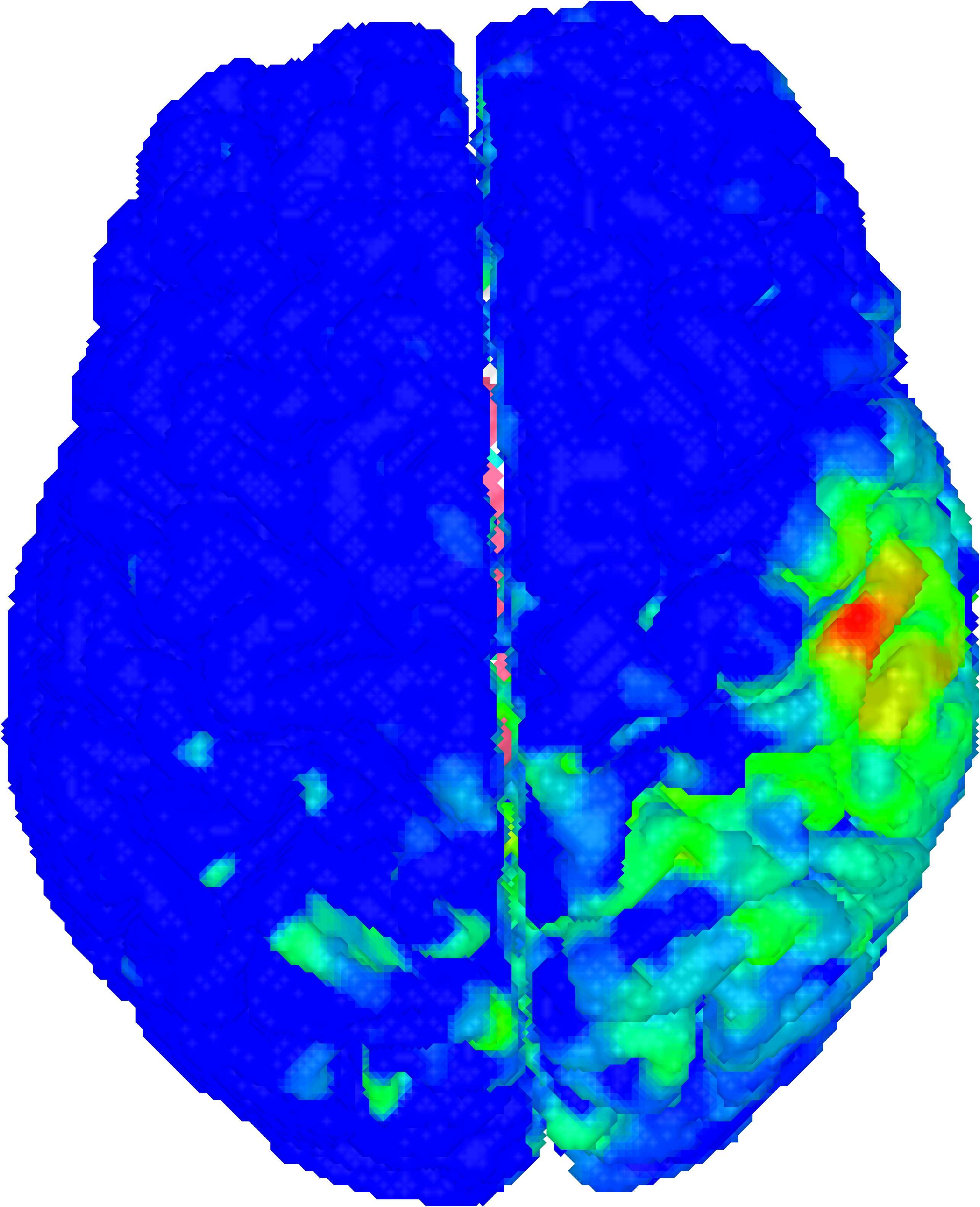}
        \caption{$\azimuth=\qty{0}{\degree}$, $\elevation=\qty{90}{\degree}$.}
        \end{subfigure}
        \hfill
        \begin{subfigure}[b]{\subFigWidth}
        \includegraphics[width=\linewidth]{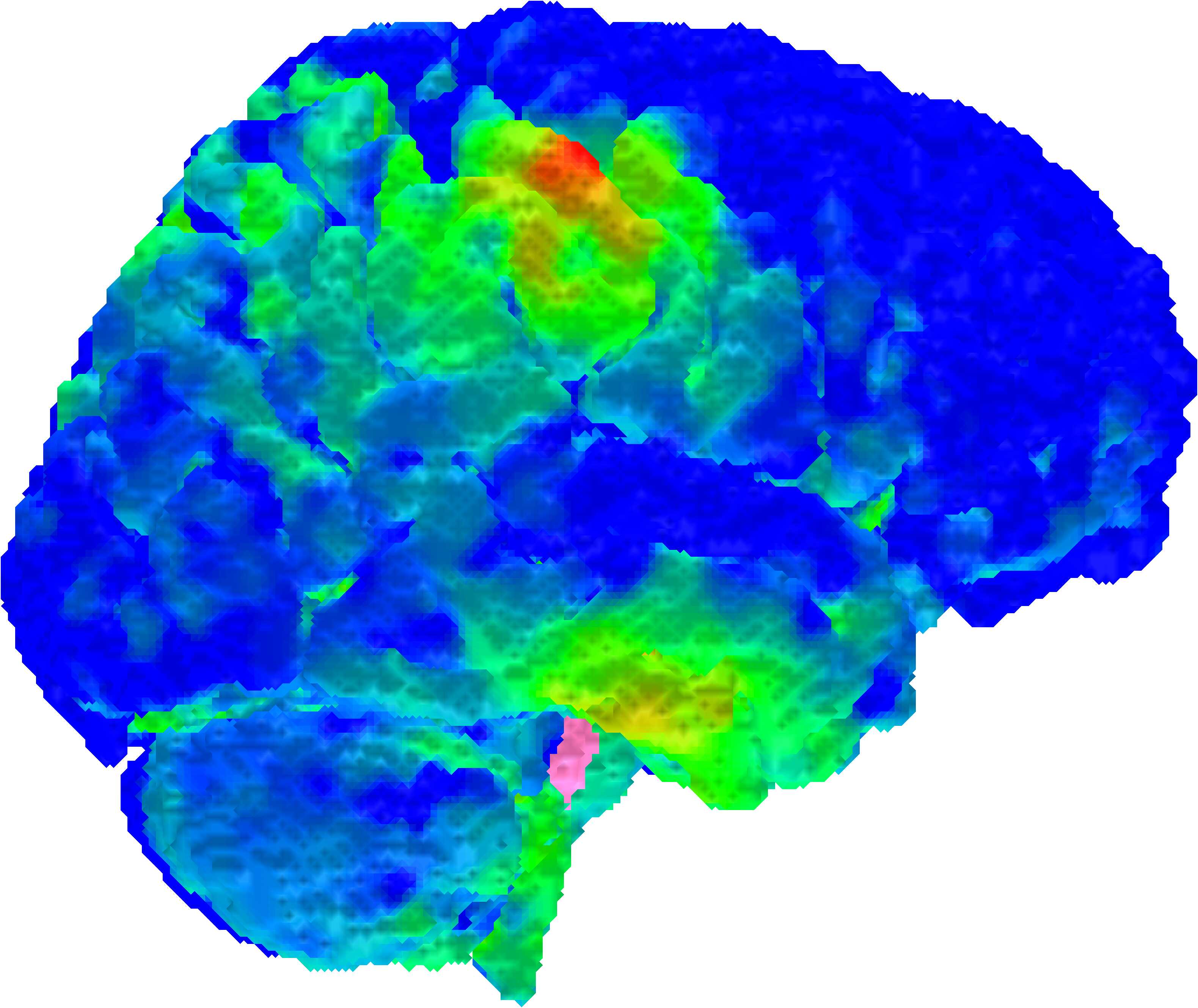}
            \caption{$\azimuth=\qty{90}{\degree}$, $\elevation=\qty{0}{\degree}$.}
        \end{subfigure}

        \def\subFigWidth{0.45\linewidth}

        \begin{subfigure}[b]{\subFigWidth}
            \includegraphics[width=\linewidth]{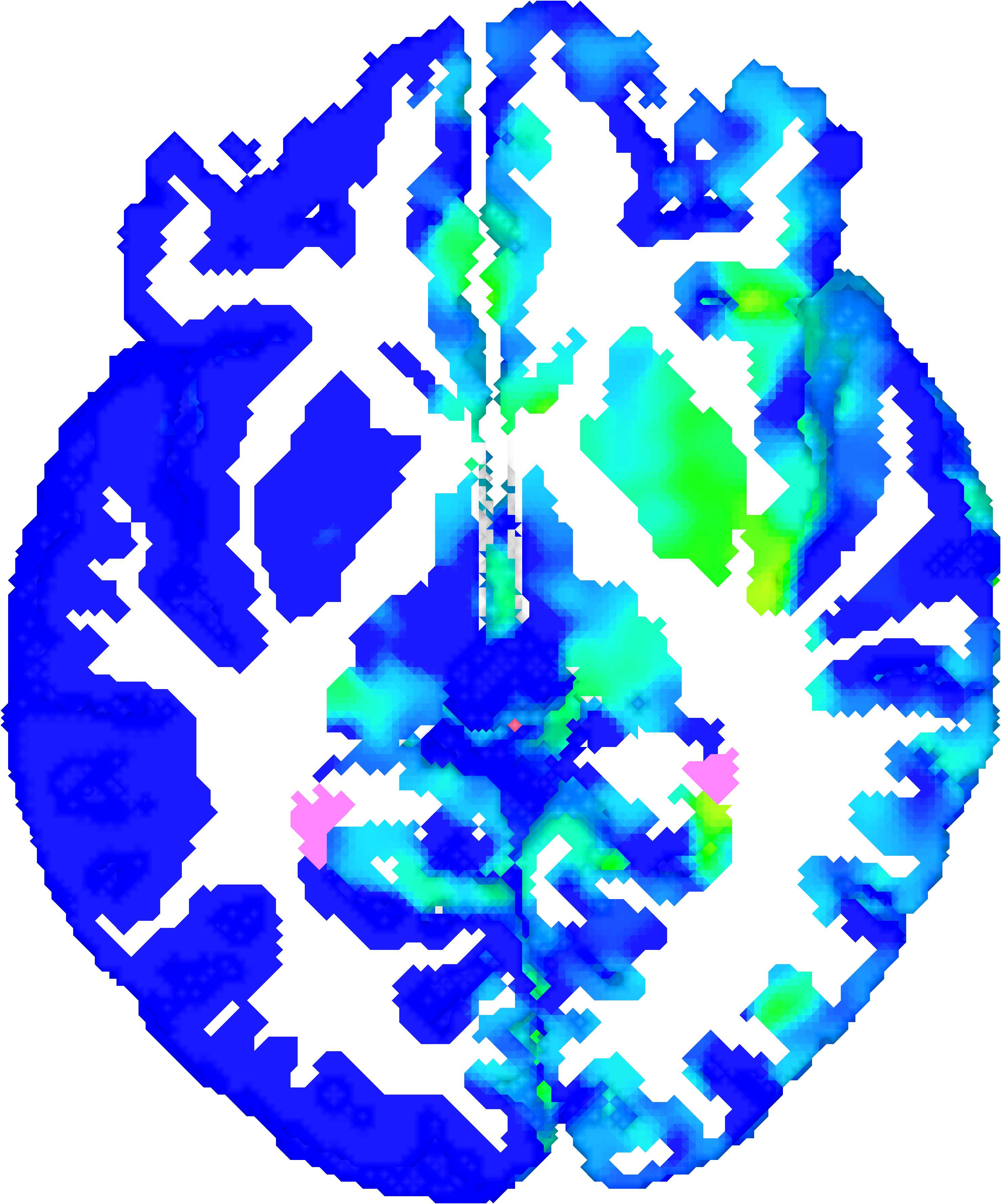}
            \caption{$\azimuth=\qty{0}{\degree}$, $\elevation=\qty{90}{\degree}$, $z=\qty{30}{\milli\meter}$.}
        \end{subfigure}
       \hfill
        \begin{subfigure}[b]{\subFigWidth}
            \includegraphics[width=\linewidth]{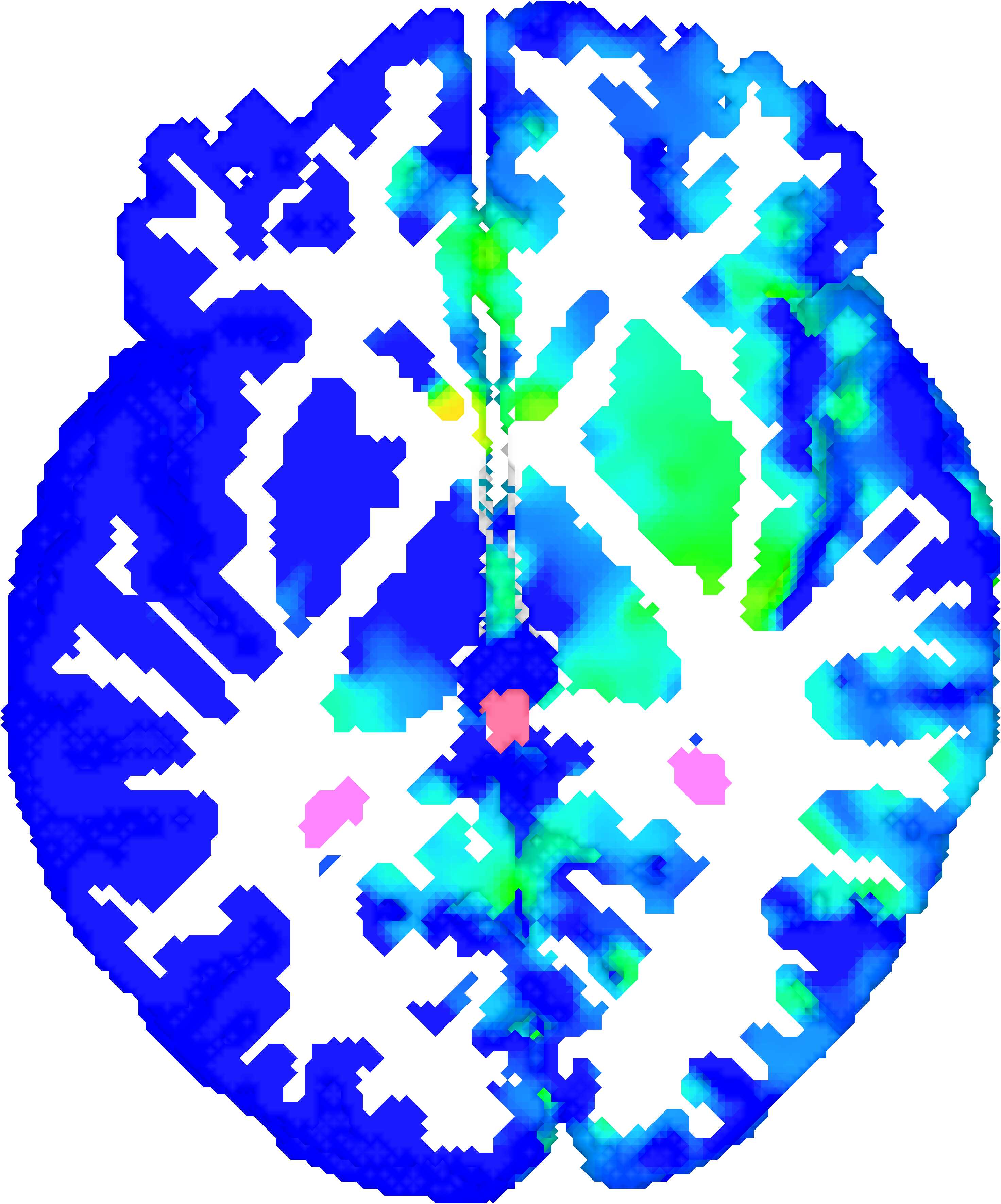}
            \caption{$\azimuth=\qty{0}{\degree}$, $\elevation=\qty{90}{\degree}$, $z=\qty{33}{\milli\meter}$.}
        \end{subfigure}

        \begin{subfigure}[b]{\subFigWidth}
            \includegraphics[width=\linewidth]{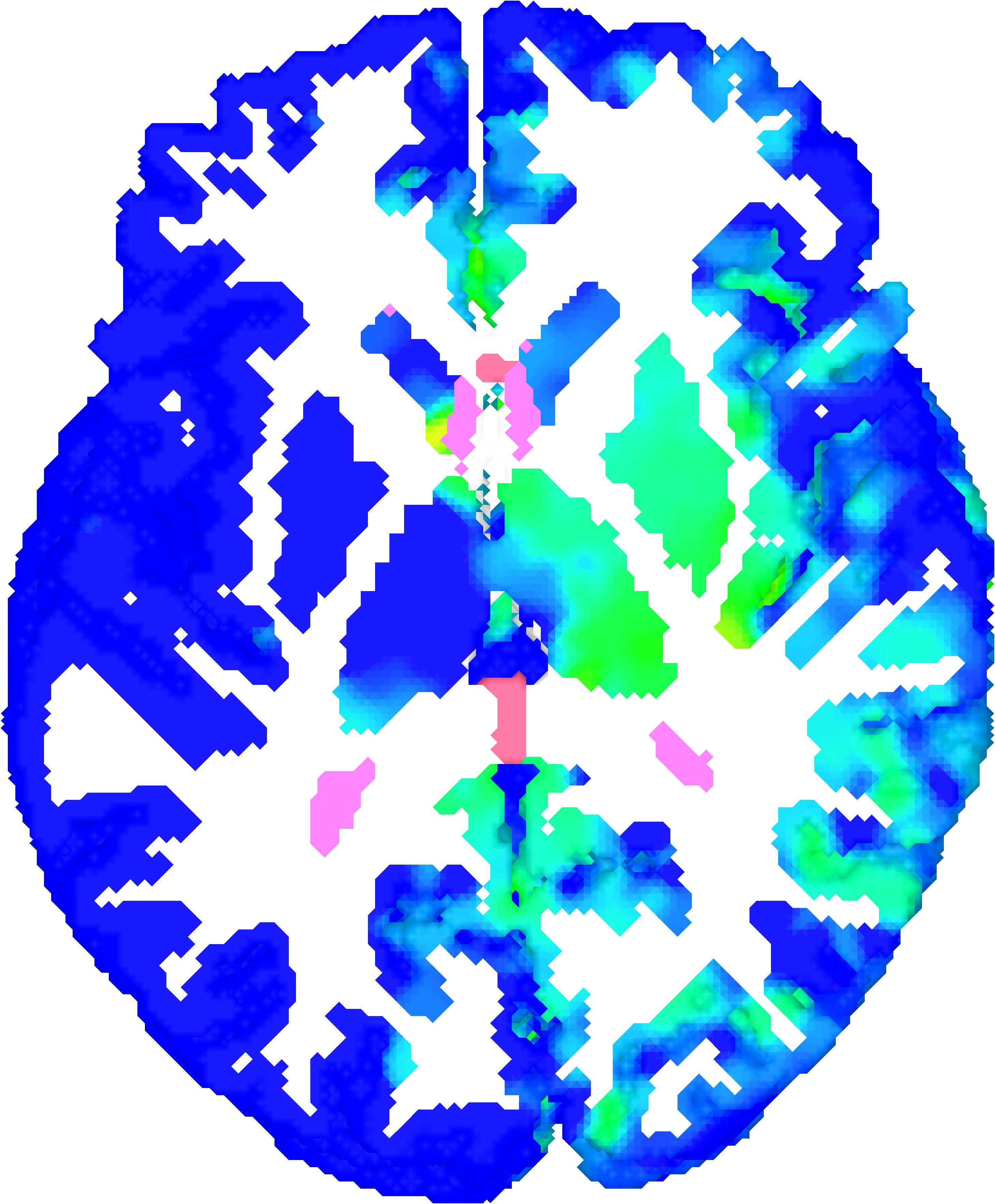}
            \caption{$\azimuth=\qty{0}{\degree}$, $\elevation=\qty{90}{\degree}$, $z=\qty{37}{\milli\meter}$.}
        \end{subfigure}
        \hfill
        \begin{subfigure}[b]{\subFigWidth}
            \includegraphics[width=\linewidth]{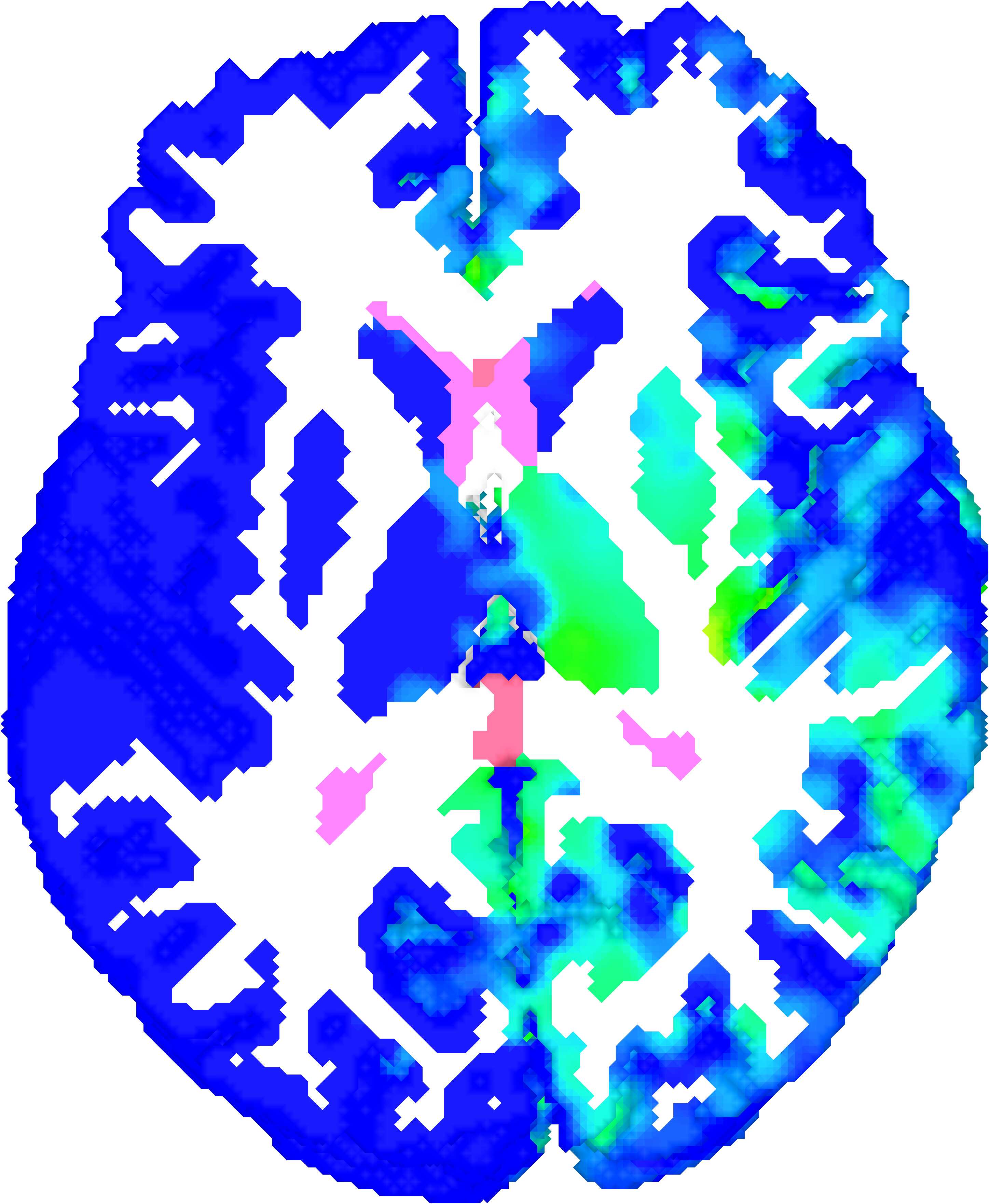}
            \caption{$\azimuth=\qty{0}{\degree}$, $\elevation=\qty{90}{\degree}$, $z=\qty{40}{\milli\meter}$.}
        \end{subfigure}
    \end{minipage}
    \hfill
    \begin{minipage}{0.09\linewidth}
    \includegraphics[width=\linewidth]{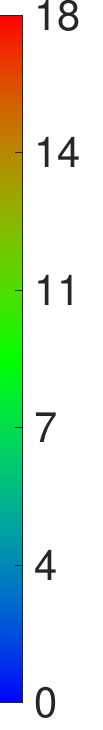}
    \end{minipage}
    \caption{Relative differences \eqref{eq:rel-diff} between the interference fields \eqref{eq:interference-field} of $\linJv$ and $\refJv$, smoothed via the function $\dBT\args{\someField/\max\someField,\qty{18}{\decibel}}$ of equation \eqref{eq:visual-log-transform} at different azimuths $\azimuth$ and elevations $\elevation$, with base stimulation frequency $\frequency=\qty{1000}{\hertz}$ and an updated $\contactResistance=\qty{5270}{\ohm}$. The current pattern was $\currentPattern=[-1.0,1.0,-1.0,1.0]\,\si{\milli\ampere}$. The maximal value of the absolute differences was roughly \qty{1.98e-6}{\milli\ampere\per\milli\meter\squared}.}
    \label{fig:linrefreldiff-Rc=5270-i=1.0-1.0}
\end{figure}

\begin{figure}[!h]
    \centering
    \begin{minipage}{0.85\linewidth}
        \def\subFigWidth{0.31\linewidth}
        \begin{subfigure}[b]{\subFigWidth}
            \includegraphics[width=\linewidth]{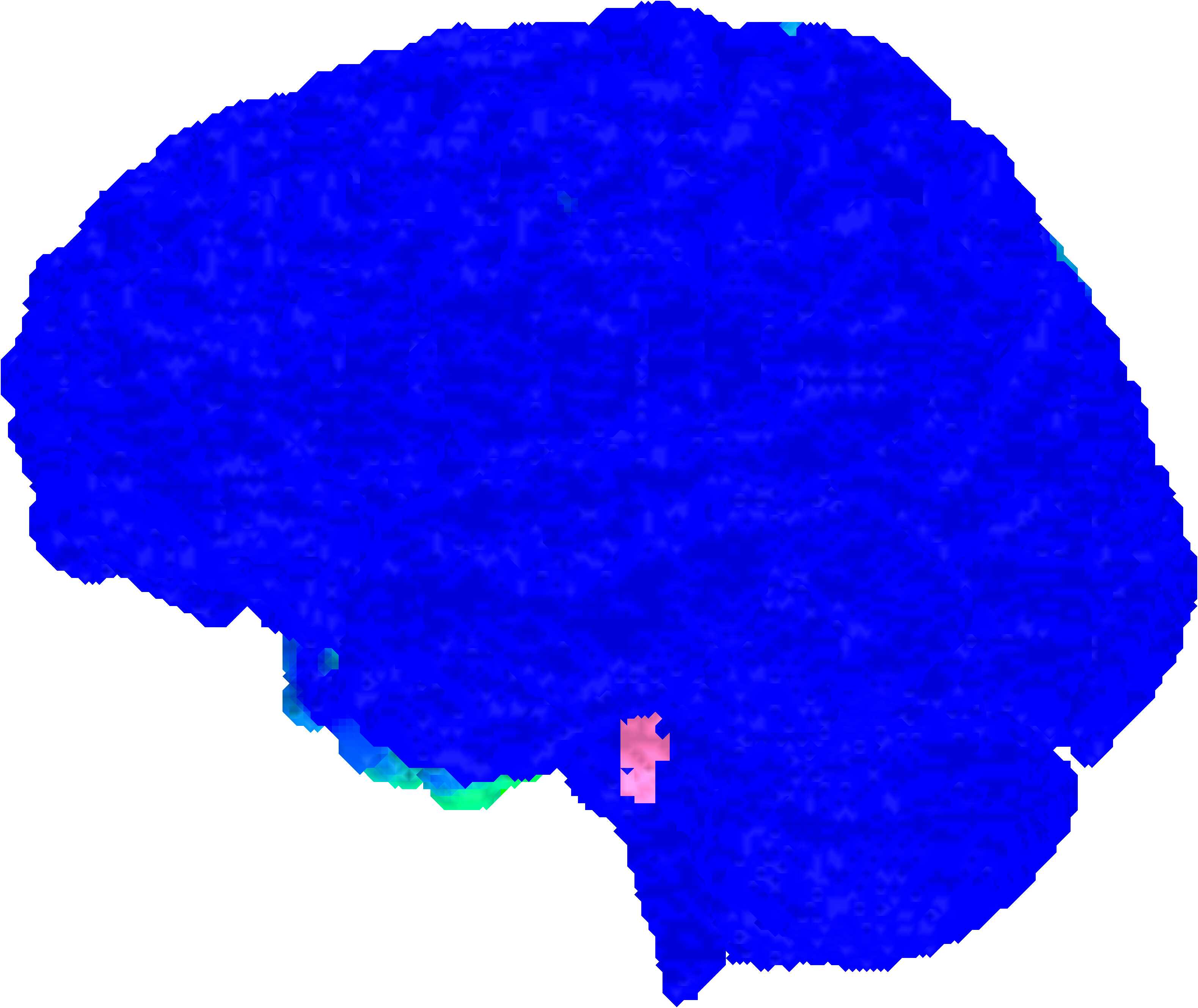}
            \caption{$\azimuth=\qty{-90}{\degree}$, $\elevation=\qty{0}{\degree}$.}
        \end{subfigure}
        \hfill
        \begin{subfigure}[b]{\subFigWidth}
        \includegraphics[width=\linewidth]{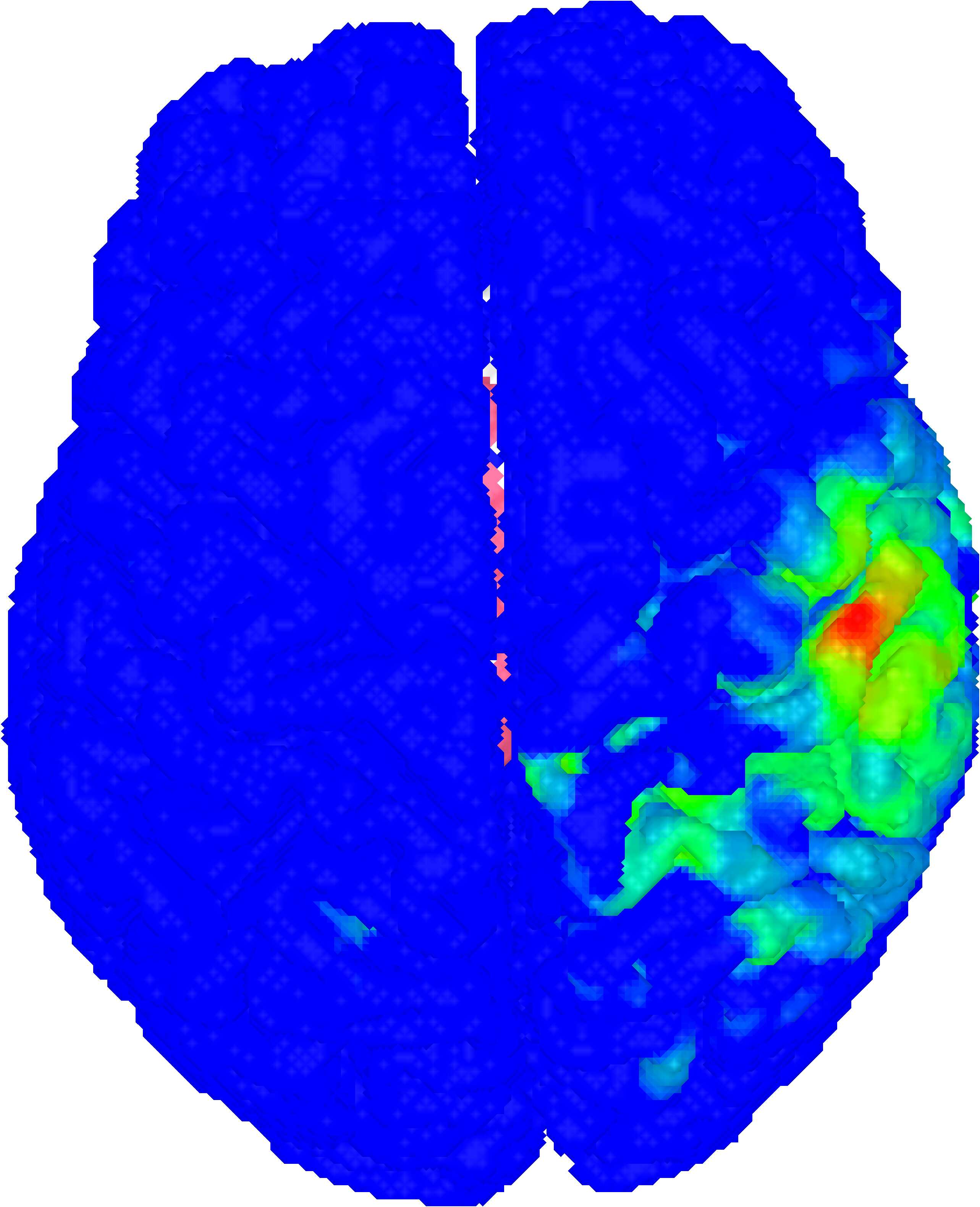}
        \caption{$\azimuth=\qty{0}{\degree}$, $\elevation=\qty{90}{\degree}$.}
        \end{subfigure}
        \hfill
        \begin{subfigure}[b]{\subFigWidth}
        \includegraphics[width=\linewidth]{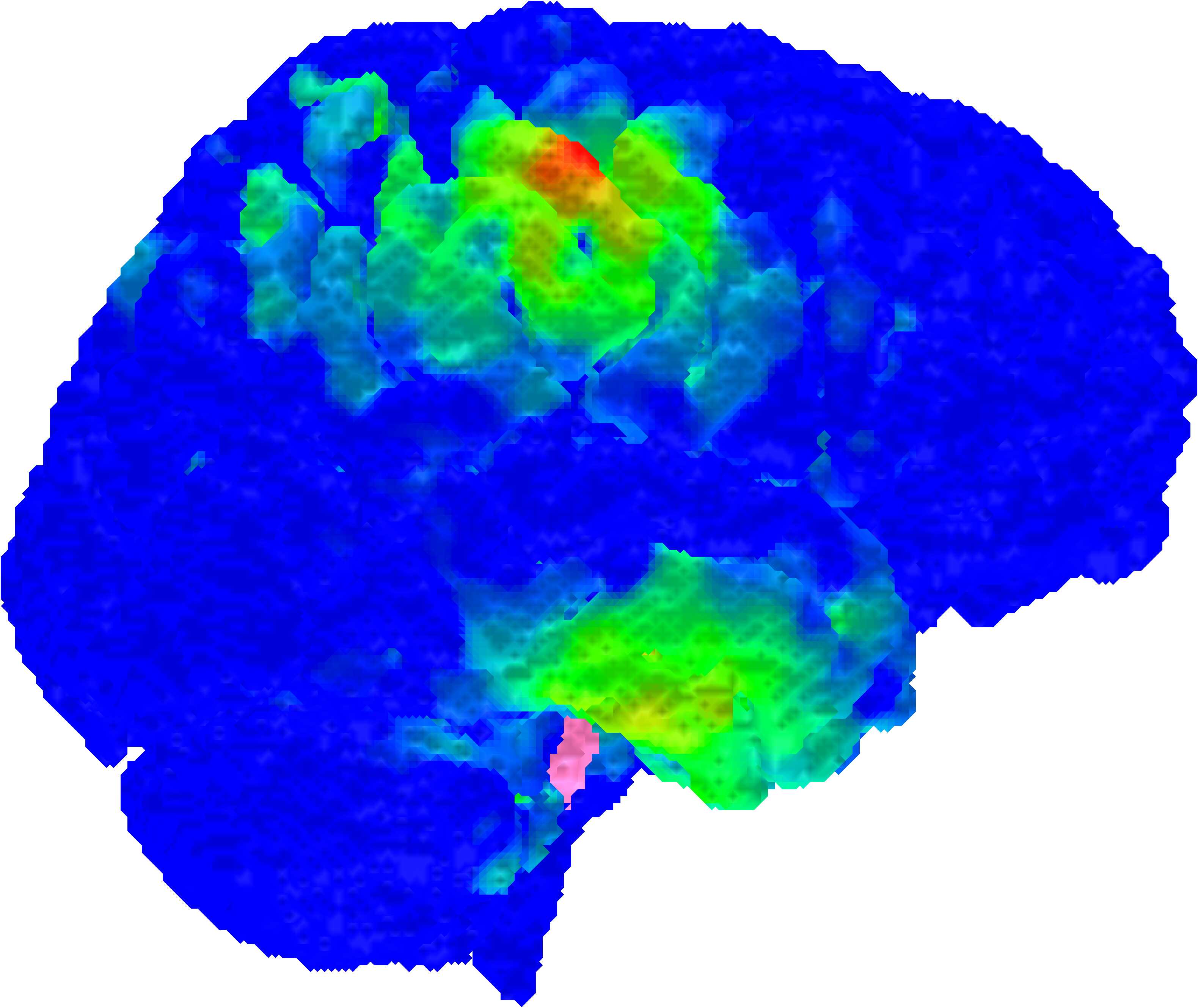}
            \caption{$\azimuth=\qty{90}{\degree}$, $\elevation=\qty{0}{\degree}$.}
        \end{subfigure}

        \def\subFigWidth{0.45\linewidth}

        \begin{subfigure}[b]{\subFigWidth}
            \includegraphics[width=\linewidth]{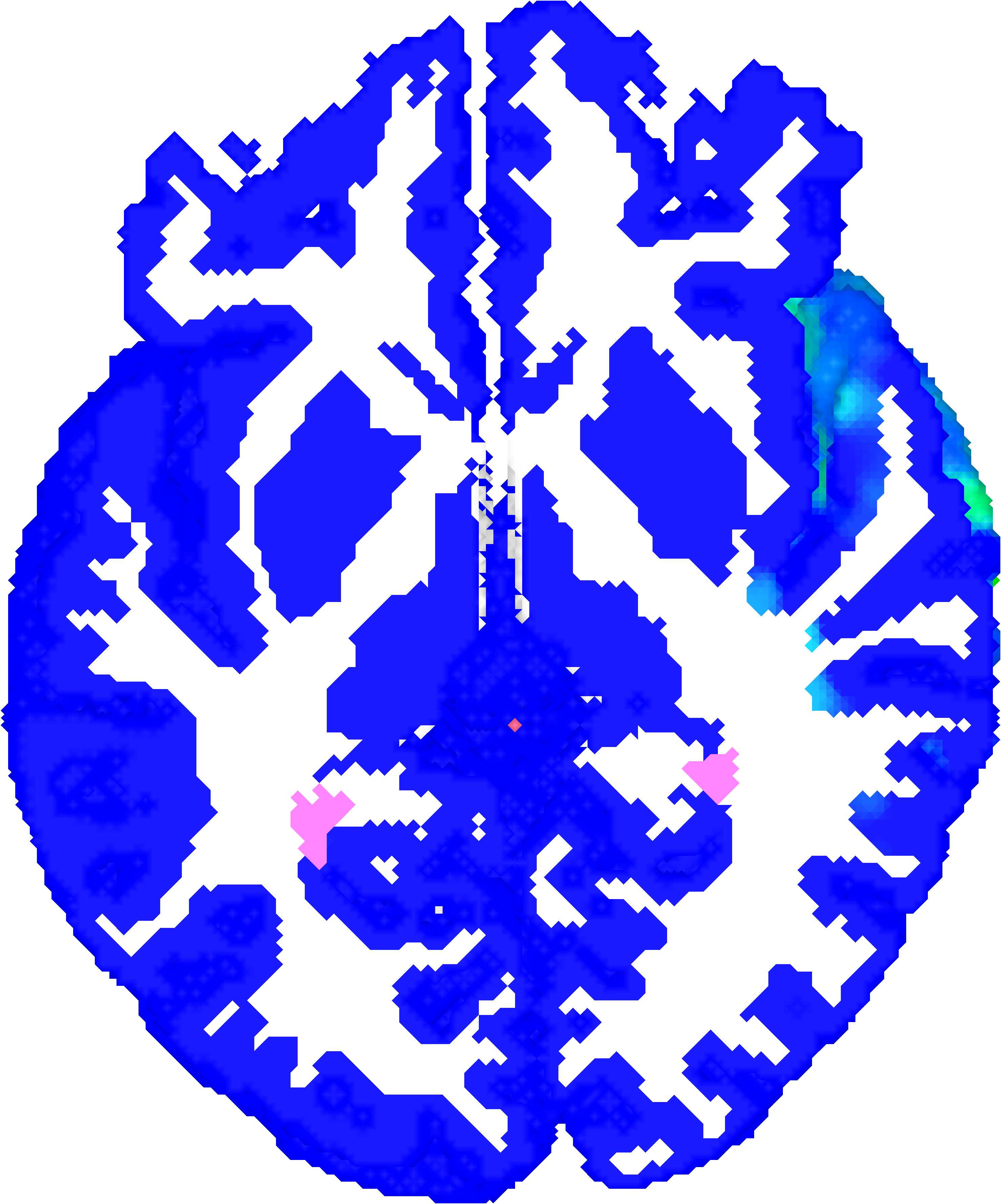}
            \caption{$\azimuth=\qty{0}{\degree}$, $\elevation=\qty{90}{\degree}$, $z=\qty{30}{\milli\meter}$.}
        \end{subfigure}
       \hfill
        \begin{subfigure}[b]{\subFigWidth}
            \includegraphics[width=\linewidth]{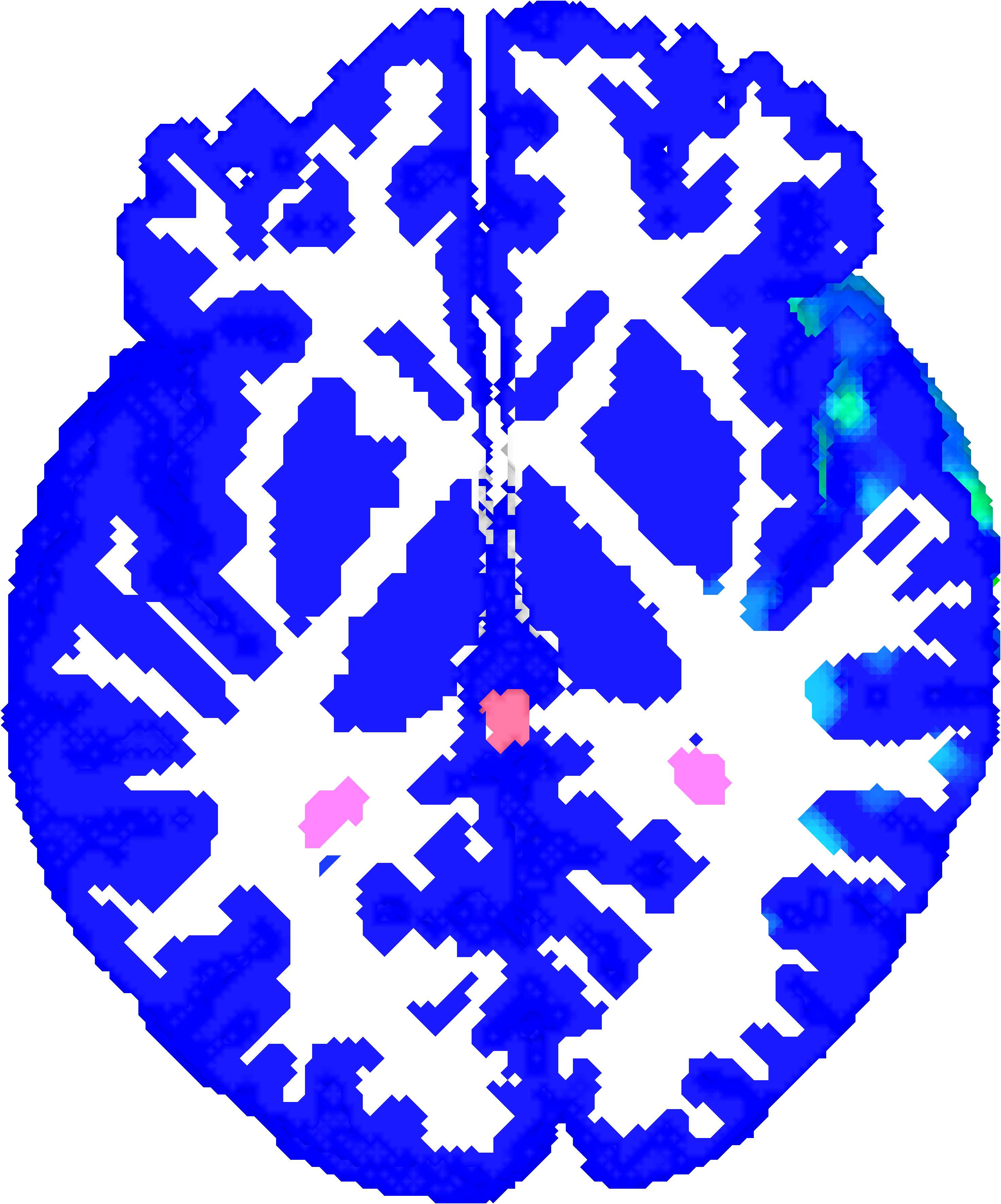}
            \caption{$\azimuth=\qty{0}{\degree}$, $\elevation=\qty{90}{\degree}$, $z=\qty{33}{\milli\meter}$.}
        \end{subfigure}

        \begin{subfigure}[b]{\subFigWidth}
            \includegraphics[width=\linewidth]{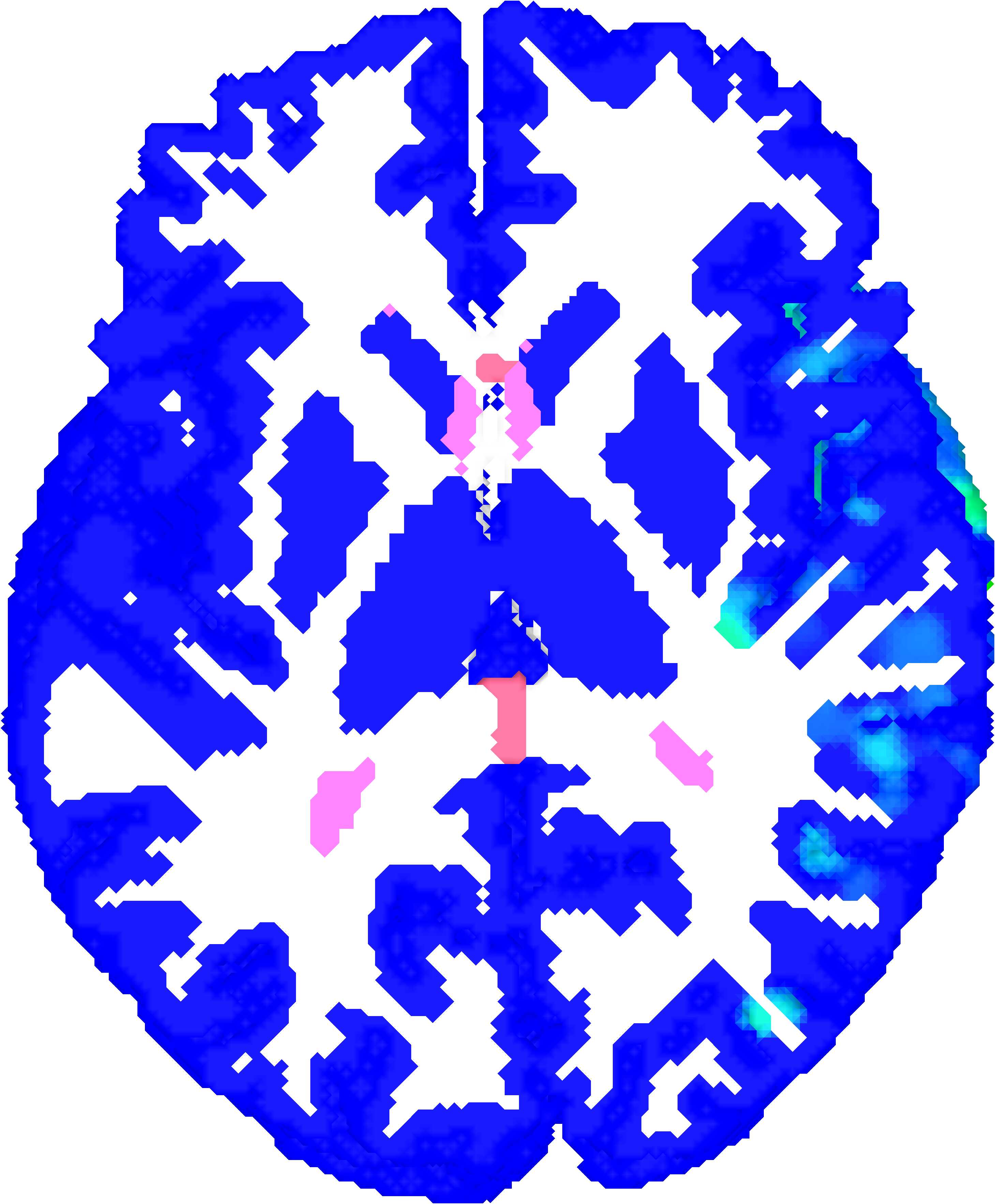}
            \caption{$\azimuth=\qty{0}{\degree}$, $\elevation=\qty{90}{\degree}$, $z=\qty{37}{\milli\meter}$.}
        \end{subfigure}
        \hfill
        \begin{subfigure}[b]{\subFigWidth}
            \includegraphics[width=\linewidth]{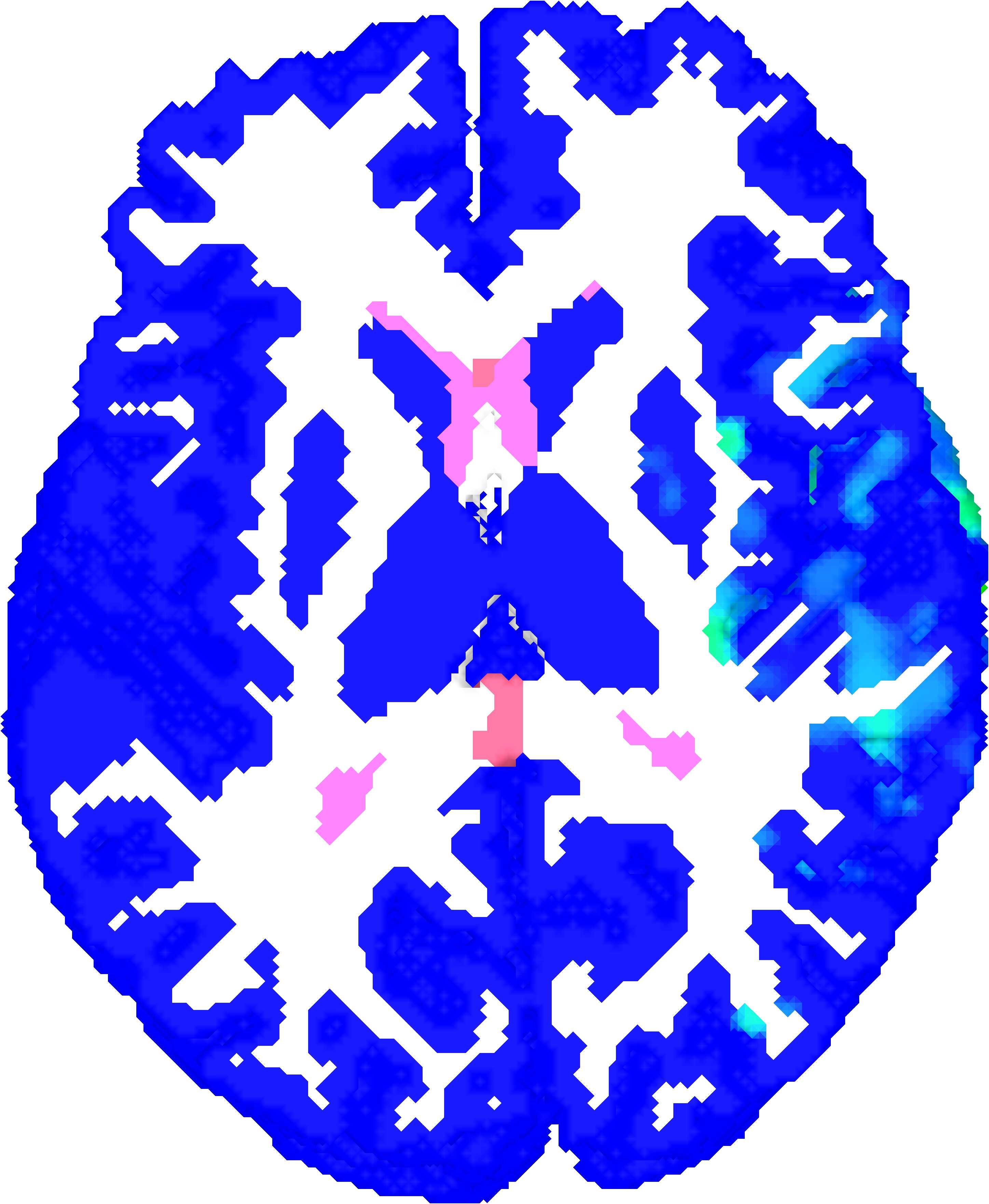}
            \caption{$\azimuth=\qty{0}{\degree}$, $\elevation=\qty{90}{\degree}$, $z=\qty{40}{\milli\meter}$.}
        \end{subfigure}
    \end{minipage}
    \hfill
    \begin{minipage}{0.09\linewidth}
    \includegraphics[width=\linewidth]{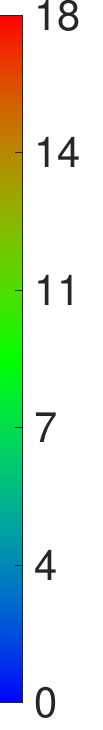}
    \end{minipage}
    \caption{Relative differences \eqref{eq:rel-diff} between the interference fields \eqref{eq:interference-field} of $\linJv$ and $\refJv$, smoothed via the function $\dBT\args{\someField/\max\someField,\qty{18}{\decibel}}$ of equation \eqref{eq:visual-log-transform} at different azimuths $\azimuth$ and elevations $\elevation$, with base stimulation frequency $\frequency=\qty{1000}{\hertz}$ and an updated $\contactResistance=\qty{1270}{\ohm}$. The current pattern was $\currentPattern=[-1.5,1.5,-0.5,0.5]\,\si{\milli\ampere}$. The maximal value of the absolute differences was roughly \qty{3.12e-6}{\milli\ampere\per\milli\meter\squared}.}
    \label{fig:linrefreldiff-Rc=1270-i=1.5-0.5}
\end{figure}

\begin{figure}[!h]
    \centering
    \begin{minipage}{0.85\linewidth}
        \def\subFigWidth{0.31\linewidth}
        \begin{subfigure}[b]{\subFigWidth}
            \includegraphics[width=\linewidth]{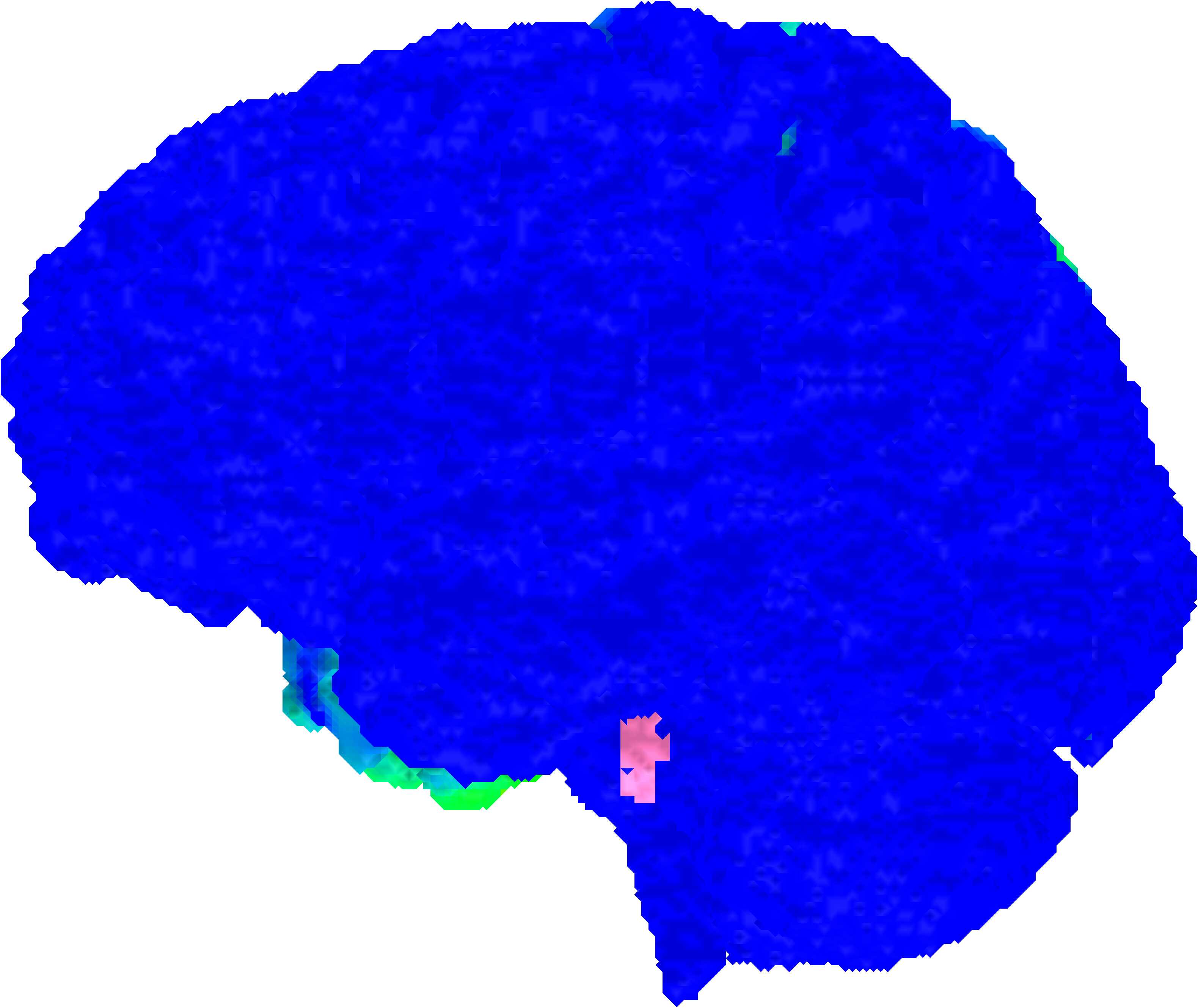}
            \caption{$\azimuth=\qty{-90}{\degree}$, $\elevation=\qty{0}{\degree}$.}
        \end{subfigure}
        \hfill
        \begin{subfigure}[b]{\subFigWidth}
        \includegraphics[width=\linewidth]{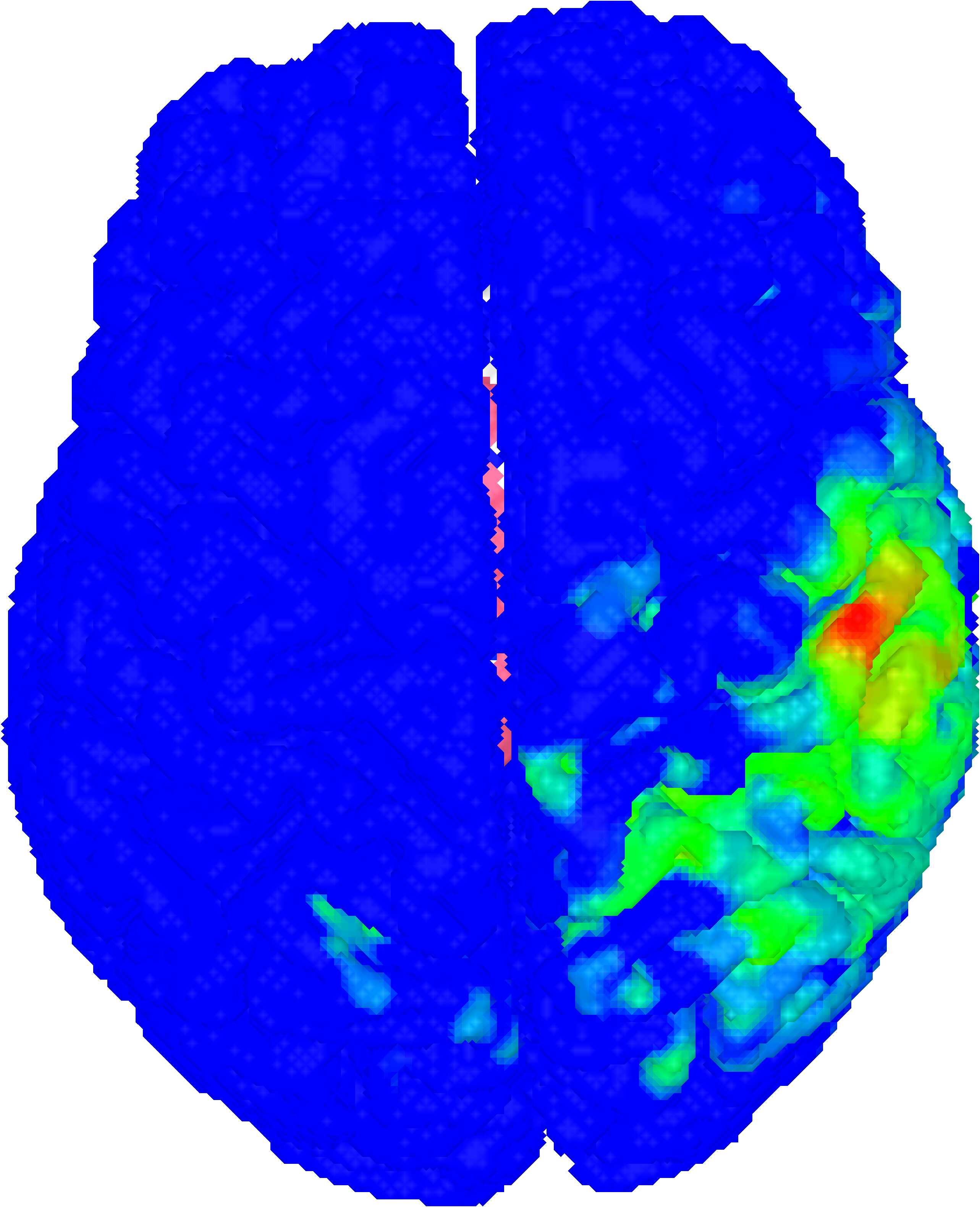}
        \caption{$\azimuth=\qty{0}{\degree}$, $\elevation=\qty{90}{\degree}$.}
        \end{subfigure}
        \hfill
        \begin{subfigure}[b]{\subFigWidth}
        \includegraphics[width=\linewidth]{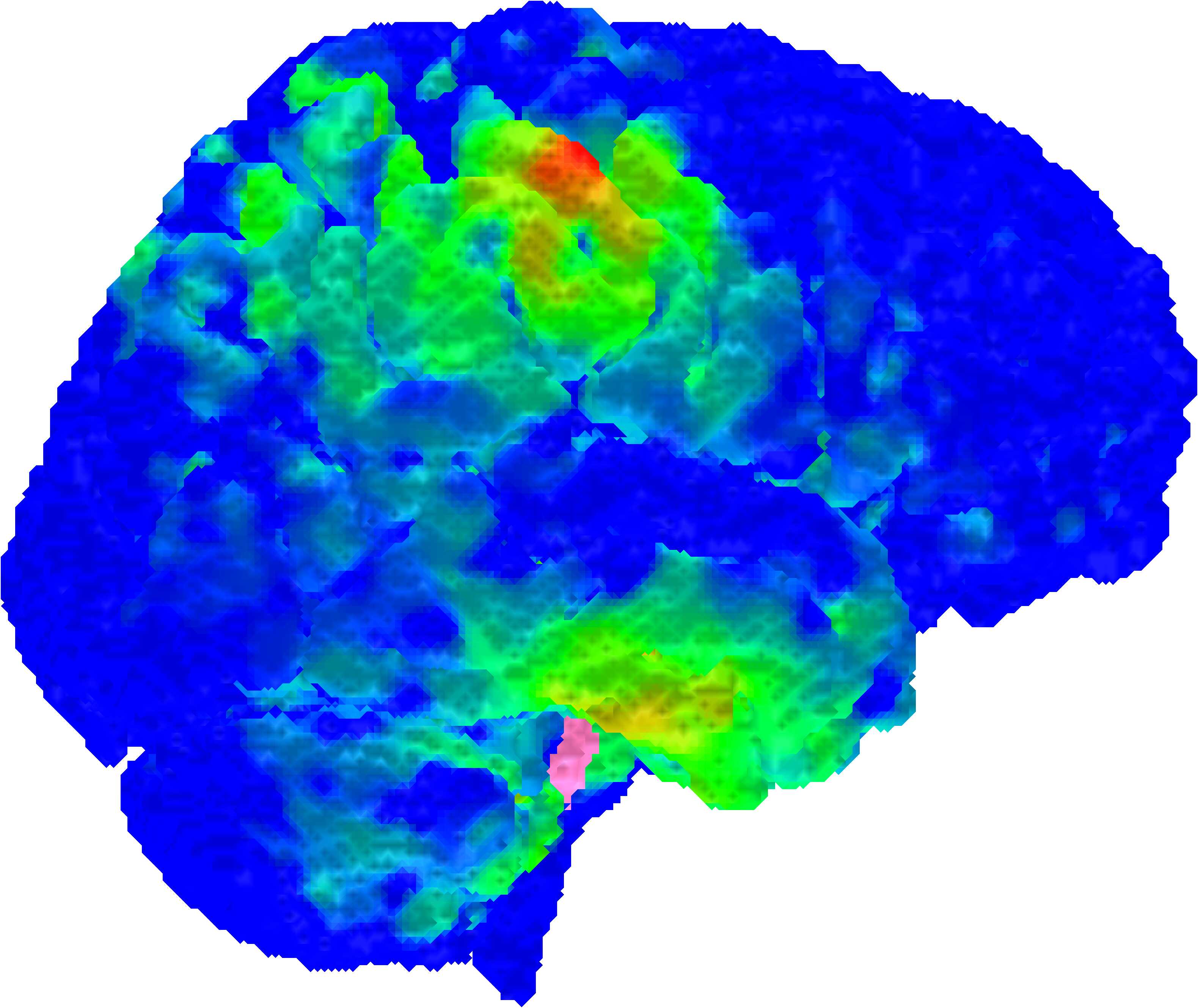}
            \caption{$\azimuth=\qty{90}{\degree}$, $\elevation=\qty{0}{\degree}$.}
        \end{subfigure}

        \def\subFigWidth{0.45\linewidth}

        \begin{subfigure}[b]{\subFigWidth}
            \includegraphics[width=\linewidth]{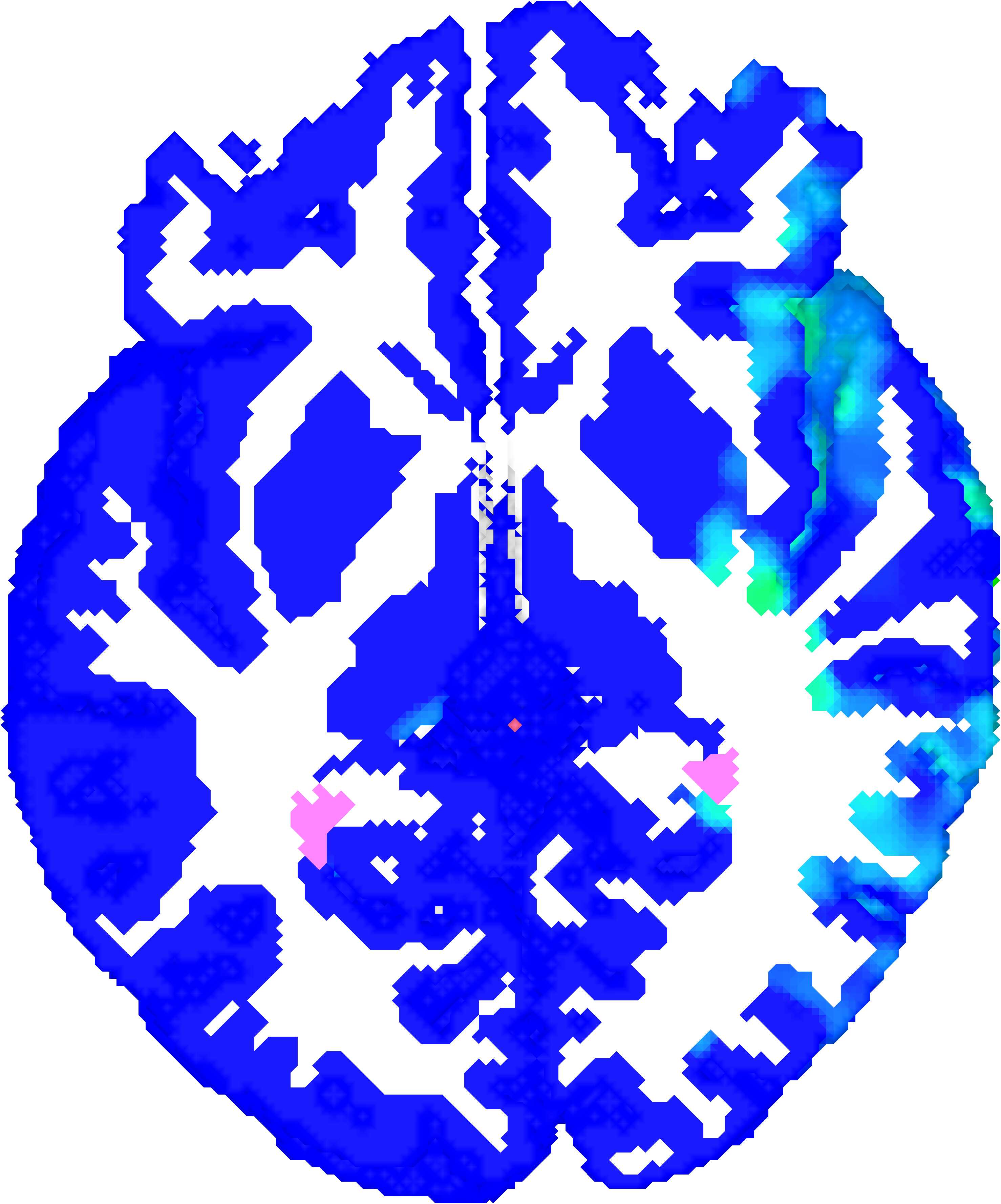}
            \caption{$\azimuth=\qty{0}{\degree}$, $\elevation=\qty{90}{\degree}$, $z=\qty{30}{\milli\meter}$.}
        \end{subfigure}
       \hfill
        \begin{subfigure}[b]{\subFigWidth}
            \includegraphics[width=\linewidth]{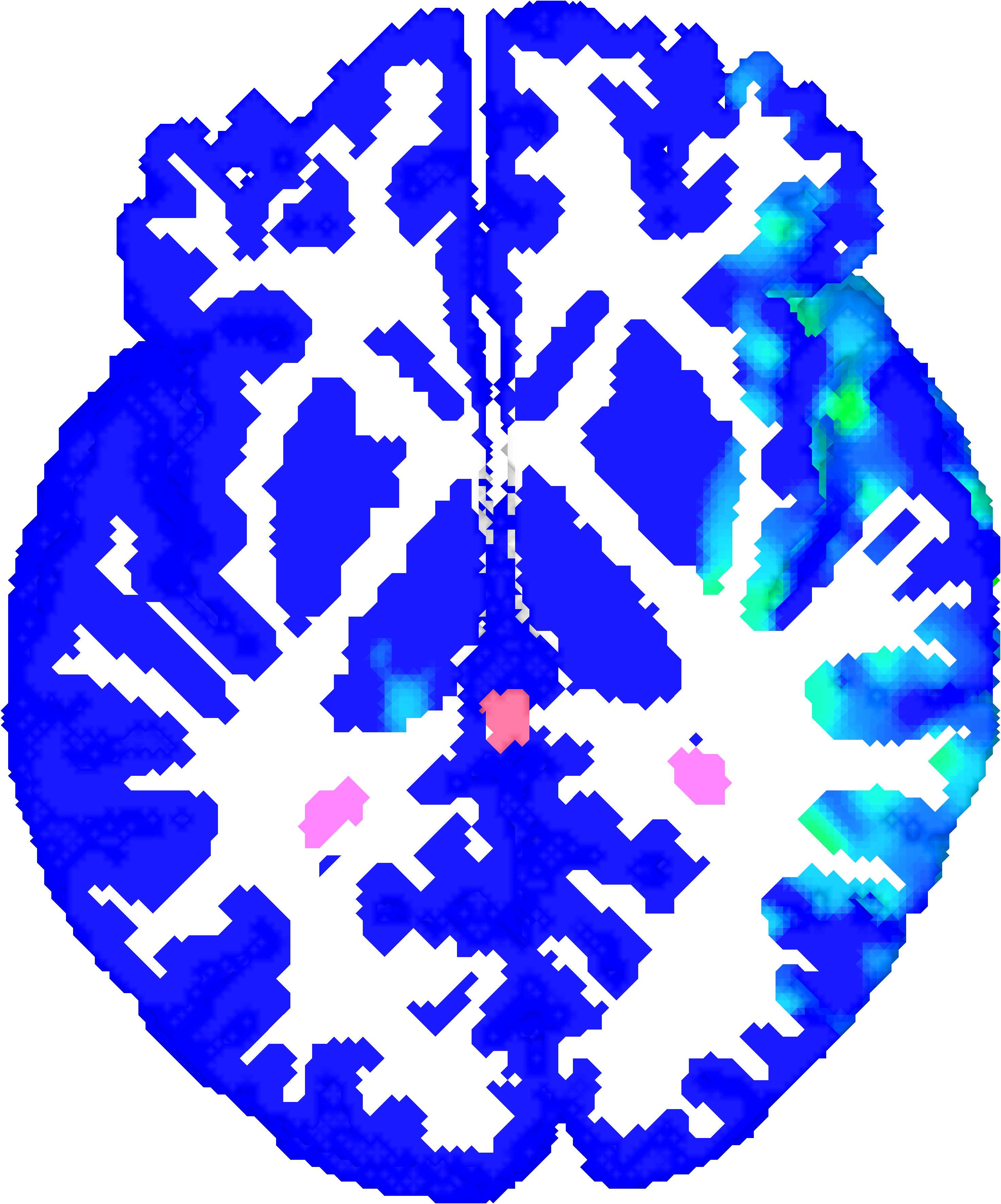}
            \caption{$\azimuth=\qty{0}{\degree}$, $\elevation=\qty{90}{\degree}$, $z=\qty{33}{\milli\meter}$.}
        \end{subfigure}

        \begin{subfigure}[b]{\subFigWidth}
            \includegraphics[width=\linewidth]{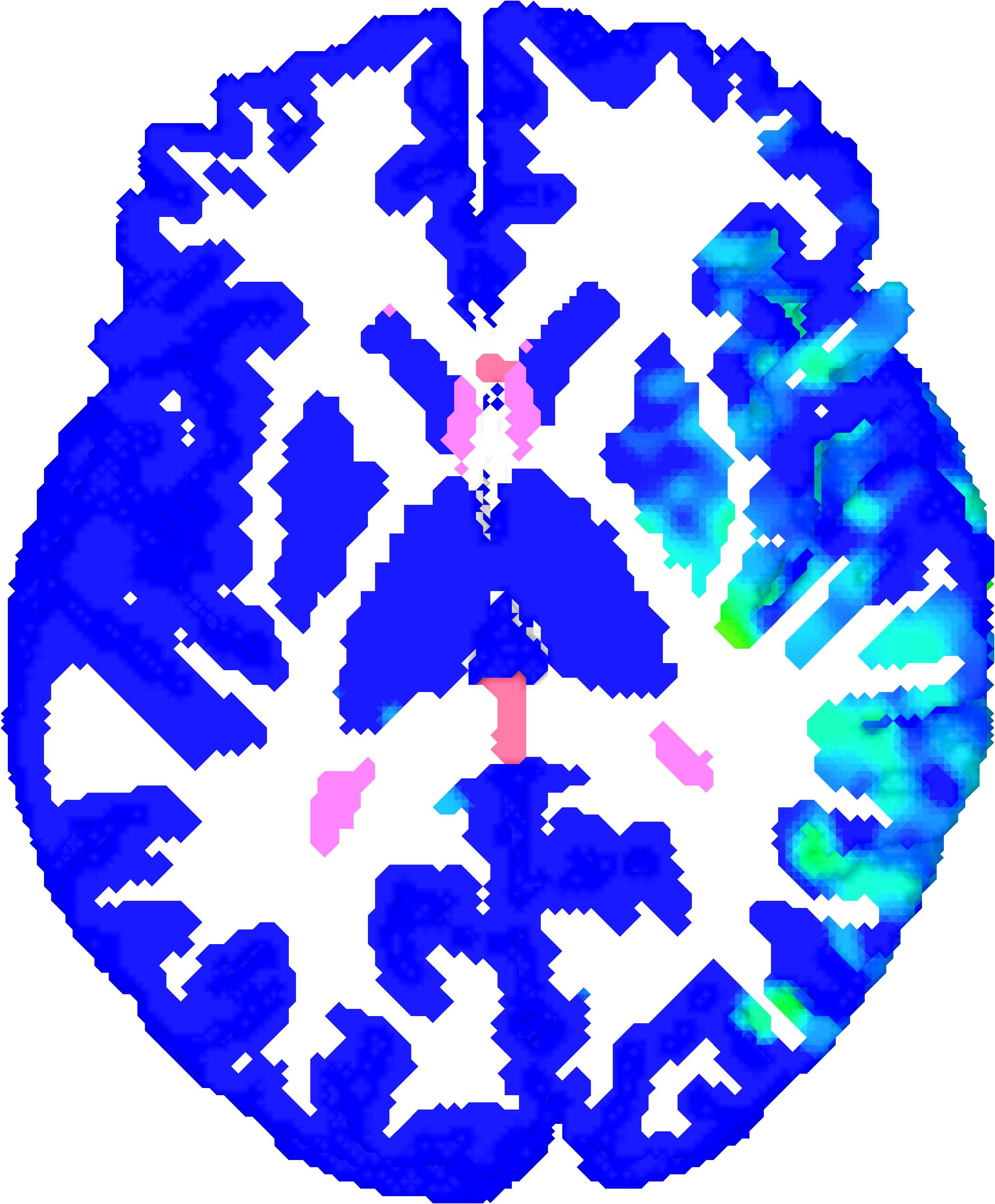}
            \caption{$\azimuth=\qty{0}{\degree}$, $\elevation=\qty{90}{\degree}$, $z=\qty{37}{\milli\meter}$.}
        \end{subfigure}
        \hfill
        \begin{subfigure}[b]{\subFigWidth}
            \includegraphics[width=\linewidth]{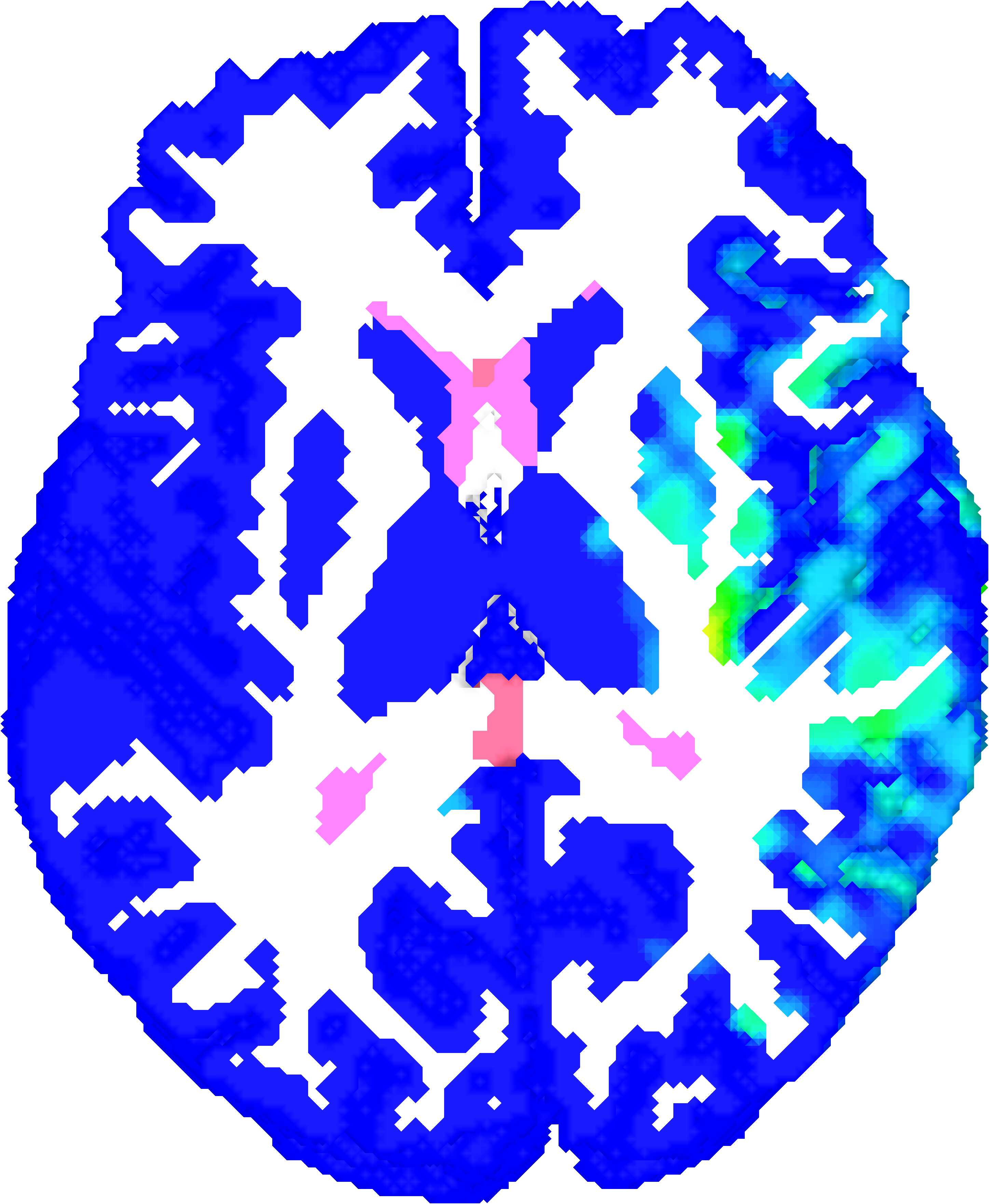}
            \caption{$\azimuth=\qty{0}{\degree}$, $\elevation=\qty{90}{\degree}$, $z=\qty{40}{\milli\meter}$.}
        \end{subfigure}
    \end{minipage}
    \hfill
    \begin{minipage}{0.09\linewidth}
    \includegraphics[width=\linewidth]{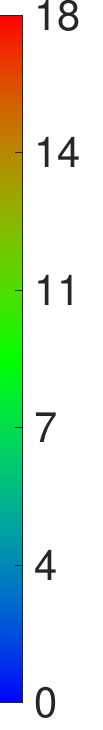}
    \end{minipage}
    \caption{Relative differences \eqref{eq:rel-diff} between the interference fields \eqref{eq:interference-field} of $\linJv$ and $\refJv$, smoothed via the function $\dBT\args{\someField/\max\someField,\qty{18}{\decibel}}$ of equation \eqref{eq:visual-log-transform} at different azimuths $\azimuth$ and elevations $\elevation$, with base stimulation frequency $\frequency=\qty{1000}{\hertz}$ and an updated $\contactResistance=\qty{5270}{\ohm}$. The current pattern was $\currentPattern=[-1.5,1.5,-0.5,0.5]\,\si{\milli\ampere}$. The maximal value of the absolute differences was roughly \qty{2.95e-6}{\milli\ampere\per\milli\meter\squared}.}
    \label{fig:linrefreldiff-Rc=5270-i=1.5-0.5}
\end{figure}

\begin{figure}[!h]
    \centering
    \begin{minipage}{0.85\linewidth}
        \def\subFigWidth{0.31\linewidth}
        \begin{subfigure}[b]{\subFigWidth}
            \includegraphics[width=\linewidth]{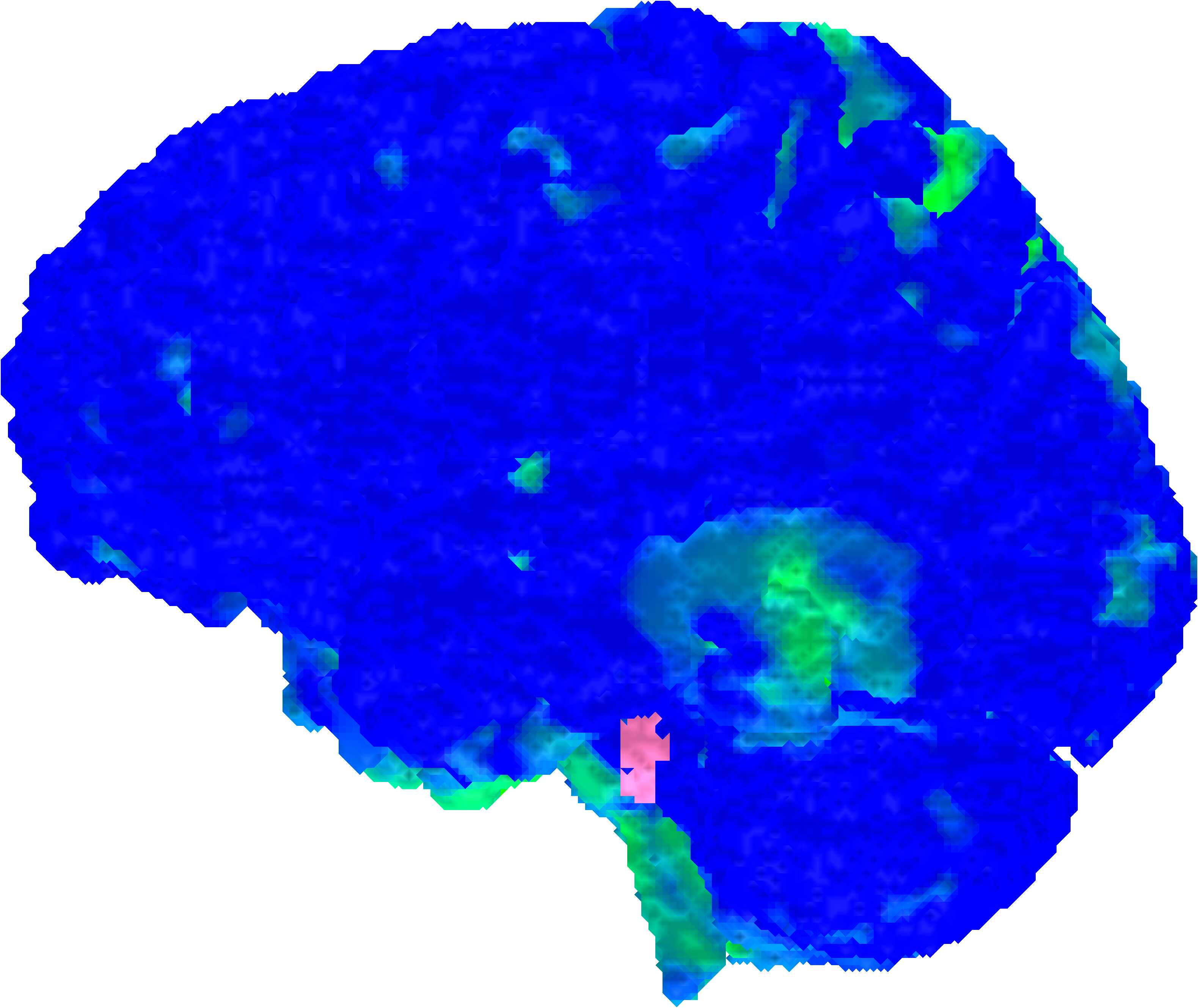}
            \caption{$\azimuth=\qty{-90}{\degree}$, $\elevation=\qty{0}{\degree}$.}
        \end{subfigure}
        \hfill
        \begin{subfigure}[b]{\subFigWidth}
        \includegraphics[width=\linewidth]{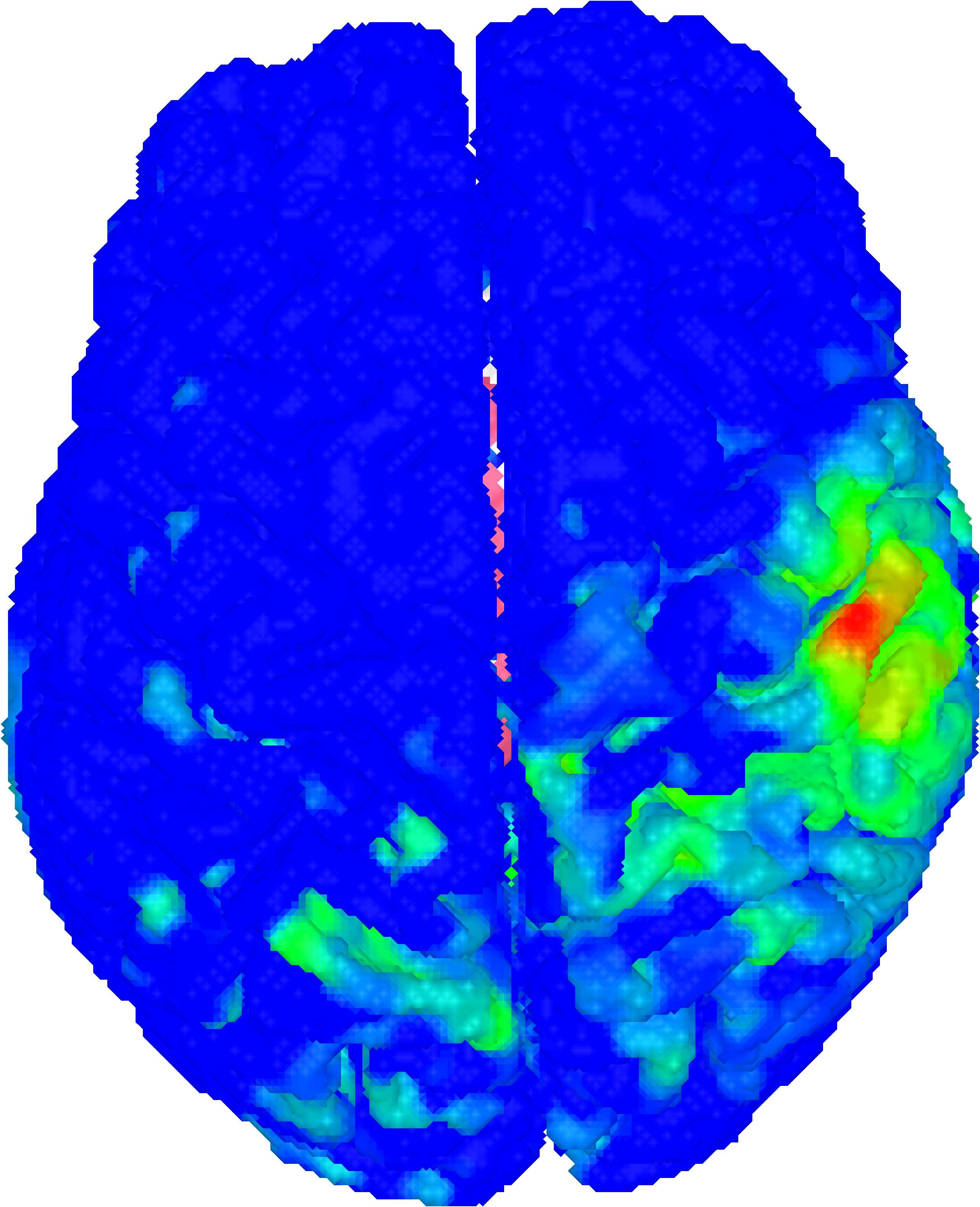}
        \caption{$\azimuth=\qty{0}{\degree}$, $\elevation=\qty{90}{\degree}$.}
        \end{subfigure}
        \hfill
        \begin{subfigure}[b]{\subFigWidth}
        \includegraphics[width=\linewidth]{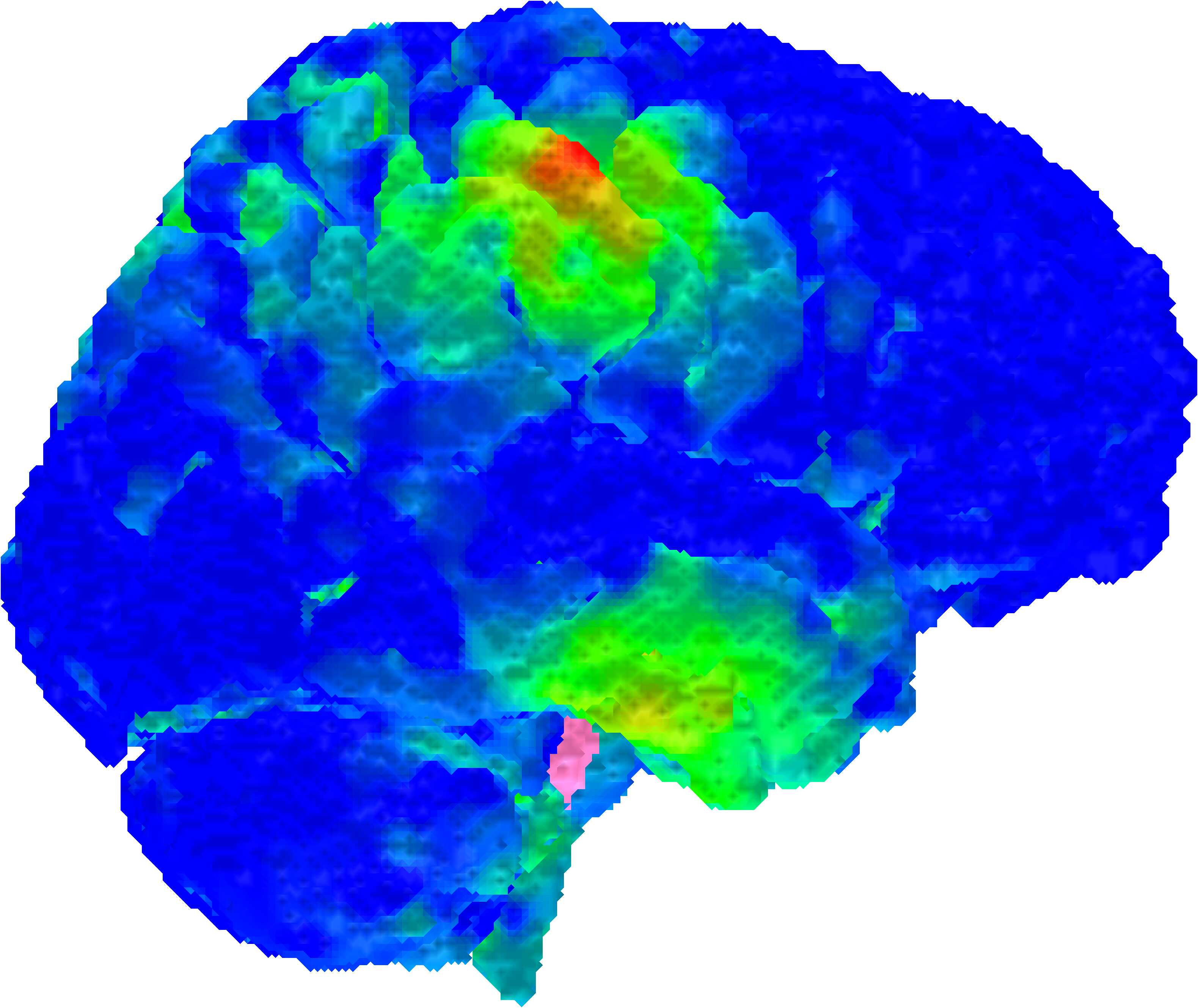}
            \caption{$\azimuth=\qty{90}{\degree}$, $\elevation=\qty{0}{\degree}$.}
        \end{subfigure}

        \def\subFigWidth{0.45\linewidth}

        \begin{subfigure}[b]{\subFigWidth}
            \includegraphics[width=\linewidth]{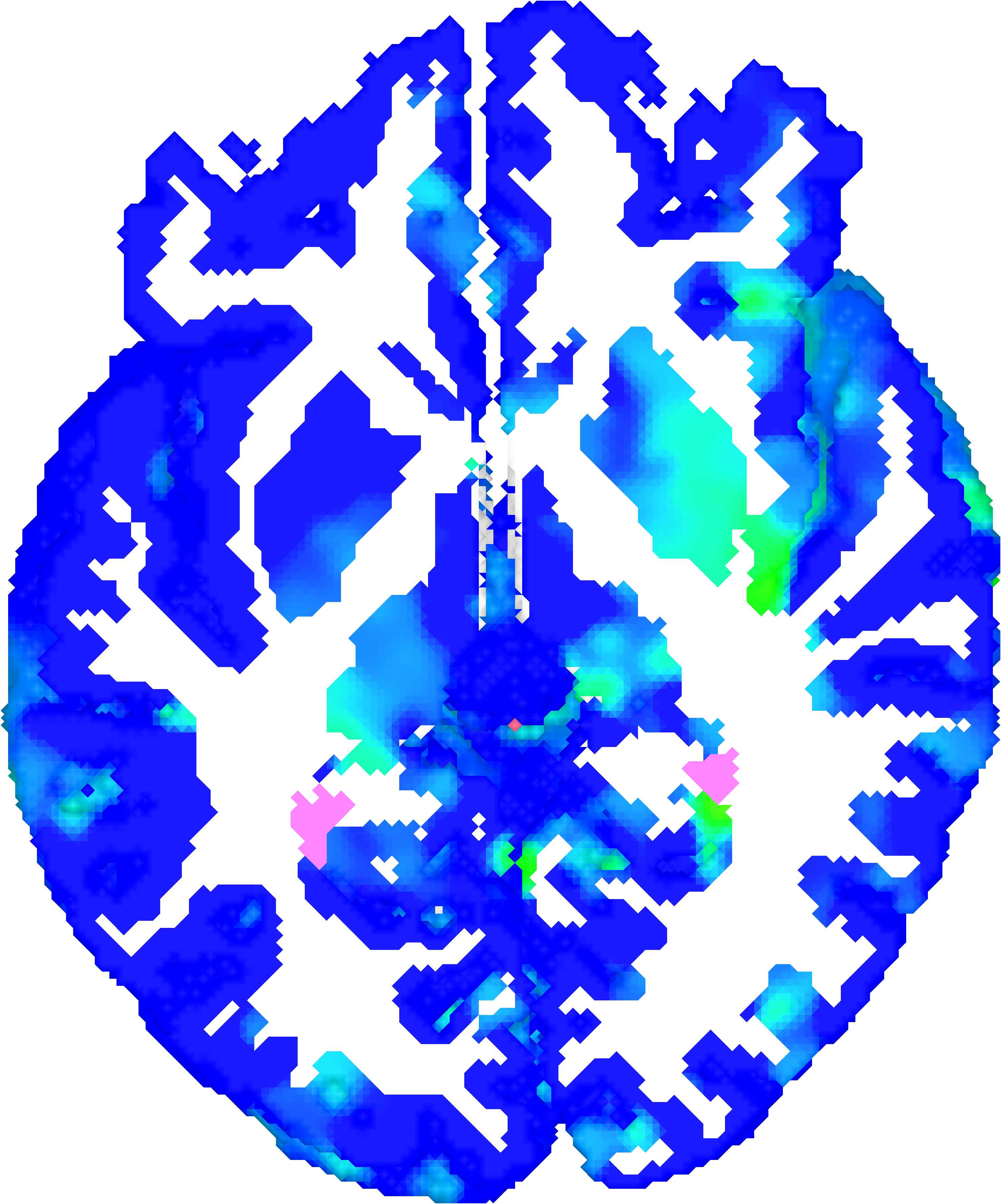}
            \caption{$\azimuth=\qty{0}{\degree}$, $\elevation=\qty{90}{\degree}$, $z=\qty{30}{\milli\meter}$.}
        \end{subfigure}
       \hfill
        \begin{subfigure}[b]{\subFigWidth}
            \includegraphics[width=\linewidth]{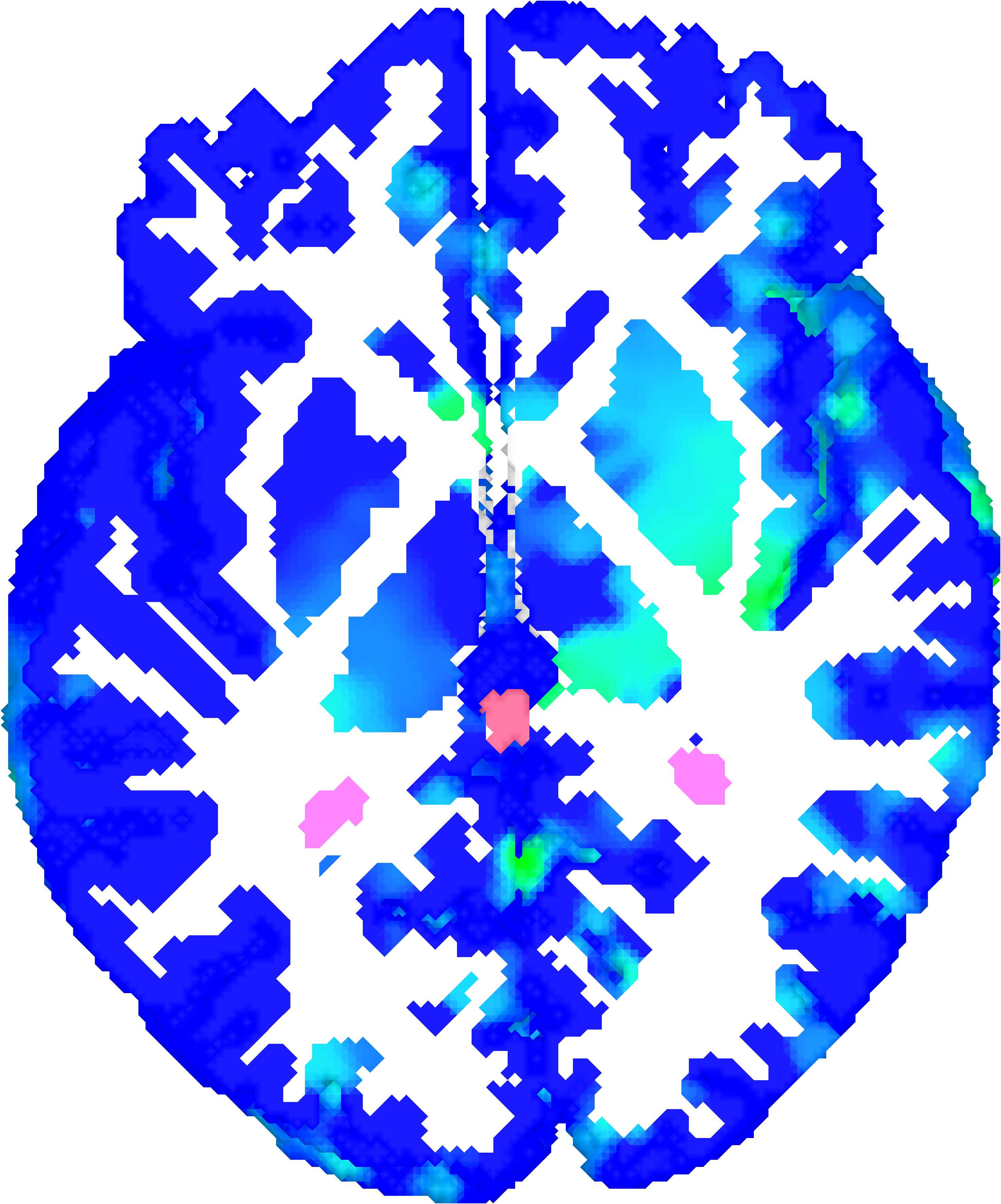}
            \caption{$\azimuth=\qty{0}{\degree}$, $\elevation=\qty{90}{\degree}$, $z=\qty{33}{\milli\meter}$.}
        \end{subfigure}

        \begin{subfigure}[b]{\subFigWidth}
            \includegraphics[width=\linewidth]{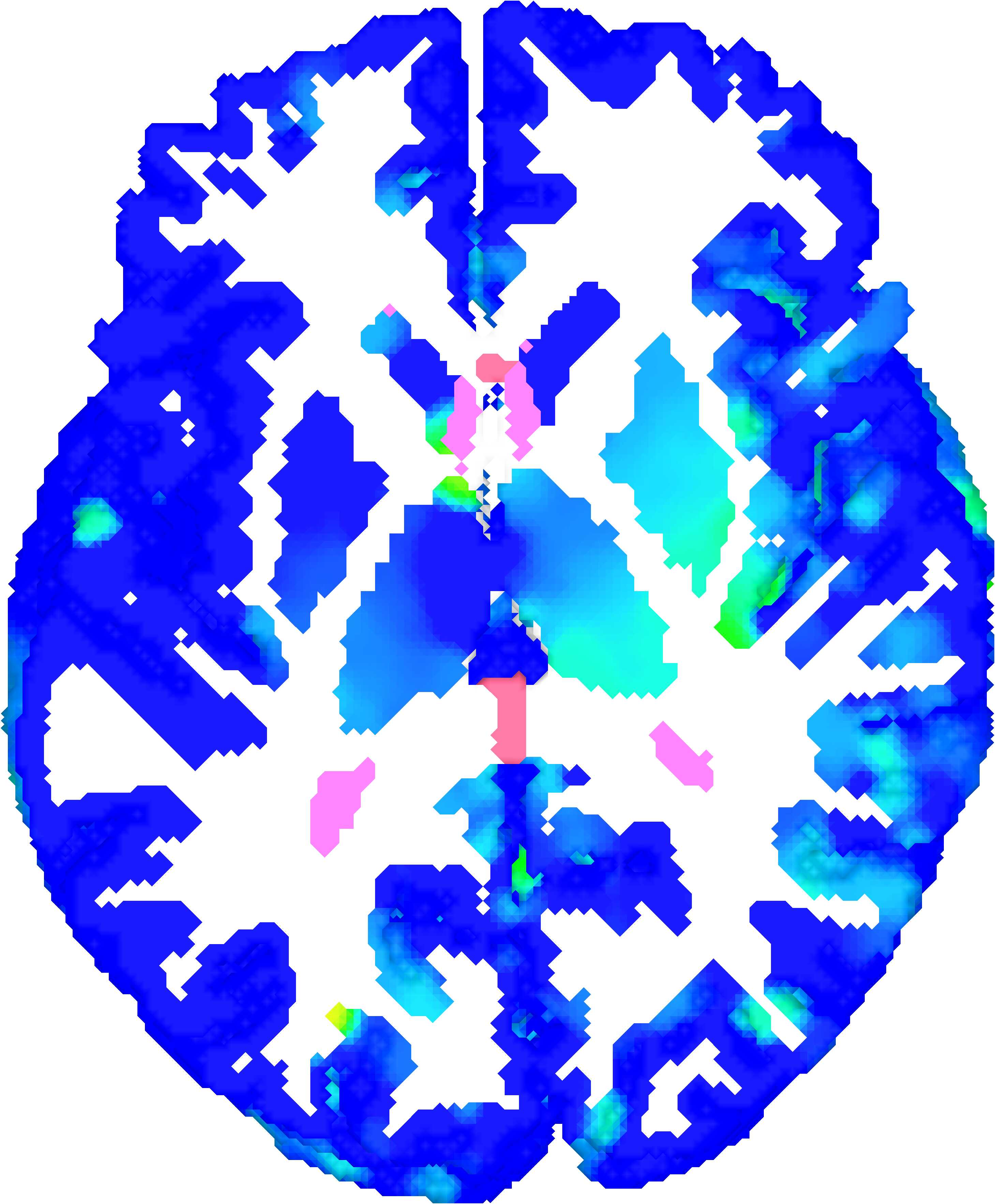}
            \caption{$\azimuth=\qty{0}{\degree}$, $\elevation=\qty{90}{\degree}$, $z=\qty{37}{\milli\meter}$.}
        \end{subfigure}
        \hfill
        \begin{subfigure}[b]{\subFigWidth}
            \includegraphics[width=\linewidth]{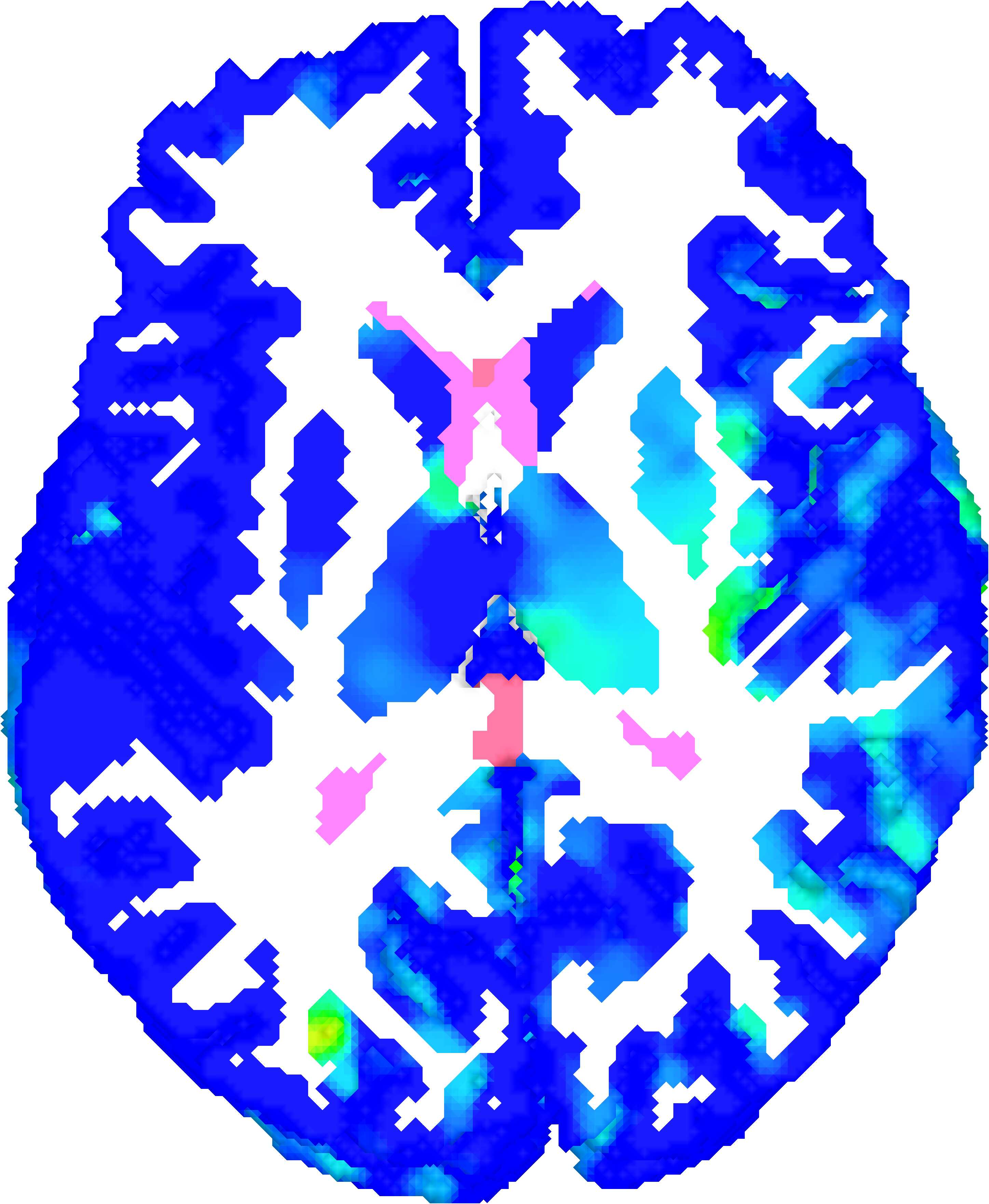}
            \caption{$\azimuth=\qty{0}{\degree}$, $\elevation=\qty{90}{\degree}$, $z=\qty{40}{\milli\meter}$.}
        \end{subfigure}
    \end{minipage}
    \hfill
    \begin{minipage}{0.09\linewidth}
    \includegraphics[width=\linewidth]{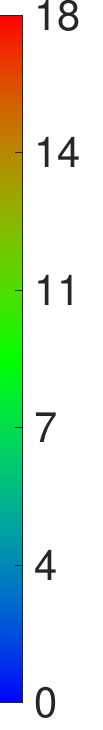}
    \end{minipage}
    \caption{Relative differences \eqref{eq:rel-diff} between the interference fields \eqref{eq:interference-field} of $\linJv$ and $\refJv$, smoothed via the function $\dBT\args{\someField/\max\someField,\qty{18}{\decibel}}$ of equation \eqref{eq:visual-log-transform} at different azimuths $\azimuth$ and elevations $\elevation$, with base stimulation frequency $\frequency=\qty{1000}{\hertz}$ and an updated $\contactResistance=\qty{1270}{\ohm}$. The current pattern was $\currentPattern=[-0.5,0.5,-1.5,1.5]\,\si{\milli\ampere}$. The maximal value of the absolute differences was roughly \qty{1.04e-6}{\milli\ampere\per\milli\meter\squared}.}
    \label{fig:linrefreldiff-Rc=1270-i=0.5-1.5}
\end{figure}

\begin{figure}[!h]
    \centering
    \begin{minipage}{0.85\linewidth}
        \def\subFigWidth{0.31\linewidth}
        \begin{subfigure}[b]{\subFigWidth}
            \includegraphics[width=\linewidth]{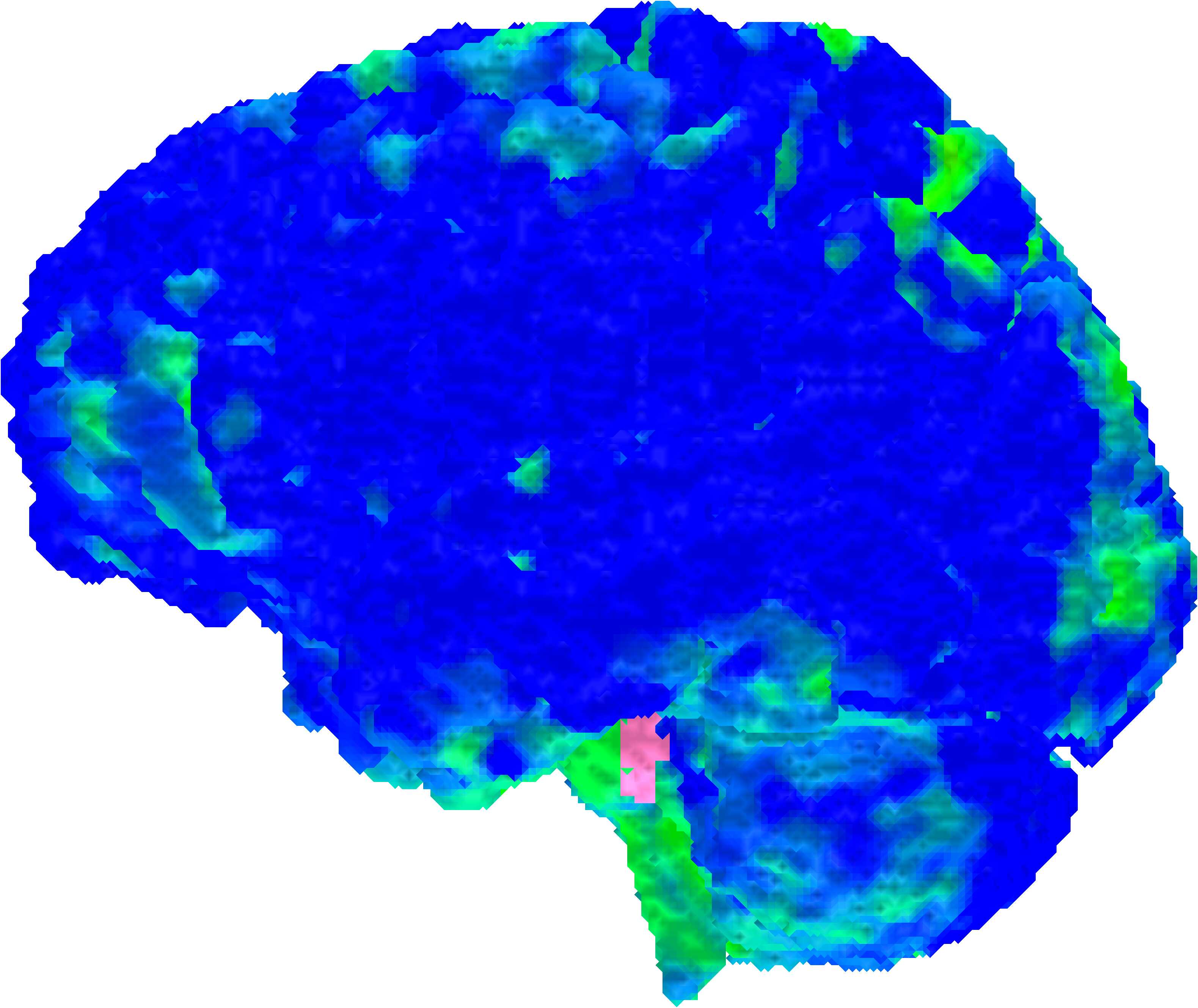}
            \caption{$\azimuth=\qty{-90}{\degree}$, $\elevation=\qty{0}{\degree}$.}
        \end{subfigure}
        \hfill
        \begin{subfigure}[b]{\subFigWidth}
        \includegraphics[width=\linewidth]{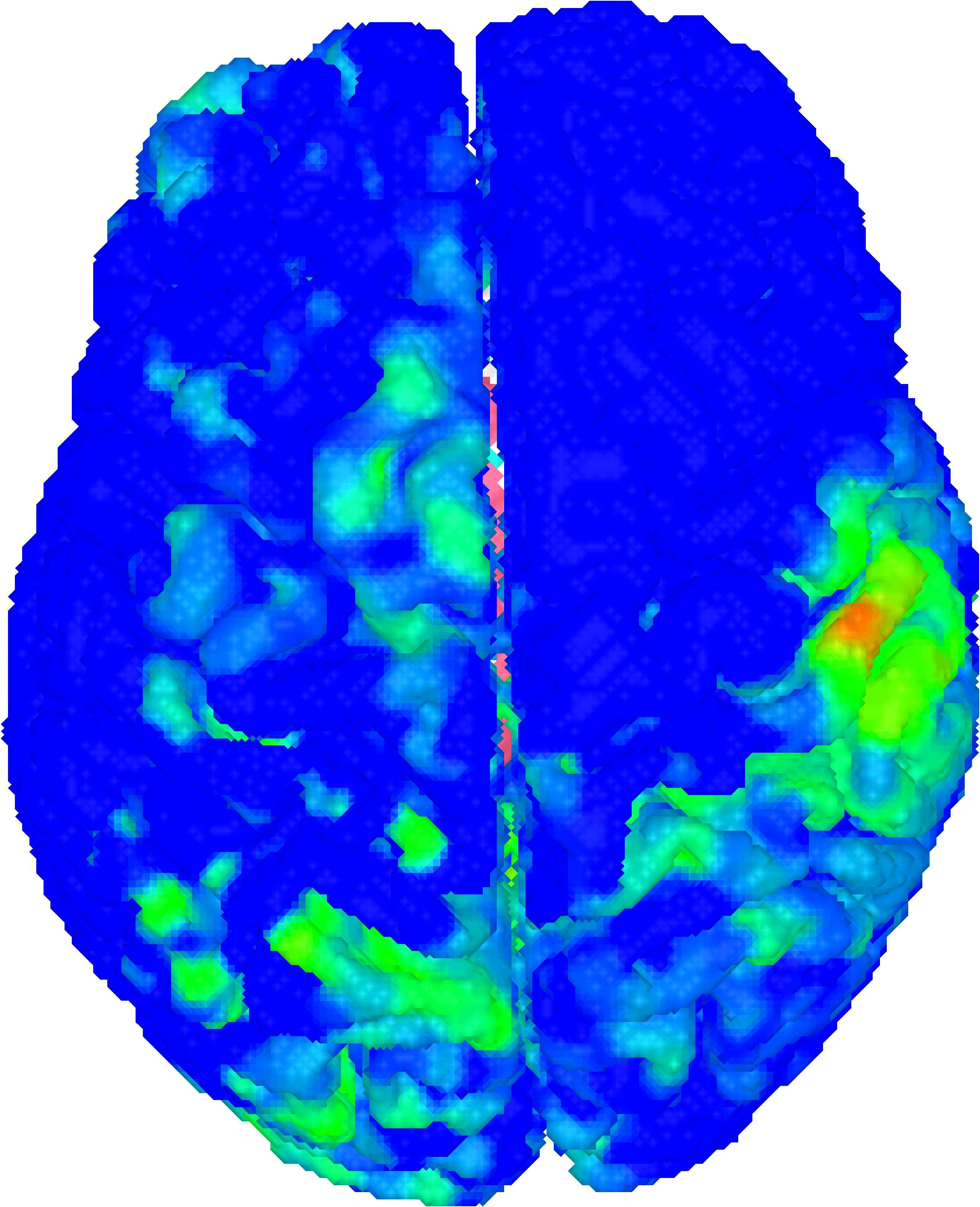}
        \caption{$\azimuth=\qty{0}{\degree}$, $\elevation=\qty{90}{\degree}$.}
        \end{subfigure}
        \hfill
        \begin{subfigure}[b]{\subFigWidth}
        \includegraphics[width=\linewidth]{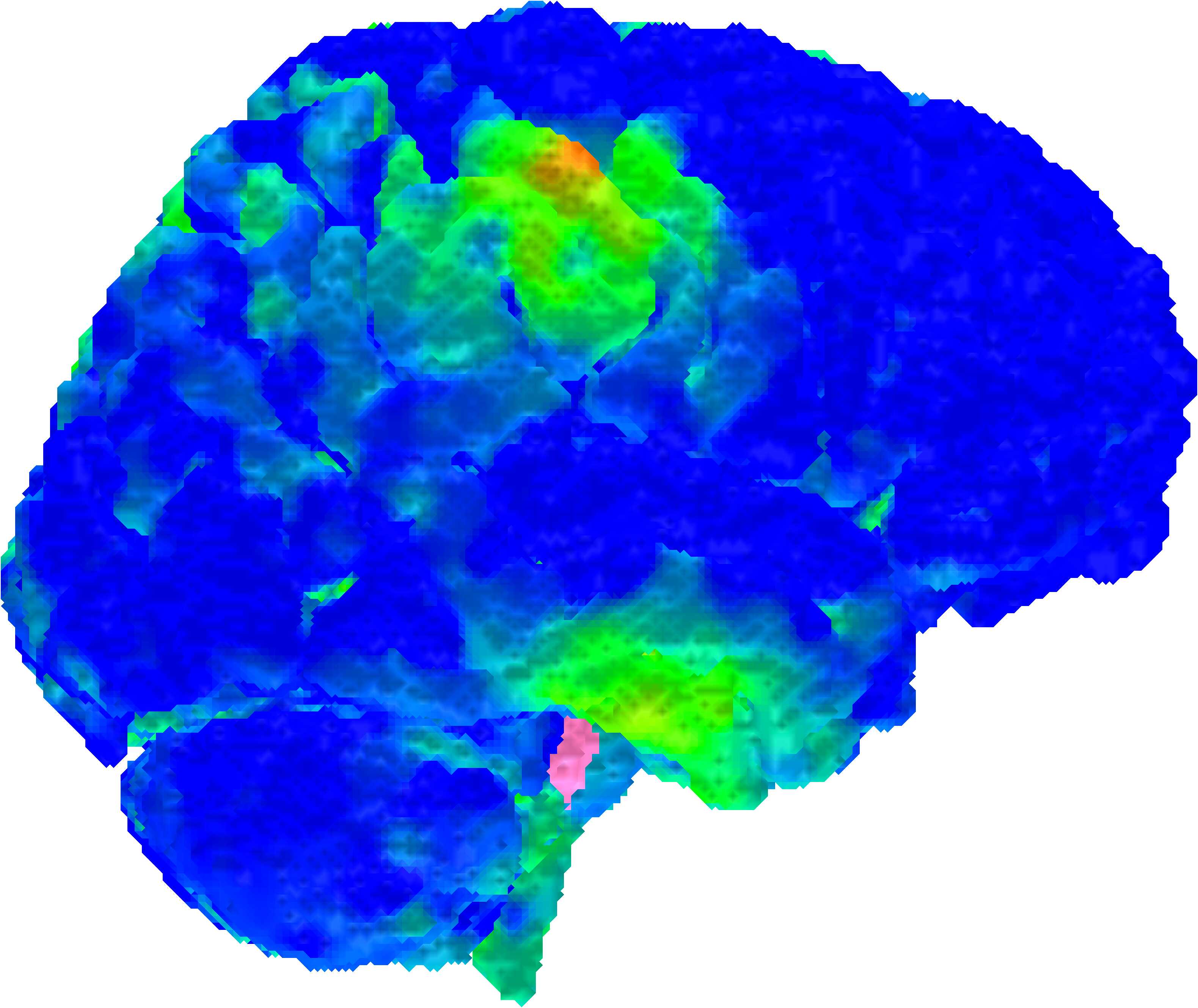}
            \caption{$\azimuth=\qty{90}{\degree}$, $\elevation=\qty{0}{\degree}$.}
        \end{subfigure}

        \def\subFigWidth{0.45\linewidth}

        \begin{subfigure}[b]{\subFigWidth}
            \includegraphics[width=\linewidth]{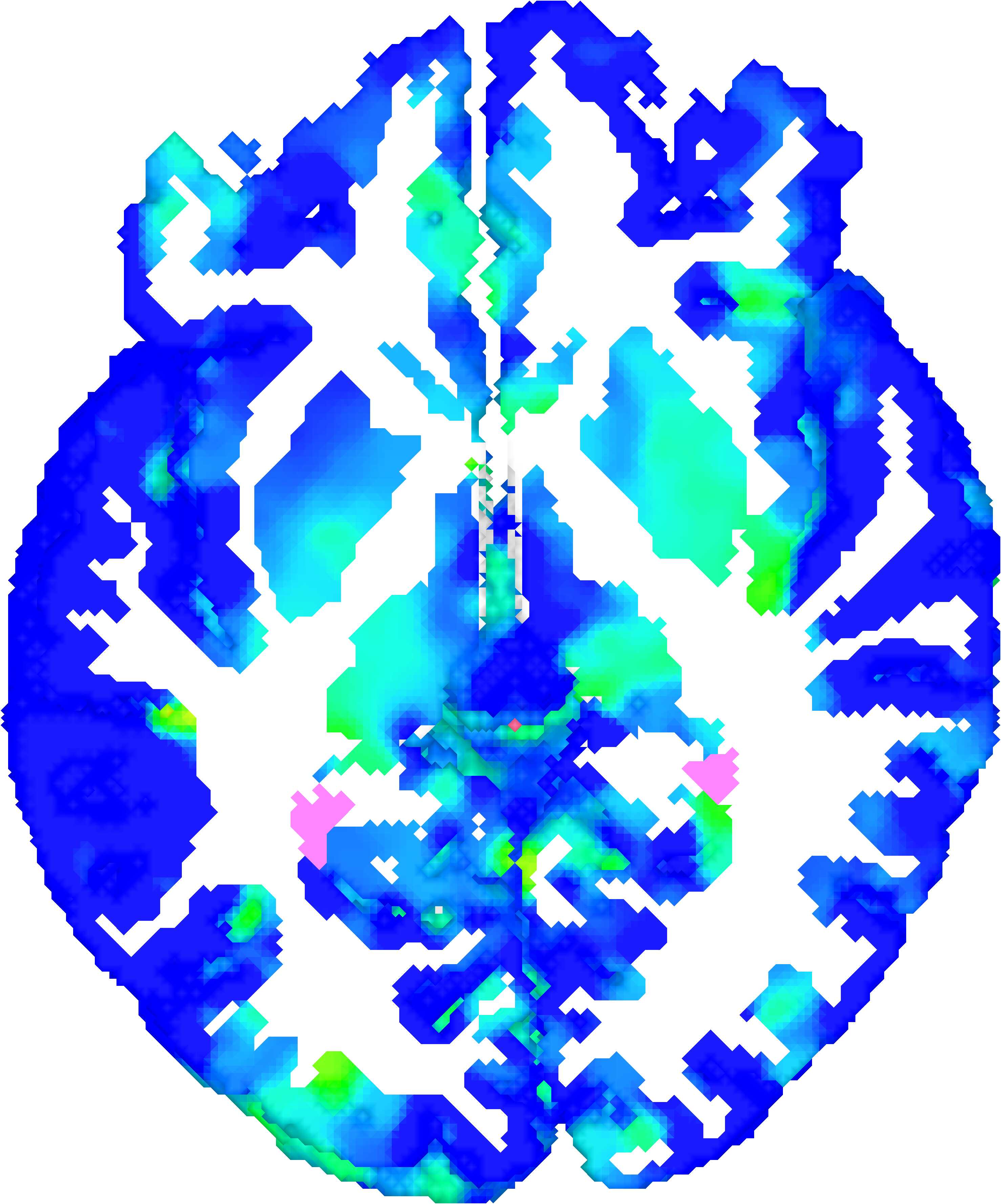}
            \caption{$\azimuth=\qty{0}{\degree}$, $\elevation=\qty{90}{\degree}$, $z=\qty{30}{\milli\meter}$.}
        \end{subfigure}
       \hfill
        \begin{subfigure}[b]{\subFigWidth}
            \includegraphics[width=\linewidth]{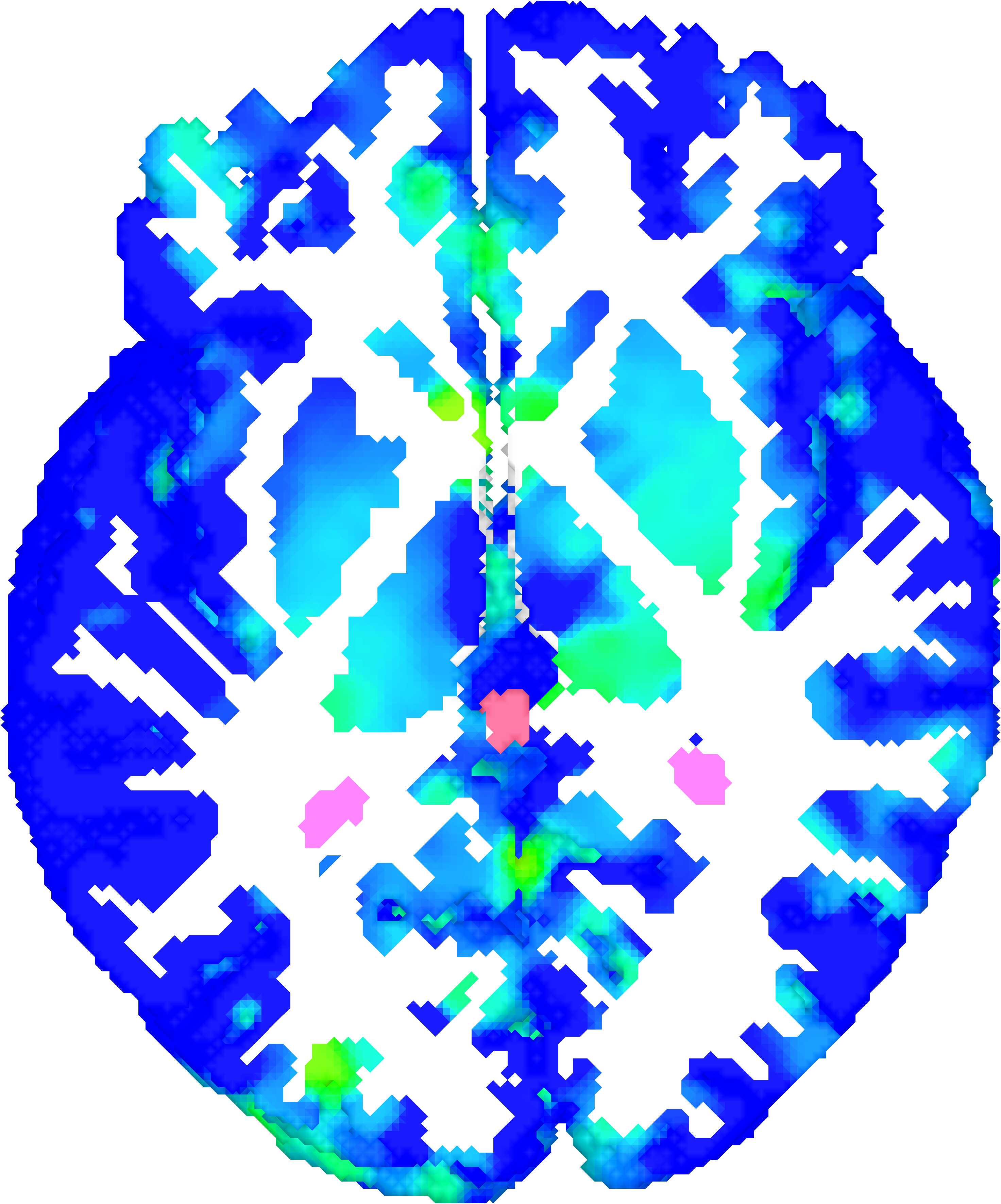}
            \caption{$\azimuth=\qty{0}{\degree}$, $\elevation=\qty{90}{\degree}$, $z=\qty{33}{\milli\meter}$.}
        \end{subfigure}

        \begin{subfigure}[b]{\subFigWidth}
            \includegraphics[width=\linewidth]{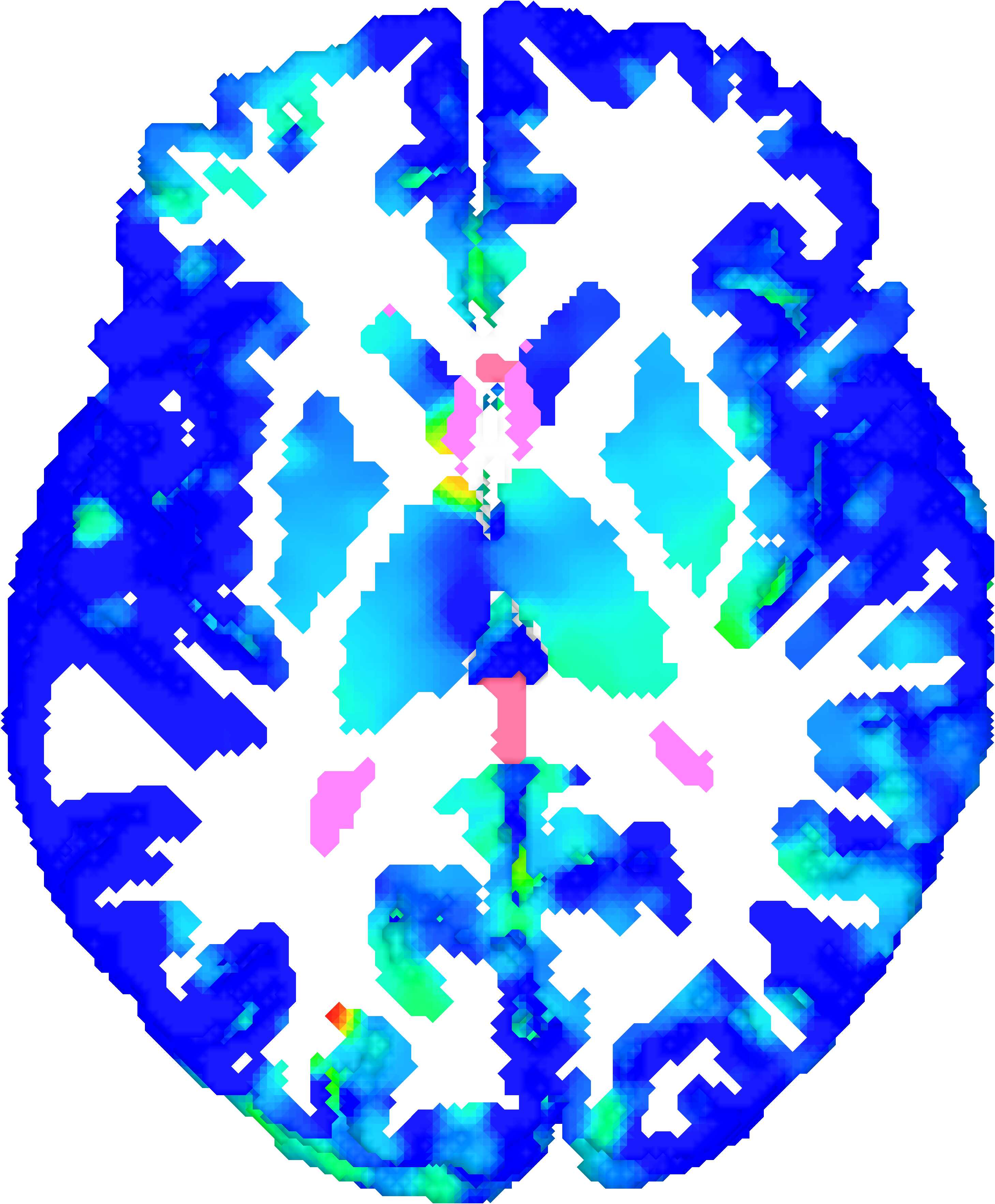}
            \caption{$\azimuth=\qty{0}{\degree}$, $\elevation=\qty{90}{\degree}$, $z=\qty{37}{\milli\meter}$.}
        \end{subfigure}
        \hfill
        \begin{subfigure}[b]{\subFigWidth}
            \includegraphics[width=\linewidth]{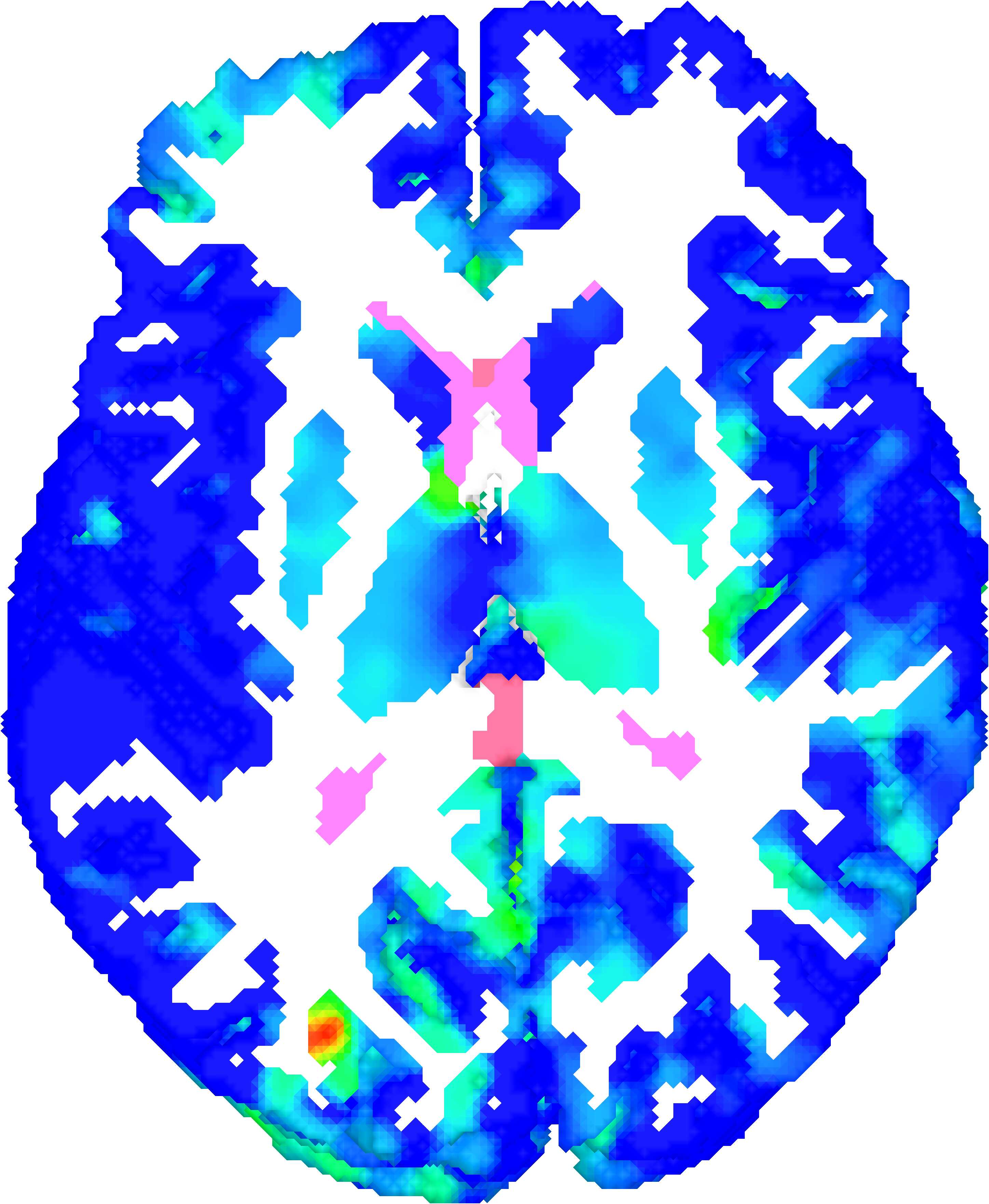}
            \caption{$\azimuth=\qty{0}{\degree}$, $\elevation=\qty{90}{\degree}$, $z=\qty{40}{\milli\meter}$.}
        \end{subfigure}
    \end{minipage}
    \hfill
    \begin{minipage}{0.09\linewidth}
    \includegraphics[width=\linewidth]{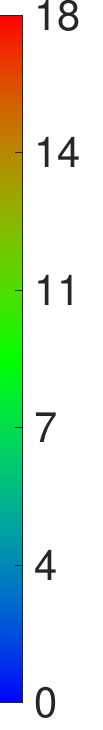}
    \end{minipage}
    \caption{Relative differences \eqref{eq:rel-diff} between the interference fields \eqref{eq:interference-field} of $\linJv$ and $\refJv$, smoothed via the function $\dBT\args{\someField/\max\someField,\qty{18}{\decibel}}$ of equation \eqref{eq:visual-log-transform} at different azimuths $\azimuth$ and elevations $\elevation$, with base stimulation frequency $\frequency=\qty{1000}{\hertz}$ and an updated $\contactResistance=\qty{5270}{\ohm}$. The current pattern was $\currentPattern=[-0.5,0.5,-1.5,1.5]\,\si{\milli\ampere}$. The maximal value of the absolute differences was roughly \qty{0.99e-6}{\milli\ampere\per\milli\meter\squared}.}
    \label{fig:linrefreldiff-Rc=5270-i=0.5-1.5}
\end{figure}

Figures~\ref{fig:linrefreldiff-Rc=1270-i=1.5-0.5}--\ref{fig:linrefreldiff-Rc=5270-i=0.5-1.5}
display relative differences~\eqref{eq:rel-diff} between
the interference fields \eqref{eq:interference-field} of
$\linJv$ and $\refJv$, with different stimulation current
patterns $\currentPattern$. A common feature regardless
of the applied current pattern seen in these comparisons
is that the greatest differences seem to occur near the
unaltered electrodes C4 and TF10, on the opposite side of
the volume conductor, with the point of reference being
the modified electrode TP9. The effect of the current
pattern on the positions of largest errors is best seen
in Figures~\ref{fig:linrefreldiff-Rc=1270-i=0.5-1.5}
and \ref{fig:linrefreldiff-Rc=5270-i=0.5-1.5}, where an
increased difference is seen in the left hemisphere, where
the maximum of the amplitude modulation is also located.
As the stimulating current pattern is altered such that the
modulation maximum moves first towards the central fissure,
as in
Figures~\ref{fig:interference-field-equal-currents} and
\ref{fig:linrefreldiff-Rc=1270-i=1.0-1.0}--\ref{fig:linrefreldiff-Rc=5270-i=1.0-1.0},
and further towards the right hemisphere, as in
Figures~\ref{fig:interference-field-weaker-right} and
\ref{fig:linrefreldiff-Rc=1270-i=1.5-0.5}--\ref{fig:linrefreldiff-Rc=5270-i=1.5-0.5},
the largest differences seem to follow along with the field
maxima.

\section{Discussion}%
\label{sec:discussion}

In this  study, we presented an extended
CEM-based~\cite{Somersalo1992, GALAZPRIETO-2022} forward
model and its linearization for tTIS~\cite{GROSSMAN-2017,
RAMPERSAD_computational_2019, Wessel_2021, Herrmann_2021,
tTIS_2022}, where both the tissue admittivities
$\admittivity = \conductivity + \iu\angFreq\permittivity$
and electrode impedances $\Zell$ are allowed to be
complex-valued. In the case of the tissue parameter, the
added imaginary part encodes the capacitive properties of
the tissue. The same is true for the electrodes, as the
usual model for the impedance of a simple double-layer
electrode only contains resistive and capacitive
components~\cite{BUTSON-2005, electrode-modelling-2017}.
CEM  allows including the electrodes as boundary patches
with properties determined by a set of boundary conditions
for the governing partial differential equation, taking
shunting effects into account~\cite{Agsten-2018}.
It has proven to be necessary in EIT of the
brain~\cite{darbas2021sensitivity, shi2023robust,
fernandez2016effects}, where just like in tTIS, a set
of alternating currents is injected transcranially.
In the spirit of Grossman et al.~\cite{GROSSMAN-2017},
modifying the ratio of these currents, while keeping
the total injected current constant, allowed us to steer
the amplitude modulation maximum from one hemisphere to
another, which we demonstrate in our numerical experiments.
Additionally, we reduced the quality of a contact electrode
and observed what kind of a deviation this causes when the
lead and interference fields are computed via a reference
FEM simulation, and a surrogate linearized model.

As our computational domain, we utilized an open-access
dataset featuring a highly detailed multicompartment
head model of a healthy male subject~\cite{Piastra_2020}.
The model was generated by extracting a triangular
surface element mesh from an MRI image using
FreeSurfer~\cite{FreeSurferWebSite,fischl2012freesurfer}
and then running the tetrahedral mesh generation algorithm
of Zeffiro Interface~\cite{galaz-prieto-etal-2023} without post
processing on the surfaces to produce a volumetric mesh.
The benefit of this approach as compared to software such
as SimNIBS~\cite{simnibs-documentation,wang-etal-2025}
is that at least by default, the finite element meshes
produced by its \texttt{charm} pipeline can consist of
elements of very differing sizes, with larger elements
being present in large regions of tissue homogeneity, with
smaller compartments consisting of smaller elements. While
this approach saves computational resources such as memory
and solver iterations due to smaller matrix sizes, there
might be large jumps between neighbouring element sizes,
which can produce modelling errors especially related
to field smoothness. Smooth transitions in element sizes
allows for modelling smoother fields more accurately. This
is why we chose a mesh with very regular element size over
computational efficiency.

The use of CEM over alternatives such as including the
electrodes as a part of the finite element mesh itself
was made due to it allowing us to quickly experiment
with electrode positions and sizes, and due to the chosen
regions of interest for the stimulation being relatively
far away from the electrodes. Embedding the electrodes
into the mesh itself would have required re-meshing, if the
position or size of the electrodes was changed. With CEM
it is simple to roughly position an electrode, and project
it to the nearest mesh node, while computing the contact
area based on the assigned electrode radius. Also, while
we simply averaged the electrode potential over the entire
electrode patches, technically nothing is preventing a more
complex potential distribution from being applied to the
nodes within the contact area of an electrode patch, via
methods such as interpolation between the electrode center
and the boundary of its contact area. Finally, including
the electrodes into the 3D domain itself would mainly be
beneficial, if one was interested in what was happening
in the immediate vicinity of and possibly within the
electrodes themselves.

Related to this, our chosen double layer electrode
capacitance was slightly small, however. We
followed the example of \cite{BUTSON-2005},
and chose our capacitance as an absolute value
$\doubleLayerCapacitance=\qty{3.3}{\micro\farad}$,
whereas relative values between \num
6--\qty{12}{\micro\farad\per\centi\meter\squared}
reported by \cite{khademi-etal-2020} might have been more
appropriate. However, the value is not orders of magnitude
different from the suggested values, and therefore this
choice of $\doubleLayerCapacitance$ was unlikely to affect
our resuls greatly.

This finite element model, along with the CEM electrode configuration
derived from the 10--20 system~\cite{Nuwer-2018},
provided a realistic framework for our numerical
simulations. The mathematical forward tTIS model
presented in this paper extends a similar model related
to transcranial electric stimulation via direct currents,
presented in~\cite{GALAZPRIETO-2022}. The electrode
capacitance was implemented assuming that the utilized
electrodes were assumed to be simple double-layer
electrodes~\cite{BUTSON-2005, electrode-modelling-2017},
where in addition to the electrode--skin contact resistance
$\contactResistance$, there was a parallel connection
between a capacitor and a resistor with respective
capacitance $\doubleLayerCapacitance$ and resistance
$\doubleLayerResistance$.

A theoretical justification for this omission of
tissue capacitance is the insignificantness of
the time derivatives in the utilized quasi-static
versions of Maxwell's equations~\cite{Hamalainen1993,
Knösche--Haueisen-2022}, which nullifies the reactive
component of the tissue~\cite{GROSSMAN-2017, BUTSON-2005}.
Note that this is in direct opposition to the advice of
\cite{BUTSON-2005}, who recommend that  where tissue volume
activation is concerned, it is tissue capacitance that
should be taken into account. However,~\cite{BUTSON-2005}
also states that with relative permittivities
$\relativepermittivity < \num{10000}$, the modelling errors
introduced by this omission seem to be insignificant.

A practical reason for the omission of the tissue
capacitance was also present: the iterative
numerical solvers we utilized failed to successfully
compute a transfer matrix $\transferMat =
\inverse\stiffMat\electrodeCurrentMat$ with a
complex-valued stiffness matrix $\stiffMat$. The attempted
solvers included our in-house PCG solver, the stabilized
bi-conjugate gradient (\biCGStab) solver of MATLAB and
our own version of \biCGStab, which in our internal tests
performed better than that of MATLAB on smaller matrices of
sizes up to $100\times100$. We therefore chose to leave the
complex tissue parameter out of these results.

While an existing linearization can be considered as a
potential surrogate to adapt the numerical system with
regard to changing electrode impedance, the linearization
of the resistance matrix~\(\resistanceMat\) can \emph{not}
be expected to generally speed up the computation of the
lead field~\(\leadFieldMat\). A speed-up only occurs, if
the stimulating electrodes \(\eell\) are very small, such
that the number of non-zero colums in the related mass
matrix \(\massMat=\sum_\ell\massMat^{\paren\ell}\) does
not reach the total number of electrodes attached to the
volume conductor $\domain$, here $4$. This is because to
solve for the total derivative~\(\der\resistanceMat\Zell\),
we need to solve for matrices such as
\(\inverse\stiffMat\pder{\stiffMat}{\Zell}\) for
each modified \(\Zell\), and the method utilized
in the inversion is the same iterative \PCG\
algorithm~\cite{sauer2018numerical} used in the
computation of the transfer matrix~\(\transferMat =
\inverse\stiffMat\electrodeCurrentMat\), where the number
of columns in $\electrodeCurrentMat$ corresponds to the
number of stimulating electrodes. The number of nodes
involved in a contact surface of an electrode grows
past the $4$ columns of $\electrodeCurrentMat$ very fast
when electrode radii are increased, especially with a
high-resolution mesh.

Our mesh was high-resolution indeed. The sparse stiffness
matrix~\(\stiffMat\in\Cset^{\Nof N\times\Nof N}\) computed
for the utilized head model (and therefore \(\massMat\))
does not fit into the memory of a modern computer, unless
implemented in a sparse matrix format. The number of
nodes in the utilized mesh was \(\Nof N = \num{4 877
338}\) and the data was stored in double-precision format,
meaning the amount of memory needed to store a single
full column of \(\electrodeCurrentMat\) or \(\massMat\)
during the inversion of \(\stiffMat\) against them was
\(\Nof N\cdot \qty{8}{\byte} = \qty{39018704}{\byte} \approx
\qty{39}{\mega\byte}\). This was much more manageable.
However, a speedup in the computation \emph{cannot} be
considered as a benefit of this linearization involving
$\massMat$.

The interference field also turned out to be steerable, as
described by~\cite{GROSSMAN-2017}, with the total current
$\currentPattern$ fed into the volume conductor kept at
a constant \qty{2}{\milli\ampere}, but with the current
coming through the left and right electrode pairs C3+TP9
and C4+TF10 of the 10--20 system~\cite{Nuwer-2018} varying
between \num{0.5}, \num{1.0} and \qty{1.5}{\milli\ampere}.
The maximum of the interference field always occurred
near a point where the contributing fields had similar
magnitudes. When equal currents were used, the maximum
appeared in and around the central fissure, with one
peak located near the front of the thalamus. With unequal
currents, the field maximum always moved towards the weaker
of the currents, which is understandable as the field from
the stronger electrodes is expected to reach further into
the volume conductor, whereas the field from the weaker
current source should attenuate faster.

Worsening the contact resistance $\contactResistance$
from \qty{270}{\ohm} to \qty{1270}{\ohm} and further
to \qty{5270}{\ohm} in the case of the reference field
$\refJv$ made no large difference in the interference
field magnitudes or the field shape itself. The maxima
still occurred where they did in the case of the
original field. Their magnitudes were changed, but only
slightly. A \emph{slightly} stronger field was observed
when $\contactResistance=\qty{1270}{\ohm}$ versus the
\qty{5270}{\ohm} case. The steerability of the field also
did not visibly suffer when compared to the case of the
original contact resistance pattern.

We then observed the relative difference of the
reference and linearized lead fields $\linL$ and
$\refL$ at two different changes in the contact
impedance $\contactResistance$ of electrode TP9, with
$\Delta\contactResistance\in\set{\num{1000},\linebreak\num{
5000}}\,\si\ohm$. These increases in $\contactResistance$
correspond to a weakening of the electrode contact,
which originally had a contact impedance suitable for
measurements: \qty{270}{\ohm}. Somewhat unsuprisingly,
with $\contactResistance=\qty{1270}{\ohm}$ at the electrode
TP9, the greatest difference in the lead fields was
produced at the electrode TP9, whose contact impedance
was increased. The differences quickly diminish when
moving further away from the electrode and into the
volume conductor. With a maximal difference value of
\qty{0.13e-3}{\per\milli\meter\squared}, the largest
differences in the observed regions were in the vicinity
of \qty{14}{\decibel}, with most of the difference range
near the electrode staying around \qty{9}{\decibel}. At
$\contactResistance=\qty{5270}{\ohm}$, we start to see a
non-zero difference in the distribution near electrodes
whose contact resistance was not altered at all as
well. However, these differences do not climb much
above \qty{7}{\decibel}, when a maximal difference value
\qty{0.16e-3}{\per\milli\meter\squared} was observed,
leaving the distribution otherwise very similar to the
$\contactResistance=\qty{1270}{\ohm}$ case.

Our final numerical experiment observed how the
reference and linearized interference current fields
$\interferenceField\paren{\volumeCurrentDensity_1,\volumeCurrentDensity_2}$
themselves differed from each
other on relative terms, when the volume current
densities $\volumeCurrentDensity$ were produced
by the lead fields $\refL$ and $\linL$, again with
$\Delta\contactResistance\in\set{\num{1000},\num{5000}}
\,\si\ohm$. Also, the visualized dynamic range
was once again \qty{18}{\decibel}, with a
maximal absolute difference in the range \num
1--\qty{3e-6}{\milli\ampere\per\milli\meter\squared},
depending on where the interference maximum was steered.
Here the largest differences occurred \emph{not} near
the electrode TP9, whose contact deteriorated, but
mostly near the electrode pair C4--TF10, which did not
experience any change in their contact properties, on
the opposite side of the volume conductor. In addition
to there always being consistent discrepancies in the
reference and linearized interference fields near the
unmodified electrodes, modifying the injected current
pattern such that the peak of the interference pattern
moved from the right lobe to the central fissure and
then towards the left hemisphere also caused larger
discrepancies to appear in those locations. However,
the largest absolute differences occurred in the case
where the amplitude modulation maximum of approximately
\qty{3e-6}{\milli\ampere\per\milli\meter\squared}
was near the unmodified electrodes C4--FT10. When
the peak was moved towards the modified electrode
TP9, the largest differences diminished to
\qty{1e-6}{\milli\ampere\per\milli\meter\squared}. Apart
from the actual maxima, much of the differences stayed in
the vicinity of \qty{-9}{\decibel}.

Our future work will include utilizing the linearization
approach in optimization applications. By incorporating
this model into optimization frameworks, it would
be possible to refine and enhance the design of tTIS
protocols, including the optimal placement of electrodes
and the selection of stimulation frequencies. This could
lead to more precise and effective stimulation strategies,
improving the outcomes in both therapeutic and research
contexts. Further experimental validation of these
optimizations will also be essential to ensure their
practical applicability and effectiveness.

Additional future work might consist of developing
an iterative solver that can handle inverting the
complex-valued tissue parameter matrix $\stiffMat$ against
the electrode current matrix $\electrodeCurrentMat$
when solving for a transfer matrix $\transferMat$. Our
in-house implemetation of the \biCGStab\ algorithm did
not implement the higher order polynomial stabilizing
step presented in~\cite{bicgstab-2022}. The use of
simple \biCGStab\ instead of $\biCGStab(L>1)$, where $L$
is the order of the stabilizing polynomial might have
partly explained the convergence issues. The randomized
shadow residual mentioned in the paper and which was
included in the implementation was not enough to tackle
the issue. Alternatively, the inversion of $\stiffMat$
might be framed using the \CtoR\ method described by
\cite{axelsson-etal-2014}, with a suitable preconditioner
to enable or speed up the convergence.

To handle more carrier frequencies than two, the method
suggested by \cite{botzanowski-etal-2025} might also be
considered in the future. Lead fields $\leadFieldMat$ for
different frequency responses might be computed, and the
suggested box-car averaging and root-mean-squaring of the
signals mapped to a region of interest performed, to find
out an approximation of the modulating signal there.

Alternative waveforms discussed by \cite{luff-etal-2024}
could also be implemented in a future study. Again, an
electrical stimulation lead field $\leadFieldMat$ is
technically a mapping of a unit stimulation current with a
given frequency to a given RoI. The mapping $\leadFieldMat$
takes no stance on what the shape of the input current
waveform is. One could then map a square wave to a region
of interest and use the Hilbert transform described in
\cite{luff-etal-2024} to extract the modulating signal
over time.

\section*{Funding and Acknowledgments}%
\label{sec:ack}

The work of Maryam Samavaki and Sampsa Pursiainen is
supported by the Research Council of Finland's (RCF's)
Centre of Excellence in Inverse Modelling and imaging
2018-2025, decision 359185, Flagship of Advanced
Mathematics for Sensing, Imaging and  Modelling, RCF's
decision 359185; Santtu Söderholm and Maryam Samavaki
have been supported by the ERA PerMed project PerEpi
(PERsonalized diagnosis and treatment for refractory focal
paediatric and adult EPIlepsy), RCF's decision 344712
and also project 359198. In 2025-2029, Maryam Samavaki
serves as an Academy Research Fellow (RCF 371055). We
thank Prof.\ Carsten H. Wolters, University of Münster,
Münster, Germany, for the fruitful discussions and support
in mathematical modelling and are grateful to DAAD (German
Academic Exchange Service) and RCF for supporting our
travels to Münster (RCF decision 354976, 367453).

\section*{Conflicts of Interest}%
\label{sec:conflict}

The authors confirm that the research utilized in this
study was entirely independent, open, and academic. They
have neither financial nor non-financial relationships,
affiliations, knowledge, or beliefs in the subject matter
or materials included in this manuscript, nor do they have
any involvement with or affiliation with any organization
or personal relationship.

\section*{Ethics statement}\label{sec:ethics-statement}

The research related to this paper did not involve any live
human subjects or other research objects, that would have
required an explicit permission from the target of research
or any other body related to the target of research. All
research data used in the study was in openly accessible
and permissively licensed form.

\begingroup
\interlinepenalty=10000 
\bibliographystyle{elsarticle-num}
\bibliography{references.bib}
\endgroup

\appendix

\section{Envelope of two interfering sinusoidal waves}%
\label{app:envelope}

Given two sinudoidal waves $w_1 = A_1
\sin\paren{\frequency_1 t + \theta_1}$ and $w_2 = A_2
\sin\paren{\frequency_2 t + \theta_2}$, their sum $w_1 +
w_2$ can be expressed via complex exponential notation or
Euler's formula as follows:
\begin{equation}
\begin{aligned}
    &
    w_1 + w_2
    \\
    &=
    A_1 \sin\paren{\frequency_1 t + \theta_1}
    +
    A_2 \sin\paren{\frequency_2 t + \theta_2}
    \\
    &=
    A_1
    \imagPart \paren*{
        \exp\paren*{\iu\theta_1}
        \exp\paren*{\iu\frequency_1 t}
    }
    +
    A_2
    \imagPart\paren*{
        \exp\paren*{\iu\theta_2}
        \exp\paren*{\iu\frequency_2 t}
    }
    \\
    &=
    \imagPart\paren*{
        A_1
        \exp\paren*{\iu\theta_1}
        \exp\paren*{\iu\frequency_1 t}
        +
        A_2
        \exp\paren*{\iu\theta_2}
        \exp\paren*{\iu\frequency_2 t}
    }
    \,.
\end{aligned}
\end{equation}
Here $\imagPart$ is the imaginary part of a complex
number. Remembering the law of cosines or the generalized
Pythagorean theorem for two complex numbers $c_1$ and
$c_2$, we have the relation
\begin{equation}
    \abs{c_1 + c_2}^2
    =
    \abs{c_1}^2
    +
    \abs{c_2}^2
    +
    2
    \abs{c_1}
    \abs{c_2}
    \cos\paren{\phi_1 - \phi_2}
    \,,
\end{equation}
where $\phi_1$ is the argument of $c_1$, and
$\phi_2$ the argument of $c_2$. Since the unit circle
$\exp\paren{\iu\phi}$ with $\phi\in\Rset$ has a modulus of
$1$, the amplitude of $w_1$ is entirely determined by $A_1$
and similarly for $w_2$ and $A_2$. Therefore
\begin{equation}
\begin{aligned}
    &
    \abs{w_1 + w_2} ^ 2
    \\
    &=
    A_1^2
    +
    A_2^2
    +
    2
    A_1
    A_2
    \cos\paren{
        \paren{\frequency_1 t + \theta_1}
        -
        \paren{\frequency_2 t + \theta_2}
    }
    \\
    &=
    A_1^2
    +
    A_2^2
    +
    2
    A_1
    A_2
    \cos\paren{
        \paren{\frequency_1 - \frequency_2}
        t
        +
        \paren{\theta_1 - \theta_2}
    }
    \\
    &=
    A_1^2
    +
    A_2^2
    +
    2
    A_1
    A_2
    \cos\paren{
        \beatF
        t
        +
        \Delta\theta
    }
    \,,
\end{aligned}
\end{equation}
where $\beatF$ is the beat frequency between the signals
and $\Delta\theta$ their phase difference. The actual
time-varying upper envelope is then
\begin{equation}\label{eq:abs-of-sum}
    \abs{w_1 + w_2}
    =
    \sqrt{
        A_1^2
        +
        A_2^2
        +
        2
        A_1
        A_2
        \cos\paren{
            \beatF
            t
            +
            \Delta\theta
        }
    }
    \,.
\end{equation}

To observe how the modulus $\abs{w_1 - w_2}$ would
behave, we note that $-w_2$ is structurally similar
to $w_2$, except it has been phase-shifted by $\pi$ or
\qty{180}{\degree}. Then by $\cos\paren{\theta + \pi} =
-\cos\theta$ we have a similar relation
\begin{equation}\label{eq:abs-of-diff}
    \begin{aligned}
        &
        \abs{w_1 - w_2}
        \\
        &=
        \sqrt{
            A_1^2
            +
            A_2^2
            +
            2
            A_1
            A_2
            \cos\paren{
                \beatF
                t
                +
                \paren{\Delta\theta + \pi}
            }
        }
        \\
        &=
        \sqrt{
            A_1^2
            +
            A_2^2
            -
            2
            A_1
            A_2
            \cos\paren{
                \beatF
                t
                +
                \Delta\theta
            }
        }
        \,,
    \end{aligned}
\end{equation}
where the only real difference to \eqref{eq:abs-of-sum} is
the minus sign in front of the mixed term in the radical.

\section{Admittivity tensor}%
\label{app:admittivity-tensor}

The distribution of electric potential resulting
from a current stimulation by solving the equation
(\ref{poisson_com_2}) in the head model $\domain$ and
considering compatibility for each frequency component by
writing the admittivity in the following diagonal form:
\begin{equation}
    \admittivity(\position, \angFreq)
    =
    \begin{bmatrix}
        \admittivity(\position, \angFreq_1)
        &
        {\bm 0}
        &
        \cdots
        &
        {\bm 0}
        \\
        {\bm 0}
        &
        \admittivity(\position,  \angFreq_2)
        &
        \cdots
        &
        {\bm 0}
        \\
        \vdots
        &
        \vdots
        &
        \ddots
        &
        \vdots
        \\
        {\bm 0}
        &
        {\bm 0}
        &
        \cdots
        &
        \admittivity(\position, \angFreq_L)
    \end{bmatrix}.
\end{equation}
For a given position $\position$ and angular frequency
$\angFreq$, the components $\admittivity_{ij}(\position,
\angFreq)$ of the two-rank tensor $\admittivity(\position,
\angFreq)$ can be expressed as follows:
\begin{equation}
    \admittivity_{ij}(\position, \angFreq)
    =
    \begin{cases}
        { \conductivity}_{ij}(\mathbf{x}, \angFreq)
        +
        \iu
        \angFreq
        \,
        \permittivity[0]
        \,
        \,\permittivity_{ij}(\mathbf{x}, \angFreq)
        \,,
        &
        \hbox{if}
        \quad i=j
        \\
        {0} ,
        &
        \text{otherwise}
        \,.
    \end{cases}
\end{equation}
In this way, the matrix form of the symmetric conductivity
distribution $\conductivity(\mathbf{x}, \angFreq):
\domain\rightarrow\mathbb{C}^{N \times N}$ is:
\begin{equation}\label{eq:conductivity-tensor}
    \conductivity(\mathbf{x}, \angFreq)
    =
    \begin{bmatrix}
        \conductivity_{11}
        \paren{\mathbf{x}, \angFreq}
        &
        \conductivity_{12}
        \paren{\mathbf{x}, \angFreq}
        &
        \cdots
        &
        \conductivity_{1 N}
        \paren{\mathbf{x}, \angFreq}
        \\
        \conductivity_{21}
        \paren{\mathbf{x}, \angFreq}
        &
        \conductivity_{22}
        \paren{\mathbf{x}, \angFreq}
        &
        \cdots
        &
        \conductivity_{2 N}
        \paren{\mathbf{x}, \angFreq}
        \\
        \vdots
        &
        \vdots
        &
        \ddots
        &
        \vdots
        \\
        \conductivity_{N 1}
        \paren{\mathbf{x}, \angFreq}
        &
        \conductivity_{N 2}
        \paren{\mathbf{x}, \angFreq}
        &
        \cdots
        &
        \conductivity_{N N}
        \paren{\mathbf{x}, \angFreq}
    \end{bmatrix}
    \,.
\end{equation}
And the matrix representation for the permittivity
$\permittivity(\mathbf{x}, \angFreq): \domain \rightarrow
\mathbb{C}^{N \times N}$ is:
\begin{equation}\label{eq:permittivity-tensor}
    \permittivity(\mathbf{x}, \angFreq)
    =
    \begin{bmatrix}
        \permittivity_{ 11}(\mathbf{x}, \angFreq)
        &
        \permittivity_{ 12}(\mathbf{x}, \angFreq)
        &
        \cdots
        &
        \permittivity_{ 1 N}(\mathbf{x}, \angFreq)
        \\
        \permittivity_{ 21}(\mathbf{x}, \angFreq)
        &
        \permittivity_{ 22}(\mathbf{x}, \angFreq)
        &
        \cdots
        &
        \permittivity_{2 N}\left(\mathbf{x}, \angFreq\right)
        \\
        \vdots
        &
        \vdots
        &
        \ddots
        &
        \vdots
        \\
        \permittivity_{ N 1}(\mathbf{x}, \angFreq)
        &
        \permittivity_{ N 2}(\mathbf{x}, \angFreq)
        &
        \cdots
        &
        \permittivity_{ N N}(\mathbf{x}, \angFreq)
    \end{bmatrix}
    \,.
\end{equation}

\subsection{Weak Form}%
\label{app:forward_model}

The numerical solution of the complete electrode
model \eqref{CEM1}--\eqref{CEM5} relies on FEM.
The FEM approximation for this model, derived
in~\cite{GALAZPRIETO-2022}, initiates with the variational
formulation of the problem. A weak form for the electric
potential field $ \potential\in\sobolevSpace^1(\domain)$
can be derived through integration by parts. In this paper,
$\sobolevSpace^1(\domain)$ represents a \textit{Sobolev
space} of square-integrable ($\int_{\domain}
\abs\potential^2 \, \diff \position_\domain < \infty$)
functions with square-integrable partial derivatives:
\[
    \sobolevSpace^1(\domain)
    =
    \set{
        \potential\in\mathrm{L}^2(\domain)
        \,|\,
        \nabla\potential\in\mathrm{L}^2(\domain )
    }\,.
\]
By multiplying the equation~\eqref{eq:admittivity} with a
 sufficiently smooth test function $\testFn\in\surface$,
 where $\surface$ is a subspace of $\sobolevSpace^1(\domain)$, and
 integrating by parts, we get the following formulation:
\begingroup\small
    \begin{equation}
        \begin{aligned}
            0
            &=
            -\int_\domain
            \nabla
            \cdot
            (\admittivity \nabla{  u})
            \,{\conj\testFn}
            \, \diff\position_\domain
            \\
            &=
            \int_\domain
            (\conductivity+\iu
            \angFreq \permittivity)
            \,(\nabla {  u}
            \cdot
            \nabla
            {\conj\testFn} )
            \, \diff\position_\domain
            -
            \int_{\electrode[\ell]}
            (\conductivity+\iu
            \angFreq \permittivity)
            \,
            \frac{\partial { {  u}}}{\partial {\bavec n}}
            \,
            {\conj\testFn}
            \,
            \diff\surface
            \\
            &=
            \int_\domain
            (\conductivity+\iu
            \angFreq\permittivity)
            \,
            (\nabla u\cdot \nabla  {\conj\testFn})
            \,
            \diff\position_\domain
            -
            \sum_{\ell = 1}^L
            \int_{\electrode[\ell]}
            (\conductivity+\iu
            \angFreq \permittivity)
            \,
            \frac{\partial {  {  u}}}{\partial {\bavec n}}
            {\conj\testFn}
            \,
            \diff\surface\,.
        \end{aligned}
        \label{eq:integration}
    \end{equation}
\endgroup

Furthermore, we can present the following equations derived
from the complex boundary conditions governing CEM:
\begin{equation}
\begin{aligned}
    -
    \sum_{\ell = 1}^L
    \int_{\electrode[\ell]}
        \paren{
            \conductivity
            +
            \iu
            \angFreq
            \permittivity
        }
        \,
        \pder
            \potential
            {\bavec n}
        {\conj\testFn}
        \,
    \diff\surface
    &=
    -
    \sum_{\ell = 1}^L
    \,
    \frac
        {\potential_\ell}
        { \Zvec_\ell\,\areaOf\electrode }
    \int_{\electrode[\ell]}
        {\conj\testFn}
        \,
    \diff\surface
    \\
    &+
    \sum_{\ell = 1}^L
    \,
    \frac
        {1}
        {\Zvec_\ell\,\areaOf\electrode}
    \int_{\electrode[\ell]}
        \potential
        {\conj\testFn}
    \diff\surface
    \,.
\end{aligned}
\label{eq:integration_bc}
\end{equation}
Here \(\areaOf\electrode\) is the contact area of electrode
\(\electrode[\ell]\). Then
\begingroup\small
\begin{equation}
\label{eq:weakform}
\begin{aligned}
    \int_\domain
        \paren{
            \conductivity
            +
            \iu
            \angFreq
            \permittivity
        }
        \,
        \paren{
            \nabla\potential
            \cdot
            \nabla {\conj\testFn}
        }
    \, \diff\position_\domain
    &=
    \sum_{\ell = 1}^L
        \frac{1}{ \Zvec_\ell \, \areaOf\electrode^2}
        \int_{\electrode[\ell]} u \, \diff\surface
        \int_{\electrode[\ell]} {\conj\testFn} \, \diff\surface
    \\
   &+
    \sum_{\ell = 1}^L   \frac{ \current_\ell}{  \areaOf\electrode}\,
        \int_{\electrode[\ell]} {\conj\testFn} \, \diff\surface
    \\
    &-
    \sum_{\ell = 1}^L\,\frac{1}{\Zvec_\ell\,\areaOf\electrode} \int_{\electrode[\ell]} u \,{\conj\testFn} \, \diff\surface\
    \end{aligned}
\end{equation}
\endgroup
for all $\testFn\in\surface$. Thus, the problem statement
under observation becomes:
\begin{quote}
    Find $u\in\surface$ such that, for all $\testFn\in\surface$
\begin{equation}
\label{weak_com}
\begin{aligned}
    a(u, \testFn)
     &=
    \int_\domain
    \paren{
        \conductivity+\iu
        \angFreq \permittivity
    }
    \paren{
        \nabla u \cdot \nabla{\conj\testFn}
    }
    \,
    \diff\position_\domain.
\end{aligned}
\end{equation}
\end{quote}
The continuous bilinear form $a:\surface\times\surface\rightarrow
\mathbb{C}$ is defined as follows:
 \begin{align*}
a(u, \testFn)&:=
    \sum_{\ell = 1}^L
        \frac{1}{ \Zvec_\ell \,\areaOf\electrode^2}
        \int_{\electrode[\ell]} u \, \diff\surface
        \int_{\electrode[\ell]} {\conj\testFn} \, \diff\surface +
    \sum_{\ell = 1}^L   \frac{ \current_\ell}{ \areaOf\electrode}\,
        \int_{\electrode[\ell]} {\conj\testFn} \, \diff\surface
    \\
    &-
    \sum_{\ell = 1}^L\,\frac{1}{\Zvec_\ell\,\areaOf\electrode} \int_{\electrode[\ell]} u \,{\conj\testFn} \, \diff\surface\,.
\end{align*}
Now, we introduce the finite element space,
\[
    {\mathbb{V}_h}
    =
    \vecspan\set [\big] {
        \basisFn_1,\dots,\basisFn_N,
        \iu\,\basisFn_1,\dots,\iu\,\basisFn_N
    }
    \subset \mathrm{H}^1(\mathbb{V})
    \,,
\]
where $\{\basisFn_j\}_{j=1}^{N}\in\surface$ are real and
periodic functions, $N$ is the number of nodes in the FE
mesh, $\basisFn_i$ are the (piecewise linear) nodal basis
functions of the mesh, and \(\iu\) is the imaginary unit.
This leads to the discrete variational problem
\begin{quote}
    Find $u_h\in {\mathbb{V}_h}$ such that, for all ${ \testFn}\in\surface$
    \begin{equation}
    \begin{aligned}
        a(u_h, { \testFn})
        =
        \int_\domain (\conductivity
        +
        \iu\angFreq \permittivity)
        \
        ,(\nabla u_h \cdot \nabla {\conj\testFn})
        \,
        \diff\position_\domain
\end{aligned}
\label{vari}
\end{equation}
\end{quote}
The potential function $u$ takes the following form:
\begin{align*}
    u(\position, \angFreq_\ell)
    =
    \sum_{j=1}^{N}
    u_j(\angFreq_\ell)
    \,
    \basisFn_j(\position)
    +
    \iu
    \sum_{j=1}^{N}
    z_j(\angFreq_\ell)\,\basisFn_j(\position)
    \,,
\end{align*}
Then, the argument \eqref{vari} is equivalent to
\begin{quote}
    Find $u_h\in {\mathbb{V}_h}$ such that, for all $\basisFn_j\in {\mathbb{V}_h}$
\begin{equation}
\label{weak_com_2}
\begin{aligned}
     a(u_h, \basisFn_j)
     =
     \int_\domain (\conductivity
     +
     \iu\angFreq \permittivity)
     \,
     (\nabla u_h \cdot \nabla \basisFn_j)
     \,
     \diff\position_\domain
     \,,
\end{aligned}
\end{equation}
with \(1 \leq j \leq N\).
\end{quote}

\subsection{Resistance Matrix}%
\label{app:res_matrix}

In terms of the actual implementation in Matlab, since
$\potential$ and $\admittivity$ are both complex-valued,
equation \eqref{poisson_com_c} cannot be solved
directly in the PDE Toolbox. To this end, we represented
\eqref{poisson_com_c} in matrix form as follows: Given
the scalar valued functions $\basisFn_1, \basisFn_2,
\ldots, \basisFn_{N} \in \mathcal{S}$, the potential
distribution $\potential$ in the domain $\domain$
can be approximated as the finite sum $u = \sum_{i =
1}^{N} u_i \basisFn_i$. Denoting by $\potentialVec =
(\potential_1,\ldots, \potential_N)\in \mathbb{C}^{N
\times 1}$ the coefficient vector of the discretized
potential, by $\potentialLossVec=(\potentialLoss_1,\ldots,
\potentialLoss_L)\in \mathbb{C}^{L \times 1}$
the potential losses across electrodes, and by
$\currentPattern=(\current_1,\ldots, \current_L)\in
\mathbb{C}^{L \times 1}$ the injected current pattern, the
weak form~\eqref{eq:weakform} is given by
\begin{equation}
\label{eq:u_system}
     \begin{pmatrix}
     \stiffMat & -\electrodeCurrentMat \\
    -\transpose{\electrodeCurrentMat} & {\electrodeImpedanceMat}
    \end{pmatrix}
    \begin{pmatrix}
    \potentialVec  \\
    \potentialLossVec
    \end{pmatrix}
    =
    \begin{pmatrix}
    {\bf 0} \\
    \currentPattern
    \end{pmatrix}
    \,.
\end{equation}
The components of the matrix system are defined as follows:
\begin{align}
    \stiffMat_{i j}
    &=
    \int_\domain
    \admittivity
    \nabla
    \conj{\basisFn_i}
    \cdot
    \nabla \basisFn_j
    \,
    \diff\position_\domain
    +
    \sum_{\ell = 1}^L
    \frac{1}{\Zell \areaOf\electrode}
    \int_\eell
    \conj{\basisFn_i}
    \basisFn_j
    \,
    \diff\surface\,,
    \label{eq:aij}
    \\
    \electrodeCurrentMat_{i\ell}
    & =
    \frac{1}{\Zell\,\areaOf\electrode}
    \int_{e_{\ell}}
    \conj{\basisFn_i }
    \,
    \diff\surface
    \,,\label{eq:bil} \\
    \electrodeImpedanceMat_{h\ell}
    & =
    \diag
    \set*{
        \frac{\int_{e_h}\, \diff\surface}{ \Zell\,\areaOf\electrode}
    }
    \,.
    \label{eq:chl}
\end{align}
\label{eq:fem_system}
Here $\stiffMat\in \mathbb{C}^{N \times N}$
with units $\unitOf\stiffMat=\si{\per\ohm}$,
$\electrodeCurrentMat\in\mathbb{C}^{N \times L}$
with $\unitOf\electrodeCurrentMat=\si{\per\ohm}$ and
$\electrodeImpedanceMat \in \mathbb{C}^{L \times L}$
with $\unitOf\electrodeImpedanceMat=\si{\per\ohm}$.
Consequently, the resistance matrix
$\resistanceMat$~\cite{GALAZPRIETO-2022} with
units $\unitOf\resistanceMat=\si{\ohm}$, satisfying
$\potentialVec = \resistanceMat\currentPattern$, can be
expressed as
\begin{equation}
\begin{aligned}
    \label{eq:resistance}
    \resistanceMat
    &=
    \inverse\stiffMat
    \electrodeCurrentMat
    \inverse{
        \paren*{
            \electrodeImpedanceMat
            -
            \transpose{\electrodeCurrentMat}
            \inverse\stiffMat
            \electrodeCurrentMat
        }
    }
    \\&=
    \transferMat\paren{{\electrodeImpedanceMat}  - \transpose{\electrodeCurrentMat} \transferMat }^{-1}
    \\&=
    \transferMat\schurMat^{-1}
    \in \mathbb{C}^{N \times L}
    \,,
\end{aligned}
\end{equation}
where $\transferMat =
\inverse\stiffMat\electrodeCurrentMat\in \Cset^{N
\times L}$ with $\unitOf\transferMat=\si{1}$
is the so-called \emph{transfer
matrix}~\cite{Knösche--Haueisen-2022}\cite{soderholm-2024},
whose factors $\stiffMat$ and $\electrodeCurrentMat$
are described in~\ref{app:admittivity-tensor}. The
matrix~$\transferMat$ can be pre-computed for a
given electrode configuration with the iterative
\PCG~\cite{sauer2018numerical} algorithm due to the
reciprocity of brain tissue~\cite{Rush-1969} and therefore
the hermitianness of \(\stiffMat\)~\cite{Potton-2004}.
The matrix \(\schurMat={\electrodeImpedanceMat}
- \transpose{\electrodeCurrentMat}
\transferMat\in\Cset^{L\times L}\)
is the \emph{Schur complement} of
\(\stiffMat\)~\cite{schur-complement-book-2005}. The
ungrounded electrode voltages $\potentialLossVec$
can be obtained by referring to the bottom row
of \eqref{eq:u_system}, i.e., $\currentPattern =
-\transpose{\electrodeCurrentMat} \potentialVec
+{\electrodeImpedanceMat}\potentialLossVec$.
The use of $\transpose\electrodeCurrentMat$
instead of $\ctranspose\electrodeCurrentMat =
\conj{\transpose\electrodeCurrentMat}$ for the
reciprocal mapping is justified by the nature
of our conductivity tensor: it contains only
symmetric or Hermitian conductivity and permittivity
components~\cite{Potton-2004}.

\subsection{Required Partial Derivatives}%
\label{sec:partial-derivatives}

To find out the required derivatives for the linearization,
we write the components of \(\resistanceMat\) depending on
\(\Zell\) as
\begin{equation}
 \begin{aligned}
    {\stiffMat}
    &=
    \tilde{{\stiffMat}}
    +
    \sum_{\ell = 1}^L
    \frac{1}{ { \Zell}\,\areaOf\electrode}
    \,
    \tilde{\massMat}^{(\ell)}
    \,,
    \\
    {\electrodeCurrentMat }
    &=
    \tilde\electrodeCurrentMat\electrodeImpedanceMat
    \,,
\end{aligned}
\label{linearization_A_B}
\end{equation}
where the following components of the sparse
symmetric matrices $\tilde{\stiffMat},
\tilde{\massMat}^{(\ell)}\in \mathbb{C}^{N \times N}$, and
${\tilde{\electrodeCurrentMat}}^{(\ell)}\in \mathbb{C}^{N
\times L}$ are not functions of \(\Zell\):
\begin{align*}
    \tilde{\stiffMat}_{ij}
    &=
    \int_\eell
    \admittivity
    \nabla
    \basisFn_i
    \cdot
    \nabla
    \basisFn_j
    \,
    \diff\surface
    \,,
    \\
    \tilde{\massMat}^{(\ell)}_{ij}
    &=
    \int_{\electrode[\ell]}
    \basisFn_i
    \basisFn_j
    \,
    \diff\surface\,,
    \\
    \tilde{\electrodeCurrentMat}_{i\ell}
    &=
    \frac 1 {\areaOf\electrode}
    \int_{\eell}\basisFn_i
    \,
    \diff\surface
    \,.
\end{align*}
The required derivative matrix via the product rule can
then be written as
\begin{equation}
    \pder
        \resistanceMat
        \Zell
    =
    \pder
        {\inverse\stiffMat}
        \Zell
    \,
    \electrodeCurrentMat
    \,
    \inverse\schurMat
    +
    \inverse\stiffMat
    \pder
        {\electrodeCurrentMat}
        {\Zell}
    \inverse\schurMat
    +
    \inverse\stiffMat
    {\electrodeCurrentMat}
        \frac
        {\partial\inverse\schurMat}
        {\partial \Zell}
    \,.
\label{R_Z}
\end{equation}
Differentiating the equation ${\bf V} {\bf
V}^{-1} = \idMat$ by $\Zell$ on both sides and
performing some algebraic manipulations, we
have~\cite{he2020zeffiro,pursiainen-kaasalainen-2016,matrix-cookbook}
\[
    \pder
        {\inverse{\mathbf V}}
        \Zell
    =
    -
    \inverse{\mathbf V}
    \pder
        {\mathbf V}
        \Zell
    \inverse{\mathbf V}
    \,.
\]
Thus, we can express the formula \eqref{R_Z} in the
following way:
\begin{equation}
 \begin{aligned}
    \pder
        \resistanceMat
        \Zell
    &=
    -
    \inverse\stiffMat
    \,
    \pder
        \stiffMat
        \Zell
    \,
    \resistanceMat
    +
    \inverse\stiffMat
    \pder
        \electrodeCurrentMat
        \Zell
    \,
    \inverse\schurMat
    -
    \resistanceMat
    \,
    \pder
        \schurMat
        \Zell
    \,
    \inverse\schurMat
    \,.
\end{aligned}
\label{dR_Z}
\end{equation}

To clarify the expressions involving the partial
derivatives with respect to $\Zell$, we apply
\eqref{linearization_A_B} to obtain the following
equations:
\begin{equation}
 \begin{aligned}
    \pder
        \stiffMat
        \Zell
    &=
    -
    \frac
    {1}
    {\Zell^2\,\areaOf{\electrode[\ell]}}
    \,
    \tilde{\massMat}^{(\ell)}
    \,,
    \\
    \frac
    {\partial {\electrodeCurrentMat}}
    {\partial \Zell}
    &=
    -
    \frac
    {1}
    {\Zell^2}
    \,
    \tilde{\electrodeCurrentMat}^{(\ell)}
    \,,
    \\
    \pder
        {\schurMat}
        {\Zell}
    &=
    \pder
        {\electrodeImpedanceMat}
        {\Zell}
    -
    \pder
        {\ctranspose\electrodeImpedanceMat}
        {\Zell}
    \transpose{\tilde{\electrodeCurrentMat}}
    \transferMat
    +
    \ctranspose{\electrodeCurrentMat}
    \inverse\stiffMat
    \,
    \pder{\stiffMat}{\Zell}
    \,
    \transferMat
    -
    \,
    \ctranspose{\electrodeCurrentMat}
    \inverse\stiffMat
    \pder
        {\electrodeCurrentMat}
        {\Zell}
    \,,
    \\
    \pder\electrodeImpedanceMat\Zell
    &=
    -
    \frac 1 {\Zell^2}
    \idMat^{(\ell)}
    \quad\text{and}\quad
    \pder
        {\ctranspose\electrodeImpedanceMat}
        \Zell
    =
    \frac
    {\abs {\Zell} - 2 \Zell}
    {\abs {\Zell}^3}
    \idMat^{(\ell)}
    \,.
\end{aligned}
\label{cdC_Z}
\end{equation}
Here $\partial\schurMat/\partial\Zell$ was derived using
the product rule and the identity $\ctranspose{(\matrixcmd
A\matrixcmd B)} = \ctranspose{\matrixcmd B}
\ctranspose{\matrixcmd A}$. This allows us to numerically
compute the derivative of \(\resistanceMat\), by defining
subroutines that produce the required partial derivative
matrices.

\end{document}